\documentstyle[titlepage,twoside,12pt]{book}
\begin{document}

\begin{titlepage}

{\LARGE \bf SHORT INTRODUCTION TO \\ \\ NONSTANDARD ANALYSIS} \\

\hspace*{3cm} {\bf ( to be completed)} \\ \\ \\

{\bf Elem\'{e}r E Rosinger \\ Department of Mathematics \\ University of Pretoria \\ Pretoria, 0002 South Africa \\
e-mail : eerosinger@hotmail.com} \\ \\ \\

\end{titlepage}

\newpage

\vspace*{10cm} \hspace*{4cm} Dedicated to Meda

\newpage

\pagenumbering{roman} \setcounter{page}{1}

{\LARGE \bf Table of Contents} \\ \\ \\

0~  Introduction  \hfill 1

\bigskip
1~ Two Versions of Nonstandard Analysis  \hfill 11

\medskip
~~~~1.1~~ A Short First Account \hfill 11

\medskip
~~~~1.2~~ Why the Need for Mathematical Logic ? \hfill 12

\medskip
~~~~1.3~~ Two Things which Tend to Put Off \\
\hspace*{1.37cm} Usual Mathematicians \hfill 16

\medskip
~~~~1.4~~ A Very Good Reason to Learn Nonstandard Analysis \hfill
20

\medskip
~~~~1.5~~ We Opt for the Robinson Approach \hfill 22

\medskip
~~~~1.6~~ A Few Historical Notes \hfill 22

\bigskip
2~ Constructing $^*{\bf R}$ \hfill 25

\medskip
~~~~2.1~~ $^*{\bf R}$ as a Field which Contains ${\bf R}$ \hfill
26

\medskip
~~~~2.2~~ Ideals and Filters \hfill 29

\medskip
~~~~2.3~~ Maximal Ideals and Ultrafilters \hfill 32

\medskip
~~~~2.4~~ Ultrafilters and Binary Valued Non-Atomic Measures
\hfill 36

\medskip
~~~~2.5~~ Filters and co-Filters \hfill 39

\medskip
~~~~2.6~~ Transferring the Total Order from ${\bf R}$ to $^*{\bf
R}$ \hfill 41

\medskip
~~~~2.7~~ Transfer of Relations and Functions \hfill 42

\medskip
~~~~2.8~~ Countable Saturation \hfill 50

\medskip
~~~~2.9~~ The Abundance of Ultrafilters \hfill 54

\bigskip
3~ The Transfer of Standard Properties \hfill 57

\medskip
~~~~3.1~~ Towards a General Transferring Mechanism \hfill 57

\medskip
~~~~3.2~~ Simple Systems and their Simple Languages \hfill 58

\medskip
~~~~3.3~~ Interpretation of Simple Languages \hfill 64

\medskip
~~~~3.4~~ The Transfer Property in its Simple Version \hfill 69

\medskip
~~~~3.5~~ Several General Results \hfill 76

\medskip
~~~~3.6~~ The Local and Global Structure of $^*{\bf R}$ \hfill 83

\medskip
~~~~3.7~~ Self-Similarity in $^*{\bf R}$ \hfill 95

\medskip
~~~~3.8~~ Sets which Are Not Transfers of Standard \\
\hspace*{1.37cm} Sets \hfill 96

\medskip
~~~~3.9~~ A Few Basic Applications to Calculus : Sequences \\
\hspace*{1.37cm} and their Limits \hfill 98

\medskip
~~~~3.10~ Proof of the Transfer Property in its Simple Version
\hfill 104

\bigskip
4~~ Short Review \hfill 105

\medskip
~~~~4.1~~ Two Stages of Transfer \hfill 105

\medskip
~~~~4.2~~ Nonstandard Approach to General Mathematics \hfill 107

\bigskip
5~ Superstructures and their Languages \hfill 113

\medskip
~~~~5.1~~ Superstructures \hfill 113

\medskip
~~~~5.2~~ Languages for Superstructures \hfill 120

\medskip
~~~~5.3~~ Interpretations in Superstructures \hfill 125

\medskip
~~~~5.4~~ Monomorphisms of Superstructures \hfill 127

\medskip
~~~~5.5~~ Ultrapower Construction of Superstructures \hfill 139

\medskip
~~~~5.6~~ Genuine Extensions, Hyperfinite Sets, Standard and \\
\hspace*{1.37cm} Nonstandard Entities, Concurrent Relations, \\
\hspace*{1.37cm} and Exhausting Sets \hfill 144

\medskip
~~~~5.7~~ Internal and External Entities, Comprehensiveness \hfill
160

\medskip
~~~~5.8~~ Permanence \hfill 170

\medskip
~~~~5.9~~ Saturation \hfill 173

\bigskip
6~ Looeb Integration and Measures \hfill 175

\medskip
~~~~6.1~~ Standard and Nonstandard Integration Structures \hfill
177

\medskip
~~~~6.2~~ The Loeb Construction of Integration Structures by
Standardization \hfill 181

\medskip
~~~~6.3~~ Loeb Measures \hfill 187

\bigskip
Bibliography \hfill 191

\newpage

\pagenumbering{arabic}
\setcounter{page}{1}
\setcounter{chapter}{-1}
\pagestyle{myheadings}
\markboth{E E Rosinger}{Nonstandard Analysis}

\chapter{Introduction}

This book, in fact, presents {\it two} introductions to Nonstandard Analysis. \\
Chapters 2 and 3 offer what may be seen as a rather minimal necessary approach in which basics of the fundamental
operation of {\it transfer} are introduced, motivated and given a few important applications. At that stage, the treatment
ventures as little as possible into the realms of {\it formal languages} in Mathematical Logic, a subject which is not
familiar to most of present day mathematicians. Such an approach was first suggested in Keisler [2], and developed in
some detail in chapter 1 of Hurd \& Loeb. \\
In chapter 5 a second and higher level introduction to Nonstandard Analysis is presented. This is based on the concept
of Superstructures and the associated formal languages, and brings the reader to the beginnings of the present day
treatment of the subject. The idea of this second way of introduction originated in Robinson \& Zakon, and has since
known a wide following. Its main aim was to minimize the extent to which Mathematical Logic becomes involved in
Nonstandard Analysis, an extent which in the first major and systematic presentation of Nonstandard Analysis, given in
Robinson [2], had indeed been considerable, even if rather elementary, in fact. \\

But now, we may as well ask : what is so special about Nonstandard Analysis ? \\

A short answer may be as follows. Given a usual mathematical theory, like for instance, Algebra, Calculus, Topology,
Measure, Probability, or Ordinary and Partial Differential Equations, and so on, with the help of Nonstandard Analysis,
one can embed such a usual theory into a significantly more {\it rich} nonstandard one. And in doing so, two advantages
are encountered. First, the corresponding more rich nonstandard theory proves to be more {\it intuitive} and thus easier
to deal with as well. Second, one obtains a general {\it transfer} mechanism between the two theories, according to
which any statement in the usual version is true, if and only if its transfer into the nonstandard theory is true. \\

What is called {\it transfer} in Nonstandard Analysis is so far a feature unique to this theory among all other
mathematical theories. And in its essence it establishes a {\it two way} interaction between the usual mathematical
structures, which we can denote by (UMS), and on the other hand, the associated and significantly more rich
nonstandard mathematical structures, (ENMS). \\
However, it is important to note from the start that the interest in, and the power of Nonstandard Analysis is not limited
to this transfer alone. Indeed, within the associated extended nonstandard mathematical structures (ENMS) a wealth of
specific nonstandard concepts, methods and results can be obtained which, as mentioned, will typically have the
advantage of being more {\it intuitive} than their correspondents within the usual mathematical structures (UMS). And
this fact alone, even without transfer, can help in better understanding and dealing with these usual mathematical
structures. \\
In fact, as seen in Note 5.7.2 in section 7, in chapter 5, transfer is actually limited between the usual entities in (UMS),
and on the other hand, the so called {\it internal} entities in the extended nonstandard mathematical structures (ENMS).
Nevertheless, transfer is a crucial way to relate the usual and their corresponding nonstandard mathematical
structures. \\

The other way for such a relationship, a way already {\it beyond} the realms of transfer, is based on certain
developments within the nonstandard theories (ENMS) which can nevertheless help in the usual theories (UMS) as well.
Examples about the interest for usual theories (UMS) in such nonstandard developments within (ENMS), developments
which take us outside of the operations of transfer, can already be found with such simple usual topological properties
of sets, like for instance, being open, closed, or compact. Indeed, nonstandard descriptions can appear in their
characterizations, and can do so in highly clarifying, relevant and intuitive ways, even if such descriptions can no
longer be transferred back into the usual terms, see Proposition 3.6.4. \\

It follows that we can summarize the above in the diagram

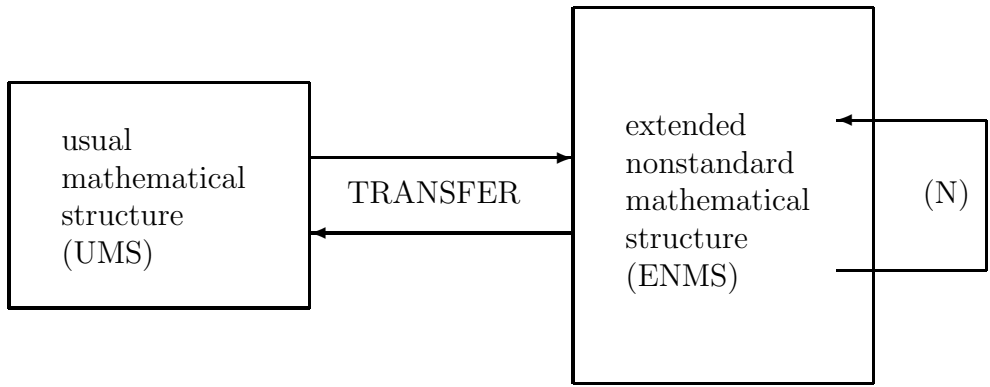
\begin{math}
\setlength{\unitlength}{1cm}
\thicklines
\begin{picture}(15,5.5)
\put(0,0.5){\line(0,1){3}}
\put(0,0.5){\line(1,0){4}}
\put(0,3.5){\line(1,0){4}}
\put(4,0.5){\line(0,1){3}}
\put(7.5,-0.5){\line(0,1){5}}
\put(7.5,-0.5){\line(1,0){4}}
\put(7.5,4.5){\line(1,0){4}}
\put(11.5,-0.5){\line(0,1){5}}
\put(4,2.5){\vector(1,0){3.5}}
\put(4.5,1.9){$\mbox{TRANSFER}$}
\put(7.5,1.5){\vector(-1,0){3.5}}
\put(0.7,2.6){$\mbox{usual}$}
\put(0.7,2.1){$\mbox{mathematical}$}
\put(0.7,1.6){$\mbox{structure}$}
\put(0.7,1.1){$\mbox{(UMS)}$}
\put(8.2, 2.8){$\mbox{extended}$}
\put(8.2,2.3){$\mbox{nonstandard}$}
\put(8.2,1.8){$\mbox{mathematical}$}
\put(8.2,1.3){$\mbox{structure}$}
\put(8.2,0.8){$\mbox{(ENMS)}$}
\put(13,3){\vector(-1,0){2}}
\put(11,1){\line(1,0){2}}
\put(13,1){\line(0,1){2}}
\put(12.15,1.9){$\mbox{(N)}$}
\end{picture}
\end{math} \\ \\

\medskip
where the arrow (N) denotes those nonstandard developments within (ENMS) which do {\it no} longer fall within the
realms of transfer. \\
Also, it should be noted that, unlike in the above diagram where for the ease of representation the structures (UMS)
and (ENMS) appear as disjoint, in reality, the usual mathematical structures (UMS) are contained as {\it particular}
cases the extended nonstandard mathematical structures (ENMS). \\

The various approaches to Nonstandard Analysis which, more than four decades after the first emergence of the
subject, still appear to compete with one another, can create a certain problem for mathematicians interested in the
subject. And one way to bring in some clarity in this regard is by trying to analyze the reasons which have led to such
a situation. \\
In the present, and in the next chapter, such an analysis is attempted. And regardless of the extent to which it may
have attained its aim, the author of this book strongly suggests the reader to take the approximately half an hour
needed in order to read it with some care. \\ \\

During the last nearly four decades, since the first systematic presentation of Nonstandard Analysis in the 1966 book
of Abraham Robinson, quite a few excellent treaties have been published on the subject, a number of them being
mentioned in the reference. \\
And yet, in spite of the clear and not seldom major advantages Nonstandard Analysis can offer in research in a large
variety of branches of usual mathematics, the fact remains that the number of mathematicians who have learned it and
applied it in their research in such branches is rather limited. \\

The only more notable exception so far is in the study of time continuous stochastic processes where, since the
introduction in 1973 by Peter A Loeb of nonstandard methods in measure and probability theory, such methods have
gained a certain ground, owing to their ability to avoid major difficulties coming from the Kolmogorov model of
probability, see for instance Albeverio, et.al. \\

Furthermore, there are mathematicians who, although familiar with Nonstandard Analysis, prefer not to use it in their
research publications, so as to avoid adding difficulties upon reviewers and editors, not to mention potential readers. \\

Quite likely, one of the main reasons for such a state of affairs is in the peculiar, and so far unique feature of
Nonstandard Analysis among various more usual branches of mathematics, to have to involve one rather inevitably in
a certain amount of Mathematical Logic, see for details sections 2 and 3, in chapter 1. \\

A further quite likely important, and closely related reason for the lack of spread in interest in Nonstandard Analysis
comes from another feature which, so far, is equally unique to Nonstandard Analysis. Namely, in this theory one
becomes essentially involved in an ongoing {\it two way interaction} between the usual mathematical structures, and on
the other hand, their {\it nonstandard extensions}. And such extensions are rather rich in structure, far more so than one
is accustomed to, when dealing with the usual extensions encountered in algebra or topology. After all, it is precisely
due to this wealth of structure of nonstandard extensions that they can often be so near to our better intuitions, and also
provide powerful methods. \\

And to add to the picture, this two way interaction between usual and nonstandard structures does at least at the
beginning require a certain amount of Mathematical Logic. Indeed, this interaction is done through the {\it transfer}
mechanism, and based on what is the center piece of Nonstandard Analysis, namely, the Transfer Property. And the
proof, as well as the use - at least at the beginning - of transfer happen in the realms of {\it formal languages}. Thus the
involvement of Mathematical Logic. \\

A good example of such a threshold in complexity, when compared to the usual mathematical approaches, can be
found in those impressive two books aimed to offer a First Course in Calculus, published by H J Keisler, in 1976. \\ \\

As it turns out so far, in the balance between advantages, and on the other hand, the difficulties coming from the
involvement of Mathematical Logic and the mentioned complexity threshold, the disadvantages proved to be more
important for large numbers of usual mathematicians, when facing Nonstandard Analysis. \\ \\

And then, the aim of the present book is to try to redress such an imbalance. And the audience aimed at are research
mathematicians involved in usual branches of mathematics. \\
In this regard, this book tries to do the following three things : \\

- offer a first contact with the subject, and as such, be short, rigorous, general enough, and also clear, \\

- cut the involvement of Mathematical Logic to a minimum, and at the same time, explain why it cannot be cut
completely, see again sections 2 and 3, in chapter 1, as well as sections 1 - 4, in chapter 3,  \\

- present an essential core of the mentioned two way interaction between the usual spaces in mathematics, and on the
other hand, their nonstandard extensions, and do so in {\it two stages}, the first, in chapters 2 and 3, and the second, in
chapter 5, in order to try to minimize the threshold in complexity which is inevitably involved. \\

A reader having gone through this book is, therefore, supposed to reach the following situation : \\

- have quite a clear understanding of what Nonstandard Analysis may in fact be about, and \\

- be ready to read rather at ease any of the excellent more detailed and voluminous treaties, in case he or she may
want to venture further in-depth into the subject. \\

And now, about the structure of this book. \\

Chapter 1 is a preliminary commentary which sets out briefly and informally some of the more important unusual
aspects one may have to face, when trying to study for the first time Nonstandard Analysis. Among such aspects, as
mentioned, is the reason one cannot do without some Mathematical Logic, and what are the further reasons usual
mathematicians may find the subject off putting. \\
Also, a few of the more important stages in the history of the subject are indicated, together with appropriate
references. \\

Chapters 2 and 3 form the {\it first stage} in the escalation of the mentioned threshold in complexity. \\

In chapter 2, the totally ordered field of nonstandard reals $^*{\bf R}$, which extends the usual reals ${\bf R}$, is
constructed {\it solely} with usual mathematical methods, such as algebra, ultrafilters and non-atomic binary valued
measures. This approach, which does {\it not} make any use of Mathematical Logic, can already allow the nonstandard
extension - through what is later, and in its general form, called {\it transfer} - of such fundamental mathematical
entities as functions, and more generally, relations. \\
Unfortunately however, within such a usual mathematical approach, one cannot come anywhere near to using the true
powers of Nonstandard Analysis, powers which rest upon a full use of the mentioned two way interaction between usual
and nonstandard mathematical structures, and interaction made by the method of {\it transfer}. \\

And then, in chapter 3, a rock bottom simple version of the {\it transfer} mechanism is presented, the mechanism which
is supposed to perform the mentioned two way interaction. This simple version was introduced in Keisler [2], and was
also used with considerable benefit in Hurd \& Loeb. With this version, however, no matter how simple it may be, one
must already employ some elements of Mathematical Logic, as one simply can no longer escape the reasons
mentioned in section 2, in chapter 2. \\

However, a special effort is made in chapter 3, in order to clarify as much as possible those few elements of
Mathematical Logic which may be less familiar to usual mathematicians. \\

The highlight of chapter 3 is a corresponding simple version of the celebrated Transfer Property in Nonstandard
Analysis, followed by some of its more simple, yet important applications. \\

Upon reading chapters 2 and 3, one may already get a certain idea about the main outlines of Nonstandard Analysis,
and as such, the reading of these two chapters may be seen as worthwhile from the point of view of the rather minimal
effort involved. A short review of the essential aspects of the constructions in chapters 2 and 3 is presented in chapter
4. This can help towards a better understanding of the general constructions based on Superstructures, introduced in
chapter 5. \\

The {\it second stage} in the escalation of the mentioned threshold in complexity starts in chapter 5, with the
introduction of Superstructures and their formal languages. Back in 1969, A Robinson and E Zakon introduced this avenue
in Nonstandard Analysis precisely in order to bring it back as much as possible to usual set theory, and thus minimize
the amount of Mathematical Logic employed. Such a return to the more usual ways of mathematics appeared to be
quite welcome, after the first major and systematic presentation of Nonstandard Analysis in Robinson [2], a
presentation which was relying quite heavily on certain elements of Mathematical Logic and Model Theory, elements
not familiar among usual mathematicians. \\
Needless to say, ever since, the method of Superstructures has proved to be particularly effective and useful. \\

Yet in this book we do not start with Superstructures, and instead leave them to the second stage. The main reason for
that is in the nature of the formal languages which are needed in order to take full advantage of the richness of usual
mathematics allowed by Superstructures. Indeed, such formal languages may appear rather complicated for those
mathematicians who are not familiar with Mathematical Logic. And then, as suggested in Keisler [2] and Hurd \& Loeb,
we opted for introducing Superstructures only in the second stage. \\

The reward in chapter 5 is in the formulation and proof of the full version of the Transfer Property, which from then on,
can be used in a variety of applications to usual mathematics. \\

In chapter 6, as one of the most important and novel applications so far, an introduction to Loeb measures is
presented. \\ \\

And at last, how about Exercises and Problems ? \\

Textbooks, and even more so introductory ones, appear to have as a compulsory component frequent exercises and
problems which pop up at various places in the main text. In some such books, one can even find a number of such
exercises and problems solved, or at least, provided with hints for solution. On occasion, several of such exercises
and problems are also used to present results of a certain theoretical interest which, however, were not included in the
main text, in order to keep to some brevity. \\

A disadvantage of such ways is that the main text with its line of thinking is often interrupted, and not seldom in a
somewhat arbitrary manner. Indeed, typically, it is not quite clear at first sight to the beginner reader how and why
many of the given exercises and problems may fit in at the precise place they happen to appear within the main line of
thinking. \\
Also, the exercises and problems which do not come with solutions or hints may disrupt by the simple fact that, when
attempted to be solved, they may take up more time than their importance deserves, or the respective reader may in
fact fail to solve them and thus gives up on them, remaining with a certain negative lingering feeling. \\
On the other hand, if one simply skips the exercises and problems, one may remain with a certain feeling of
superficiality. \\

In this book, instead of the mentioned ways, we chose to adopt what may be seen as an "informal" approach to
exercises and problems. Namely, the whole text is dedicated to the main line of thinking only, without any interruptions
or disruptions given by the usual and frequent explicit lists of exercises and problems. \\
And yet, the text often has certain limited and calculated gaps in the proofs or various statements. These gaps,
however, can be filled in based on a satisfactory understanding of the text which precedes them. \\
In this way it is hoped that two objectives can be attained at the same time : the main line of thinking is followed
continually, and at the same time, at various places, the reader is expected to apply what he or she is supposed to
have learned and understood already from the earlier part of the text. \\ \\

\chapter{Two versions of Nonstandard Analysis}

{\bf 1. A Short First Account} \\ \\

Nonstandard Analysis is {\it not} only about extending the set of usual real numbers ${\bf R}$ into the set of
nonstandard reals $^*{\bf R}$. \\
In fact, every usual mathematical theory can have a similar nonstandard extension. And the extension of ${\bf R}$ to
$^*{\bf R}$ has only been the one which, as they say, started the ball rolling ... \\

The first full version of Nonstandard Analysis was presented in 1966 by A Robinson. This version relied on a certain
familiarity with Mathematical Logic, and in particular, Model Theory, and as such, it proved to be less than an easy
reading for most of usual mathematicians. Nevertheless, this approach has the major, and so far, unique advantage to
offer a {\it constructive} presentation of the various nonstandard realms, and consequently, it has known important
follow up, such as for instance in Machover \& Hirschfeld, Keisler [1,2], Stroyan \& Luxemburg, Davis, Goldblatt. \\
More recently, by placing this constructive approach to Nonstandard Analysis in the framework of what is called
Superstructures, much of the technical aspects related to Mathematical Logic have been dispensed with, see Hurd \&
Loeb, Albeverio et.al., Cutland, Kursaev \& Kutateladze. \\

However, as seen in section 2, a certain involvement of Mathematical Logic in Nonstandard Analysis is {\it essential},
and thus unavoidable. And in fact, Nonstandard Analysis is, so far, the first major branch of mathematics, a branch
which in its themes and aims is clearly outside of Mathematical Logic, yet  it appears to need to contain a certain
amount of Mathematical Logic among its methods. \\

A second and alternative approach to Nonstandard Analysis was presented in 1977 by Edward Nelson. The striking
aspect of this approach is that it is {\it axiomatic}. More precisely, by adding only one new predicate of one variable,
and three new axioms to usual set theory, it manages to recover just about all of Nonstandard Analysis. Therefore, from
the start, and ever after, this approach does not take one outside of the customary ways of everyday mathematics,
except for the fact that the mentioned additional new predicate and three axioms appear in formulations which involve
what in Mathematical Logic is called Predicate Calculus. And as such, they are somewhat complicated, and certainly
less than intuitive at first. Therefore, their subsequent use, which of course is essential in order to get into Nonstandard
Analysis, does require a certain amount of preliminary time and exercise, until one may indeed become familiar and at
ease with them, see Lutz \& Goze, or Diener \& Diener. \\ \\

{\bf 2. Why the Need for Mathematical Logic ?} \\ \\

From the start, the embedding of the set of real numbers ${\bf R}$ into the set of nonstandard reals $^*{\bf R}$ brings
up several fundamental facts, which are new and unprecedented in usual mathematics. \\
Indeed, in this construction the interest is not merely to embed the set ${\bf R}$ into some large set $^*{\bf R}$, since
this could clearly be done in many trivial ways. \\

Instead, the interest is : \\

-~~~ to extend ${\bf R}$ as a {\it totally ordered field} - that is, as a structure \\
      \hspace*{0.5cm} $( {\bf R},~ + ~,~ . ~,~ \leq )$ - into a {\it larger} totally ordered field $( ^*{\bf R},~ + ~,~ . ~,~ \leq )$ \\

and what is of {\it outmost concern} in addition : \\

-~~~ to have lots of important properties of ${\bf R}$ preserved in $^*{\bf R}$, as a \\
      \hspace*{0.5cm} consequence of a strong and clear connection between such \\
      \hspace*{0.5cm} properties of ${\bf R}$ and $^*{\bf R}$. \\

And certainly, the property of being a totally ordered field does go over - or in nonstandard terms, does {\it
transfer} - from ${\bf R}$ to $^*{\bf R}$. \\

Thus the question is : which are {\it all} those properties of ${\bf R}$ which do transfer to $^*{\bf R}$ as well ? \\
Furthermore, how can we find a clear and rigorous {\it general} mechanism for such a transfer ? \\
Also, which are those properties of ${\bf R}$ which do {\it not} transfer to $^*{\bf R}$ ? \\
And certainly, there are some important such properties. For instance, ${\bf R}$ is Dedekind order complete and
Archimedean, while $^*{\bf R}$ is neither. \\

And related to the above, we should note the important objective - typical for science, and even more so for
mathematics - that such a transfer mechanism should rather function as a {\it wholesale industry}, and not merely as a
piece by piece artisanship. \\

However, here already, certain basic problems arise, since such a transfer is far from being trivial. \\
For instance, as is well known, ${\bf R}$ is the only totally ordered field which is complete and also Archimedean.
Therefore, $^*{\bf R}$ being a lager totally ordered field, thus different form ${\bf R}$, it can {\it no} longer be
Archimedean. \\
There also are other properties of ${\bf R}$ which do no longer hold in $^*{\bf R}$. For instance, and as mentioned,
${\bf R}$ is Dedekind order complete, that is, each bounded subset in ${\bf R}$ has a lower bound and upper bound in
${\bf R}$, while as we shall see, this property does not hold in $^*{\bf R}$. \\

Now, in usual mathematical theories, we have been busy establishing certain properties {\it inside} the framework of
one or another given mathematical structure. And therefore, we have not yet encountered the situation when the issue
arose to be able as well to {\it transfer} such properties from one structure to another, let alone do so {\it systematically}
on a massive scale, and not merely upon a tedious and elaborate case by case approach. \\

Indeed, one should not fail immediately to note that in Nonstandard Analysis such a transfer does {\it no} longer take
place in the usual manner between two mathematical structures, where a certain limited number of properties in the two
given structures become related, like for instance it is done by algebraic homomorphism, topological
homeomorphisms, and so on, or in terms of category theory, by morphisms, functors or natural transformations. \\
Instead, such a transfer is now establishing a general and systematic correspondence which ranges over a very large
variety of often very different {\it properties} in the respective two structures. \\

Therefore, the entities which are supposed to be transferred this time are very large classes of {\it rigorous
mathematical sentences} which express those properties. \\

And then, quite naturally, the {\it rigorous} construction of such a transfer of large classes of sentences will by
necessity involve a certain amount of {\it formalism}, typical for Mathematical Logic. \\

It is, in this way, in Nonstandard Analysis that we face for the first time this specific kind of problem of transfer between
two usual mathematical structures. And indeed, here we are facing two such structures, since what is of interest to us is
not only how much ${\bf R}$ and $^*{\bf R}$ are similar through transfer, but also how much they are different, in spite
of all the transfer. After all, there would not be much interest in $^*{\bf R}$, if it were just about the same with
${\bf R}$. \\

So then, after we manage to construct $^*{\bf R}$ as a totally ordered field which contains strictly ${\bf R}$, a
construction which in fact is not so difficult, a main issue arises, which as mentioned, among all usual mathematical
theories so far is typical for Nonstandard Analysis only. \\
Namely, to {\it formalize rigorously} the mechanism of {\it transfer}. \\
And so far, the only way we know how to perform a corresponding rigorous formalization is by involving a certain amount
of Mathematical Logic. \\

In this regard, it should be noted that none of mentioned two approaches to Nonstandard Analysis, namely, those
originated by A Robinson and E Nelson, respectively, can do completely without Mathematical Logic. \\

After all, more than three centuries ago, Leibniz had quite a clear intuition of infinitesimals, and used them appropriately
in developing the first stages of Calculus. What prevented for just about three centuries, and until our own times, the
emergence of Nonstandard Analysis, however, has been the lack of a rigorous enough mathematical theory which
would be able to deal with the transfer of large classes of mathematical sentences that express mathematical
properties of interest. \\
The Mathematical Logic of the 20th century, in particular, Model Theory, which owes much to A Robinson, proved to
provide for the first time such a possibility. \\

In conclusion, {\it transfer} will among others have to mean the following two essential things : \\

(TE)~:~~ the {\it transfer} of various {\it mathematical entities} from ${\bf R}$ to $^*{\bf R}$, \\
         \hspace*{1.45cm} the first and simplest case of it being the field homomorphism \\
         \hspace*{1.45cm} $^*(~)$, see (2.1.4), \\

(TS)~:~~ the {\it transfer} of {\it sentences} $\Phi$ from a given {\it language} ${\cal L}_{\bf R}$ about \\
         \hspace*{1.45cm} mathematical entities in ${\bf R}$, to corresponding sentences $^*\Phi$ in \\
         \hspace*{1.45cm} a language ${\cal L}_{^*{\bf R}}$ about mathematical entities in $^*{\bf R}$. \\

Furthermore, the transfer of sentences is a {\it two way} process. Namely, when we want to prove a standard property of
${\bf R}$, we formulate it as a sentence $\Phi$ in the language ${\cal L}_{\bf R}$, then we transfer it to the sentence
$^*\Phi$ in the language ${\cal L}_{^*{\bf R}}$, and prove this transferred sentence $^*\Phi$ in $^*{\bf R}$, after
which we deduce by reverse transfer that the initial sentence $\Phi$ is true about ${\bf R}$. \\

What is particularly important to note and keep in mind with respect to the above two stages (TE) and (TS) of transfer
is that (TS) is far more {it powerful} and {\it systematic} than (TE). The stage (TE), however, has the advantage that
it can be performed by employing usual mathematics. On the other hand, the stage (TS) needs a few basic concepts and
methods from Mathematical Logic. Chapter 2 is an example of performing the stage (TE) of transfer, while in chapter 3
both stages (TE) and (TS) are employed, with an accent on the latter. \\
Starting with chpater 5, one does no longer care much to differentiate between such stages. \\ \\

{\bf 3. Two Things which Tend to Put Off Usual Mathematicians} \\ \\

As mentioned, in the last two or three decades, both approaches to Nonstandard Analysis have managed to reduce
significantly the amount of Mathematical Logic involved. And yet, there is a significant reluctance among most of the
usual mathematicians in embracing the methods of Nonstandard Analysis, although such methods are often strikingly
powerful and also close to one's intuition, and can be used with great success in a large variety of usual mathematical
theories, as seen in most, if not in fact, all of the items in the References relating to the subject. \\

An answer to such a state of affairs may be as follows : there appear to be a few important stumbling blocks. \\
Two of them, which one encounters from the start, when trying to enter into the realms of Nonstandard Analysis, will be
mentioned here. \\

There are, however, at lest three other ones as well, which will emerge, starting with chapter 5. Let us mention them
only in passig now, in order to make the reader aware, when their time comes. \\
The use of Superstructures leads to the need to go beyond the customary ways in mathematics, where we only deal
with the following three, or at most four levels, namely, elements in some given set, subsets of such elements, sets of
such subsets, and perhaps, sets of sets of such subsets. On the other hand, in Superstructures one has to iterate
countably infinite times the transition from a set to the set of all of its subsets. Of course, one usually is only involved
in a few such transitions at a time, however, unlike in usual mathematics, one has to keep in one's awareness the fact
that in a variety of important nonstandard constructions and arguments, a countably infinite number of such transition
are implicitly involved. \\
Then there is the richness, and consequent complexity of the nonstandard concepts and situations, which is beyond
what one is accustomed to usually in mathematics. Indeed, one has to deal with such concepts as {\it internal} versus
{\it external}, {\it standard} versus {\it nonstandard}, or {\it hyperfinite}, {\it overflow}, {\it underflow}, {\it saturation}, and
so on. \\
But perhaps above all, one has to get accustomed to the idea that, unlike our good old and unique real line ${\bf R}$,
its nonstandard extensions are {\it far from} being unique. And this in itself is far from being a weakness of
Nonstandard Analysis. On the contrary, it places clearly and firmly in front of us the fact that all what we are dealing
with - be it ${\bf R}$, or the variety of its nonstandard extensions $^*{\bf R}$ - are rather {\it models} of a deeper and
quite hard to fathom entity, which we may call {\it one dimensional continuum}. After all, it is not for nothing that the
Axiom of the Continuum in Set Theory is {\it independent} of the other axioms. And it may only be due to our long time
habit with the real line ${\bf R}$, a habit which made us believe that we have already perfectly well understood what a
one dimensional continuum is, that we may now be surprised, and perhaps even put off, when having to face the idea
that, so far, we have only been dealing with models. \\

And now, let us return to what appear to be the first two important stumbling blocks. \\

First, we should note that in usual mathematics, one is assumed to have acquired a natural rigorous way of thinking,
as far as its {\it logical correctness} is concerned. Consequently, one's thinking is assumed to function quite perfectly
and very much freely and at ease, even if rather informally, intuitively and automatically, whenever its logical aspects
are involved. In fact, one is {\it never ever } supposed to get concerned about the logical correctness of one's thinking,
since it can simply be taken for granted. \\
Indeed, when doing usual mathematics, there are two {\it dividing lines} involved. \\
First, usual mathematics is quite outside of Mathematical Logic. \\
Second, there is a certain separation between {\it rigour} and {\it formalism}. Namely, {\it rigour} is, of course, assumed
to apply to everything we do, including the logical aspects of out thinking. On the other hand, {\it formalism} is nearly
exclusively limited to the final expression of the specific mathematical structure studied and developed. And when it
comes to one's logical thinking used in the process, one simply assumes that one already knows quite well how to do
it rigorously, and thus one can feel pretty free and informal about it, without any risk to rigour.

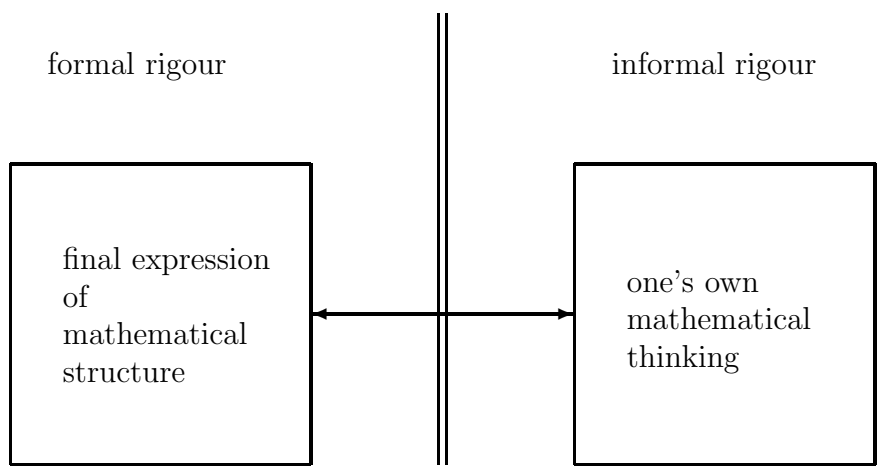
\begin{math}
\setlength{\unitlength}{1cm}
\thicklines
\begin{picture}(15,7)
\put(1,0){\line(0,1){4}}
\put(1,0){\line(1,0){4}}
\put(1,4){\line(1,0){4}}
\put(5,0){\line(0,1){4}}
\put(8.5,0){\line(0,1){4}}
\put(8.5,0){\line(1,0){4}}
\put(8.5,4){\line(1,0){4}}
\put(12.5,0){\line(0,1){4}}
\put(5,2){\vector(1,0){3.5}}
\put(6,2){\vector(-1,0){1}}
\put(6.7,0){\line(0,1){6}}
\put(6.8,0){\line(0,1){6}}
\put(1.5,5.2){$\mbox{formal rigour}$}
\put(9,5.2){$\mbox{informal rigour}$}
\put(1.7,2.6){$\mbox{final expression}$}
\put(1.7,2.1){$\mbox{of}$}
\put(1.7,1.6){$\mbox{mathematical}$}
\put(1.7,1.1){$\mbox{structure}$}
\put(9.2, 2.3){$\mbox{one's own}$}
\put(9.2,1.8){$\mbox{mathematical}$}
\put(9.2,1.3){$\mbox{thinking}$}
\end{picture}
\end{math} \\ \\

Now, when it comes to the issue of transfer in Nonstandard Analysis, we have to face the unprecedented fact in usual
mathematics that the above two customary dividing lines are somewhat shifted, as suggested in the next illustration.
And consequently, it can easily appear as a strange, if not in fact, intolerable imposition to have, from the start,
subjected one's own mathematical thinking to certain criteria formulated in terms of Mathematical Logic, criteria which
are needed in order to perform a proper, that is, rigorous transfer. \\

Of course, after a certain period of familiarization, which is needed in both approaches to Nonstandard Analysis, one
will fully recover one's feeling of natural ease ... \\
This time, however, it is not merely about getting familiar with one more usual mathematical theory, but about getting
familiar with the situation suggested in :

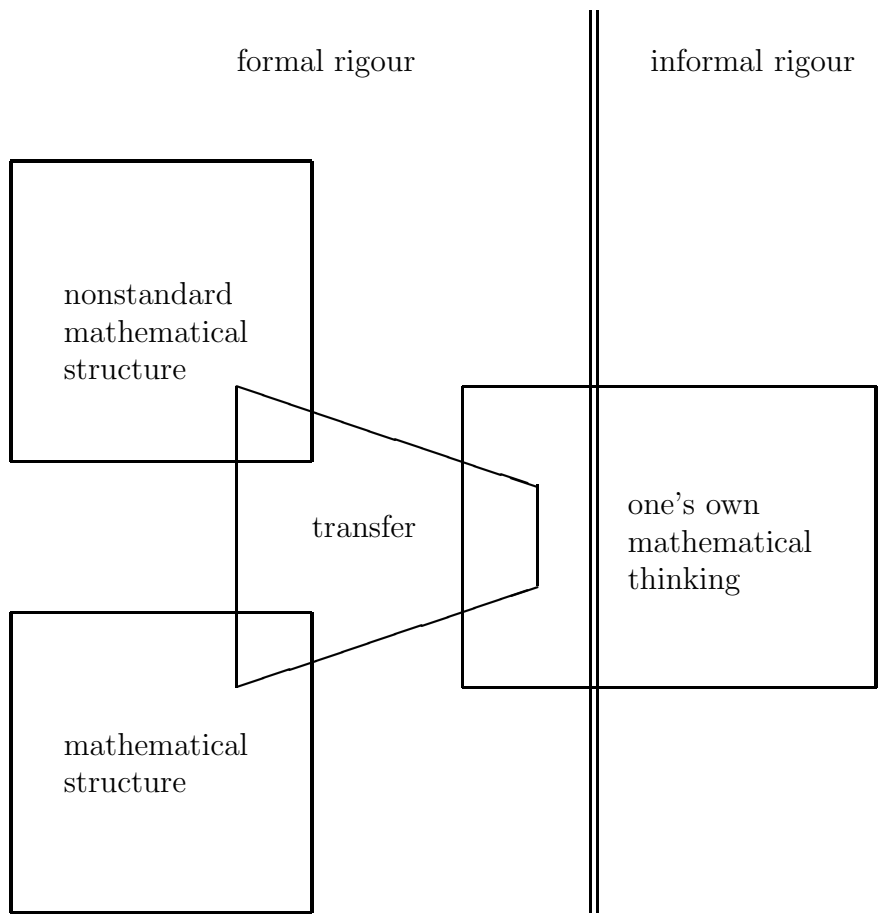
\begin{math}
\setlength{\unitlength}{1cm}
\thicklines
\begin{picture}(15,13.)
\put(1,6){\line(0,1){4}}
\put(1,6){\line(1,0){4}}
\put(1,10){\line(1,0){4}}
\put(5,6){\line(0,1){4}}
\put(1.7,8.1){$\mbox{nonstandard}$}
\put(1.7,7.6){$\mbox{mathematical}$}
\put(1.7,7.1){$\mbox{structure}$}
\put(1,0){\line(0,1){4}}
\put(1,0){\line(1,0){4}}
\put(1,4){\line(1,0){4}}
\put(5,0){\line(0,1){4}}
\put(7,3){\line(0,1){4}}
\put(7,3){\line(1,0){5.5}}
\put(7,7){\line(1,0){5.5}}
\put(12.5,3){\line(0,1){4}}
\put(8.7,0){\line(0,1){12}}
\put(8.8,0){\line(0,1){12}}
\put(4,11.2){$\mbox{formal rigour}$}
\put(9.5,11.2){$\mbox{informal rigour}$}
\put(1.7,2.1){$\mbox{mathematical}$}
\put(1.7,1.6){$\mbox{structure}$}
\put(9.2, 5.3){$\mbox{one's own}$}
\put(9.2,4.8){$\mbox{mathematical}$}
\put(9.2,4.3){$\mbox{thinking}$}
\put(4,3){\line(0,1){4}}
\put(4,3){\line(3,1){4}}
\put(4,7){\line(3,-1){4}}
\put(8,4.35){\line(0,1){1.35}}
\put(5,5){$\mbox{transfer}$}
\end{picture}
\end{math} \\ \\

In other words, it is in Nonstandard Analysis that for the first time in usual mathematical theories, the entities
subjected to {\it rigorous formal} mathematical  concern are {\it not only} the usual mathematical structures, {\it but
also} the rigorous mathematical sentences about such objects. Thus, it is for the first time that {\it formal}
mathematical concern touches {\it both} upon usual mathematical structures, {\it and} at least a good part of our thinking
about them. \\
As a consequence, when starting with Nonstandard Analysis, at least at the beginning, we can no longer enjoy thinking
all the time freely and quite informally, based alone on our well trained usual mathematical intuitive and informal
rigour. \\
And this certainly can put off ... \\

A second thing which can put off is the following. We are far better at learning more and more about a game which we
are already familiar with, than to learn a rather new game, and become fast enough sufficiently proficient at it. \\
Now, Nonstandard Analysis is quite a new game, not least because of the mentioned ways it intrudes upon our usual
mathematical thinking. Added to it are two more rather off putting facts. One is that a usual mathematician would learn
Nonstandard Analysis not so much for its own sake, but for its powerful applications in the field of mathematics of one's
specific interest. And whatever can be proved in usual mathematics by methods of Nonstandard Analysis, can also be
proved, even if often in far more difficult and far less intuitive ways, by usual mathematical methods. \\

And yet ... there is ... \\ \\

{\bf 4. A Very Good Reason to Learn Nonstandard Analysis} \\ \\

And there are, in fact, several good reasons ... \\

First and foremost is quite likely that, once in a decade or so, one should anyhow spend some time and learn the
basics of yet another mathematical theory. And the best way to do that is from suitable short enough books.
Unfortunately however, there have not been many like that written in more recent times ... \\
It is indeed but a result of a rather sorry state of too early and too narrow specialization that one does not even think of
engaging voluntarily and with joy in such a venture. And yet, our knowledge of and interest in mathematics should not
look like a ... Delta function ... which is extremely high over a very narrow range, and quite zero everywhere else ... \\
Instead, our knowledge and interest should rather look like a Gauss bell, where the ranges of lower knowledge, but not
necessarily interest as well, may spread quite widely ... \\
And if we are at Gauss bells, then why not even more general such curves ? \\
Like for instance, curves with more than one maximum ? \\

Of course, such a less than high level knowledge is not meant for pursuing research in the respective fields. However,
keeping up the interest in a larger variety of mathematical fields, and doing so rather as an matter of culture, can have
in the longer run any number of advantages, including for one's own specialized research. \\
And by far the best way in order to understand such advantages is, of course, by practicing a ... non-Delta function ...
approach ... \\

And then, Nonstandard Analysis is a very good such mathematical theory to learn. Indeed, here are some of the
reasons : \\

-~~~ It is a relatively new theory, being started as such in the 1960s. \\

-~~~ It is, among usual mathematical theories, quite special due, among \\
      \hspace*{0.5cm} others, to its unique features mentioned in sections 2 and 3. \\

-~~~ It has massive applications in a large variety of usual mathematical \\
      \hspace*{0.5cm} theories. \\

-~~~ Due to the introduction of infinitely small and infinitely large \\
      \hspace*{0.5cm} elements, it enriches the usual concept of sets in wonderful, and \\
      \hspace*{0.5cm} also useful ways. And based on that, it offers our mathematical \\
      \hspace*{0.5cm} intuition a well grounded and fruitful freedom to operate. \\

Related to the last point above, let us recall a rather instructive event from the history of mathematics which can show
quite clearly the power and advantages of using an {\it intuition} based on {\it infinitesimals}, and on the other hand,
the limitations and possible errors one may suffer when not doing so. \\

During the middle of the 19th century, in connection with integration, Cauchy formulated and used a theorem which
stated that a sequence of continuous functions on a bounded and closed real interval will converge uniformly to a
continuous function, if it converges pointwise. As it happens, if we consider this statement within the usual reals ${\bf
R}$, that is, {\it without} infinitesimals, then simple examples show it to be false. \\
However, in a letter to Abel written in 1853, Cauchy speaks about what in nonstandard terms amount to points in $^*{\bf
R}$, and shows that in such a context the statement is indeed true, see ?? \\ \\

{\bf 5. We Opt for the Robinson Approach} \\ \\

In this book we shall follow the approach to Nonstandard Analysis due to A Robinson, and in its simplified version
based on Superstructures. The main reason for that is in the clear constructive ways which are typical for that
approach. \\
Once, however, one got familiar and at ease with such an approach, it is not so difficult to go over to the approach of
Internal Set Theory, originated by E Nelson. \\

The converse process, namely, starting first with Internal Set Theory, and then going over to Superstructures is also
possible, of course. It may appear, however, that a first familiarization with Nonstandard Analysis may be more easy
with the constructive approach, that is, of Superstructures. \\ \\

{\bf 6. A Few Historical Notes} \\ \\

The origin of modern Nonstandard Analysis is in the 1961 paper of A Robinson, which was soon followed up in
Luxemburg [1]. \\

Not much earlier, in the late 1950s, in papers by D Laugwitz and C Schmieden, extensions of ${\bf R}$ were
constructed, based on a certain generalization of the well known construction of the real numbers by Cantor. These
extensions, however, did not lead to fields, since they gave algebras with zero divisors, and as such, they proved to
be somewhat cumbersome in applications. \\
The earliest traceable modern publication with implications for a nonstandard theory appears to be a paper of T A
Skolem, published in 1934. \\

A variety of related, well documented, instructive, and in fact, rather comprehensive details on what may be seen as
attempts prior to the 1960s at the creation of a nonstandard theory can be found in Robinson [2], at pages 260-282. Also
of interest in this regard are the Preface and Introduction to Stroyan \& Luxemburg, and pages 3,4 in Albeverio, et.al. \\

The first major and systematic presentation of Nonstandard Analysis, including of a number of applications of interest,
was in the 1966 book of A Robinson. A difficulty of this book for usual mathematicians was the systematic use of
certain elements of Mathematical Logic, and in particular, Model Theory, among others, of the theory of types, the
compactness property of first order predicate calculus, and so on. \\

Soon after, in their 1969 paper, A Robinson and E Zakon introduced Superstructures and their formal languages,
managing in this way to reduce Nonstandard Analysis to set theory, and with only a minimal and rather elementary
involvement of Mathematical Logic. \\

The next major departure was with the 1977 paper of E Nelson which created a parallel and equivalent avenue for
Nonstandard Analysis, through what is called Internal Set Theory. \\

And yet, if one happens to be familiar with the amount of Mathematical Logic - an amount by no means significant -
which is used in the 1966 book of A Robinson, one can easily conclude that that book is by far the best, shortest and
clearest introduction to Nonstandard Analysis. \\
In fact, given the nature of the subject, a nature which is novel within usual mathematics, one may as well conclude
that, simply, there cannot be a better way of introduction. \\

In this way, one may as well say that there are in fact {\it three} approaches to Nonstandard Analysis : the presently
used two approaches, namely, by Superstructures, respectively, by Internal Set Theory, and the third one - which was
originally the first, and by now it is neglected, although it is {\it naturally} by far the best - in the 1966 book of A
Robinson. \\ \\

\chapter{Constructing $^*{\bf R}$}

This and chapter 3, next, present a {\it first stage} into Nonstandard Analysis. And as such, the aim in these two
chapters is to keep as much to usual mathematics as possible, and only make use of a rock bottom minimal amount of
Mathematical Logic, and in particular, of formal languages. \\
In this respect, the present chapter only employs usual mathematics, and the mentioned first encounter with formal
languages, as they are defined in Mathematical Logic, will happen in chapter 3. \\

One of the features of this and the next chapter is that, in spite of the maximum simplicity they pursue, they
nevertheless give an idea about the outlines of Nonstandard Analysis, as well as about some of its main features, such
as for instance, {\it transfer}. Needless to say, for a better understanding of the subject, {\it stage two}, which starts with
chapter 5, is needed, and in addition to the present book, one should also study some of the excellent and more
voluminous treaties mentioned in the reference. \\
In chapter 4, a short review of the main ideas and results in chapters 2 and 3 is presented as a stepping stone to the
general approach to Nonstandard Analysis introduced in chapter 5. \\

Let us start here by pointing out the {\it two main features} which underlie all the subsequent constructions, and
which should always be kept in mind. Here we present these two features in the particular case of the relationship
between ${\bf R}$ and $^*{\bf R}$, however, they are the same as well in the general case introduced in chapter 5.
Namely : \\

(TE)~:~~ the {\it transfer} of various {\it mathematical entities} from ${\bf R}$ to $^*{\bf R}$, \\
         \hspace*{1.45cm} the first and simplest case of it being the field homomorphism \\
         \hspace*{1.45cm} $^*(~)$ in (2.1.4) below, \\

(TS)~:~~ the {\it transfer} of {\it sentences} $\Phi$ from a given {\it language} ${\cal L}_{\bf R}$ about \\
         \hspace*{1.45cm} mathematical entities in ${\bf R}$, to corresponding sentences $^*\Phi$ in \\
         \hspace*{1.45cm} a language ${\cal L}_{^*{\bf R}}$ about mathematical entities in $^*{\bf R}$. \\

Furthermore, the transfer of sentences is a {\it two way} process. Namely, when we want to prove a standard property of
${\bf R}$, we formulate it as a sentence $\Phi$ in the language ${\cal L}_{\bf R}$, then we transfer it to the sentence
$^*\Phi$ in the language ${\cal L}_{^*{\bf R}}$, and prove this transferred sentence $^*\Phi$ in $^*{\bf R}$, after
which we deduce by reverse transfer that the initial sentence $\Phi$ is true about ${\bf R}$. \\

In this chapter we shall step by step and by employing only ususal mathematics implement the above stage (TE) of
transfer. In chapter 3, we shall present a first and simple way to implement stage (TS) of transfer as well. This however
will have to involve certain elements of formal languages in Mathematical Logic, elements which will be kept to a
minimum. \\ \\

{\bf 1. $^*{\bf R}$ as a Field which Contains ${\bf R}$} \\ \\

It is possible in just {\it two} simple steps to obtain $^*{\bf R}$ as a totally ordered field which contains ${\bf R}$ strictly,
and also, as a totally ordered subfield. The first step is to embed ${\bf R}$ into the vastly {\it larger} algebra
${\bf R}^{\bf N}$, as done in (2.1.1). Then in the second step, this much larger algebra ${\bf R}^{\bf N}$ is {\it reduced} by
a quotient construction with the use of certain kind of maximal ideals, as in (2.1.3). \\

Let us note here that the general form of such two step constructions,  a form which reaches far beyond such algebraic
realms, is called {\it reduced powers}, or even more generally, {\it reduced products}, and it is one of the basic tools in
Model Theory, which is a modern branch of Mathematical Logic, see Bell \& Slomson, Marker, or Jech. \\
However, we shall only use a particular form of reduced powers, namely, {\it ultrapowers}, and in doing so, in this
chapter we shall be able to keep to usual mathematics. \\

Let us now proceed to the first step and construct the {\it algebra homomorphism} given by the {\it injective} mapping

\bigskip
(2.1.1) \quad $ {\bf R} \ni r ~~\longmapsto~~ < r, r, r, ~.~.~.~ >~ \in {\bf R}^{\bf N} $

\medskip
where we observe that ${\bf R}^{\bf N}$ is indeed a commutative algebra over ${\bf R}$, when considered with the
termwise operations on its elements which are all the sequences $s ~=~ < s_1, s_2, s_3, ~.~.~.~ >$ of real numbers
$s_i \in {\bf R}$, with $i \in {\bf N}$. \\
However, ${\bf R}^{\bf N}$ is obviously not a field, since it has zero divisors, for instance

\bigskip
(2.1.2) \quad $ \begin{array}{l}
                                    < 1, 0, 1, 0, 1, 0, ~.~.~.~ > . < 0, 1, 0, 1, 0, 1, ~.~.~.~ > ~=~  \\ \\
                                        ~~~~~~~~~~~~~~~~ = ~< 0, 0, 0, 0, 0, 0, ~.~.~.~ >
                      \end{array} $

\medskip
It follows that ${\bf R}^{\bf N}$ {\it cannot} further be embedded into a field. Therefore, we can only proceed as follows in
our second step. We take a {\it maximal ideal} ${\cal M}$ in ${\bf R}^{\bf N}$ and obtain the {\it canonical quotient
algebra homomorphism} given by the {\it surjective} mapping

\bigskip
(2.1.3) \quad $ {\bf R}^{\bf N} \ni s ~~\longmapsto~~ [ s ] ~=~ s + {\cal M} \in {\bf R}^{\bf N} ~/~ {\cal M} $

\medskip
Here we note that ${\bf R}^{\bf N} ~/~ {\cal M}$ is a {\it field}, since ${\cal M}$ is supposed to be a maximal ideal, which
as is well known, gives the necessary and sufficient condition to obtain the respective quotient as a field. \\

And now, bringing the above together, we obtain the commutative diagram of algebra homomorphisms

\begin{math}
\setlength{\unitlength}{1cm}
\thicklines
\begin{picture}(15,7)
\put(2,5){${\bf R}$}
\put(2.7,5.15){\vector(1,0){6}}
\put(5.5,5.4){$^*(~)$}
\put(9,5){$^*{\bf R} ~=~ {\bf R}^{\bf N} ~/~ {\cal M}$}
\put(6,1.5){${\bf R}^{\bf N}$}
\put(2.7,4.8){\vector(1,-1){3}}
\put(6.8,2){\vector(1,1){2.7}}
\put(0,3){$(2.1.4)$}
\end{picture}
\end{math}

where the {\it injective field homomorphism}, thus {\it field extension} $^*(~)$ is given by

\bigskip
(2.1.5) \quad $ \begin{array}{l}
                                    {\bf R} \ni r ~~\longmapsto~~ ^* r ~=~ [~ < r, r, r, ~.~.~.~ > ~] ~=~ \\ \\
                                       ~~~~~~~~~~~~~ =~ < r, r, r, ~.~.~.~ > + {\cal M} \in ^*{\bf R} ~=~ {\bf R}^{\bf N} ~/~ {\cal M}
                     \end{array} $

\medskip
between the two fields ${\bf R}$ and $~^*{\bf R} ~=~ {\bf R}^{\bf N} ~/~ {\cal M}$. \\ \\

{\bf Note 2.1.1} \\

From (2.1.4) it is obvious that the cardinal of $~^*{\bf R}$ {\it cannot} be larger than that of the {\it continuum}. Indeed, the
cardinal of ${\bf R}^{\bf N}$ is the cardinal of the continuum, therefore, when it is subjected to a quotient by the ideal
${\cal M}$, in order to give $~^*{\bf R}$, the cardinal cannot increase.

\hfill $\Box$ \\ \\

As for the {\it total order} on $~^*{\bf R} ~=~ {\bf R}^{\bf N} ~/~ {\cal M}$, let us start by noting that the algebra ${\bf R}^{\bf
N}$ inherits form ${\bf R}$ the following termwise order relation, which however, is only a partial order. Given $s =
< s_1, s_2, s_3, ~.~.~.~ >$ and $t = < t_1, t_2, t_3, ~.~.~.~ >$ in ${\bf R}^{\bf N}$, we define

\bigskip
(2.1.6) \quad $ s ~\leq~ t ~~~\Longleftrightarrow~~ s_n ~\leq~ t_n,~~~\mbox{with}~ n \in {\bf N} $

\medskip
and this partial order is compatible with the algebra structure of ${\bf R}^{\bf N}$. \\

It follows that there are four questions to clarify : \\

1. How to construct maximal ideals ${\cal M}$ in ${\bf R}^{\bf N}$ ? \\

2. When is the algebra homomorphism $^*(~)$ in (2.1.4) and (2.1.5) injective ? \\

3. When is the field $~^*{\bf R} ~=~ {\bf R}^{\bf N} ~/~ {\cal M}$ strictly larger than ${\bf R}$ ? \\

4. Which is the total order on $~^*{\bf R} ~=~ {\bf R}^{\bf N} ~/~ {\cal M}$ ? \\

In order to answer these four questions it will be useful to make a short detour into filters, ultrafilters and binary valued
non-atomic measures. Indeed, filters on ${\bf N}$ are closely related to ideals in ${\bf R}^{\bf N}$, while ultrafilters and
binary valued non-atomic measures on ${\bf N}$ are closely related to maximal ideals ${\bf R}^{\bf N}$. \\ \\

{\bf 2. Ideals and Filters} \\ \\

Since some of the following results hold not only for ${\bf N}$, but also for any other infinite set, we shall present them
in the general case, as such a presentation does not involve additional complications. Furthermore, this general
apporach will be useful in chapter 5, when studying Superstructures. \\
Let us therefore take any infinite set $I$ and consider the corresponding algebra ${\bf R}^I$, which obviously can be
identified with the algebra on ${\bf R}$

$$ {\cal A} ( I ) ~= ~ \{~ s ~|~ s : I ~\longrightarrow~ {\bf R} ~\} ~=~ {\bf R}^I $$

\medskip
of all the real valued functions on $I$. In order to understand the structure of the ideals in ${\cal A} ( I )$, which therefore
means, of those in ${\bf R}^I$ as well, we proced as follows. \\

With each such function $s \in {\cal A} ( I )$, let us associate its {\it zero set} given by $Z ( s ) ~=~ \{~ i \in I ~|~ s ( i ) ~=~ 0
~\}$, which is a subset of $I$. \\

Further, it is useful to introduce the following concept. A family ${\cal F}$ of subsets of $I$, that is, a subset ${\cal F}
\subseteq {\cal P} ( I )$, is called a {\it filter} on $I$, if and only if it satisfies the following three conditions

\bigskip
(2.2.1) \quad $ \begin{array}{l} 1.~~ \phi \notin {\cal F} \neq \phi \\ \\
                                                  2.~~ J,~ K \in {\cal F} ~~\Longrightarrow~~ J \bigcap K \in {\cal F} \\ \\
                                                  3.~~ I \supseteq K \supseteq J \in {\cal F} ~~\Longrightarrow~~ K \in {\cal F}
                         \end{array} $

\medskip
Given an ideal ${\cal I}$ in ${\bf R}^I$, let us associate with it the set of its zero sets, namely

\bigskip
(2.2.2) \quad $ {\cal F}_{\cal I} ~=~ \{~ Z ( s ) ~~|~~ s \in {\cal I} ~\} ~\subseteq~ {\cal P} ( I ) $

\bigskip
Then

\bigskip
(2.2.3) \quad $ {\cal F}_{\cal I} ~~\mbox{is a filter on}~ I $

\medskip
Indeed, ${\cal I} \neq \phi$, thus ${\cal F}_{\cal I} \neq \phi$. Further, assume that $Z ( s ) = \phi$, for a certain
$s \in {\cal I}$. Then $s ( i ) \neq 0$, for $i \in I$. Therefore we can define $t : I ~\longrightarrow~ {\bf R}$, by
$t ( i ) = 1 / s ( i )$, with $i \in I$. Then however $t . s = 1 \in {\bf R}^I$, hence ${\cal I}$ cannot be an ideal in
${\bf R}^I$, which contradicts the hypothesis. In this way condition 1 in (2.2.1) holds for ${\cal F}_{\cal I}$. \\
Let now $s,~ t \in {\cal I}$, then clearly $s^2 + t^2 \in {\cal I}$, and $Z ( s^2 + t^2 ) = Z ( s ) \bigcap Z ( t )$,
thus condition 2 in  (2.2.1) is also satisfied by ${\cal F}_{\cal I}$. \\
Finally, let $s \in {\cal I}$ and $K \subseteq I$, such that $K \supseteq Z ( s )$. Let $t$ be the characteristic
function of $I \setminus K$. Then $t . s \in {\cal I}$, since ${\cal I}$ is an ideal. Now obviously $Z ( t . s ) = K$,
which shows that ${\cal F}_{\cal I}$ satisfies as well condition 3 in (2.2.1). \\

There is also the {\it converse} construction. Namely, let ${\cal F}$ be any filter on $I$, and let us associate with it the
set of functions

\bigskip
(2.2.4) \quad $ {\cal I}_{\cal F} ~=~ \{~ s : I ~\longrightarrow~ {\bf R} ~~|~~ Z ( s ) \in {\cal F} ~\} ~\subseteq~ {\bf R}^I $

\medskip
Then

\bigskip
(2.2.5) \quad $ {\cal I}_{\cal F} ~~\mbox{is an ideal in}~ {\bf R}^I $

\medskip
Indeed, we have $Z ( s + t ) \supseteq Z ( s ) \bigcap Z ( t )$, thus $s,~ t \in {\cal I}_{\cal F}$ implies that $s + t
\in {\cal I}_{\cal F}$. Also $Z ( s . t ) \supseteq Z ( s )$, therefore $s \in {\cal I}_{\cal F},~ t \in {\bf R}^I$
implies that $s . t \in {\cal I}_{\cal F}$. Further we note that $Z ( \lambda s ) = Z ( s )$, for $\lambda \in {\bf R},~
\lambda\neq 0$. Finally, it is clear that ${\cal I}_{\cal F} \neq {\bf R}^I$, since $s \in {\cal I}_{\cal F}
~\Longrightarrow~ Z ( s ) \neq \phi$, as ${\cal F}$ satisfies condition 1 in (2.2.1).Therefore (2.2.5) does indeed
hold. \\

Let now ${\cal I},~ {\cal J}$ be two ideals in ${\bf R}^I$, while ${\cal F},~ {\cal G}$ are two filters on $I$. Then it
is easy to see that

\bigskip
(2.2.6) \quad $ \begin{array}{l}
                                       {\cal I} ~\subseteq~ {\cal J} ~~~\Longrightarrow~~ {\cal F}_{\cal I} ~\subseteq~ {\cal F}_{\cal J} \\ \\
                                       {\cal F} ~\subseteq~ {\cal G} ~~~\Longrightarrow~~ {\cal I}_{\cal F} ~\subseteq~ {\cal I}_{\cal G}
                         \end{array} $

\medskip
We can also note that, given an ideal ${\cal I}$ in ${\bf R}^I$ and a filter ${\cal F}$ on $I$, we have by iterating the
above constructions in (2.2.2) and (2.2.4)

\bigskip
(2.2.7) \quad $ \begin{array}{l}
                              {\cal I} ~~~\longrightarrow~~~ {\cal F}_{\cal I} ~~~\longrightarrow~~~ {\cal I}_{{\cal F}_{\cal I}} ~=~ {\cal I} \\ \\
                              {\cal F} ~~~\longrightarrow~~~ {\cal I}_{\cal F} ~~~\longrightarrow~~~ {\cal F}_{{\cal I}_{\cal F}} ~=~ {\cal F}
                        \end{array} $

\medskip
Indeed, in view of (2.2.2), (2.2.4), we have for $s \in {\bf R}^I$ the equivalent conditions

$$ s \in {\cal I} ~~~\Longleftrightarrow~~ Z ( s ) \in {\cal F}_{\cal I} ~~~\Longleftrightarrow~~~ s \in {\cal I}_{{\cal F}_{\cal I}} $$

\medskip
Further, for $J \subseteq I$, we have the equivalent conditions

$$ J \in {\cal F}_{{\cal I}_{\cal F}} ~~~\Longleftrightarrow~~~ J ~=~ Z ( s ), ~~\mbox{for some}~ s \in {\cal I}_{\cal F} $$

\medskip
But for $s \in {\bf R}^I$, we also have the equivalent conditions

$$ s \in {\cal I}_{\cal F} ~~~\Longleftrightarrow~~~ Z ( s ) \in {\cal F} $$

\medskip
and the proof of (2.2.7) is completed. \\

In view of (2.2.7), it follows that every ideal in ${\bf R}^I$ is of the form ${\cal I}_{\cal F}$, where ${\cal F}$ is a certain
filter on $I$. Also, every filter on $I$ is of the form ${\cal F}_{\cal I}$, where ${\cal I}$ is a certain ideal in ${\bf R}^I$. \\ \\

{\bf 3. Maximal Ideals and Ultrafilters} \\ \\

Now  we can obtain the answer to Question 1 at the end of section 1, as follows. \\

First we note in general that, given a filter ${\cal F}$ on $I$, then

\bigskip
(2.3.1) \quad $ {\cal I}_{\cal F} ~~\mbox{is a maximal ideal in}~ {\bf R}^I ~~~\Longleftrightarrow~~~
                                                                                                  {\cal F} ~~\mbox{is an ultrafilter on}~ I $

\medskip
In other words, all {\it maximal ideals} ${\cal M}$ in ${\bf R}^I$ are of the form

\bigskip
(2.3.2) \quad $ {\cal M} ~=~ {\cal I~}_{\cal U} $

\medskip
where ${\cal U}$ are {\it ultrafilters} on $I$. \\

Here we recall that a filter ${\cal U}$ on $I$ is an {\it ultrafilter}, if and only if for every filter ${\cal F}$ on $I$, we have

\bigskip
(2.3.3) \quad $ {\cal U} ~\subseteq~ {\cal F} ~~~\Longleftrightarrow~~~ {\cal U} ~=~ {\cal F} $

\medskip
It follows that the equivalence in (2.3.1) is a direct consequence of (2.2.2),  (2.2.4),  (2.2.6) and (2.3.3). \\

Relation (2.3.2), in the particular case when $I ~=~ {\bf N}$, gives the answer to Question 1 at the end of section 1. \\

Coming now to Question 3 at the end of section 1, we have to distinguish between {\it two} kind of ultrafilters. Namely,
we call an ultrafilter ${\cal U}$ on $I$ to be {\it fixed}, if and only if for a certain given $i \in I$, we have $\{~ i ~\} \in
{\cal U}$, which is equivalent with

\bigskip
(2.3.4) \quad $ {\cal U} ~=~ \{~ J ~\subseteq~ I ~|~ i \in J ~\} $

\medskip
All other ultrafilters on $I$ are called {\it free}. \\

Here it should be noted that the existence of free ultrafilters - which are precisely those we shall need in the
construction of $^*{\bf R}$ - is equivalent with the Axiom of Choice, see Appendix 1. \\
Since ultrafilters play a fundamental role in proving the {\it existence} of nonstandard extensions of ${\bf R}$, or as will
be seen in chapter 5, of nonstandard extensions of general mathematical structures, we recall in section 9 at the end of
this chapter a well known result showing that there are plenty of such ultrafilters on every infinite set. \\

Let us see what happens if in (2.1.4) and (2.1.5) we use a maximal ideal ${\cal M}$ generated by a fixed ultrafilter
${\cal U}$ on ${\bf N}$. Then in view of (2.3.4), we have for a certain $n \in {\bf N}$ the relation ${\cal U} ~=~ \{~ L
~\subseteq~ {\bf N} ~~|~~ n \in L ~\}$. Thus  (2.2.4) results in

$$ {\cal M} ~=~ \{~ s ~=~ < s_1, s_2, s_3, ~.~.~.~ >~ \in {\bf R}^{\bf N} ~~|~~ s_n ~=~ 0 ~\} $$

\medskip
hence

$$ {\bf R}^{\bf N} ~/~ {\cal M} ~=~ {\bf R} $$

\medskip
and we failed to obtain in (2.1.4) and (2.1.5) a quotient field ${\bf R}^{\bf N} ~/~ {\cal M}$ which is larger than ${\bf R}$. \\

It follows that our only chance left to obtain in (2.1.4) and (2.1.5) a quotient field ${\bf R}^{\bf N}~/~ {\cal M}$ which is
larger than ${\bf R}$ is to use a {\it free} ultrafilter ${\cal U}$ on ${\bf N}$, when we generate the corresponding maximal
ideal ${\cal M}$. And as we can see next, this will always work. \\

First, we need an important and characteristic {\it dichotomy} type property of ultrafilters which will be proved in
Appendix 1. \\

A filter ${\cal F}$ on $I$ is an ultrafilter, if and only if for every $J \subseteq I$, we have

\bigskip
(2.3.5) \quad $ \mbox{either}~~J \in {\cal U}, ~~~\mbox{or}~~ I \setminus J \in {\cal U} $ \\ \\

{\bf Note 2.3.1} \\

Let us note the following about the dichotomy type property (2.3.5) of ultrafilters, which plays a fundamental role in
Nonstandard Analysis, and in general, in ultrapowers used in Model Theory. \\

Any ultrafilter ${\cal U}$ on $I$ is a set of subsets of $I$, that is, ${\cal U} \subset {\cal P} ( I )$. Thus given an
arbitrary subset $J \subseteq I$, there are {\it three} mutually exclusive logical possibilities, namely \\

1. $J \in {\cal U}$, \\

2. $I \setminus J \in {\cal U}$, \\

3. $J,~ I \setminus I \notin {\cal U}$ \\

since both $J$ and $I \setminus J$ cannot simultaneously be in ${\cal U}$, which satisfies the filter conditions 1
and 2 in (2.2.1). \\
Now what property (2.3.5) says is that for ultrafilters, the above alternative 3 is {\it not} possible either, therefore, we
always have only {\it two mutually exclusive logical possibilities}, namely, either 1, or 2 above. \\

As we shall see, this {\it elimination} of the third possibility 3 in the case of ultrafilters has fundamental consequences
regarding the properties of $^*{\bf R}$. \\
One such example is in section 6, where the total order on $^*{\bf R}$ is defined. \\

In fact, the dichotomy property (2.3.5) leads to the following alternative characterization of ultrafilters ${\cal U}$ on $I$

\bigskip
(2.3.6) \quad $ \begin{array}{l}
                             \forall~~ J_1, ~.~.~.~ , J_n ~\subseteq~ I ~: \\ \\
                             ~~~~J_1 \bigcup ~.~.~.~ \bigcup J_n \in {\cal U} ~~~\Longrightarrow~~ J_i \in {\cal U},~~~\mbox{for some}~~
                                              1 \leq i \leq n
                      \end{array} $

\medskip
Indeed, assume that $J_i \notin {\cal U}$, with $1 \leq i \leq n$. Then according to (2.3.5), we have $I \setminus J_i \in
{\cal U}$, for $1 \leq i \leq n$, hence $( I \setminus J_1 ) \bigcap ~.~.~.~ ( I \setminus J_n ) \in {\cal U}$. But this means
that $I \setminus ( J_1 \bigcup ~.~.~.~ \bigcup J_n ) \in {\cal U}$, and in view of (2.3.5), this contradicts the hypothesis.

\hfill $\Box$ \\ \\

Now we can return and show that ${\bf R}^{\bf N} ~/~ {\cal M}$ is indeed larger than ${\bf R}$, whenever ${\cal M}$ is
obtained through (2.2.4) from a free ultrafilter ${\cal U}$ on ${\bf N}$. \\
Let us take the sequence $~\omega ~=~ < 1, 2, 3, ~.~.~.~ >~ \in {\bf R}^{\bf N}$, then $[ \omega ] ~=~ \omega + {\cal M} \in
{\bf R}^{\bf N} ~/~ {\cal M}$. Now we show that $[ \omega ] \neq~ ^*r$, for every $r \in {\bf R}$. Indeed, let us assume that
for a certain $r \in {\bf R}$ we have $[ \omega ] ~=~ ^*r$, then $~\omega ~-~ < r, r, r, ~.~.~.~ >~ \in {\cal M}$, therefore

$$ \{~ n \in {\bf N} ~~|~~ n ~=~ r ~\} \in {\cal U} $$

\medskip
However, this set is either void, or it contains only one element. And both of these cases contradict the assumption that
${\cal U}$ is a free ultrafilter on ${\bf N}$. \\

Finally, we also obtain the answer to Question 2 at the end of section 1. Indeed, the algebra homomorphism $^*(~)$ in
(2.1.4) and (2.1.5) is always injective, when ${\cal M}$ is an ideal in ${\bf R}^{\bf N}$, even when it is not a maximal
ideal. To show that, let us assume that for a certain $r \in {\bf R},~ r \neq 0$, we have $^*r ~=~ 0 \in {\bf R}^{\bf N} ~/~
{\cal M}$. Then obviously $< r, r, r, ~.~.~.~ >~ \in {\cal M}$. This however contradicts the assumption that ${\cal M}$ is an
ideal in ${\bf R}^{\bf N}$, since $r \neq 0$. \\

We can now conclude that the construction in (2.1.4) and (2.1.5) does indeed lead to a {\it larger} field

\bigskip
(2.3.7) \quad $ ^*{\bf R} ~=~ {\bf R}^{\bf N} ~/~ {\cal M} $

\medskip
which has ${\bf R}$ as a {\it strict subfield}, if and only if ${\cal M}$ is the ideal corresponding through (2.2.4) to a {\it free}
ultrafilter ${\cal U}$ on ${\bf N}$, in which case ${\cal M}$ is a {\it maximal} ideal. \\ \\

{\bf Convention 2.3.1} \\

Form now on, unless specified otherwise, we shall assume that a free ultrafilter ${\cal U}$ on ${\bf N}$ was chosen,
together with the corresponding maximal ideal ${\cal M}$. \\ \\

{\bf Note 2.3.2} \\

If we consider a compact topology on the infinite set $I$, and instead of the algebra ${\bf R}^I$ of all real valued
functions on $I$, we only consider its subalgebra ${\cal C} ( I )$ of all continuous real valued functions on $I$, then the
{\it maximal} ideals in ${\cal C} ( I )$ are given by

$$ {\cal M}_i ~=~ \{~ s \in {\cal C} ( I ) ~~|~~ s ( i ) ~=~ 0 ~\},~~~\mbox{with}~~ i \in I $$

\medskip
therefore, contrary to what happens in the algebra ${\bf R}^I$, they correspond to the {\it fixed} ultrafilters on $I$,
namely

$$ {\cal M}_i ~=~ {\cal I~}_{{\cal U~}_i} ~\bigcap~ {\cal C} ( I ), ~~~\mbox{with}~~ i \in I $$

\medskip
where ${\cal U~}_i ~=~ \{~ J \subseteq I ~~|~~ i \in J ~\} $. This result is due in essence to Gelfand and Kolmogorov,
and was established in the 1930s, see Gillman \& Jerison, or Walker. \\ \\

{\bf Note 2.3.3} \\

It is obvious that on a finite set $I$ all ultrafilters are fixed. Therefore, as seen above, in such a case we would always
obtain

$$ {\bf R}^I ~/~ {\cal M} ~=~ {\bf R} $$

\medskip
for maximal ideals ${\cal M}$ in ${\bf R}^I$. It follows that in order to obtain ${\bf R}^I ~/~ {\cal M}$ as a genuine
extension of ${\bf R}$, even when ${\bf R}$ is considered merely as a set, that is, with no any other structure, we must
choose an infinite set $I$. \\
In this chapter, and well as in chpater 3, we choose the smallest such infinite set, namely $ I ~=~ {\bf N}$. \\ \\

{\bf 4. Ultrafilters and Binary Valued Non-Atomic Measures} \\ \\

Since we shall be dealing a lot with the quotient structure $^*{\bf R} ~=~ {\bf R}^{\bf N} ~/~ {\cal M}$, it is useful
to introduce a {\it measure} on ${\bf N}$, in terms of which the equivalence relation on ${\bf R}^{\bf N}$ with respect to
the maximal ideal ${\cal M}$ can be expressed as an {\it almost everywhere} property on ${\bf N}$. \\

Again, such a construction is valid on every infinite set $I$, and thus we present it here in general, without introducing
by that any complications. Let ${\cal U}$ be a {\it free} ultrafilter on $I$. Then we associate with it a {\it finitely additive
binary valued} measure on $I$

\bigskip
(2.4.1) \quad $ \mu~_{\cal U} : {\cal P} ( I ) ~~\longrightarrow~~ \{~ 0, 1 ~\} $

\medskip
which is defined on all the subsets $J \subseteq I$, according to

\bigskip
(2.4.2) \quad $ \mu~_{\cal U}~ ( J ) ~=~ \Bigg | \begin{array}{l}
                                                                        ~~~ 1 ~~~\mbox{if}~~ J \in {\cal U} \\

                                                                        ~~~ 0 ~~~\mbox{if}~~ J \notin {\cal U}
                                                                      \end{array} $

\medskip
A meaning of such a measure is that sets $J \in {\cal U}$ are {\it large}, or {\it relevant}, while sets $J \notin {\cal U}$
are {\it small}, or {\it negligible}, see section 5, next, for a further elaboration of such an interpretation. \\

An important property of this measure $\mu~_{\cal U}$ is that it is {\it non-atomic}, namely

\bigskip
(2.4.3) \quad $ \mu~_{\cal U}~ ( i ) ~=~ 0, ~~~\mbox{with}~~ i \in I $

\medskip
and it is also {\it nontrivial}, since obviously

\bigskip
(2.4.4) \quad $ \mu~_{\cal U}~ ( I ) ~=~ 1 $

\medskip
Indeed, let us show the finite additivity of $\mu~_{\cal U}$. If $J,~ K \subseteq I,~ J \bigcap K = \phi$, then the
only case we have to consider is when $\mu~_{\cal U} ( I ) ~=~ \mu~_{\cal U} ( K ) ~=~ 1$. However, this means that
$J,~ K \in {\cal U}$, which is absurd since $J \bigcap K = \phi$ and ${\cal U}$ is a filter. The fact that
$\mu~_{\cal U}$ is non-atomic results from ${\cal U}$ being a free ultrafilter, for in case $\mu~_{\cal U}~ ( i )
\neq 0$ for a certain $i \in I$, then $\mu~_{\cal U}~ ( i ) ~=~ 1$, hence $\{~ i ~\} \in {\cal U}$, which is absurd. \\

We also have the {\it converse} construction. Namely, given a {\it finitely additive, non-atomic, nontrivial binary
valued} measure defined on all the subsets of an infinite set $I$

\bigskip
(2.4.5) \quad $ \mu : {\cal P} ( I ) ~~\longrightarrow~~ \{~ 0, 1 ~\} $

\medskip
then we can associate with it the {\it free} ultrafilter on $I$, given by

\bigskip
(2.4.6) \quad $ {\cal U~}_\mu ~=~ \{~ J \subseteq I ~~|~~ \mu ( J ) ~=~ 1 ~\} $

\medskip
Let us indeed check that ${\cal U~}_\mu$ is a free ultrafilter on $I$. Clearly $\phi \notin {\cal U~}_\mu$, since $\mu
( \phi ) = 0$. Also ${\cal U~}_\mu \neq \phi$, since $\mu ( I ) = 1$. Assume now $J,~ K \in {\cal U~}_\mu$ and $J
\bigcap K \notin {\cal U~}_\mu$. Then $\mu ( J \bigcap K ) = 0$. But $J \setminus K,~ J \bigcap K$ and $K \setminus J$
are pairwise disjoint, thus $\mu ( J \bigcup K ) = \mu ( J \setminus K ) + \mu ( J \bigcap K ) + \mu ( K \setminus J ) =
1 + 0 + 1 = 2$, which is absurd. In view of (2.3.5), ${\cal U~}_\mu$ is an ultrafilter, since $\mu$ can only have the
values 0 and 1. Since $\mu$ is non-atomic, it follows that ${\cal U~}_\mu$ is a free ultrafilter. \\

Now similar with (2.2.7), one can easily obtain the iterated constructions

\bigskip
(2.4.7) \quad $ \begin{array}{l}
                              {\cal U} ~~~\longrightarrow~~~ \mu~_{\cal U} ~~~\longrightarrow~~~ {\cal U~}_{\mu~_{\cal U}} ~=~ {\cal U} \\ \\
                              \mu ~~~\longrightarrow~~~ {\cal U~}_\mu ~~~\longrightarrow~~~ \mu~_{{\cal U~}_\mu} ~=~ \mu
                        \end{array} $

\medskip
Our interest in the above constructions is as follows. We choose a free ultrafilter ${\cal U}$ on ${\bf N}$, and
associate with it the maximal ideal ${\cal M}_{\cal U}$ in ${\bf R}^{\bf N}$ in order to obtain the quotient structure

$$ ^*{\bf R} ~=~ {\bf R}^{\bf N} ~/~ {\cal M}_{\cal U} $$

\medskip
giving the {\it nonstandard} field which extends the field ${\bf R}$ of real numbers. And when working with this
nonstandard extension $^*{\bf R}$, we shall have to deal with the quotient structure defined by ${\cal M}_{\cal U}$,
that is, with the {\it equivalence modulo} ${\cal M}_{\cal U}$. \\
Fortunately, we can easily translate this algebraic equivalence into the measure theoretic concept of {\it almost
everywhere}, by using the measure $\mu~_{\cal U}$, associated to the free ultrafilter ${\cal U}$. And the basic pattern
in this regard is as follows. \\

Let $s ~=~ < s_1, s_2, s_3, ~.~.~.~ >,~ t ~=~ < t_1, t_1, t_3, ~.~.~.~ >~ \in {\bf R}^{\bf N}$, then obviously

\bigskip
(2.4.8) \quad $ \begin{array}{l}
                         t - s \in {\cal M}_{\cal U} ~~~\Longleftrightarrow~~~ \{~ n \in {\bf N} ~|~ s_n = t_n ~\} \in {\cal U}
                     \end{array} $

\medskip
where upon convenience, first we interpret $s$ and $t$ as sequences in ${\bf R}^{\bf N}$, and then as functions $s,~ t :
{\bf N} ~\longrightarrow~ {\bf R}$. Now by definition, we shall write

\bigskip
(2.4.9) \quad $ t ~=~ s ~~~~\mu~_{\cal U}~~\mbox{almost everywhere on}~ {\bf N} $

\medskip
if and only if (2.4.8) holds. \\
Clearly, this sort of notation can be further extended to relations, functions, and so on. For instance, we can write

\bigskip
(2.4.10) \quad $ \begin{array}{l}
                             t ~\leq~ s ~~~~\mu~_{\cal U}~~\mbox{almost everywhere on}~ {\bf N} ~~~\Longleftrightarrow~~~ \\ \\
                              ~~~~~~~~~~~~\Longleftrightarrow~~~  \{~ n \in {\bf N} ~|~ s_n \leq t_n ~\} \in {\cal U}
                           \end{array} $ \\ \\

{\bf Convention 2.4.1} \\

As mentioned in Convention 2.3.1, we have assumed given an ultrafilter ${\cal U}$ on ${\bf N}$. And for brevity of
notation, we shall not longer specify it as an index when we use the associated maximal ideal ${\cal M}$, or the
associated measure $\mu$. \\
Also, instead of writing "$\mu~\mbox{almost everywhere on}~ {\bf N}$", we shall only write "$( a. e. )$", as customary in
measure theory. \\
When it is convenient, we shall in addition specify the integer variable, say $n \in {\bf N}$, with respect to which the
"$( a. e. )$" property holds, by writing "$( a. e. ) ~~\mbox{in}~ n \in {\bf N}$". \\ \\

{\bf 5. Filters and co-Filters} \\ \\

There is yet another way to approach the concept of a filter ${\cal F}$ on a set $I$. Namely, we can see the elements
$J \in {\cal F}$ of the filter, which therefore are subsets $J \subseteq I$, as being {\it large} in $I$. It follows then
obviously that $\phi \notin {\cal F}$, while on the contrary, $I \in {\cal F}$. Furthermore, if $J \in {\cal F}$ and $J
\subseteq K \subseteq I$, then also $K \in {\cal F}$. \\
What is the only less than trivially obvious condition we ask about {\it large} subsets of $I$ is that the intersection of
two large sets is still large, in other words, that $J, K \in {\cal F} \Longrightarrow J \bigcap K \in {\cal F}$. \\
And with these conditions we are back to the definition of a filter on $I$, given in (2.2.1). \\

In view of the above, it is natural to talk {\it dually} about sets ${\cal G}$ of {\it small} subsets $J \subseteq I$. And then
the conditions which they should satisfy are the following

\bigskip
(2.5.1) \quad $ \begin{array}{l}
                        1.~~ I \notin {\cal G} \neq \phi \\ \\
                        2.~~ J, K \in {\cal G} ~~~\Longrightarrow~~ J \bigcup K \in {\cal G} \\ \\
                        3.~~ J ~\subseteq~ K \in {\cal G} ~~~\Longrightarrow~~ J \in {\cal G}
                 \end{array} $

\medskip
We call {\it co-filter} every set ${\cal G}$ of subsets of $I$ which satisfies (2.5.1). \\

Now the connection between filters and co-filters, that is, between large and small subsets of $I$ is obvious. Namely,
we have the {\it duality}

\bigskip
(2.5.2) \quad $ {\cal F} ~~\mbox{filter on}~~ I ~~~\Longleftrightarrow~~~ {\cal G} ~=~ \{~ I \setminus J ~~|~~ J
                                                                                                               \in {\cal F} ~\} ~~\mbox{co-filter on}~~ I $

\medskip
Further, let us recall that a filter ${\cal U}$ on $I$ is an {\it ultrafilter}, if and only if for every $J \subseteq I$, we have
either $J \in {\cal U}$, or $I \setminus J \in {\cal U}$. \\
Thus it follows that, {\it dually}, a cofilter ${\cal V}$ on $I$ is {\it prime}, if and only if for every $J \subseteq I$, we
have either $J \in {\cal V}$, or $I \setminus J \in {\cal V}$. \\
Consequently, just like ultrafilters among filters, see (2.3.3), {\it prime} co-filters are also {\it maximal} among co-filters,
namely

\bigskip
(2.5.3) \quad $ {\cal V} ~\subseteq~ {\cal G} ~~~\Longrightarrow~~~ {\cal V} ~=~ {\cal G} $

\medskip
for every prime co-filter ${\cal V}$ and co-filter ${\cal G}$ on $I$. Furthermore, if ${\cal U}$ and ${\cal V}$ are a dual pair
of ultrafilter and prime co-filter, respectively, then

\bigskip
(2.5.4) \quad $ {\cal V} ~=~ {\cal P} ( I ) ~\setminus~ {\cal U},~~~ {\cal U} ~=~ {\cal P} ( I ) ~\setminus~ {\cal V} $

\medskip
In view of (2.5.4), the existence of ultrafilters, see Appendix 1, implies the existence of prime co-filters. \\

A simple and useful example of a pair of a dual filter and co-filter, in case of {\it infinite} sets $I$, can be obtained as
follows

\bigskip
(2.5.5) \quad $ \begin{array}{l}
                          {\cal P}_\infty ( I ) ~=~ \{~ J ~\subseteq~ I ~~|~~ I ~\setminus~ J ~~\mbox{finite} ~\} ~~\mbox{is a filter on}~~ I \\ \\
                          {\cal P}_F ( I ) ~=~ \{~ J ~\subseteq~ I ~~|~~ J ~~\mbox{finite} ~\} ~~\mbox{is a co-filter on}~~ I
                         \end{array} $

\medskip
and we note that ${\cal P}_\infty ( I )$ is called the {\it Fr\'{e}chet} filter on $I$. Obviously, an ultrafilter ${\cal U}$ on $I$
is {\it free}, if and only if ${\cal P}_\infty ( I ) \subseteq {\cal U}$. \\ \\

{\bf 6. Transferring the Total Order from ${\bf R}$ to $^*{\bf R}$} \\ \\

Now as a second step in accomplishing the stage (TE) of transfer, we answer Question 4 at the end of section 1, by
transferring from ${\bf R}$ to $^*{\bf R}$ the property of being a {\it totally ordered field}. And we do so again with the
help of usual mathematics, since we have not yet set up the general mechanism of transfer, due to be presented in
chapter 3, in its first and simpler form. We shall show here at the same time that ${\bf R}$ is embedded into $^*{\bf R}$
as a totally ordered subfield. \\

Let therefore $^*x,~ ^*y  \in~ ^*{\bf R} ~=~ {\bf R}^{\bf N} ~/~ {\cal M}$. Then $^*x = [ s ],~ ^*y = [ t ]$, for certain suitable
$s =~ < s_1, s_2, s_3, ~.~.~.~ >,~ t =~ < t_1, t_2, t_3, ~.~.~.~ >~ \in {\bf R}^{\bf N}$. \\
Now we define

\bigskip
(2.6.1) \quad $ \begin{array}{l}
                                     ^*x ~\leq~ ^*y ~~~\Longleftrightarrow~~~ \{~ n \in {\bf N} ~~|~~ s_n ~\leq~ t_n ~\} \in {\cal U}
                                      ~~~\Longleftrightarrow \\ \\
                                      ~~~~~~~~~~\Longleftrightarrow~~~ s ~\leq~ t~~~~ ( a. e. )
                         \end{array} $

\medskip
Here $\mu$ and ${\cal M}$ are the measure and the maximal ideal, respectively, which correspond to the free ultrafilter
${\cal U}$ on ${\bf N}$, which is considered given once and for all. And as mentioned, for the sake of brevity, we no
longer indicate this ultrafilter in the notation of $\mu$ and ${\cal M}$. \\

In view of (2.4.8), it is easy to see that the above definition (2.6.1) is correct, since it does not depend on the sequences
$s$ and $t$. Let us therefore show that the binary relation $\leq$ which it defines on  $^*{\bf R}$ is indeed a {\it total
order}. \\

First we note that $\leq$ is reflexive. Certainly, given $^*x  \in~ ^*{\bf R}$, with  $^*x = [ s ]$, for a certain  $s =~
< s_1, s_2, s_3, ~.~.~.~ >~ \in {\bf R}^{\bf N}$, we obviously have $ s = s~~~ ( a. e. )$. Next, $\leq$ is transitive. Indeed, if
$^*x,~ ^*y,~ ^*z  \in~ ^*{\bf R}$ and $^*x \leq ^*y ~\leq ^*z$, with $^*x = [ s ],~ ^*y = [ t ],~ ^*z = [ v ]$, for some $s =~
< s_1, s_2, s_3, ~.~.~.~ >,~ t =~ < t_1, t_2, t_3, ~.~.~.~ >~, v =~ < v_1, v_2, v_3, ~.~.~.~ >~ \in {\bf R}^{\bf N}$, then $ s = t~~~
( a. e. ),~ t = v~~~ ( a. e. )$, hence $s = v~~~ ( a. e. )$, which means that $^*x \leq ^*z$. Also, $\leq$ is antisymmetric, since
$s \leq t~~~ ( a. e. ),~ t \leq s~~~ ( a. e. )$ obviously imply $s = t~~~ ( a. e. )$. \\

Corresponding to the order $\leq$ in (2.6.1), we define as usual the strict order $<$ by

\bigskip
(2.6.2) \quad $ ^*x ~<~ ^* y ~~~\Longleftrightarrow~~~ ^*x ~\leq~ ^*y ~~\mbox{and}~~ ^* x ~\neq~ ^*y $

\medskip

Let us now show that $\leq$ is a total order. For that, let $^*x,~ ^*y  \in~ ^*{\bf R}$, and $^*x = [ s ],~ ^*y = [ t ]$, for some
$s =~ < s_1, s_2, s_3, ~.~.~.~ >,~ t =~ \\=~ < t_1, t_2, t_3, ~.~.~.~ >~ \in {\bf R}^{\bf N}$. Now, in case $^*x ~\leq~ ^*y$ does not
hold, then

$$ \{~ n \in {\bf N} ~~|~~ s_n \leq t_n ~\} \notin {\cal U} $$

\medskip
thus according to (2.3.5) we must have

$$ \{~ n \in {\bf N} ~~|~~ s_n > t_n ~\} \in {\cal U} $$

\medskip
which means that $^*x ~>~ ^*y$. \\

The fact that the total order $\leq$ on $^*{\bf R}$ satisfies the required compatibility conditions with the field structure
results in a similar manner. \\

Finally, let us show that ${\bf R}$ is embedded into $^*{\bf R}$ as a totally ordered subfield. Let $q,~ r \in {\bf R}$. Then
according to (2.1.5), we have $^*q ~=~ \\=~ [~ < q, q, q, ~.~.~.~ > ~],~ ^*r ~=~ [~ < r, r, r, ~.~.~.~ > ~] \in ^*{\bf R}$. Thus in view of
(2.6.1), it is obvious that

$$ q ~\leq~ r ~~\mbox{in}~~ {\bf R} ~~~\Longleftrightarrow~~~ ^*q ~\leq~ ^*r ~~\mbox{in}~~ ^*{\bf R} $$ \\ \\

{\bf 7. Transfer of Relations and Functions} \\ \\

So far we managed to transfer from ${\bf R}$ to $^*{\bf R}$, and do so one by one, the algebra and total order structures,
and in the process we only employed usual mathematics, and not the more general transfer mechanism (TS) which will
be presented in chapter 3, and which by necessity must employ certain elementary methods from Mathematical
Logic. \\
These transfers, which we did in the previous sections, belong to the stage (TE), and as such, concerned binary
functions corresponding to addition and multiplication, as well as a binary relation giving the total order, all them
defining the usual structure of ${\bf R}$. \\
However, there are many more functions and relations on ${\bf R}$ which are of mathematical interest, and which
should therefore be transferred to $^*{\bf R}$. Here in this section we show the way to do that in general within the
stage (TE) of transfer, and again, we shall only employ usual mathematics in the process, and not refer to the more
general and powerful transfer mechanism (TS) to be introduced in chapter 3, and then in its full power, in chapter 5. \\ \\

{\bf Definition 2.7.1} \\

Let $P$ be an n-ary relation on ${\bf R}$. Its {\it transfer} to an n-ary relation $^*P$ on $^*{\bf R}$ is given by the set of
all n-tuples $( s^1, ~.~.~.~ , s^n ) \in ( ^*{\bf R} )^n$ such that, if $s^i = [~ < s^i_1, s^i_2, s^i_3, ~.~.~.~ > ~]$, with $1 \leq i \leq n$,
then

\bigskip
(2.7.1) \quad $ P < s^1_j, ~.~.~.~ , s^n_j > ~~~~( a. e. )~~ \mbox{in}~ j \in {\bf N} $

\medskip
which. as seen easily, is equivalent with

\bigskip
(2.7.2) \quad $ \{~ j \in {\bf N} ~~|~~ P < s^1_j, ~.~.~.~ , s^n_j > ~\} \in {\cal U} $

\hfill $\Box$ \\ \\

Let us show that this is a correct definition, that is, it does not depend on the representations $s^i =
[~ < s^i_1, s^i_2, s^i_3, ~.~.~.~ > ~]$, with $1 \leq i \leq n$. Indeed, assume that we have another representation $s^i =
[~ < t^i_1, t^i_2, t^i_3, ~.~.~.~ > ~]$, with $1 \leq i \leq n$. Then (2.6.2) gives

$$ A ~=~ \{~ j \in {\bf N} ~~|~~ P < s^1_j, ~.~.~.~ , s^n_j > ~\} \in {\cal U} $$

\medskip
On the other hand, due to the assumed alternative equivalent representations, we also have for $1 \leq i \leq n$

$$ B^i ~=~ \{~ j \in {\bf N} ~~|~~  s^i_j ~=~ t^i_j ~\} \in {\cal U} $$

\medskip
therefore $A \bigcap B^1 \bigcap ~.~.~.~ \bigcap B^n \in {\cal U}$, which means that

$$ P < t^1_j, ~.~.~.~ , t^n_j > ~~~~( a. e. )~~ \mbox{in}~ j \in {\bf N} $$

\medskip
and the independence from representation is proved. \\

In order to transfer functions, we deal with them as particular cases of relations. Let therefore $f : D \subseteq {\bf R}^n
\longrightarrow {\bf R}$ be a function of n variables defined on the domain $D$. Then $f$ is equivalent with the
(n + 1)-ary relation $P$, given for $r_1, ~.~.~.~ , r_n, r_{n + 1} \in {\bf R}$, by

\bigskip
(2.7.3) \quad $ P < r_1, ~.~.~.~ , r_n, r_{n + 1} > ~~~\Longleftrightarrow~~~ f ( r_1, ~.~.~.~ , r_n ) ~=~ r_{n + 1} $

\medskip
Let now $( s^1, ~.~.~.~ , s^n, s^{n + 1} ) \in ( ^*{\bf R} )^{n + 1}$ be such that $s^i = \\=[~ < s^i_1, s^i_2, s^i_3, ~.~.~.~ > ~]$, with
$1 \leq i \leq n + 1$, then in view of (2.7.1) - (2.7.3), we have the equivalence

\bigskip
(2.7.4) \quad $ \begin{array}{l}
                                            ^*f ( s^1, ~.~.~.~ , s^n ) ~=~ s^{n + 1} ~~~\Longleftrightarrow~~ f ( s^1_j, ~.~.~.~ , s^n_j ) ~=~
                                                                                             s^{n + 1}_j \\ \\
                                                        ~~~~~~~~~~~~~~~~~~~~~~~~~~~~~~~~~~~~~~~~~~~~~~~~~~~~~~~~~( a. e. )~~ \mbox{in}~ j \in {\bf N}
                         \end{array} $ \\ \\

{\bf Convention 2.7.1} \\

For the sake of convenience, we shall introduce a few simplifying notational conventions. \\
Since the mapping $^*(~)$ in (2.1.5) is an injective field homomorphism which preserves the total order, it is convenient
to see ${\bf R}$ as a {\it subset} of $^*{\bf R}$. For this reason, we shall identify every $r \in {\bf R}$ with its
image $^*r \in~ ^*{\bf R}$, in other words, we shall assume that

\bigskip
(2.7.5) \quad $ {\bf R} \ni r ~=~ ^*r \in~ ^*{\bf R} $

\medskip
Further, the addition "$+$" and the multiplication "$.$" on ${\bf R}$ do transfer to corresponding binary functions
"$~^*+$" and "$~^*.$" on $^*{\bf R}$. However, for convenience, we shall omit the "$~^*~$", and write them as usual,
that is, "$+$" and "$.$", even when they operate on $^*{\bf R}$. The same we shall do with the binary relation of total
order "$\leq$", when transferred to $^*{\bf R}$.

\hfill $\Box$ \\ \\

Any subset $A \subseteq {\bf R}$ can be seen as a unary relation $P$ on ${\bf R}$, defined by $P < r >
~\Longleftrightarrow~ r \in A$, for $r \in {\bf R}$. Therefore, according to (2.7.1), its transfer is given by

\bigskip
(2.7.6) \quad $ \begin{array}{l}
                            ^*A ~=~ \{~ s = [~ < s_1, s_2, s_3, ~.~.~.~ > ~] \in~ ^*{\bf R} ~~|~~ \\ \\

                            ~~~~~~~~~~~~~~~~~~~~~~~~~~~~~~~~~~~~~~~s_j \in A ~~~~( a. e. )~~ \mbox{in}~ j \in {\bf N} ~\}
                \end{array} $

\medskip
On the other hand, in view of the above identification (2.7.5), we may often be interested in the subset $A \subseteq
{\bf R}$ as seen in $^*{\bf R}$, namely, as the image of the mapping $^*(~)$ in (2.1.5). For that purpose, sometime it may
be convenient to use the notation

\bigskip
(2.7.7) \quad $ A_* ~=~ \{~ ^* r ~=~ r ~~|~~ r \in A ~\} ~\subset~ ^*{\bf R} $

\medskip
although we shall often make the identification $A_* ~=~ A$. \\

It is easy to see that

\bigskip
(2.7.8) \quad $ \begin{array}{l}
                            1)~~~ A_* ~\subseteq~ ^*A \\ \\
                            2)~~~ A_* ~=~ ^*A ~~~\Longleftrightarrow~~ A ~~\mbox{is a finite subset of}~~ {\bf R}
                         \end{array} $

\medskip
Indeed, 1) is immediate. For 2), let us assume $A_* ~=~ ^*A$, which means that $^*A \subseteq A_*$. Then (2.7.6) gives
for every $s =\\= [~ < s_1, s_2, s_3, ~.~.~.~ > ~] \in~ ^*{\bf R}$, with $s_j \in A ~~~~( a. e. )~~ \mbox{in}~ j \in {\bf N}$, a certain
$r \in A$, such that $s_ j = r ~~~( a. e. )~~ \mbox{in}~ j \in {\bf N}$. Therefore, $A$ must indeed be finite, since ${\cal U}$ is
a free ultrafilter. \\
Conversely, let $A = \{~ a_1, ~.~.~.~ , a_n ~\}$ be finite, and let any $s =\\= [~ < s_1, s_2, s_3, ~.~.~.~ > ~] \in~ ^*{\bf R}$, with
$s_j \in A ~~~~( a. e. )~~ \mbox{in}~ j \in {\bf N}$. Then $\{~ j \in {\bf N} ~|~ s_j \in A ~\} \in {\cal U}$, hence \\$\{~ j \in {\bf N} ~|~
s_j = a_1 ~\} \bigcup ~.~.~.~ \bigcup \{~ j \in {\bf N} ~|~ s_j = a_n ~\} \in {\cal U}$, therefore (2.3.6) gives for a certain $1 \leq
i \leq n$, the relation $\{~ j \in {\bf N} ~|~ s_j = a_i ~\} \in {\cal U}$. But this obviously means that $s = a_i \in A_*$. \\

In order further to clarify (2.7.8), let us show how for any {\it infinite} set $A \subseteq {\bf R}$, we can always find
elements in $^*A$ which are {\it not} in $A_*$. \\

Indeed, since $A$ is infinite, we can take an infinite sequence of elements \\$a_1, a_2, a_3, ~.~.~.~ \in A$ which are {\it
pairwise different}. But according to (2.7.6), we obviously have $a ~=~ [~ < a_1, a_2, a_3, ~.~.~.~ > ~] \in~ ^*A$. If now we
assume that $a \in~ A_*$, then there exists $r \in A$, such that $a_j = r ~~~( a. e. )~~ \mbox{in}~ j \in {\bf N}$. However,
since the $a_j$-s are pairwise different, there can at most be one single $j \in {\bf N}$, such that $a_j = r$. And this
contradicts the assumption that the ultrafilter ${\cal U}$ is free. Thus we obtained

\bigskip
(2.7.9) \quad $ a ~=~ [~ < a_1, a_2, a_3, ~.~.~.~ > ~] \in~ ^*A \setminus A_* ~=~ ^*A \setminus A $

\medskip
As an immediate application, we obtain the existence of {\it nonstandard integers}, since (2.7.8), (2.7.9) imply

\bigskip
(2.7.10) \quad $ \omega ~=~ [~ < 1, 2, 3, ~.~.~.~ > ~] \in ~^*{\bf N} \setminus {\bf N}_* ~=~  ~^*{\bf N} \setminus {\bf N} $

\medskip
Another useful example occurs in the case of {\it bounded} intervals of real numbers $A = [ a, b ] \subset {\bf R}$. Then
(2.7.6) gives $^*A =~ ^*[ ^*a, ^*b ] = \{~ x \in~ ^*{\bf R} ~|~ ^*a \leq x \leq~ ^*b ~\} = [ a, b ]$, where the last interval is
considered in $^*{\bf R}$. Similar relations hold for half open or open bounded real intervals. \\

Here we should note that, since $- \infty,~ + \infty \notin {\bf R}$, the above transfer does not apply to them. However, if
for instance, we have the {\it unbounded} real interval $A = [ a, + \infty ) \subset {\bf R}$, then it can be written as $A =
\{~ x \in {\bf R} ~|~ a \leq x ~\}$, thus (2.7.6) will give $^*A = \{~ x \in~ ^*{\bf R} ~|~ ^*a \leq x ~\} = [ a, + \infty )$, where again,
the latter interval is considered in $^*{\bf R}$. Therefore, $+\infty$ in this latter interval means the positive end of
$^*{\bf R}$, which is far larger than the positive end of ${\bf R}$ would be, when ${\bf R}$ is considered as a subset of
$^*{\bf R}$. By the way, this so called positive end of ${\bf R}$, with ${\bf R}$ considered as a subset of $^*{\bf R}$,
does not even exist, since as noted earlier, $^*{\bf R}$ is not Dedekind order complete. \\

The identification in (2.7.5) has the important consequence that the transfer $^*P$ to $^*{\bf R}$ of any relation $P$ on
${\bf R}$ can be seen as an {\it extension} of $P$ to $^*{\bf R}$. And in particular, the transfer $^*f$ to $^*{\bf R}$ of any
function $f$ on a subset $A \subseteq {\bf R}$ can be seen as an {\it extension} of $f$ to $^*A \subseteq~ ^*{\bf R}$.
Indeed, let $P$ be an n-ary relation on ${\bf R}$ and $< r_1, ~.~.~.~ , r_n > \in {\bf R}^n$, then (2.7.5), (2.7.1) give

\bigskip
(2.7.11) \quad $ P  < r_1, ~.~.~.~ , r_n > ~~~\Longleftrightarrow~~~ ^*P  < r_1, ~.~.~.~ , r_n >  $

\medskip
In particular, if $f : D \subseteq {\bf R}^n \longrightarrow {\bf N}$ is a function of n variables defined on the domain $D$,
then for any  $< r_1, ~.~.~.~ , r_n, r_{n + 1} > \in {\bf R}^{n + 1}$, we have

\bigskip
(2.7.12) \quad $ f ( r_1, ~.~.~.~ , r_n ) ~=~ r_{n + 1} ~~~\Longleftrightarrow~~~ ^*f ( r_1, ~.~.~.~ , r_n ) ~=~ r_{n + 1} $

\medskip
Let us further note that any n-ary relation $P$ can be identified with its {\it characteristic} function $\chi_P$, defined by

\bigskip
(2.7.13) \quad $ \chi_P ~ ( \alpha_1, ~.~.~.~ , \alpha_n ) ~=~ \Bigg | \begin{array}{l} ~~ 1 ~~~\mbox{if}~~ P < \alpha_1, ~.~.~.~ , \alpha_n > \\
                                                                                         ~~ 0 ~~~\mbox{otherwise}
                                                                \end{array} $

\medskip
Then for any n-ary relation $P$ on ${\bf R}$ we have

\bigskip
(2.7.14) \quad $ ^* ( \chi_P ) ~=~ \chi_{(~^*P)} $

\medskip
Indeed, let us start with computing the left hand side in (2.7.14). According to (2.7.3), there exists an (n + 1)-ary
relation $Q \subseteq {\bf R}^{n + 1}$, such that for $r_1, ~.~.~.~ , r_n, r_{n + 1} \in {\bf R}$, we have

$$ \chi_P~ ( r_1, ~.~.~.~ , r_n ) ~=~ r_{n + 1} ~~~\Longleftrightarrow~~~
                                               Q < r_1, ~.~.~.~ , r_n, r_{n + 1} > $$

\medskip
It follows that $Q <  r_1, ~.~.~.~ , r_n, r_{n + 1} > $ is equivalent with

$$ (~ P < r_1, ~.~.~.~ , r_n > ~\mbox{and}~ r_{n + 1} = 1 ~) ~~~\mbox{or}~~~
   (~ P~^\prime < r_1, ~.~.~.~ , r_n > ~\mbox{and}~ r_{n + 1} = 0 ~) $$

\medskip
where $P~^\prime$ is the negation of $P$, namely, $P~^\prime = {\bf R}^n \setminus P$. \\
Now (2.7.4) gives for $( s^1, ~.~.~.~ , s^n, s^{n + 1} ) \in ( ^*{\bf R} )^{n + 1}$, where \\
$s^i = [~ < s^i_1, s^i_2, s^i_3, ~.~.~.~ > ~]$, with $1 \leq i \leq n + 1$, the relation

$$ ^*( \chi_P )~ ( s^1, ~.~.~.~ , s^n ) ~=~ s^{n + 1} ~~~\Longleftrightarrow~~~
                                          ^*Q < s^1, ~.~.~.~ , s^n, s^{n + 1} > $$

\medskip
On the other hand, the right hand side of (2.7.14) follows from (2.7.13) applied to $^*P$, namely

$$ \chi_{(~^*P)}~ ( s^1, ~.~.~.~ , s^n ) ~=~ 1 ~~~\Longleftrightarrow~~~
                                          ^*P < s^1, ~.~.~.~ , s^n > $$

And now if we compare the above computations of the two sides of (2.7.14), and take into account (2.7.1), the proof of
(2.7.14) is completed. \\

Similar with the above, one can easily establish the following relations. Given $A_1, ~.~.~.~ , A_n \subseteq {\bf R}$, then

\bigskip
(2.7.15) \quad $ \begin{array}{l}
                                  ^*( A_1 \bigcup ~.~.~.~ \bigcup A_n ) ~=~ ^*A_1 \bigcup ~.~.~.~ \bigcup~ ^*A_n \\ \\
                                   ^*( A_1 \bigcap ~.~.~.~ \bigcap A_n ) ~=~ ^*A_1 \bigcap ~.~.~.~ \bigcap~ ^*A_n
                          \end{array} $

\medskip
Indeed, first we note that (2.7.6) obviously gives

\bigskip
(2.7.16) \quad $ A ~\subseteq~ B ~\subseteq~ {\bf R} ~~~\Longrightarrow~~~ ^*A ~\subseteq~ ^*B ~\subseteq~ ^*{\bf R} $

\medskip
thus we have the inclusion "$\supseteq$" in the first, and the inclusion "$\subseteq$" in the second of the relations
(2.7.15). Let us now show the inclusion "$\subseteq$" in the first of the relations (2.7.15). According to (2.7.6), if $s =
[~ < s_1, s_2, s_3, ~.~.~.~ > ~] \in ^*( A_1 \bigcup ~.~.~.~ \bigcup A_n )$, then

$$ s_j \in A_1 \bigcup ~.~.~.~ \bigcup A_n ~~~~\mbox{( a. e. )}~~ j \in {\bf N} $$

\medskip
thus

$$ \bigcup_{1 \leq i \leq n}~ \{~ j \in {\bf N} ~~|~~ j \in A_i  ~\} ~=~ \{~ j \in {\bf N} ~~|~~ j \in A_1 \bigcup ~.~.~.~
                                                                                        \bigcup A_n ~\} \in {\cal U} $$

\medskip
which in view of (2.3.6) gives $\{~ j \in {\bf N} ~~|~~ j \in A_i  ~\}  \in {\cal U}$, for a certain $1 \leq i \leq n$. Thus (2.7.6)
yields $s \in A_i$, and the first relation in (2.7.15) is proved. \\
Let us show now that the inclusion "$\supseteq$" holds in the second of the relations (2.7.15). If $s \in~ ^*A_1 \bigcap
~.~.~.~ \bigcap~ ^*A_n$ then (2.7.6) gives for every $1 \leq i \leq n$ the relation $\{~ j \in {\bf N} ~~|~~ j \in A_i  ~\}  \in {\cal U}$.
And since ${\cal U}$ is a filter, it follows that

$$ \{~ j \in {\bf N} ~~|~~ j \in A_1 \bigcap ~.~.~.~ \bigcap A_n ~\} ~=~ \bigcap_{1 \leq i \leq n}~ \{~ j \in {\bf N} ~~|~~ j \in A_i  ~\}
                                                                                                                            \in {\cal U} $$

\medskip
and thus in view of (2.7.6), the proof of the second relation in (2.7.15) is completed. \\

In the next chapter, as one of the immediate applications of the more general transfer method ot (TS) type presented
there, we shall reestablish (2.7.14) and (2.7.15), and also give a number of other similar basic results. \\ \\

{\bf Convention 2.7.2} \\

As we have seen, in view of convention (2.7.5), we have the {\it extension} property (2.7.11), and its particular case
(2.7.12). In this way, when we {\it transfer} a relation $P$, or a function $f$, from ${\bf R}$ to the corresponding relation
$^*P$, respectively, function $^*f$ on $^*{\bf R}$, we in fact obtain extensions of the entities transferred. \\

In view of that, and unless explicitly mentioned otherwise, we shall often use the term {\it extension} of sets,
functions, or relations, with the above meaning, namely, that they were obtained through {\it transfer}. \\ \\

{\bf Note 2.7.1} \\

So far, we managed to extend the standard reals ${\bf R}$ into a strictly larger totally ordered field $^*{\bf R}$, and in
the process we could also transfer relations and functions from ${\bf R}$ to $^*{\bf R}$. However, we should not forget
that all these constructions depend on the {\it choice} of the free ultrafilter ${\cal U}$ which through the corresponding
maximal ideal ${\cal M}$ gives the construction of $^*{\bf R}$ in (2.1.4). \\

And then the question arises to what extent does the whole nonstandard theory depend on the choice of the respective
free ultrafilter ? Is it possible that there is no such a dependence ? Are the various totally ordered fields $^*{\bf R}$
which correspond to different choices of free ultrafilters isomorphic ? \\

As mentioned in chapter 5, this issue is related to certain additional assumptions of what is called saturation. And
under such assumptions, different free ultrafilters do actually lead to one and the same nonstandard totally ordered
field $^*{\bf R}$ in (2.1.4), see Kreisler [2] for further details. Otherwise, as specified later, there can be a large variety
of different nonstandard totally ordered field extensions  $^*{\bf R}$ of ${\bf R}$. \\ \\

{\bf 8. Countable Saturation} \\ \\

Saturation is an important and typical concept in Nonstandard Analysis, and it will be dealt with in chapter 5. Here we
indicate a particular case of it, called {\it countable saturation}, see Proposition 2.8.1 below. And as we can see, this
highly nontrivial nonstandard porperty can already be obtained with the rather basic tools developed so far. \\

Let us start with a simple example. If we consider the standard intervals $( 0, 1/n) \subset {\bf R}$, with $n \in {\bf N},~
n \geq 1$, then obviously they have the {\it finite intersection} property, namely, the intersection of any finite number of
them is not void. However, obviously

$$ \bigcap_{1 \leq n < \infty}~ ( 0, 1/n ) ~=~ \phi $$

\medskip
On the other hand, if we consider the corresponding intervals $^*( 0, 1/n) \\ \subset~ ^*{\bf R}$, with $n \in {\bf N},~ n \geq
1$, they also have the {\it finite intersection} property, and in addition, and unlike in the standard case, we have

$$ \bigcap_{1 \leq n < \infty}~ ^*( 0, 1/n ) ~=~ mon ( 0 ) ~\cap~  \{~ s \in~ ^*{\bf R} ~~|~~ s  ~>~ 0 ~\} ~\neq~  \phi $$

\medskip
This phenomenon is a particular case of a more general property. Let us clarify that by first introducing the following
concept. \\ \\

{\bf Definition 2.8.1} \\

A subset $A \subseteq~ ^*{\bf R}$ is called {\it internal}, if and only if there exists a sequence of subsets $A_i \subseteq
{\bf R}$, with $i \in {\bf N}$, such that

\bigskip
(2.8.1) \quad $ A ~=~ \{~ [~ < s_1, s_2, s_3, ~.~.~.~ > ~] ~~|~~ \{~ i \in {\bf N} ~~|~~ s_i \in A_i ~\} \in {\cal U} ~\} $

\hfill $\Box$ \\

The general case of the concept of {\it internal} entity will be given in chapter 5. \\

Obviously, if $A ~=~ ^*B$, for some $B ~\subseteq~ {\bf R}$, then $A$ is internal. In particular, the above intervals $^*( 0,
1/n )$ are internal. \\ \\

{\bf Proposition 2.8.1} \\

Let

\bigskip
(2.8.2) \quad $ ^*{\bf R} \supseteq A_0 ~\supseteq~ A_1 ~\supseteq~ A_2 ~\supseteq~ ~.~.~.~ ~\supseteq~ A_n
                                                                 ~\supseteq~ ~.~.~.~ $

\medskip
be a decreasing sequence of nonvoid internal subsets. Then

\bigskip
(2.8.3) \quad $ \bigcap_{n \in {\bf N}}~ A_n ~\neq~ \phi ~ $ \\

{\bf Proof} \\

Let us take for each $A_n$ a sequence of subsets $A_{n,~i} \subseteq {\bf R}$, with $i \in {\bf N}$, which satisfies (2.8.1).
Since we assumed that $A_n \neq \phi$, it follows that

$$ \{~ i \in {\bf N} ~~|~~ A_{n,~i} ~\neq~ \phi ~\} \in {\cal U} $$

\medskip
Furthermore, the inclusion $A_n \supseteq A_{n + 1}$ gives

$$ \{~ i \in {\bf N} ~~|~~ A_{n,~i} ~\supseteq~ A_{n + 1,~i} ~\} \in {\cal U} $$

\medskip
Hence for every $m \in {\bf N}$, we have

$$ B_m ~=~ \{~ i \in {\bf N} ~~|~~ A_{n,~0} ~\supseteq~ ~.~.~.~ ~\supseteq~ A_{n,~m} ~\neq~ \phi ~\} \in {\cal U} $$

\medskip
Also obviously

$$ B_0 ~\supseteq~ B_1 ~\supseteq~ B_2 ~\supseteq~ ~.~.~.~ $$

\medskip
Let us now construct a sequence $s_i \in {\bf R}$, with $i \in {\bf N}$, as follows. For $i \in B_0$, we define

$$ m_i ~=~ \max~ \{~ m \in {\bf N} ~~|~~ m ~\leq~ i,~~ i \in B_m ~\} $$

\medskip
and since we took $i \in B_0$, it follows obviously that the above is a correct definition, since the set in the right hand
term is not void. Furthermore, obviously $i \in B_{m_i}$. Hence we can take $s_i \in A_{i,~0} \cap ~.~.~.~ \cap
A_{i,~m_i}$. \\
If $i \notin B_0$, we take $s_i$ arbitrary in ${\bf R}$. \\

Let us prove now that for every $n \in {\bf N}$, we have

\bigskip
(2.8.4) \quad $ \{~ i \in {\bf N} ~~|~~ n ~\leq~ i ~\} \bigcap B_n ~\subseteq~ \{~ i \in {\bf N} ~~|~~ s_i \in A_{n,~ i} ~\} $

\medskip
Indeed, if $n \leq i$ and $i \in B_n$, then $n \leq m_i$, thus the above choice of $s_i$ gives $s_i \in A_{i,~n}$. \\

Now it follows that $s = < s_0, s_1, s_2, ~.~.~.~ >$ has the property that, for every $n \in {\bf N}$

\bigskip
(2.8.5) \quad $ \{~ i \in {\bf N} ~~|~~ s_i \in A_{n,~i} ~\} \in {\cal U} $

\medskip
Indeed, since ${\cal U}$ is a free ultrafilter on ${\bf N}$, it follows that $\{~ i \in {\bf N} ~~|~~ n ~\leq~ i ~\} \in {\cal U}$. Thus
(2.8.4) and the fact shown above, namely that $B_n \in {\cal U}$, will give indeed (2.8.5). However, that relation means
precisely that $[ s ] \in  \bigcap_{n \in {\bf N}}~ A_n$.

\hfill $\Box$ \\

The {\it countable saturation} property in (2.8.2), (2.8.3) has important consequences. Some of them, however, create
difficulties, like for instance the following one in Measure Theory. \\ \\

{\bf Corollary 2.8.1} \\

Let ${\cal S}$ be a $\sigma$-algebra of internal subsets of $^*{\bf R}$, and let $E_n \in {\cal S}$, with $n \in {\bf N}$, be
an increasing sequence of sets. Then there exists $m \in {\bf N}$, such that

\bigskip
(2.8.6) \quad $ \bigcup_{n \in {\bf N}}~ E_n ~=~ E_m $

\medskip
Therefore $E_n = E_m$, for $n \in {\bf N},~ n \geq m$. \\

{\bf Proof} \\

Let $E = \bigcup_{n \in {\bf N}}~ E_n$, then by the assumption, we have $E \in {\cal S}$, thus $E$ is internal. Let us now
take $A_n = E \setminus E_n$, with $n \in {\bf N}$. It is easy to see that each $A_n$ is again internal, furthermore

$$ A_0 ~\supseteq~ A_1 ~\supseteq~ A_2 ~\supseteq~ ~.~.~.~ $$

\medskip
Now in case $A_n \neq \phi$, for $n \in {\bf N}$, then Proposition 2.8.1 gives

$$ \bigcap_{n \in {\bf N}}~ A_n ~\neq~ \phi ~ $$

\medskip
On the other hand, by the above construction, we have

$$ \bigcap_{n \in {\bf N}}~ A_n ~=~ \bigcap_{n \in {\bf N}}~ ( E \setminus E_n ) ~=~ \phi ~ $$

\medskip
The contradiction obtained proves (2.8.6).

\hfill $\Box$ \\

Obviously, property (2.8.6) {\it trivializes} $~\sigma$-algebras in $^*{\bf R}$, as it prevents the existence of genuine
ones, in which one could have {\it strictly} increasing countably {\it infinite} sequences of sets. \\
The theory of Loeb measures and integration, presented in short in chapter 6, was originally introduced precisely in
order to deal with the kind of difficulties in the above corollary. \\

Let us now mention a few properties of internal sets which follow immediately from (2.8.1), and which we shall use in
the sequel. If $A,~ B$ are internal, then the same hold for $A \cup B,~ A \cap B$ and $A \setminus B$. \\

Regarding the {\it size} of internal sets, we have the following {\it dichotomy}. \\ \\

{\bf Corollary 2.8.2} \\

If $A \subseteq~ ^*{\bf R}$ is internal, then $A$ is either finite, or it is uncountably infinite. \\

{\bf Proof} \\

Let us assume that $A$ is infinite and countable, namely $A = \{~ a_0, a_1, a_2, ~.~.~.~ ~\}$. In view of (2.8.1), it is easy
to see that $A \setminus B$ is internal, whenever $A,~ B$ are internal. Also finite sets are internal. Therefore, the sets

$$ A_n ~=~ A \setminus \{~ a_0, ~.~.~.~ , a_n ~\} ~=~ \{~ a_{n + 1}, a_{n + 2}, a_{n + 3}, ~.~.~.~ ~\},~~~ n\in {\bf N} $$

\medskip
are internal. However $A_0 \supseteq A_1 \supseteq A_2 \supseteq ~.~.~.~$ and $\bigcap_{n \in {\bf N}}~ A_n = \phi$, thus
(2.8.3) is contradicted. \\ \\

{\bf 9. The Abundance in Ultrafilters} \\

The method of nonstandard extension is not only relatively new, but it is also not widely enough familiar among
mathematicians. And as seen even during the first steps in this chapter in the construction of the nonstandard
extension $^*{\bf R}$ of ${\bf R}$, somewhat unusual mathematical entities such as {\it ultrafilters}, for instance, are
essentially involved. \\
It is, therefore, important to note that the construction of nonstandard extensions is not constrained in any way to, or
dependent on certain rather limited or singular situations, or scarce entities. Indeed, as seen later in chapter 5, and
specifically related to the issue of saturation, there is a large variety of such nonstandard extensions which can be
constructed for general mathematical structures, and not only for the real numbers ${\bf R}$. \\

Here, as further pointing to the large supply of ingredients which are available in the construction of nonstandard
extensions, we present a well known result related to the {\it abundance} of ultrafilters on arbitrary infinite sets. The
proof of this result can be found in Jech [p. 75]. \\

Let be given an infinite set $I$ which has the cardinal $\kappa$. As every filter ${\cal F}$ on $I$ is obviously a
subset of ${\cal P} ( I )$, there cannot be more than $2^{2^\kappa}$ ultrafilters ${\cal U}$ on $I$. \\
The interesting thing is that, in fact, one does always have this {\it maximum} number of ultrafilters. Furthermore, as
seen next, one has even a stronger result. \\

Indeed, let us consider the class of ultrafilters given by \\ \\

{\bf Definition 2.9.1} \\

If $I$ is an infinite set with cardinal $\kappa$ and ${\cal U}$ is an ultrafilter on $I$, we say that ${\cal U}$ is {\it
uniform}, if and only if we have $car J = \kappa$, for every $J \in {\cal U}$.

\hfill $\Box$ \\

Obviously, only {\it free} ultrafilters can be uniform. And the number of uniform ultrafilters is given by  \\ \\

{\bf Theorem 2.9.1 ( Pospisil )} \\

On every infinite set $I$ of cardinal $\kappa$, there exist $2^{2^\kappa}$ uniform ultrafilters. \\ \\

\chapter{The Transfer of Standard Properties}

{\bf 1. Towards a General Transferring Mechanism} \\ \\

As mentioned in section 2, in chapter 1, when doing mathematics in ${\bf R}$, we are interested in {\it large} classes of
properties, and we would like to know, whether or not, it is possible to transfer them to $^*{\bf R}$. Furthermore, we
would like to have a {\it general} enough transfer mechanism, so that  we can avoid doing transfer in the tedious, one
by one way, by finding and using each time a specific particular proof. \\
However, as the properties of ${\bf R}$ of being Archimedean or Dedekind complete, for instance, show it, and which
are {\it no} longer valid for $^*{\bf R}$ as well, such a transfer is not a trivial task. \\

Therefore, after having constructed in chapter 2 the nonstandard totally ordered field extension $^*{\bf R} ~=~ {\bf R}^{\bf
N} ~/~ {\cal M}$ of ${\bf R}$, now as a second step, we have to set up a {\it wholesale industry} for performing the
mentioned kind of transfers of large classes of properties from ${\bf R}$ to $^*{\bf R}$. \\

But before we start setting up in this chapter a general enough method for transfer, let us recall that in chapter 2, we
already managed to transfer in an ad-hock manner from ${\bf R}$ to $^*{\bf R}$ the property of being a {\it totally ordered
commutative field}. And for that purpose, all what we used were a few well known results in algebra related to
factorization with a maximal ideal, plus some basic property of ultrafilters. On that occasion, it also followed that ${\bf
R}$ is a totally ordered subfield of $^*{\bf R}$. Similarly, in chapter 2, we transferred relations and functions from ${\bf
R}$ to $^*{\bf R}$. \\

However, the fact is that there are far more properties of ${\bf R}$ which we would like to transfer to $^*{\bf R}$, and
such a case by case approach as used in chapter 2 will not lead us far enough. Therefore, as mentioned in section 2, in
chapter 1, we shall have to make use as well of a certain amount of Mathematical Logic. \\

In order better to clarify the purpose of this chapter, let us recall the way we classified the operations of transfer,
namely : \\

(TE)~:~~ the {\it transfer} of various {\it mathematical entities} from ${\bf R}$ to $^*{\bf R}$, \\
         \hspace*{1.45cm} the first and simplest case of it being the field homomorphism \\
         \hspace*{1.45cm} $^*(~)$ in (2.1.4) next \\

(TS)~:~~ the {\it transfer} of {\it sentences} $\Phi$ from a given {\it language} ${\cal L}_{\bf R}$ about \\
         \hspace*{1.45cm} mathematical entities in ${\bf R}$ to corresponding sentences $^*\Phi$ in \\
         \hspace*{1.45cm} a language ${\cal L}_{^*{\bf R}}$ about mathematical entities in $^*{\bf R}$. \\

In these terms, the aim of the present chapter is to introduce one simple version of the stage (TS) of transfer. \\ \\

{\bf 2. Simple Systems and their Simple Languages} \\ \\

For the convenience of those not familiar with Mathematical Logic, we shall set up the mechanism of transfer in {\it two}
stages, starting first with an easier and limited version of it, which only involves most simple tools of Mathematical
Logic. These tools center around the concepts of {\it simple systems} and {\it simple languages} which are introduced
in this section, see Keisler [1,2] and Hurd \& Loeb. \\
Later, in chapter 5, stage two will be presented with the introduction of the method of Superstructures and the
corresponding more evolved formal languages. \\

Since there is no complication involved, and in fact, it may make things more clear, we shall begin with a more general
setup than that of ${\bf R}$. \\ \\

{\bf Definition 3.2.1} \\

A {\it simple system} is any structure of the form

\bigskip
(3.2.1) \quad $ {\cal S} ~=~ ( S,~ ( P_i ~|~ i \in I ),~ ( f_j ~|~ j \in J ) ) $

\medskip
where $S$ is a nonvoid set, $P_i \subseteq S^{n_i}$ are $n_i$-ary relations on $S$, while $f_j : D_j \longrightarrow S$
are functions of $m_j$ variables defined on their respective domains $D_j \subseteq S^{m_j}$, and with values in
$S$.

\hfill $\Box$ \\ \\

We shall of course be interested in the two simple systems

\bigskip
(3.2.2) \quad $ \begin{array}{l}
                                 {\cal R} ~=~  ( ~{\bf R},~ (~ =,~  \leq~ ),~ (~ +,~ .~ )~ ) \\ \\
                                 ^*{\cal R} ~=~ (~ ^*{\bf R},~ (~ =, ~\leq~ ),~ ( ~+,~ . ~)~ )
                       \end{array} $

\medskip
however, as specified later, see (3.3.1) and (3.4.1), we shall consider a far large number of relations and functions, both
on ${\bf R}$ and $^*{\bf R}$. \\

For the moment, nevertheless, let us pursue the general form of simple systems given in (3.2.1). \\

Given an n-ary relation $P \subseteq S^n$, and $s_1, ~.~.~.~ , s_n \in S$, it will on occasion be convenient to use one or
another of the following two equivalent notations

$$ ( s_1, ~.~.~.~ , s_n ) \in P ~~~\Longleftrightarrow~~~ P < s_1, ~.~.~.~ , s_n > $$

\medskip
As already mentioned, functions are particular cases of relations. Namely, an n variable function can be seen as an
( n + 1 )-ary relation. Indeed, if we have an n variable function $f : D \longrightarrow S$, with domain $D \subseteq S^n$,
then for $s_1, ~.~.~.~ , s_n, s_{n + 1} \in S$, with $( s_1, ~.~.~.~ , s_n ) \in D$, we have

$$ f ( s_1, ~.~.~.~ , s_n ) ~=~ s_{n + 1} ~~~\Longleftrightarrow~~~ ( s_1, ~.~.~.~ , s_n, s_{n + 1} ) \in f $$

\medskip
And what is specific for functions among relations is that

$$ ( s_1, ~.~.~.~ , s_n, s_{n + 1} ),~ ( s_1, ~.~.~.~ , s_n, t_{n + 1} ) \in f ~~~\Longrightarrow~~~ s_{n + 1} ~=~ t_{n + 1} $$ \\

{\bf Definition 3.2.2} \\

With each simple system ${\cal S} ~=~ ( S,~ ( P_i ~|~ i \in I ),~ ( f_j ~|~ j \in J ) )$ we associate a {\it simple language}
${\cal L}_{\cal S}$ as follows. \\

The {\it alphabet} of ${\cal L}_{\cal S}$ is made up from the next four categories of symbols : \\

1. Two logical connectives : $\wedge$ and $\longrightarrow$ which are interpreted as "and" and "implies". \\

2. One quantifier symbol : $\forall$ which is interpreted as "for all". \\

3. Six parentheses : [ , ] , ( , ) , $<$ ,  $>$ . \\

4. Variable symbols : a countable set x, y, ~.~.~.~ \\

These symbols do not depend on the simple system ${\cal S} ~=~\\=~ ( S,~ ( P_i ~|~ i \in I ),~ ( f_j ~|~ j \in J ) )$, while the
next three categories do, namely : \\

5. Constant symbols : for each $s \in S$ we associate its name $\underline{s}$. \\

6. Relation symbols : for each $P_i$, with $i \in I$, we associate its name $\underline{P_i}$. \\

7. Function symbols : for each $f_j$, with $j \in J$, we associate its name $\underline{f_j}$. \\

Let us note here that {\it variable symbols} are {\it not} associated names in the simple languages ${\cal L}_{\cal S}$,
and they appear in these languages directly, as themselves. After all, the variable symbols belong to that part of the
alphabet of simple languages ${\cal L}_{\cal S}$ which do not depend on the respective simple systems ${\cal S}$. \\ \\

{\bf Convention 3.2.1} \\

For convenience, in the case of the better known and more often used elements, relations and functions, we shall not
underline them when writing their names. For instance, we shall not underline constants such as 0, 1, $\pi$, and so on,
or relations like $\leq$, $<$, or functions given by +,~ .~, ~ $\sin$, and so on.

\hfill $\Box$ \\ \\

It is important to note that such a {\it differentiation} as at 5 , 6 and 7 above, between, on the one hand, usual
mathematical entities such as elements in a set, relations on a set, or functions on a set, and on the other hand,
symbols which name them is fundamental in Mathematical Logic. In particular, the components of the simple languages
${\cal L}_{\cal S}$ which we are constructing now, are {\it not} entities within any of the sets in the simple systems
${\cal S} ~=~ ( S,~ ( P_i ~|~ i \in I ),~ ( f_j ~|~ j \in J ) )$ which they may be talking about, but on the contrary, they are entities
outside of such sets. And the {\it only} connections between such simple languages ${\cal L}_{\cal S}$ and various
simple system ${\cal S} ~=~ ( S,~( P_i ~|~ i \in I ),~ ( f_j ~|~ j \in J ) )$ is through {\it interpretation}, that is

\bigskip
\begin{math}
\setlength{\unitlength}{1cm}
\thicklines
\begin{picture}(15,2)
\put(1,1){${\cal L}_{\cal S}$}
\put(2.4,1.4){$\mbox{interpretation}$}
\put(2,1.1){\vector(1,0){3.5}}
\put(6,1){${\cal S} ~=~ ( S,~( P_i ~|~ i \in I ),~ ( f_j ~|~ j \in J ) )$}
\end{picture}
\end{math}

which will be defined below. \\

One of the main reasons for such a differentiation is that a simple language ${\cal L}_{\cal S}$ may have {\it different}
interpretations ${\cal S} ~=~ ( S,~ ( P_i ~|~ i \in I ),\\~ ( f_j ~|~ j \in J ) )$. \\
Also, the same element, relation or function may be given different names. In this way, the correspondences

\bigskip
\begin{math}
\setlength{\unitlength}{1cm}
\thicklines
\begin{picture}(15,2)
\put(1,1){${\cal L}_{\cal S}$}
\put(2,1.3){\vector(1,0){3.5}}
\put(5.5,1){\vector(-1,0){3.5}}
\put(6,1){${\cal S} ~=~ ( S,~( P_i ~|~ i \in I ),~ ( f_j ~|~ j \in J ) )$}
\end{picture}
\end{math}

can be defined in {\it multiple} ways. \\

In short, on the one hand, the entities of usual mathematics, and on the other hand, the languages we use when we
think, write or talk about them {\it cannot} be claimed to belong to the same realm, although in their everyday
mathematical activity usual mathematicians have a rather strong, even if automatic and less than conscious tendency
to do so. \\
Further, the ways we associate names to entities, of entities to names are far from being unique or canonical. \\ \\

{\bf Definition 3.2.2 (continuation)} \\

The words, called {\it terms}, of the simple language ${\cal L}_{\cal S}$ are defined inductively as follows : \\

8. Each constant or variable symbol is a {\it term}. \\

9. If $\underline{f}$ is the name of a function of n variables and $\tau^1, ~.~.~.~ , \tau^n$ are terms, then
$\underline{f} ( \tau^1, ~.~.~.~ , \tau^n)$ is again a {\it term}. \\

10. A term which does not contain variable symbols is called a {\it constant term}. \\

At last, we get to the sentences of the simple language ${\cal L}_{\cal S}$, and they are called {\it simple sentences},
having one the following two forms : \\

11. {\it Atomic sentences}, given by $\underline{P} < \tau^1, ~.~.~.~ , \tau^n >$, where $\underline{P}$ is the name of an
n-ary relation, while $\tau^1, ~.~.~.~ , \tau^n$ are {\it constant terms}. \\

12. {\it Compound sentences}, given by

$$ ( \forall~ x_1 ) ~.~.~.~ ( \forall~ x_n )~ \Big [ ~ \bigwedge_{1 \leq i \leq k}~ \underline{P_i} < \xi_i > ~\longrightarrow~
                                                                     ~ \bigwedge_{1 \leq j \leq l}~ \underline{Q_j} < \chi_j > ~ \Big ] $$

\medskip
Here

$$ \bigwedge_{1 \leq i \leq k}~\underline{P_i} < \xi_i > $$

\medskip
denotes $\underline{P_1} < \xi_1 > \wedge ~.~.~.~ \wedge \underline{P_k} < \xi_k >$. Further, $\underline{P_i}$ is the
name of an n$_i$-ary relation, while $\xi_i$ are n$_i$ terms $\tau^1_i, ~.~.~.~ , \tau^{n_i}_i$ which may be constant, or
contain some, or all of the variable symbols $x_1, ~.~.~.~ , x_n$, but do not contain any other variable symbols.
Similarly for $\bigwedge_{1 \leq j \leq l}~ \underline{Q_j} < \chi_j >$. \\ \\

{\bf Note 3.2.1} \\

In the context of formal languages in Mathematical Logic, the concept of {\it sentence} is fundamental. And as also seen
above in 11 and 12, which define the simple sentences of the simple language ${\cal L}_{\cal S}$, the caharcteristic
feature of sentences is that they do {\it not} contain variable symbols which are not in the range of a quantifier. Let us
give some examples which may clarify this issue. \\
A formula in a formal language which contains a variable symbol not in the range of a quantifier, such as for instance
$x$ in

$$ x + x^2 = 2 $$

\medskip
makes a statement about the possible values of that variable symbol, in this case that $x \in \{ -2, 1 \}$, and as such,
that is, as a formula, it is not supposed to be either true, or false. \\
On the ohter hand, a sentence, which by definition is a formula without variable symbols not in the range of a quantifier,
may fail to make a statement about the values of any of the variable symbols which it may contain. However, it can
make a statement about the theory to which that sentence may belong. In particular, such a sentence may be true or
false in the respective theory, or alternatively, it may be independent of that theory. The latter case means that the
respective sentence, as much as its negation, cannot be proved within the means of that theory alone. \\
For instance, the sentence

$$ ( \forall~ x )~ [~ x^2 + 1 > 0 ~] $$

\medskip
is true in ${\bf R}$, but false in ${\bf C}$, and in either case, it does not say anything about the particular values the
variable symbol $x$ may take, except that it can range over all of ${\bf R}$, respectively, ${\bf C}$.

\hfill $\Box$ \\ \\

As far as examples are concerned, there will be plenty of them given, once we start applying systematically the method
of transfer. Until then, we shall pursue the theoretical presentation. In this regard, next we deal with the issue of
interpretation. \\ \\

{\bf 3. Interpretation of Simple Languages} \\ \\

A main aim of setting up formal languages, such as for instance the simple languages ${\cal L}_{\cal S}$, is to obtain a
rigorous and systematic method for identifying and studying those sentences which are {\it true}. \\
The aim of {\it interpretation} of a formal language is to give one of the possible ways for finding out which of its
sentences are true. \\
Of course, the quality of being {\it true} only applies to the {\it sentences} of a formal language, and not also to its other
auxiliary entities, such as symbols, terms, and so on, which are used in building up sentences. These auxiliary entities
are only {\it interpreted}, without the issue of their truth arising. \\

A most important point, typical for Mathematical Logic, is the following. Prior to the dichotomy "{\it true - false}", there is
the dichotomy "{\it interpretable - non interpretable}", to which the auxiliary entities of formal languages are subjected.
Namely, in our case of simple languages, {\it constant symbols} and {\it constant terms} are subjected only to the
dichotomy \\

\bigskip
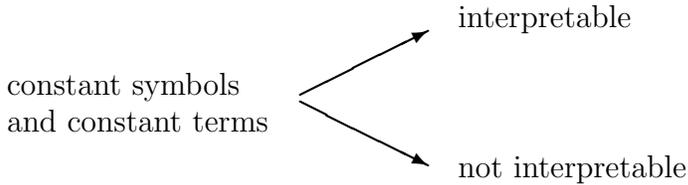
\begin{math}
\setlength{\unitlength}{1cm}
\thicklines
\begin{picture}(15,2.5)
\put(0,1.5){$\mbox{constant symbols}$}
\put(0,1){$\mbox{and constant terms}$}
\put(3.9,1.5){\vector(2,1){1.7}}
\put(6.,2.4){$\mbox{interpretable}$}
\put(3.9,1.4){\vector(2,-1){1.7}}
\put(6,0.4){$\mbox{not interpretable}$}
\end{picture}
\end{math}

\medskip
while {\it simple sentences}, that is, {\it atomic sentences} and {\it compound sentences}, are subjected to the
successive two dichotomies

\bigskip
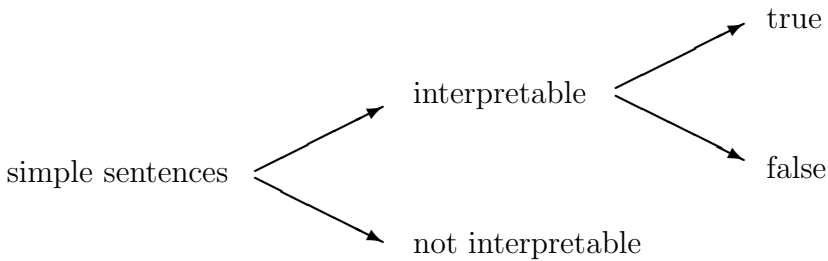
\begin{math}
\setlength{\unitlength}{1cm}
\thicklines
\begin{picture}(15,4.5)
\put(0,1.35){$\mbox{simple sentences}$}
\put(3.3,1.5){\vector(2,1){1.7}}
\put(5.4,2.4){$\mbox{interpretable}$}
\put(8.1,2.6){\vector(2,1){1.7}}
\put(10.1,3.4){$\mbox{true}$}
\put(8.1,2.5){\vector(2,-1){1.7}}
\put(10.1,1.4){$\mbox{false}$}
\put(3.3,1.4){\vector(2,-1){1.7}}
\put(5.4,0.4){$\mbox{not interpretable}$}
\end{picture}
\end{math}

\medskip
and {\it only} those simple sentences which are interpretable can further be subjected to the dichotomy "true -
false". \\

As far as {\it variable symbols} or terms with such symbols are concerned, they are not interpreted, this being
precisely one of the essential points in the conept of 'variable' which, just as in algebra, is only supposed to be a
place holder. \\ \\

{\bf Definition 3.3.1} \\

The {\it interpretation} of the simple language ${\cal L}_{\cal S}$ {\it within} the simple system ${\cal S} ~=~
( S,~ ( P_i ~|~ i \in I ),~ ( f_j ~|~ j \in J ) )$ is done inductively, as follows : \\

1. A constant symbol $\underline{s}$ which names an element $s \in S$ has as {\it interpretation} this element $s$. \\
We note that a constant symbol is always interpretable. \\

2. A term $\underline{f} ( \tau^1, ~.~.~.~ , \tau^n )$ is {\it interpretable}, if and only if each of the terms $\tau^i$ has an interpretation
as some element $s^i \in S$, and in addition, $( s^1, ~.~.~.~ , s^n )$ is in the domain of definition of $f$, where
$\underline{f}$ is the name of the function of n variables $f$. Then the {\it interpretation} of the term $\underline{f}
( \tau^1, ~.~.~.~ , \tau^n)$ is the element $f ( s^1, ~.~.~.~ , s^n) \in S$. \\
It follows that not all such terms are interpretable. \\

3. An atomic sentence $\underline{P} < \tau^1, ~.~.~.~ , \tau^n >$ is {\it interpretable}, if and only if each of the terms $\tau^i$  has an
interpretation as some element $s^i \in S$. \\
Further, the atomic sentence $\underline{P} < \tau^1, ~.~.~.~ , \tau^n >$ is {\it true}, if and only if  it is interpretable, and in addition
$P < s^1, ~.~.~.~ , s^n >$ is true, where $\underline{P}$ is the name of the n-ary relation $P$. \\
Again, it follows that not every atomic sentence is interpretable, let alone, true. And an atomic sentence {\it fails} to be
true in one of the following two cases : either it is not interpretable, or it is interpretable, but the corresponding $P
< s^1, ~.~.~.~ , s^n >$ is false. \\

At last \\

4. A compound sentence

$$ ( \forall~ x_1 ) ~.~.~.~ ( \forall~ x_n )~ \Big [ ~ \bigwedge_{1 \leq i \leq k}~ \underline{P_i} < \xi_i > ~\longrightarrow~
                                                                     ~ \bigwedge_{1 \leq j \leq l}~ \underline{Q_j} < \chi_j > ~ \Big ] $$

\medskip
is {\it interpretable}, if and only if each $\underline{P_i} < \xi_i >$ and each $\underline{Q_j} < \chi_j >$ is interpretable for certain
substitutions in them of the {\it variable symbols} $x_1, ~.~.~.~ ,  x_n$ with constant symbols. \\
Further, the above compound sentence is {\it true}, if and only if whenever it is interpretable, and each of the corresponding
$\underline{P_i} < \xi_i >$ are true for certain substitutions of the variable symbols, then each of $\underline{Q_j}
< \chi_j >$ with the same substitutions of variable symbols are also true.

\hfill $\Box$ \\ \\

Here we should note that the above definition of truth can allow a certain {\it indirect} specification of the range of a
universal quantifier $\forall$ in a compound sentence, a specification influenced by some of the sentences in that
compound sentence, to which the quantifier may refer. \\

Let us illustrate that effect in the case of usual real numbers in ${\bf R }$. \\

For that purpose, first we enrich the simple system ${\cal R}$ in (3.2.2) in the following form which we shall use from
now on, namely

\bigskip
(3.3.1) \quad $ {\cal R} ~=~ (~ {\bf R},~ {\cal P},~ {\cal F} ~) $

\medskip
where ${\cal P}$ is the set of all possible relations $P$ of finite arity on ${\bf R}$, while ${\cal F}$ is the set of all
possible functions $f$ of finitely many variables in ${\bf R }$, and with values in ${\bf R }$. \\
Of course, we shall also be interested in the corresponding simple system

\bigskip
(3.3.2) \quad $ ^*{\cal R} ~=~ (~ ^*{\bf R},~ ^*{\cal P},~ ^*{\cal F} ~) $

\medskip
specified later in (3.4.1), (3.4.2). \\

Now for instance, the compound sentence in ${\cal L}_{\cal R}$

$$ ( \forall~ x )~ \Big [~ ( \sqrt x > -1 ) ~\longrightarrow~ ( \sqrt x \geq 0 ) ~\Big ] $$

\medskip
is true in ${\cal R}$, since $\sqrt x$ is only interpretable for $x \geq 0$, in which case both sentences $\sqrt x > -1$ and
$\sqrt x \geq 0$ are true. In this way an indirect specification of the range of $( \forall~x )$ got implemented. \\

One of the advantages of simple systems and simple languages is that in such a framework it is so much easier to
prove the fundamental Transfer Property, see section 4, next. \\
However, the same simplicity can create difficulties as well, since in mathematics we often deal with sentences which
have a more complex structure than those allowed in simple languages. For instance, let us consider the true sentence
in ${\cal R}$

$$ ( \forall~x ) ( \exists~ y )~ \Big [~ (~ \underline{{\bf R}} < x > \wedge~ x \neq 0 ~) ~\longrightarrow~ (~  x . y = 1 ~) ~\Big ] $$

\medskip
Here $\underline{{\bf R}} < x >$ denotes the name in the simple language ${\cal L}_{\cal R}$ of the unary relation $x \in
{\bf R}$, that is, of the property of $x$ of being a usual real number. Obviously, this sentence does {\it not} belong to the
simple language ${\cal L}_{\cal R}$, in view of the existential quantifier $\exists$ which is present in it. \\
Also, simple languages only allow the two logical connectives $\wedge$ and $\longrightarrow$, while there are several
other such connective and their use can be particularly convenient. \\
Nevertheless, with some manipulations, such shortcomings of simple languages can be overcome. Moreover, we do
not have to dwell too much on such issues, since we only use simple languages as a stepping stone towards the
method of Superstructures and their languages, which do no longer suffer from any such shortcomings. \\

Here, for illustration, we show in an example the way the lack of the existential quantifier $\exists$ in simple languages
can be overcome with the use of the so called {\it Skolem functions}, see Keisler [1,2] and Hurd \& Loeb. \\
Let us consider the above sentence which is true in ${\cal R}$, but as we have noted, it is not a simple sentence in the
simple language ${\cal L}_{\cal R}$. \\
In order to deal with it, we can assume the existence of a function

$$ \sigma : (~ {\bf R} \setminus \{~ 0 ~\} ~) ~~\longrightarrow~~ {\bf R} $$

\medskip
such that $\sigma ( x ) ~=~ 1 / x$, for $x \in {\bf R},~ x \neq 0$. Then the above sentence is equivalent with the following
simple sentence from ${\cal L}_{\cal R}$, namely

$$ ( \forall~x )~ \Big [~ (~ \underline{{\bf R}} < x > \wedge~ x \neq 0 ~) ~\longrightarrow~
                                                                             (~  x . \underline{\sigma} ( x ) = 1 ~) ~\Big ] $$

\medskip
Any function $\sigma$ performing such a task of substituting for the existential quantifier $\exists$ is called a {\it
Skolem function} in Mathematical Logic. \\

Connected with {\it quantifiers} in the simple languages ${\cal L}_{\cal S}$ associated with simple systems ${\cal S} ~=~
( S,~ ( P_i ~|~ i \in I ),~ ( f_j ~|~ j \in J ) )$ we can note the following very important {\it limitation}. According to 12 in
Definition 3.2.2, {\it variable symbols} $x$ can appear in the quantifications $( \forall~ x)$ that are allowed in compound
sentences. However, in view of the way the interpretation of compound sentences is done, see 4. in Definition 3.3.1,
such variable symbols can {\it only} be replaced with {\it constant symbols}. And as 5 in Definition 3.2.2 shows it,
constant symbols can only be associated with {\it elements} in the sets $S$ of the respective simple systems.
Therefore, such variable symbols {\it cannot} range over subsets of $S$, sets of subsets of $S$ and so on, or over
relations or functions on $S$, let alone, over sets of such relations or functions, and so on.
And yet, in usual mathematics, variable symbols which appear under quantification do range not only over the elements
of a fixed set. \\

This very important limitation of simple languages to what is called in Mathematical Logic a First Order Predicate
Calculus, is the reason why we shall have to go to Superstructures and the more evolved languages  associated with
them, languages which do no longer suffer from such a limitation. A price paid, however, for such an increase in
generality is the the corresponding full version of the Transfer Property is more difficult to prove, see chapter 5. \\

To illustrate the important and ubiquitous role played in usual mathematics by quantifications which cannot be
included in simple languages, let us formulate here the property of ${\bf R}$ of being Dedekind order complete. And for
simplicity, let us only formulate the part which says that every subset in ${\bf R}$ which is bounded form above, has an
upper bound in ${\bf R}$, namely

$$ \begin{array}{l}
           ( \forall~ A \subseteq {\bf R} )~ (~ (~ ( \exists~ M \in {\bf R} )~ ( \forall~ a \in A )~ (~ a ~\leq~ M ~)~) ~~\longrightarrow~~ \\ \\
           ~~~~~~ (~ ( \exists~ \overline{a} \in {\bf R} )~ (~ (~ ( \forall~ a \in A )~ (~ a ~\leq~ \overline{a} ~) ~) ~\wedge~ \\ \\
           ~~~~~~~~~~ (~ ( \forall~ b \in {\bf R} )~ (~ (~ ( \forall~ a \in A )~ (~ a ~\leq~ b ~) ~) ~~\longrightarrow~~
                                  (~ \overline{a} ~\leq~ b ~) ~) ~) ~) ~) ~)
      \end{array} $$

\medskip
Clearly, in this sentence which is true in ${\bf R}$, the first quantification, namely, $( \forall~ A \subseteq {\bf R} )$ is
over all possible {\it subsets} of ${\bf R}$, and this places the sentence outside of the simple language
${\cal L}_{\cal R}$. \\ \\

{\bf 4. The Transfer Property in its Simple Version} \\ \\

Before we give examples of transfer of sentences, let us specify here the enriched form of the simple system
$^*{\cal R}$ in (3.3.2), a form which we shall use from now on, and which corresponds to the simple system ${\cal R}$
in its enriched form in (3.3.1), namely

\bigskip
(3.4.1) \quad $ ^*{\cal R} ~=~ (~ ^*{\bf R},~ ^*{\cal P},~ ^*{\cal F} ~) $

\medskip
where, see (2.7.1), (2.7.4)

\bigskip
(3.4.2) \quad $ ^*{\cal P} ~=~ \{~ ^*P ~~|~~ P \in {\cal P} ~\},~~~~ ^*{\cal F} ~=~ \{~ ^*f ~~|~~ f \in {\cal F} ~\} $

\medskip
in other words, the relations in $^*{\cal R}$ are those obtained by transfer from relations in ${\cal R}$, and similarly, the
functions in $^*{\cal R}$ are those obtained by transfer from functions in ${\cal R}$. \\

Now in the case of our interest, namely, the two simple systems

$$ {\cal R} ~=~ (~ {\bf R},~ {\cal P},~ {\cal F} ~),~~~ ^*{\cal R} ~=~ (~ ^*{\bf R},~ ^*{\cal P},~ ^*{\cal F} ~) $$

\medskip
we shall associate to them, according to sections 2 and 3, their respective simple languages ${\cal L}_{\cal R}$ and
${\cal L}_{~^*{\cal R}}$. \\

We can now give the {\it rigorous} and {\it general} procedure to {\it transfer} simple sentences from ${\cal L}_{\cal R}$
into simple sentences in ${\cal L}_{~^*{\cal R}}$. \\
This procedure obviously goes parallel with the steps 5 - 12 in Definition 3.2.2 of simple languages, particularized this
time to the simple system ${\cal R} = (~ {\bf R},~ {\cal P},~ {\cal F} ~)$. \\ \\

{\bf Definition 3.4.1} \\

1. {\bf Transfer of constant symbols from ${\cal R}$}. If $\underline{r}$ is the name in ${\cal L}_{\cal R}$ of $r \in
{\bf R}$, then it will be kept as well as the name in ${\cal L}_{~^*{\cal R}}$ of $r =~ ^*r \in~ ^*{\bf R}$, see (2.1.5)
and (2.7.5). \\
Let us note that in view of Definitions 3.2.2 and 3.3.1, {\it variable symbols} are {\it not} subjected to being given
names or interpretation. Furthermore, they remain the same under transfer as well. \\

2. {\bf Transfer of relations from ${\cal R}$}. If $\underline{P}$ is the name in ${\cal L}_{\cal R}$ of the relation
$P$ in ${\bf R}$, then $\underline{^*P}$ will be the name in ${\cal L}_{~^*{\cal R}}$ of $^*P$, see (2.7.1). \\

3. {\bf Transfer of functions from ${\cal R}$}. If $\underline{f}$ is the name in ${\cal L}_{\cal R}$ of the function
$f$ in ${\bf R}$, then $\underline{^*f}$ will be the will be the name in ${\cal L}_{~^*{\cal R}}$ of $^*f$, see
(2.7.4). \\

It is useful to present in a diagram the above transfer of constant symbols, relations and functions, namely

\bigskip
\begin{math}
\setlength{\unitlength}{1cm}
\thinlines
\begin{picture}(15,4)

\put(0,2){${\cal R},$}
\put(0.7,2){$^*{\cal R}$}
\put(0,0.5){${\cal L}_{\cal R}$,}
\put(1,0.5){${\cal L}_{~^*{\cal R}}$}

\put(3,3.5){${\bf R}$}
\put(4.5,3.5){${^*{\bf R}}$}
\put(3,2){$r ~~~=~~~ ^*r$}
\put(3,0.5){$\underline{r} ~~~=~~~ \underline{^*r}$}
\put(3.1,1.8){\vector(0,-1){0.8}}
\put(4.8,1.8){\vector(0,-1){0.8}}

\put(6.5,3.5){${\bf R}$}
\put(8,3.5){${^*{\bf R}}$}
\put(6.5,2){$P$}
\put(7.1,2.1){\vector(1,0){0.6}}
\put(8,2){$^*P$}
\put(6.5,0.5){$\underline{P}$}
\put(7.1,0.6){\vector(1,0){0.6}}
\put(8,0.5){$\underline{^*P}$}
\put(6.6,1.8){\vector(0,-1){0.8}}
\put(8.3,1.8){\vector(0,-1){0.8}}

\put(10,3.5){${\bf R}$}
\put(11.5,3.5){${^*{\bf R}}$}
\put(10,2){$f$}
\put(10.6,2.1){\vector(1,0){0.6}}
\put(11.5,2){$^*f$}
\put(10,0.5){$\underline{f}$}
\put(10.6,0.6){\vector(1,0){0.6}}
\put(11.5,0.5){$\underline{^*f}$}
\put(10.1,1.8){\vector(0,-1){0.8}}
\put(11.8,1.8){\vector(0,-1){0.8}}

\end{picture}
\end{math}

Here we should note that, in fact, we have

$$ \underline{^*r} = \underline{( ^*r )} = r,~~~ \underline{^*P} = \underline{( ^*P )},~~~ \underline{^*f} = \underline{( ^*f )} $$

\medskip
this being the meaning of the above three commutative diagrams, while $^* ( \underline{r} ),~ ^*( \underline{P} ),~
^*( \underline{f} )$ do not yet make sense, since the operation of transfer $^*(~)$ has not yet been defined on
${\cal L}_{\cal R}$, this being only done below, startig immediately next. \\

4. {\bf Transfer of terms from ${\cal L}_{\cal R}$}. If $\tau$ is a variable or constant term in ${\cal L}_{\cal R}$, then its
transfer to ${\cal L}_{~^*{\cal R}}$ is given by $^*\tau = \tau$. Note that this extends 1 above to variable symbols. \\
If $\tau = \underline{f} ( \tau_1, ~~.~.~.~ , \tau_n )$ in ${\cal L}_{\cal R}$, then its transfer to ${\cal L}_{~^*{\cal R}}$ is given
by $^*\tau = \underline{^*f} ( ^*\tau_1, ~.~.~.~ , ^*\tau_n )$. \\

5. {\bf Transfer of atomic sentences}. If

$$ \Phi ~=~ \underline{P} < \tau_1, ~~.~.~.~ , \tau_n > $$

\medskip
is an atomic sentence in ${\cal L}_{\cal R}$, then its transfer to ${\cal L}_{~^*{\cal R}}$ is given by

$$ ^*\Phi ~=~ \underline{^*P} <~ ^*\tau_1, ~~.~.~.~ , ^*\tau_n > $$ \\

6. {\bf Transfer of compund sentences}. Finally, if

$$ \Phi ~=~  ( \forall~ x_1 ) ~.~.~.~ ( \forall~ x_n )~ \Big [ ~ \bigwedge_{1 \leq i \leq k}~ \underline{P_i} < \xi_i >
                                        ~\longrightarrow~ \bigwedge_{1 \leq j \leq l}~ \underline{Q_j} < \chi_j > ~ \Big ] $$

\medskip
is a compound sentence in ${\cal L}_{\cal R}$, then its transfer to ${\cal L}_{~^*{\cal R}}$ is given by

$$ ^*\Phi ~=~  ( \forall~ x_1 ) ~.~.~.~ ( \forall~ x_n )~ \Big [ ~ \bigwedge_{1 \leq i \leq k}~ \underline{^*P_i} <~ ^*\xi_i >
                                        ~\longrightarrow~ \bigwedge_{1 \leq j \leq l}~ \underline{^*Q_j} <~ ^*\chi_j > ~ \Big ] $$

\medskip
where $^*\xi_i = <~ ^*\tau^1_i, ~.~.~.~ , ^*\tau^{n_i}_i >$, if $\xi_i = < \tau^1_i, ~.~.~.~ , \tau^{n_i}_i >$, and similarly for
$^*\chi_j$. \\ \\

{\bf Convention Summary 3.4.1} \\

As an extension of Convention 2.7.1, in the case of the binary functions "~+~" and "~.~", as well as the binary relation
"$\leq$" on ${\bf R}$ and $^*{\bf R}$, we shall in all the situations use them just as they are written usually, that
is, without underlining them even when they represent names, or without placing a star "$^*$" to their left when they
are subjected to transfer. \\
As a further simplification, let us also recall Convention 3.2.1, according to which the better known and more often
used elements, relations and functions will not be underlined, when their names are employed. However, upon their
transfer, such relations and functions may still have a star "$^*$" placed to their left. \\

Consequently, the transfer of any simple sentence $\Phi$ from ${\cal L}_{\cal R}$ into a simple sentence $^*\Phi$ in
${\cal L}_{~^*{\cal R}}$ is done by placing a star "~*~" to the left of every relation and function, or name of
relation and function in $\Phi$, except for "$\leq$", "~+~", "~.~", and those specified above.

\hfill $\Box$ \\

And now, let us give two examples which illustrate the {\it ease} and {\it generality}, and at the same time, the {\it
nontriviality} of transfer of simple sentences. \\

The simple sentence in ${\cal L}_{\cal R}$

$$ ( \forall~ x ) ( \forall~ y ) \Big [~ (~ \underline{{\bf R}} < x > \wedge~ \underline{{\bf R}} < y > ~)
                                                                                                    ~\longrightarrow~ (~ x + y = y + x ~) ~\Big ] $$

\medskip
where as before, $\underline{{\bf R}} < x >$ denotes the name of the unary relation $x \in {\bf R}$, that is, of the property
of $x$ of being a usual real number, is true in the simple system ${\cal R}$, since it expresses the commutativity of
addition in that system. \\
According to the above Definition 3.4.1, the transfer of that simple sentence to ${\cal L}_{~^*{\cal R}}$ is given by

$$ ( \forall~ x ) ( \forall~ y ) \Big [~ (~ \underline{^*{\bf R}} < x >~\wedge~ \underline{^*{\bf R}} < y > ~)
                                                                                                    ~\longrightarrow~ (~ x + y = y + x ~) ~\Big ] $$

\medskip
where $\underline{^*{\bf R}}$ is the transfer of the unary relation $\underline{{\bf R}}$, therefore, according to (2.7.1), it
is easy to see that we have $\underline{^*{\bf R}} < s >$, for every $s \in~ ^*{\bf R}$. In this way, the above transferred
simple sentence is true in the simple system $^*{\cal R}$, as it expresses the commutativity of addition in that
system. \\

Let us turn now to the property of ${\bf R}$ of being Archimedean. This can be formulated as follows

$$ ( \exists~ u ) ( \forall~ x ) ( \exists~ n) \Big[~ (~ \underline{{\bf R_+}} < u >~ \wedge~ \underline{{\bf R_+}} < x > \wedge~
                                  \underline{{\bf N}} < n > ~) ~\longrightarrow~ (~ x \leq n . u ~) ~\Big] $$

\medskip
where $\underline{{\bf R_+}} < y >$ is the name of the unary relation $( y \in {\bf R} ) \wedge ( y > 0 )$, while
$\underline{{\bf N}} < n >$ is the name of the unary relation $n \in {\bf N}$, that is, of $n$ being a standard natural
number, and in particular, $n \geq 1$. Obviously, this sentence, although true in the simple system ${\cal R}$, is
nevertheless not a simple sentence in ${\cal L}_{\cal R}$, owing to the double presence of the existential quantifier
$\exists$. However, we can overcome that as follows. First, we can replace $u$ with the constant 1, and thus get in
${\cal R}$ the equivalent true sentence

$$ ( \forall~ x ) ( \exists~ n) \Big[~ (~  \underline{{\bf R_+}} < x > \wedge~
                                  \underline{{\bf N}} < n > ~) ~\longrightarrow~ (~ x \leq n . 1 ~) ~\Big] $$

\medskip
Then for the elimination of the existential quantifier left, we can use a Skolem function

$$ \sigma : ( 0, \infty ) ~~\longrightarrow~~ {\bf N} $$

\medskip
similar with the way done at the end of section 3. Thus we obtain in the simple system ${\cal R}$ the equivalent true
sentence

$$ ( \forall~ x ) \Big[~ (~  \underline{{\bf R_+}} < x >  ~) ~\longrightarrow~ (~ x ~\leq~ \underline{\sigma} ( x ) . 1 ~) ~\Big] $$

\medskip
which by now is obviously also a simple sentence in ${\cal L}_{\cal R}$. Then according to Definition 3.4.1, the transfer
to ${\cal L}_{^*{\cal R}}$ of this sentence is

$$ ( \forall~ x ) \Big[~ (~  \underline{^*{\bf R_+}} < x >  ~) ~\longrightarrow~
                                                      (~ x ~\leq~ \underline{^*{\sigma}} ( x ) . 1 ~) ~\Big] $$

\medskip
and as we can see from the Transfer Property next, this simple sentence in ${\cal L}_{^*{\cal R}}$ is true in the simple
system $^*{\cal R}$, since it is the transfer of a simple sentence in ${\cal L}_{\cal R}$ which is true in the simple
system ${\cal R}$. \\

However, it is important to note that this latest sentence is {\it not} about the Archimedean property of $^*{\bf R}$.
Indeed, it is easy to see that, according to (2.7.4), when we transfer the Skolem function $\sigma$ into the function
$^*\sigma$, the values of this latter function will no longer be in ${\bf N}$, but in the much larger $^*{\bf N}$. And this is
not what the Archimedean property is about, which in the inequality

$$ x ~\leq~ n . 1 $$

\medskip
does essentially require a standard natural number $n \in {\bf N}$. \\

This example shows clearly that, although transfer is easy and general, it is nevertheless not trivial. Indeed, we may
start with a true sentence in ${\bf R}$, transfer it to a true sentence in  $^*{\bf R}$, but the result of the transfer may no
longer present an interest from the point of view of ${\bf R}$. \\ \\

{\bf Transfer Theorem ( the simple version)} \\

If $\Phi$ is a simple sentence in ${\cal L}_{\cal R}$ and it is true in the simple system ${\cal R}$, then its transfer
$^*\Phi$ is a simple sentence in ${\cal L}_{^*{\cal R}}$ and it is true in $^*{\cal R}$.

\hfill $\Box$ \\

The proof of this theorem will be given in section 10. Meanwhile, we shall present a number of application of this
theorem which illustrate that {\it transfer} is one of the most powerful instruments in Nonstandard Analysis. \\

{\bf Convention 3.4.2} \\

For conveniece, instead of saying that a sentence $\Phi$ in the language ${\cal L}_{\cal R}$ is true in the simple
system ${\cal R}$, we shall often say that $\Phi$ is true ${\cal L}_{\cal R}$. The same we shall do with sentences in the
language ${\cal L}_{^*{\cal R}}$. \\

Of course, the {\it truth} of a sentence in a language as we defined it so far is a {\it semantic} concept, that is, it depends
on the specific {\it interpretation} of the resepctive language. However, in case that interpretation is clear from the
context, then for brevity, we shall not mention it when talking about the truth of a sentence, and instead, we shall only
mention the language to which the respective sentence belongs. \\

Here we can recall that in Mathematical Logic the {\it truth} of a sentence can also be a {\it syntactic} concept, that is,
depending only on the structure of the respective sentence, and independent of any particular interpretation, see Bell
\& Slomson, or Marker. \\ \\

{\bf Note 3.4.1} \\

It is {\it important} to note the following two facts. First, when stated in its full generality, see 5 in Definition 5.4.2, as
well as Theorem 5.5.2, the Transfer Property acts {\it both ways} between true sentences in ${\cal L}_{\cal R}$ and
${\cal L}_{^*{\cal R}}$, namely

\bigskip
(3.4.3) \quad $ \Phi~~ \mbox{true in}~ {\cal L}_{\cal R} ~~~\Longleftrightarrow~~~ ^*\Phi~~ \mbox{true in}~
                                                                                       {\cal L}_{^*{\cal R}} $

\medskip
Indeed, let us assume that we are given a sentence $\Phi$ in ${\cal L}_{\cal R}$ which is interpretable but false. And
let us assume for the moment - what is not actually the case, owing to the rudimentary structure of ${\cal L}_{\cal R}$
- that the negation non-$\Phi$ of this sentence is also in ${\cal L}_{\cal R}$. Then obviously non-$\Phi$ is
interpretable and true in ${\cal L}_{\cal R}$. Therefore, according to the Transfer Property, its transfer $^*($non-$\Phi) =$
non-$^*\Phi$ is true in ${\cal L}_{^*{\cal R}}$. \\
It follows that in case we have a sentence $\Phi$ which is interpretable in ${\cal L}_{\cal R}$,  and its transfer
$^*\Phi$ is true in ${\cal L}_{^*{\cal R}}$, then $\Phi$ itself must be true in ${\cal L}_{\cal R}$. \\
However, in order to obtain such a two way transfer (3.4.3), we must go to Superstructures and their languages,
see chapter 5. Indeed, if $\Phi$ is a {\it compund sentence} in ${\cal L}_{\cal R}$, see 12 in Definition 3.2.2, then clearly
non-$\Phi$ does {\it not} belong to ${\cal L}_{\cal R}$, since it contains the quantifier $\exists$ which is not allowed in
${\cal L}_{\cal R}$. Therefore, in the above Transfer Theorem, instead of the full generality of transfer given in (3.4.3)
and which, as mentioned, will be obtained in chapter 5, we can only have

$$ \Phi~~ \mbox{true in}~ {\cal L}_{\cal R} ~~~\Longrightarrow~~~ ^*\Phi~~ \mbox{true in}~ {\cal L}_{^*{\cal R}} $$

\medskip
Second, even in the case of the full transfer property (3.4.3), it is essential to note that, when studying ${\bf R}$, we {\it
start} with sentences $\Phi$ in  ${\cal L}_{\cal R}$. Then by {\it transfer}, we associate with them sentences $^*\Phi$ in
${\cal L}_{^*{\cal R}}$. And we {\it cannot} do it the other way round, namely, start with arbitrary sentences in ${\cal
L}_{^*{\cal R}}$, and try to associate with them sentences in${\cal L}_{\cal R}$, and then still hope for relating their truth
to one another. Indeed, in ${\cal L}_{^*{\cal R}}$ there are far {\it more} sentences than those which correspond through
transfer to sentences in ${\cal L}_{\cal R}$. \\ \\

{\bf 5. Several General Results} \\ \\

Let us start with a few of the more simple and immediate applications of transfer. In this respect, a review of, and a
certain addition to the results in section 7, in chapter 2, is appropriate. \\
Before doing so, however, let us introduce the following \\ \\

{\bf Convention 3.5.1} \\

The sentence $\Sigma$

$$ ( \forall~ x_1 ) ~.~.~.~ ( \forall~ x_n )~ \Big [ ~ \bigwedge_{1 \leq i \leq k}~ \underline{P_i} < \xi_i >
                                 ~\longleftrightarrow~ \bigwedge_{1 \leq j \leq l}~ \underline{Q_j} < \chi_j > ~ \Big ] $$

which in view of Definition 3.2.2 is {\it not} in the simple language ${\cal L}_{\cal R}$, will denote the following {\it two}
sentences

$$ ( \forall~ x_1 ) ~.~.~.~ ( \forall~ x_n )~ \Big [ ~ \bigwedge_{1 \leq i \leq k}~ \underline{P_i} < \xi_i > ~\longrightarrow~
                                                                     ~ \bigwedge_{1 \leq j \leq l}~ \underline{Q_j} < \chi_j > ~ \Big ] $$

$$ ( \forall~ x_1 ) ~.~.~.~ ( \forall~ x_n )~ \Big [ ~ \bigwedge_{1 \leq j \leq l}~ \underline{Q_j} < \chi_j > ~\longrightarrow~
                                                                    ~ \bigwedge_{1 \leq i \leq k}~ \underline{P_i} < \xi_i > ~ \Big ] $$

\medskip
both of which are obviously in the simple language ${\cal L}_{\cal R}$. In this way, we shall say that the sentence
$\Sigma$ is also in the simple language ${\cal L}_{\cal R}$. \\
Clearly, a similar notation can be used in the simple language ${\cal L}_{~^*{\cal R}}$ as well. \\

Furthermore, the sentence $\Sigma$ in the simple language ${\cal L}_{\cal R}$ can be transferred to a sentence
$^*\Sigma$ in the simple language ${\cal L}_{^*{\cal R}}$, by transferring each of the two sentences which $\Sigma$
represents. And then in view of the Transfer Property, if $\Sigma$ is true in ${\cal L}_{\cal R}$, so will be $^*\Sigma$ in
${\cal L}_{^*{\cal R}}$. \\ \\

{\bf Proposition 3.5.1} \\

If $A_1, ~.~.~.~ , A_n \subseteq {\bf R}$, then

\bigskip
(3.5.1) \quad $ \begin{array}{l}
                                  ^*( A_1 \bigcup ~.~.~.~ \bigcup A_n ) ~=~ ^*A_1 \bigcup ~.~.~.~ \bigcup~ ^*A_n \\ \\
                                   ^*( A_1 \bigcap ~.~.~.~ \bigcap A_n ) ~=~ ^*A_1 \bigcap ~.~.~.~ \bigcap~ ^*A_n
                          \end{array} $ \\

{\bf Proof} \\

Let us denote $A = A_1 \bigcap ~.~.~.~ \bigcap A_n$, then we have the unary relations in ${\cal R}$

$$ A < x >,~ A_1 < x >, ~.~.~.~ , A_n < x > $$

\medskip
which express respectively the properties $x \in A,~ x \in A_1, ~.~.~.~ , x \in A_n$. Obviously, we have in
${\cal L}_{\cal R}$ the following true sentence

$$ ( \forall~ x )~ \Big [~ \underline{A} < x >  ~\longleftrightarrow~ (~ \underline{A_1} < x > ~\bigwedge~.~.~.~
                                                                                      \bigwedge~ \underline{A_n} < x > ~) ~\big ] $$

\medskip
which according to the above Convention 3.5.1, represents two sentences in ${\cal L}_{\cal R}$. Therefore, by transfer
we will have in ${\cal L}_{~^*{\cal R}}$ two true sentences, which in view of the same convention can be written as

$$ ( \forall~ x )~ \Big [~ \underline{^*A} < x >  ~\longleftrightarrow~ (~ \underline{^*A_1} < x > ~\bigwedge~.~.~.~
                                                                                      \bigwedge~ \underline{^*A_n} < x > ~) ~\big ] $$

\medskip
thus the second relation in (3.5.1) is proved. \\

For the first relation in (3.5.1) we use the De Morgan formula

$$ x \in A_1 \bigcup ~.~.~.~ \bigcup A_n ~~~\Longleftrightarrow~~~ x \in (~ A_1~^\prime ~\bigcap ~.~.~.~
                                                                                                                    \bigcap A_n~^\prime ~)~^\prime $$

\medskip
where $B~^\prime$ denotes the complementary in ${\bf R}$ of the subset $B \subseteq {\bf R}$. In this way, we reduce
the problem to what has just been proved.

\hfill $\Box$ \\ \\

The above direct and simple proof of (3.5.1) does {\it not} make use of any explicit reference to ultrafilters, maximal
ideals, or to the ( a. e. ) relations. Thus, when compared with the earlier proof given to the same result in (2.7.15), one
can see a rather typical instance of the {\it power} of the Transfer Property. \\ \\

{\bf Proposition 3.5.2} \\

If $P$ is an n-ary relation on ${\bf R}$ and $\chi_P$ is its characteristic function, then \\

1.~~ $^*P$ is an extension of $P$. \\

2.~~ $^*( \chi_P ) ~=~ \chi_{(~ ^*P )}$, and \\

3.~~ $^*( P~^\prime ) ~=~ ( ^*P )~^\prime$ \\

{\bf Proof} \\

1. Let $r_1, ~.~.~.~ , r_n \in {\bf R}$, then according to 3 in Definition 3.3.1, $P < r_1, ~.~.~.~ , r_n >$, if and only if the
atomic sentence $\underline{P} < r_1, ~.~.~.~ , r_n >$ is true in ${\cal L}_{\cal R}$. However, if this sentence is true in
${\cal L}_{\cal R}$ then by transfer it follows that the sentence $\underline{^*P} < r_1, ~.~.~.~ , r_n >$ is true in
${\cal L}_{~^*{\cal R}}$, and this completes the proof of 1. \\

2. In view of (2.7.3) and (2.7.13), for $r_1, ~.~.~.~ , r_n \in {\bf R}$, we have

$$ \chi_P~ ( r_1, ~.~.~.~ , r_n ) ~=~ 1 ~~~\Longleftrightarrow~~ P <  r_1, ~.~.~.~ , r_n > $$

$$ \chi_P~ ( r_1, ~.~.~.~ , r_n ) ~=~ 0 ~~~\Longleftrightarrow~~ P~^\prime <  r_1, ~.~.~.~ , r_n > $$

\medskip
where as before $P~^\prime$ denotes the complementary of $P$. Thus the sentences

$$ ( \forall~ x_1 ) ~.~.~.~ ( \forall~ x_n )~ \Big [~  \underline{\chi_P}~ ( x_1, ~.~.~.~ , x_n ) ~=~ 1 ~~~\longleftrightarrow~~
                                                                                                                  \underline{P} <  x_1, ~.~.~.~ , x_n > ~\Big ] $$

$$ ( \forall~ x_1 ) ~.~.~.~ ( \forall~ x_n )~ \Big [~  \underline{\chi_P}~ ( x_1, ~.~.~.~ , x_n ) ~=~ 0 ~~~\longleftrightarrow~~
                                                                                                           \underline{P~^\prime} <  x_1, ~.~.~.~ , x_n > ~\Big ] $$

\medskip
are true in ${\cal L}_{\cal R}$. It follows by transfer that the sentences

$$ ( \forall~ x_1 ) ~.~.~.~ ( \forall~ x_n )~ \Big [~  \underline{^*( \chi_P} )~ ( x_1, ~.~.~.~ , x_n ) ~=~ 1 ~~~\longleftrightarrow~~
                                                                                                                    \underline{^*P} <  x_1, ~.~.~.~ , x_n > ~\Big ] $$

$$ ( \forall~ x_1 ) ~.~.~.~ ( \forall~ x_n )~ \Big [~  \underline{^*( \chi_P} )~ ( x_1, ~.~.~.~ , x_n ) ~=~ 0 ~~~\longleftrightarrow~~
                                                                                                      \underline{^*( P~^\prime )} <  x_1, ~.~.~.~ , x_n > ~\Big ] $$

\medskip
are true in ${\cal L}_{^*{\cal R}}$. Thus equivalently, for $s^1, ~.~.~.~ , s^n \in~ ^*{\bf R}$, we also have in $~^*{\cal R}$

$$ ^*( \chi_P )~ ( s^1, ~.~.~.~ , s^n ) ~=~ 1 ~~~\Longleftrightarrow~~ ^*P <  s^1, ~.~.~.~ , s^n > $$

$$ ^*( \chi_P )~ ( s^1, ~.~.~.~ , s^n ) ~=~ 0 ~~~\Longleftrightarrow~~ (~ ^*P )~^\prime <  s^1, ~.~.~.~ , s^n > $$

\medskip
On the other hand, when applied to $^*P$, the corresponding version of (2.7.13) gives in $^*{\cal R}$ for $s^1, ~.~.~.~ ,
s^n \in~ ^*{\bf R}$

$$ \chi_{(~ ^*P )}~ ( s^1, ~.~.~.~ , s^n ) ~=~ 1 ~~~\Longleftrightarrow~~ ^*P <  s^1, ~.~.~.~ , s^n > $$

$$ \chi_{(~ ^*P )}~ ( s^1, ~.~.~.~ , s^n ) ~=~ 0 ~~~\Longleftrightarrow~~ ( ^*P )~^\prime <  s^1, ~.~.~.~ , s^n > $$

\medskip
and these two relations, together with the previous two ones, complete the proof of 2. \\

3. In view of (2.7.13), we have the true sentence in  ${\cal L}_{\cal R}$

$$ ( \forall~ x_1 ) ~.~.~.~ ( \forall~ x_n )~ \Big [~  \underline{\chi_P}~ ( x_1, ~.~.~.~ , x_n ) ~=~ 0
                                                                      ~\longleftrightarrow~ \underline{P~^\prime}~ < x_1, ~.~.~.~ , x_n > ~\Big ] $$

\medskip
which by transfer becomes the true sentence in ${\cal L}_{^*{\cal R}}$

$$ ( \forall~ x_1 ) ~.~.~.~ ( \forall~ x_n )~ \Big [~  \underline{^*( \chi_P )}~ ( x_1, ~.~.~.~ , x_n ) ~=~ 0
                                                       ~\longleftrightarrow~ \underline{^*( P~^\prime ~)}~ < x_1, ~.~.~.~ , x_n > ~\Big ] $$

\medskip
In other words, for $s^1, ~.~.~.~ , s^n \in~ ^*{\bf R}$ we have

$$ ^*( \chi_P )~ ( s^1, ~.~.~.~ , s^n ) ~=~ 0 ~~~\Longleftrightarrow~~~ ^*( P~^\prime ~) < s^1, ~.~.~.~ , s^n > $$

\medskip
thus in view of 2, it follows that

$$ \chi_{(~^*P )}~ ( s^1, ~.~.~.~ , s^n ) ~=~ 0 ~~~\Longleftrightarrow~~~ ^*( P~^\prime ~) < s^1, ~.~.~.~ , s^n > $$

\medskip
and then clearly

$$ \chi_{((~^*P )~^\prime~)}~ ( s^1, ~.~.~.~ , s^n ) ~=~ 1 ~~~\Longleftrightarrow~~~ ^*( P~^\prime ~) < s^1, ~.~.~.~ , s^n > $$

\medskip
But the corresponding version of (2.7.13) gives

$$ \chi_{((~^*P )~^\prime~)}~ ( s^1, ~.~.~.~ , s^n ) ~=~ 1 ~~~\Longleftrightarrow~~~ ( ^*P) ~^\prime ~ < s^1, ~.~.~.~ , s^n > $$

\medskip
thus the last two relations end the proof of 3.

\hfill $\Box$ \\ \\

{\bf Proposition 3.5.3} \\

If $f$ is a function of n variables defined on a domain in ${\bf R}^n$ and with values in ${\bf R}$, see (2.7.3), (2.7.4), then
$^*f$ is a function of n variables defined on a domain in $(~ ^*{\bf R})^n$ and values in $^*{\bf R}$, furthermore \\

1.~~ $^*f$ is an extension of $f$ \\

2.~~ $domain~ ^*f ~=~ ^*( domain~ f),~~~~~ range~ ^*f ~=~ ^*( range~ f )$ \\

{\bf Proof} \\

1. It follows from (2.7.3) and 1 in Proposition 3.5.2. \\

2. Let $P$ be the (n + 1)-ary relation on ${\bf R}$ associated to $f$, according to  (2.7.3). If $r_1, ~.~.~.~ , r_n \in {\bf R}$,
then

$$ ( r_1, ~.~.~.~ , r_n ) \in domain~ f ~~~\Longleftrightarrow~~ P < r_1, ~.~.~.~ , r_n, f ( r_1, ~.~.~.~ , r_n ) > $$

\medskip
thus we have in  ${\cal L}_{\cal R}$ the true sentence

$$ \begin{array}{l}
                   ( \forall~ x_1 ) ~.~.~.~ ( \forall~ x_n )~ \Big [~ \underline{domain~ f}  < x_1, ~.~.~.~ , x_n >
                                                       ~\longleftrightarrow~ \\ \\
                               ~~~~~~~~~~~~~~~~~~\longleftrightarrow~ \underline{P~} < x_1, ~.~.~.~ , x_n, f ( x_1, ~.~.~.~ , x_n )> ~\Big ]
      \end{array} $$

\medskip
which through transfer becomes the true sentence in ${\cal L}_{^*{\cal R}}$

$$ \begin{array}{l}
                   ( \forall~ x_1 ) ~.~.~.~ ( \forall~ x_n )~ \Big [~ \underline{^*( domain~ f)}  < x_1, ~.~.~.~ , x_n >
                                                       ~\longleftrightarrow~ \\ \\
                               ~~~~~~~~~~~~~~~~~~\longleftrightarrow~ \underline{^*P} < x_1, ~.~.~.~ , x_n, ^*f ( x_1, ~.~.~.~ , x_n )> ~\Big ]
      \end{array} $$

\medskip
This obviously means that, for $s^1, ~.~.~.~ , s^n \in~ ^*{\bf R}$, the following is true in $^*{\cal R}$

$$ ( s^1, ~.~.~.~ , s^n ) \in~ ^*( domain~ f ) ~~~\Longleftrightarrow~~ ^*P < s^1, ~.~.~.~ , s^n,~ ^*f ( s^1, ~.~.~.~ , s^n ) > $$

\medskip
On the other hand, in view of (2.7.1) - (2.7.4), we have

$$ ( s^1, ~.~.~.~ , s^n ) \in domain~ ^*f ~~~\Longleftrightarrow~~ ^*P < s^1, ~.~.~.~ , s^n,~ ^*f ( s^1, ~.~.~.~ , s^n ) > $$

\medskip
and thus indeed $^*( domain~ f)~=~ domain~ ^*f$. \\

We turn now to the property of the respective function ranges, and we start with the particular case when n = 1. Let us
define a Skolem function $\psi : range ~f ~\longrightarrow~ domain ~f$, such that $f ( \psi ( b ) ) = b$, for
$b \in range ~f$, then we have in ${\cal L}_{\cal R}$ the two true sentences

$$ ( \forall~ x )~( \forall~ y ) ~\Big [~ \underline{f} ( x ) ~=~ y ~\longrightarrow~ \underline{( range ~f )} < y > ~\Big ] $$

$$ ( \forall~ y ) ~\Big [~ \underline{( range ~f )} < y > ~\longrightarrow~
                                                            \underline{f} (~ \underline{\psi} ( y )~ ) ~=~ y ~\Big ] $$

\medskip
Thus by transfer we have in  ${\cal L}_{^*{\cal R}}$ the two true sentences

$$ ( \forall~ x )~( \forall~ y ) ~\Big [~ \underline{^*f} ( x ) ~=~ y ~\longrightarrow~
                                                          \underline{(^*( range ~f ))} < y > ~\Big ] $$

$$ ( \forall~ y ) ~\Big [~ \underline{(^*( range ~f ))} < y > ~\longrightarrow~
                                                                       \underline{^*f} (~ \underline{^*\psi} ( y )~) ~=~ y ~\Big ] $$

\medskip
Now the first of them gives in $^*{\cal R}$ the inclusion $range~ ^*f ~\subseteq~ ^*( range~ f )$, while the
second gives the converse inclusion. \\

For arbitrary n $>$ 1, it is easy to extend the above argument.

\hfill $\Box$ \\ \\

The proofs of above propositions can give a good enough idea about the power of transfer, and the ways to use it. The
next results can be obtained in similar ways, and we give them here without proofs. \\ \\

{\bf Proposition 3.5.4} \\

We have $~~~^*\phi ~=~ \phi$. \\

Further, given an infinite family of subsets $A_i \subseteq {\bf R}^n$, with $i \in I$, then

$$ \bigcup_{i \in I}~ ^*A_i ~\subseteq~ ^*(~ \bigcup_{i \in I}~ A_i ~) ~~~~~~~ \bigcap_{i \in I}~ ^*A_i ~\supseteq~
                                                                     ^*(~ \bigcap_{i \in I}~ A_i ~) $$ \\ \\

{\bf Proposition 3.5.5} \\

Given two functions $f$ and $g$ of n variables on ${\bf R}$, then on \\ $( domain~ ^*f  ) ~\bigcap~
( domain~ ^*g )$ we have

$$ ^*(~ f ~+~ g ~) ~=~ ^*f ~+~ ^*g ~~~~~~~~~~~~ ^*(~ f ~.~ g ~) ~=~ ^*f ~.~ ^*g $$

Also, on $domain~ ^*f$ we have

$$ ^* |~ f ~| ~=~ |~ ^*f ~| $$ \\ \\

{\bf 6. The Local and Global Structure of $^*{\bf R}$} \\ \\

As we shall see in this section, going from ${\bf R}$ to $^*{\bf R}$ involves {\it two expansions}, namely, a {\it local}
one, at each point $x \in {\bf R}$, as well as a {\it global} one, at both infinite ends of ${\bf R}$. \\

A particularly astute way to illustrate it, see Keisler [2], is by saying that, when we want to go from ${\bf R}$ to
$^*{\bf R}$, we need {\it two} instruments which can give us views outside of ${\bf R}$ and into the not yet seen, and not
yet even known worlds situated beyond the confines of the standard realm of ${\bf R}$. Namely, we need a {\it
microscope} in order to see the {\it monads} which give the new and additional local structure in $^*{\bf R}$, and we
also need a {\it telescope} for being able to look at {\it galaxies} giving the new and additional global structure in
$^*{\bf R}$. \\

In order to understand the structure of $^*{\bf R}$, the following three things are therefore enough : \\

- to keep in mind that $^*{\bf R}$ is a {\it totally} ordered field, \\

- to understand how the {\it monads} create the {\it local} structure of $^*{\bf R}$, and in this respect, it is enough to
understand how they create the local structure of the {\it galaxy} of $0 \in {\bf R}$, which can be seen as the central
galaxy, \\

and at last, \\

- to understand how by uncountably many translations to the right and left, the galaxy of $0 \in {\bf R}$, that is, the
central galaxy, creates the {\it global} structure of $^*{\bf R}$. \\

And now, let us make the above mathematically precise. \\ \\

{\bf Definition 3.6.1} \\

A number $s \in~ ^*{\bf R}$ is called {\it infinitesimal}, if and only if $| s | \leq r$, for every $r \in {\bf R},~ r > 0$. We denote by
$mon ( 0 )$ the set of infinitesimal numbers, and call that set the {\it monad} of $0 \in {\bf R}$. \\

For $s \in~ ^*{\bf R}$ it will be convenient to denote $mon ( s ) = s + mon ( 0 )$. \\

A number $s \in~ ^*{\bf R}$ is called {\it finite}, if and only if $| s | \leq r$, for some $r \in {\bf R},~ r > 0$. We denote by $Gal ( 0 )$
the set of all finite numbers , and call that set the {\it galaxy} of $0 \in {\bf R}$. \\

For $s \in~ ^*{\bf R}$ it will be convenient to denote $Gal ( s ) = s + Gal ( 0 )$. \\

A number $s \in~ ^*{\bf R}$ is called {\it infinite}, if and only if $| s | \geq r$, for every $r \in {\bf R},~ r > 0$.

\hfill $\Box$ \\ \\

It is easy to see that

\bigskip
(3.6.1) \quad $ \epsilon ~=~ [~ < 1, 1/2, 1/3, ~.~.~.~ , 1/n, ~ .~.~.~ > ~] \in mon ( 0 ) \setminus {\bf R} $

\medskip
is a nonzero and positive infinitesimal, while

\bigskip
(3.6.2) \quad $ \omega ~=~ 1 / \epsilon ~=~ [~ < 1, 2, 3, ~.~.~.~ , n, ~.~.~.~ > ~] \in ~^*{\bf R} \setminus Gal  ( 0 ) $

\medskip
is a positive infinite number. \\

However, let us note that in the above, we may have gone somewhat ahead of ourselves. Indeed, in Definition 3.6.1
use is made of the absolute value of arbitrary elements of $^*{\bf R}$. \\
Now the standard absolute value is of course a function $| ~.~ | : {\bf R} \longrightarrow {\bf R}$, and as such, it can be
transferred to a function $^*| ~.~ | : ~^*{\bf R} \longrightarrow ~^*{\bf R}$, according to Proposition 3.5.3. And then, we have
to see how, through this transfer, the usual properties of the standard absolute value function are preserved. This will
be done in Lemma 3.6.1, related to the proof of Proposition 3.6.1 below. So far, that is, in Definition 3.6.1 and with the
two nonstandard numbers $\epsilon$ and $\omega$ mentioned above, we only used the existence of the transferred
absolute value function $^*| ~.~ |$, and did not refer to any of its properties. \\

By the way, in Definition 3.6.1, we used the customary notation $| ~.~ |$ for the transferred absolute value function
$^*| ~.~ |$. This is in line with the earlier simplifying notational conventions, and we shall keep to it from here on, as
well. \\ \\

{\bf Proposition 3.6.1} \\

We have \\

1.~~ $mon ( 0 ) ~\bigcap~ {\bf R} ~=~ \{~ 0 ~\},~~~~ mon ( 0 ) ~\bigcup~ {\bf R}~\subset~ Gal ( 0 ) ~\subset~ ^*{\bf R}$ and
\hspace*{0.7cm} $^*{\bf R} ~\setminus~ Gal ( 0 )$~~ is the set of infinite numbers in $^*{\bf R}$, furthermore, \\
\hspace*{0.7cm} the mapping

$$ (~ ^*{\bf R} ~\setminus~Gal ( 0 ) ~) \ni s ~\longmapsto~ 1 / s \in (~ mon ( 0 ) \setminus \{~ 0 ~\} ~) $$

\medskip
\hspace*{0.7cm} is a {\it bijection}. \\

2.~~ $mon ( 0 )$~ and ~$Gal ( 0 )$~ are subrings in $^*{\bf R}$, and also algebras over  ${\bf R}$. \\

3.~~ $mon ( 0 )$~ is an ideal in~ $Gal ( 0 )$. \\

4.~~ If $x,~ y \in {\bf R}$ and $x \neq y$, then

$$ (~ x ~+~ mon ( 0 ) ~) ~\bigcap~ (~ y ~+~ mon ( 0 ) ~) ~=~ \phi $$

\medskip
5.~~ We have the representation given by a union of pairwise disjoint \\
\hspace*{0.7cm} sets

\bigskip
(3.6.3) \quad $ Gal ( 0 ) ~=~ \bigcup_{~r \in {\bf R}}~ (~ r ~+~ mon ( 0 ) ~)$ \\

\medskip
6.~~ There exists a {\it unique surjective algebra homomorphism}~ \\
\hspace*{0.7cm} $st : Gal ( 0 ) \longrightarrow {\bf R}$, called the {\it standard part} mapping, such that

$$ s - st ( s ) \in mon ( 0 ),~~~ s \in Gal ( 0 ) $$

\medskip
\hspace*{0.7cm} In particular, for $s \in Gal ( 0 )$, we have

$$ st ( s ) ~=~ s ~~~\Longleftrightarrow~~ s \in {\bf R} $$

$$ st ( s ) ~=~ 0 ~~~\Longleftrightarrow~~ s \in mon ( 0 ) $$

\medskip
\hspace*{0.7cm} Furthermore, for $s, t \in Gal ( 0 )$, we have

$$ mon ( s ) ~=~ mon ( t ) ~~~\Longleftrightarrow~~~ st ( s ) ~=~ st ( t ) $$

\medskip
7.~~ $Gal ( 0 ) ~=~ mon ( 0 ) ~+~ {\bf R}$, and the mapping

$$ {\bf R} \ni r ~~\longmapsto~~ r + mon ( 0 ) \in Gal ( 0 ) / mon ( 0 ) $$

\medskip
\hspace*{0.7cm} is a {\it ring isomorphism}, thus~ ${\bf R}$~ and ~$Gal ( 0 ) / mon ( 0 )$~ are \\
\hspace*{0.7cm} {\it isomorphic} rings, and thus {\it isomorphic} fields as well. \\

8.~~ Let be given any Hausdorff topology on~ $Gal ( 0 )$~ in which all open \\
\hspace*{0.7cm} intervals

$$ ( a, b ) ~=~ \{~ s \in ~^*{\bf R} ~~|~~ a < s < b ~\} $$

\medskip
\hspace*{0.7cm} with $a,~ b \in Gal ( 0 ),~ a <  b$, are open subsets. Then such a \\
\hspace*{0.7cm} topology, when restricted to ${\bf R}$, gives the {\it discrete} topology \\
\hspace*{0.7cm} on ${\bf R}$. \\ \\

{\bf Note 3.6.1} \\

The use of the term {\it monad} is inspired by Leibniz who first employed infinitesimals when, parallel with and
independently of Newton, started the development of Calculus in the late 1600s. \\

The fact that

\bigskip
(3.6.4) \quad $ mon ( 0 ) ~\bigcap~ {\bf R} ~=~ \{~ 0 ~\} $

\medskip
which means that zero is the {\it only} standard real which is infinitesimal, created, starting with Leibniz, all sort of
difficulties when, prior to the creation of modern Nonstandard Analysis, one tried to deal with infinitesimals. \\
Similarly, the fact that

\bigskip
(3.6.5) \quad $ ^*{\bf R} ~\setminus~ Gal ( 0 ) $

\medskip
is the set of infinite numbers in $^*{\bf R}$, thus there are {\it no} standard reals which are infinite, brought with it
difficulties when dealing with infinitely large numbers, prior to the modern theory of Nonstandard Analysis. \\

The {\it local} structure of $Gal ( 0 )$ is presented in (3.6.3), and as seen in (3.6.9), that will give an understanding
of the local structure of $^*{\bf R}$ as well. \\

The {\it order reversing bijection}

\bigskip
(3.6.6) \quad $ (~ ^*{\bf R} ~\setminus~Gal ( 0 ) ~) \ni s ~\longmapsto~ 1 / s \in (~ mon ( 0 ) \setminus \{~ 0 ~\} ~) $

\medskip
is fundamental in understanding the {\it connection} between the local and global structure of $^*{\bf R}$. Indeed, it
shows that the local structure of $^*{\bf R}$ {\it mirrors} its global structure, and vice versa. \\
It also establishes the link between Keisler's microscope and telescope, the former letting us see into $mon ( 0 )$,
while the latter allowing us to look out into $^*{\bf R} ~\setminus~Gal ( 0 )$. \\

The properties 1 - 7 express the {\it local} structure of $^*{\bf R}$, more precisely, within its part given by $Gal ( 0 )$,
that is, the structure of the finite numbers. This local structure consists of a sort of infinitesimal neighbourhood around
every standard $x \in {\bf R}$ given by the translate ~$x ~+~ mon ( 0 )$~ of ~$mon ( 0 )$. And in view of (3.6.4), this local
structure cannot be seen from ${\bf R}$, thus we have to use Keisler's microscope. \\

Property 8 illustrates the fact that topological type consideration on $^*{\bf R}$ need to be pursued rather carefully. \\ \\

{\bf Proof} of Proposition 3.6.1. \\

1 - 3 follow easily from Definition 3.6.1 and Lemma 3.6.1. \\
4. It follows from  $mon ( 0 ) ~\bigcap~ {\bf R} ~=~ \{~ 0 ~\}$ at 1. \\
5. It is an immediate consequence of 4. \\
6. It results form 5. \\
7. It is implied by 6. \\
8. Let us take any $r \in {\bf R}$. Then by the hypothesis, the open interval $( r - \epsilon, r + \epsilon )$, see (3.6.1), is
a neighbourhood of $r$ in $^*{\bf R}$. However, in view of 4, when restricted to ${\bf R}$, this interval reduces to the
point $r$ alone. \\ \\

{\bf Lemma 3.6.1} \\

The nonstandard absolute value function $| ~.~ | : ~^*{\bf R} ~\longrightarrow~ ^*{\bf R}$, which is the transfer of the usual
absolute value function  $| ~.~ | : {\bf R} ~\longrightarrow~ {\bf R}$, has the following properties : \\

1.~~~ $|~ s ~| ~\geq~ 0,~~~  s \in~ ^*{\bf R}$ \\

2.~~~  For $s \in~ ^*{\bf R}$~ we have \\

~~~~~~ $|~ s ~| ~=~ 0 ~~~\Longleftrightarrow~~~ s ~=~ 0$ \\

~~~~~~ $|~ s ~| ~=~ s ~~~\Longleftrightarrow~~~ s ~\geq~ 0$ \\

~~~~~~ $|~ s ~| ~=~ - s ~~~\Longleftrightarrow~~~ s ~\leq~ 0$ \\

3.~~~ $|~ s + t ~| ~\leq~ |~ s ~| + |~ t ~|,~~~ |~ s . t ~| ~=~ |~ s ~| . |~ t ~|,~~~~s,~ t \in~ ^*{\bf R}$ \\

{\bf Proof} \\

Obviously 1 follows form 2. In order to prove 2, we recall that the sentence

$$ ( \forall~ r ) ~\Big [~ |~ r ~| ~=~ 0 ~~\longleftrightarrow~~ r ~=~ 0 ~\Big ] $$

\medskip
is true in ${\cal L}_{\cal R}$. Then through transfer, we obtain the first equivalence in 2. The other two equivalences in 2
result in a similar way. The same goes for 3.

\hfill $\Box$ \\ \\

For a better understanding of the structure of $^*{\bf R}$ it is useful to introduce the following two equivalence relations
on it, see 2 in Proposition 3.6.1 \\ \\

{\bf Definition 3.6.2} \\

Let $s,~ t \in ~^*{\bf R}$. We say that $s$ and $t$ are {\it infinitesimally close}, if and only if $t - s \in mon ( 0 )$, in which case we
use the notation $s \approx t$. \\
And we say that $s$ and $t$ are {\it finitely close}, if and only if $t - s \in Gal ( 0 )$, and then we denote $s \sim t$.

\hfill $\Box$ \\ \\

Clearly

\bigskip
(3.6.7) \quad $ \begin{array}{l}
                                         mon ( 0 ) ~=~ \{~ s \in ~^*{\bf R} ~~|~~ s \approx 0 ~\} \\ \\
                                         Gal ( 0 ) ~=~ \{~ s  \in ~^*{\bf R} ~~|~~ s \sim 0 ~\}
                        \end{array} $

\medskip
Also ~$s \approx t ~~\Longrightarrow~~ s \sim t$,~ for $s,~ t \in ~^*{\bf R}$,~ see 1 in Proposition 3.6.1. \\

In view of 2 in the same proposition, it follows that ~$\approx$~ and ~$\sim$~ are compatible with the ring structure
of $^*{\bf R}$. \\
Namely, if $s,~ t,~ u,~ v \in ~^*{\bf R}$, with $s \approx t$ and $u \approx v$, then $s + u \approx t + v$, $s - u \approx t -
v$, and $s . u \approx t . v$. \\
Furthermore, if $u \not\approx 0$, thus also $v \not\approx 0$, then $s / u \approx t / v$. \\
And we have the similar properties for~ $\sim$. \\

Here, let us also note that one can add to the properties of the standard mapping $st$ in 6 in Proposition 3.6.1, the
following. If $s,~ t \in Gal ( 0 )$, then

\bigskip
(3.6.8) \quad $ \begin{array}{l} st ( t ) \neq 0 ~~\Longrightarrow~~ st ( s / t ) ~=~ st ( s ) / st ( t ) \\ \\
                                                  s ~\leq~ t ~~\Longrightarrow~~ st ( s ) ~\leq~ st ( t )
                         \end{array} $

\medskip
In view of the above, and especially in view of the properties of~ $mon ( 0 )$, and of the restriction to $Gal ( 0 )$ of the
equivalence relation $\approx$, the structure of $Gal ( 0 )$, and in particular, the relationship between ${\bf R}$ and
$Gal ( 0 )$ is by now quite clear. \\

Concerning the {\it global} structure of $^*{\bf R}$, the equivalence relationship $\sim$ shows that $^*{\bf R}$ is
obtained from copies of $Gal ( 0 )$ translated to the right and to the left of~ $Gal ( 0 )$. This in view of (3.6.3),
also explains the {\it local} structure of $^*{\bf R}$. \\
A further information about the {\it global} structure of $^*{\bf R}$ is given by the {\it bijection} (3.6.6) \\

Therefore, one of the issues left about the {\it global} structure of $^*{\bf R}$ is {\it how many} copies of $Gal ( 0 )$ does
it contain ? \\

The answer to this question is given in \\ \\

{\bf Proposition 3.6.2} \\

$^*{\bf R}$ is a union of a number of pairwise disjoint copies of $Gal ( 0 )$ which is {\it uncountable}, but it has
the cardinal not greater than that of the {\it continuum}. In other words

\bigskip
(3.6.9) \quad $ ^*{\bf R} ~=~ \bigcup_{\lambda \in \Lambda}~ (~ s^\lambda + Gal ( 0 ) ~) $

\medskip
where $\Lambda$ is uncountable, but of cardinal at most of the continuum, while $s^\lambda \in~ ^*{\bf R}$, and $s^\mu
- s^\lambda \notin Gal ( 0 )$, for $\lambda,~ \mu \in \Lambda,~ \lambda \neq \mu$. \\

{\bf Proof} \\

Let us assume that the set of nonnegative numbers in $^*{\bf R}$ is made up of only countably many pairwise disjoint
copies of $Gal ( 0 )$. If we take one single number in each such copy, we obtain a sequence $s^1, s^2, s^3, ~.~.~.~ \geq
0$, and correspondingly, the decomposition

\bigskip
(3.6.10) \quad $ \{~ s \in ~^*{\bf R} ~~|~~ s ~\geq~ 0 ~\} ~=~ \bigcup_{n \in {\bf N}}~ (~ s^n + Gal ( 0 ) ~) $

\medskip
Let us assume now that $s^n = [~ < s^n_1, s^n_2, s^n_3, ~.~.~.~ > ~] \in~^*{\bf R}$, with $n \in {\bf N}$. Then by the
following diagonal procedure, we can take $s = [~ < s_1, s_2, s_3, ~.~.~.~ > ~] \in ~^*{\bf R}$, such that

$$ \begin{array}{l}
                         s_1 \geq 1 + s^1_1 \\ \\
                         s_2 \geq 2 + s^1_2 + s^2_2 \\ \\
                         s_3 \geq 3 + s^1_3 + s^2_3 + s^3_3 \\ \\
                         ..................................................................... \\ \\
                         s_n \geq n + s^1_n + ~.~.~.~ + s^n_n \\ \\
                         ............................................................................
       \end{array} $$

\medskip
However, the resulting $s\in ~^*{\bf R},~ s ~\geq~ 0$, cannot belong to any of the terms in the above union in (3.6.10).
Indeed, let us assume that $s \in (~ s^n + Gal ( 0 ) ~)$, for a certain $n \in {\bf N}$. Then $s - s^n \in Gal ( 0 )$, hence
there exists $r \in {\bf R},~ r >0$, such that $| s - s^n | \leq r$. This means that

\bigskip
(3.6.11) \quad $| s_j - s^n_j | \leq r,~~~~ ( a. e. ) ~~ \mbox{in}~ j \in {\bf N} $

\medskip
However, from the above construction of $s$ it follows that for $j \geq n$, we have $| s_j - s^n_j | \geq j$, and we
obtained a contradiction with (3.6.11), since the ultrafilter which we use in defining "$( a. e. )$" is free. \\

The fact that $^*{\bf R}$ is the union of a number of copies of $Gal ( 0 )$ which does not exceed the cardinal of the
continuum follows from Note 2.1.1.

\hfill $\Box$ \\ \\

In this way, the only issue which remains to be clarified regarding the {\it global} structure of $^*{\bf R}$ is the {\it
manner} in which it is in (3.6.9) an uncountable union of pairwise disjoint copies of $Gal ( 0 )$. And in this regard we
should note that $^*{\bf R}$ being totally ordered implies that the uncountable family of copies of $Gal ( 0 )$ whose
union generates it in (3.6.9) is also totally ordered. However, this need not mean that the {\it manner} in which this
union is totally ordered is obvious. \\
A first result about the way copies of $Gal ( 0 )$ fill up $^*{\bf R}$, namely, in a {\it dense} manner, is given next. This
result will later be improved in Proposition 5.6.2, based on the concept of {\it hyperfinite}, and thus also {\it
internal} sets. \\ \\

{\bf Proposition 3.6.3} \\

In the uncountable disjoint union (3.6.9) which gives $^*{\bf R}$, between any two copies of $Gal ( 0 )$ there are
countably many other disjoint copies of $Gal ( 0 )$. \\

{\bf Proof} \\

We shall prove that between any two copies of $Gal ( 0 )$ there is another disjoint copy of $Gal ( 0 )$, and from this
obviously follows the desired result, in view of the total order of the copies of $Gal ( 0 )$ in the representation (3.6.9). \\

Let $\lambda,~ \mu \in \Lambda$ be such that $s^\lambda < s^\mu$. Clearly we have $s = ( s^\lambda + s^\mu ) / 2 \in~
^*{\bf R}$, therefore $s \in s^\nu + Gal ( 0 )$, for a suitable $\nu \in \Lambda$. Thus $| s - s^\nu | \leq r$, for
some $r \in {\bf R},~ r > 0$. \\
But $s^\nu \notin s^\lambda + Gal ( 0 )$. Indeed, otherwise we would have $| s^\nu - s^\lambda | \leq r_1$, for a
certain $r_1 \in {\bf R},~ r_1 > 0$, and this would give

$$ | s^\mu - s^\lambda | ~=~ 2 | ( s^\lambda + s^\mu ) / 2 - s^\lambda | ~=~ 2 | s - s^\lambda | ~\leq~
                      2 ( | s - s^\nu | + | s^\nu - s^\lambda | ) ~\leq~ 2 ( r + r_1 ) $$

\medskip
which would contradict the hypothesis that the terms in the union (3.6.9) are disjoint. \\
Similarly $s^\nu \notin s^\mu + Gal ( 0 )$. On the other hand, obviously $s^\lambda < s^\nu < s^\mu$. Thus $s^\nu +
Gal ( 0 )$ is indeed between $s^\lambda + Gal ( 0 )$ and $s^\mu + Gal ( 0 )$, and it is disjoint from both of
them. \\ \\

{\bf Corollary 3.6.1} \\

1.~~ In $^*{\bf R}$, there is {\it no} term in the disjoint union (3.6.9), in other \\
\hspace*{0.7cm} words, there is no copy of $Gal ( 0 )$, which would come immediately \\
\hspace*{0.7cm} to the right, or to the left of the set $Gal ( 0 )$ of all nonstandard \\
\hspace*{0.7cm} finite numbers. \\

2.~~ In $^*{\bf R}$, the set $Gal ( 0 )$ of all nonstandard finite numbers does {\it not} \\
\hspace*{0.7cm} have an infimum or a supremum. \\

3.~~ In $^*{\bf R}$, the set $mon ( 0 )$ of infinitesimals does {\it not} have an infimum \\
\hspace*{0.7cm} or a supremum. \\

4.~~ The set $^*{\bf Z}$ of all nonstandard integers is {\it uncountable}. \\

{\bf Proof} \\

1 follows from Proposition 3.6.3. \\
2 results from 1. \\
3 is a consequence of 2 and the bijection (3.6.6). \\
4 We note that, given any $s \in~ ^*{\bf R}$, there exists $t \in~ ^*{\bf Z}$, such that $s - t \in Gal ( 0 )$. Indeed,
assume that $s = [~ < s_1, s_2, s_3, ~.~.~.~ . ~]$. Then we can take $t = [~ < t_1, t_2, t_3, ~.~.~.~ . ~]$, where $t_n$
is the integer part of $s_n$, with $n \in {\bf N}$. \\
It follows that, for $\lambda \in \Lambda$, there exists $t^\lambda \in~ ^*{\bf Z}$, such that $t^\lambda \in s^\lambda
+ Gal ( 0 )$. But the uncountable union in (3.6.9) is made up from pairwise disjoint sets, thus the proof of 4 is
completed.

\hfill $\Box$ \\

Let us give a few basic {\it topological} results which use the concepts of {\it monad} and {\it standard part}. \\ \\

{\bf Proposition 3.6.4} \\

Let $A \subseteq {\bf R}$, then : \\

1.~~~ $A ~~~\mbox{is open in}~~{\bf R} ~~~\Longleftrightarrow~~~ mon ( r ) ~\subseteq~ ^*A,~~~\mbox{for}~~ r \in A$ \\

2.~~~ $A ~~~\mbox{is closed in}~~{\bf R} ~~~\Longleftrightarrow~~~ st ( s ) \in A, ~~~\mbox{for}~~ s \in~ ^*A \bigcap Gal ( 0 ) $ \\

3.~~~ $A ~~~\mbox{is compact in}~~{\bf R} ~~~\Longleftrightarrow~~~
                                       \left (~ \begin{array}{l}
                                                            ^*A ~\subseteq~ Gal ( 0 ), ~~~\mbox{and}~~~ \\ \\
                                                            st ( s ) \in A, ~~~\mbox{for}~~ s \in~ ^* A
                                                   \end{array} \right ) $ \\

{\bf Proof} \\

1. Assume that $A$ is not open. Then for some $r \in A$, every neighbourhood $V \subseteq A$ of $r$ will intersect ${\bf
R} \setminus A$. Thus for every $n \in {\bf N}$, we can take $r_n \in {\bf R} \setminus A$, such that $| r - s_n | < 1 /
( n + 1 )$. It follows that $s = [~ < s_1, s_2, s_3, ~.~.~.~ > ~] \not\in~ ^*A$. However, obviously $x \in mon ( r )$. \\
Conversely, if $A$ is open, let $r \in A$. Then for every neighbourhood $V \subseteq A$ of $r$ we have $mon ( r )
\subseteq~ ^*V \subseteq~ ^*A$. \\
2. Let $s = [~ < s_1, s_2, s_3, ~.~.~.~ > ~] \in~ ^*A \bigcap Gal ( 0 )$, such that $r = st ( s ) \not\in A$. Then $r - s \in mon
( 0 )$. However, $s_n \in A$, for $n \in {\bf N}$, hence $r$ is an acummulation point of $A$, which means that $A$ is not
closed. \\
Conversely, we prove that ${\bf R} \setminus A$ is open. According to 1 above, it is enough to show that for $r \in {\bf R}
\setminus A$, we have $mon ( r ) \subseteq~ ^*( {\bf R} \setminus A) =~ ^*{\bf R} \setminus~ ^*A$. Or in other words, we
have $mon ( r ) \bigcap~ ^*A = \phi$, for $r \in {\bf R} \setminus A$. \\
Let us assume that $s \in mon ( r ) \bigcap~ ^*A$, for a certain $r \in {\bf R} \setminus A$. Then by the assumption, we
have $st ( s ) \in A$. However $s \in mon ( r )$ implies that $st ( s) = r$, thus we get the contradiction that $r \in A$. \\
3. We note that for $A \subseteq {\bf R}$ we obviously have $^*A \subseteq Gal ( 0 )$, if and only if $A$ is bounded. The
rest then follows from 2 above. \\ \\

{\bf 7. Self-Similarity in $^*{\bf R}$} \\ \\

The local and global structures of $^*{\bf R}$, as we have seen in section 6, are closely related. Here we shall point
to a {\it self-similar} aspect of this interrelation which may remind one of a typical feature of fractals. \\

First we recall that, see (3.6.6), we have the order reversing bijection

\bigskip
(3.7.1) \quad $ (~ ^*{\bf R} ~\setminus~Gal ( 0 ) ~) \ni s ~\longmapsto~ 1 / s \in (~ mon ( 0 ) \setminus \{~ 0 ~\} ~) $

\medskip
Now the global structure of $^*{\bf R}$ is given by, see (3.6.9)

\bigskip
(3.7.2) \quad $ ^*{\bf R} ~=~ \bigcup_{\lambda \in \Lambda}~ (~ s^\lambda + Gal ( 0 ) ~) $

\medskip
while its local structure is described by, see (3.6.3)

\bigskip
(3.7.3) \quad $ Gal ( 0 ) ~=~ \bigcup_{~r \in {\bf R}}~ (~ r ~+~ mon ( 0 ) ~)$ \\

\medskip
In this way we obtain the {\it self-similar order reversing bijections}, which are expressed in terms of $mon ( 0 )$,
namely \\

\bigskip
(3.7.4) \quad $ \begin{array}{l}
                (~ \bigcup_{~r \in {\bf R},~ \lambda \in \Lambda_0}~ (~ r + s^\lambda + mon ( 0 ) ~) ~) \ni s ~\longmapsto~ \\ \\
                  ~~~~~~~~~~~~~~~~~~\longmapsto~ 1 / s \in (~ mon ( 0 ) \setminus \{~ 0 ~\} ~)
                \end{array} $ \\

\medskip
and

\bigskip
(3.7.5) \quad $ \begin{array}{l}
                   (~ mon ( 0 ) \setminus \{~ 0 ~\} ~) \ni s ~\longmapsto~ \\ \\
             ~~~~~~~~~~\longmapsto~ 1 / s \in (~ \bigcup_{~r \in {\bf R},~ \lambda \in \Lambda_0}~ (~ r + s^\lambda + mon ( 0 ) ~) ~)
                       \end{array} $ \\

\medskip
where $\Lambda_0 = \Lambda \setminus \{~ \lambda_0 ~\}$, with $\lambda_0 \in \Lambda$ such that $s^{\lambda_0} \in
Gal ( 0 )$. \\

The above bijections in (3.7.4), (3.7.5) are given by the very simple algebraic, explicit, and order reversing mapping
$s \longmapsto 1 / s$ which involves what is essentially a field operation, namely, division. And these two bijections
take the place of the much simpler order reversing bijections in the case of the standard reals ${\bf R}$, namely

\bigskip
(3.7.6) \quad $ (~ {\bf R} \setminus ( -1 , 1 ) ~) \ni r ~\longmapsto~ 1 / r \in (~ [ - 1, 1 ] \setminus \{~ 0 ~\} ~) $

\bigskip
(3.7.7) \quad $ (~ [ - 1, 1 ] \setminus \{~ 0 ~\} ~) \ni r ~\longmapsto~ 1 / r \in (~ {\bf R} \setminus ( -1 , 1 ) ~) $ \\ \\

{\bf 8. Sets which Are Not Transfers of Standard Sets} \\ \\

In section 7, in chapter 2, we have seen how to transfer relations and functions from ${\bf R}$ to $^*{\bf R}$. In
particular, in (2.7.6), the way subsets $A \subseteq {\bf R}$ are transferred to subsets $^*A \subseteq~ ^*{\bf R}$ was
presented. In this respect it is important to note that by far most subsets of $^*{\bf R}$ are {\it not} obtainable as
such transfers of subsets of ${\bf R}$. \\
Later, in chapter 5, we shall see how this issue is related to the dichotomy between the concepts of {\it internal} and
{\it external} entities in $^*{\bf R}$, where the internal entities generalize those which are elements of subsets $^*A
\subseteq~ ^*{\bf R}$ that are transfers of subsets $A \subseteq {\bf R}$. \\

Here, as a simple illustration we give several more important examples of subsets in $^*{\bf R}$ which are {\it not}
transfers of subsets in ${\bf R}$. \\
In this regard it is important to remember that, as mentioned from the beginning, transfer - no matter how crucial in
Nonstandard Analysis - is nevertheless not the only powerful instrument available. Indeed, the wealth and intuitive
nature of nonstandard concepts, methods and results which are beyond the realms transfer operates are part of the
interest in, and power of Nonstandard Analysis. And to mention only one such concept, let us recall that the {\it
infinitesimals} in $^*{\bf R}$ - intuited already by Leibniz - constitute the set $mon ( 0 )$ which, as seen next, {\it
cannot} be obtained by transfer from any subset $A \subseteq {\bf R}$. \\
Needless to say, as seen in the sequel, there are many important subsets in $^*{\bf R}$ which, similarly, cannot be
obtained by transfer from subsets of ${\bf R}$. And in general, as seen in chapters 5 and 6, there are many important
entities in Nonstandard Analysis which are not obtainable from standard entities by transfer. \\

And now, let us give a few examples of subsets in $^*{\bf R}$ which are not transfers of subsets in ${\bf R}$. \\

1.~~ Let us note that ${\bf R}$ itself, when seen as a subset of $^*{\bf R}$, is {\it not} of the form $^*A$, for any
subset $A \subseteq {\bf R}$. \\

Indeed, assume that ${\bf R} =~ ^*A$ for a certain $A \subseteq {\bf R}$. Now if $A$ is bounded from above, then we
have $A \subseteq \{~ r \in {\bf R} ~~|~~ r \leq a ~\}$, in other words, the sentence

$$ ( \forall~ x )~ \Big [~ \underline{A} < x > ~\longrightarrow~ x \leq a ~\Big ] $$

\medskip
is true in ${\cal L}_{\cal R}$, where $A < x >$ is the unary relation $x \in A$. Therefore by transfer,the sentence

$$ ( \forall~ x )~ \Big [~ \underline{^*A} < x > ~\longrightarrow~ x \leq a ~\Big ] $$

\medskip
is true in ${\cal L}_{^*{\cal R}}$. However, as we assumed that ${\bf R} =~ ^*A$, this implies ${\bf R} \subseteq
( \infty, a ]$, which is of course absurd. \\
Now in the alternative case when $A$ is not bounded from above, then we can choose $a_n \in A$, with $n \in {\bf N}$,
such that $\lim_{n \to \infty}~ a_n = \infty$. In this case (2.7.6) gives $a = [~ < a_1, a_2, a_3, ~.~.~.~ > ~] \in
^*A$. However, obviously $a \notin Gal ( 0 )$, thus $a \notin {\bf R}$. \\

2.~~ We have similarly that ${\bf N}$ and ${\bf Z}$ are {\it not} of the form $^*A$, for any subset $A \subseteq {\bf R}$. \\

3.~~ Also, $mon ( 0 )$ and $Gal ( 0 )$ are {\it not} of the form $^*A$, for any subset $A \subseteq {\bf R}$. \\

Indeed, let us assume that $mon ( 0 ) =~ ^*A$ for a certain $A \subseteq {\bf R}$. Then 1) in (2.7.8) gives $A \subseteq
mon ( 0 )$, and in view of (3.6.4) it follows that $A = \{~ 0 ~\}$, which according to 2) in (2.7.8) implies $^*A = A = \{~ 0 ~\}$,
thus a contradiction. \\
The case of $Gal ( 0 )$ can be dealt with in a way similar to that of ${\bf R}$. \\

4. Further, {\it none} of the sets $^*{\bf N} \setminus {\bf N}$,~~ $^*{\bf Z} \setminus {\bf Z}$,~~ or $^*{\bf R} \setminus
{\bf R}$,~~ are of the form $^*A$, for any subset $A \subseteq {\bf R}$. \\

Let us assume that $^*{\bf N} \setminus {\bf N} =~ ^*A$ for a certain $A \subseteq {\bf R}$. Let us take any $a \in A$.
Then 1) in (2.7.8) gives $a \in~ ^*A = ^*{\bf N} \setminus {\bf N}$, which is absurd, since $a \in Gal ( 0 )$. The proofs for
the other two sets $^*{\bf Z} \setminus {\bf Z}$,~~ and $^*{\bf R} \setminus {\bf R}$~~are similar. \\ \\

{\bf 9. A Few Basic Application to Calculus :  Sequences \\ \hspace*{0.45cm} and their Limits} \\ \\

At this stage we can already present a few  of the important application of Nonstandard Analysis to Calculus. Needless
to say, here the limitation to such applications comes exclusively from the rather rudimentary nature of the simple
languages in which we obtained so far the Transfer Property. \\
Nevertheless, even this limited framework can offer several rather relevant examples about the power and advantages
of Nonstandard Analysis. \\
Once, however, we introduce Superstructures and their corresponding languages, all such limitations will be done
away with. \\

Let $(~ x_n ~)_{n \in {\bf N}}$ be any given usual sequence of real numbers in ${\bf R}$. It will be conveninet to look at
this sequence as a usual real valued function $x : {\bf N} \longrightarrow {\bf R}$, with $x ( n ) = x_n$, for $n \in
{\bf N}$. \\
Now in view of Proposition 3.5.3, the transfer of this function is

$$ ^*x :~ ^*{\bf N} \longrightarrow~ ^* {\bf R} $$

\medskip
and this is an extension of the given usual sequence $x : {\bf N} \longrightarrow {\bf R}$. In view of that, it will be
convenient to denote $^*x$ simply by $x$, even upon its extended range $^*{\bf N}$. In particular, instead of writing
$^*x ( n )$, we shall write $x_n$, even when $n \in~  ^*{\bf N} \setminus {\bf N}$. \\ \\

{\bf Proposition 3.9.1} \\

Given a real number $l \in {\bf R}$, we have in the usual sense $\lim_{n \to \infty}~ x_n = l$, if and only if

\bigskip
(3.9.1) \quad $ x_n ~\approx~ l,~~~ \mbox{for}~~ n \in~ ^*{\bf N} \setminus {\bf N} $ \\ \\

{\bf Note 3.9.1} \\

The above {\it nonstandard} convergence characterization for the usual sequence $(~ x_n ~)_{n \in {\bf N}}$ is both
surprising and edifying. Indeed \\

- it says {\it nothing} explicitly about the way the usual terms $x_n$, with $n \in {\bf N}$, of the sequence are supposed
to behave, \\

and on the ohter hand \\

- it requires the very {\it strong} and {\it intuitively} obvious condition on the nonstandard terms $x_n$, with $n \in~
^*{\bf N} \setminus {\bf N}$, to be {\it infinitely near} to the limit $l$. \\

A similar kind of {\it dichotomy} in which, on one hand, the usual terms in a sequence $(~ x_n ~)_{n \in {\bf N}}$ are not
subjected to explicit conditions, while on the other hand, the nonstandard terms are required strong and intuitively
obvious conditions can be seen in the next few results as well. \\

Here however it is important to note that the {\it lack} of explicitly stated conditions on the usual terms in a sequence
$(~ x_n ~)_{n \in {\bf N}}$ does {\it not} at all mean that there are no conditions whatsoever on these usual terms.
Indeed, it is precisely these usual terms alone which define the transferred sequence $(~ x_n ~)_{n \in~  ^*{\bf N}}$.
Thus any conditions on the nonstandard terms of the transferred sequence $(~ x_n ~)_{n \in~  ^*{\bf N}}$ do in fact mean
{\it implicit} conditions on the usual terms of the sequence $(~ x_n ~)_{n \in {\bf N}}$. \\ \\

{\bf Proof} of Proposition 3.9.1. \\

Let us assume that we have $\lim_{n \to \infty}~ x_n = l$. As we know from Calculus, then

$$ \begin{array}{l}
                      \forall~ \epsilon \in {\bf R},~ \epsilon > 0 ~: \\
                      \exists~ m \in {\bf N} ~: \\
                      \forall~ n \in {\bf N},~ n \geq m ~: \\
                      ~~~ |~ l - x_n ~| \leq \epsilon
     \end{array} $$

\medskip
Now let us fix any $\epsilon > 0$. Then we obtain a certain $m \in {\bf N}$, such that the sentence

$$ ( \forall n ) ~\Big [~ (~ (~ \underline{{\bf N}} < n > ~) \wedge (~ n \geq m ~) ~)
                                         ~\longrightarrow~ (~ |~ l - x_n ~| \leq \epsilon ~) ~\Big ] $$

\medskip
is true in ${\cal L}_{\cal R}$. Therefore by transfer, the sentence

$$ ( \forall n ) ~\Big [~ (~ (~ \underline{^*{\bf N}} < n > ~) \wedge (~ n \geq m ~) ~)
                                         ~\longrightarrow~ (~ |~ l - \underline{^*x} ( n ) ~| \leq \epsilon ~) ~\Big ] $$

\medskip
is true in ${\cal L}_{^*{\cal R}}$. In particular, it follows that

$$ |~ l - x_n ~| \leq \epsilon,~~~ n \in~ ^*{\bf N} \setminus {\bf N} $$

\medskip
And since $\epsilon \in {\bf R},~ \epsilon > 0$ is arbitrary, it means that (3.9.1) does indeed hold. \\

Conversely, assume that, in the usual sense, we do not have $\lim_{n \to \infty}~ x_n = l$. Then there exists $\epsilon
\in {\bf R},~ \epsilon> 0$, such that

$$ \begin{array}{l}
                 \forall~ m \in {\bf N} ~: \\
                 \exists~ n \in {\bf N},~ n \geq m ~: \\
                 ~~~ |~ l - x_n ~| \geq \epsilon
     \end{array} $$

\medskip
Therefore, there exists a Skolem function $\psi : {\bf N} \longrightarrow {\bf N}$ with the property $\psi ( m ) \geq m$, for
$m \in {\bf N}$, and such that in addition the sentence

$$ ( \forall m ) ~\Big [~ {\bf N} < m > ~\longrightarrow~  (~ |~ l - x_{\psi ( m )} ~| \geq \epsilon ~) ~\Big ] $$

\medskip
is true in ${\cal L}_{\cal R}$. Therefore by transfer, the sentence

$$ ( \forall m ) ~\Big [~ \underline{^*{\bf N}} < m > ~\longrightarrow~
                                     (~ |~ l - \underline{^*x} ( \underline{^*\psi ( m )} ) ~| \geq \epsilon ~) ~\Big ] $$

\medskip
is true in ${\cal L}_{^*{\cal R}}$. This means that

$$  |~ l - x_{\psi ( m )} ~| \geq \epsilon,~~~ m \in~ ^*{\bf N} $$

\medskip
and in particular (3.9.1) cannot hold. \\ \\

{\bf Proposition 3.9.2} \\

The sequence $(~ x_n ~)_{n \in {\bf N}}$ is Cauchy in the usual sense, if and only if

\bigskip
(3.9.2) \quad $ x_n ~\approx~ x_m,~~~ \mbox{for}~~ n,~ m \in~ ^*{\bf N} \setminus {\bf N} $ \\

{\bf Proof} \\

We can follow the proof of Proposition 3.9.1 with a few modifications. Here we only indicate them in proving that (3.9.2)
is necessary for the sequence $(~ x_n ~)_{n \in {\bf N}}$ to be Cauchy. The modifications for obtaining the converse are
quite obvious. \\
We recall from Calculus that the necessary and sufficient condition for the sequence $(~ x_n ~)_{n \in {\bf N}}$ to be
Cauchy is that

$$ \begin{array}{l}
                    \forall~ \epsilon \in {\bf R},~ \epsilon > 0 ~: \\
                    \exists~ k \in {\bf N} ~: \\
                    \forall~ n,~m \in {\bf N},~ n,~ m \geq k ~: \\
                    ~~~ |~ x_n - x_m ~| \leq \epsilon
      \end{array} $$

\medskip
Let us assume that the sequence $(~ x_n ~)_{n \in {\bf N}}$ is Cauchy and take any fixed $\epsilon \in {\bf R},~
\epsilon > 0$. Then for the resulting $k \in {\bf N}$, the sentence

$$ \begin{array}{l}
                ( \forall n )~ ( \forall~ m ) ~\Big [~ (~ (~ \underline{{\bf N}} < n > ~) \wedge (~ \underline{{\bf N}} < m > ~) \wedge
                        (~ n \geq k ~) \wedge (~ m \geq k ~) ~)  ~\longrightarrow~ \\ \\
                ~~~~~~~~~~~~~~~\longrightarrow ~(~ |~ x_n - x_m ~| \leq \epsilon ~) ~\Big ]
      \end{array} $$

\medskip
is true in ${\cal L}_{\cal R}$. Therefore by transfer, the sentence

$$ \begin{array}{l}
                ( \forall n )~ ( \forall~ m ) ~\Big [~ (~ (~ \underline{^*{\bf N}} < n > ~) \wedge (~ \underline{^*{\bf N}} < m > ~) \wedge
                        (~ n \geq k ~) \wedge (~ m \geq k ~) ~)  ~\longrightarrow~ \\ \\
                ~~~~~~~~~~~~~~~\longrightarrow~ (~ |~ x_n - x_m ~| \leq \epsilon ~) ~\Big ]
      \end{array} $$

\medskip
is true in ${\cal L}_{^*{\cal R}}$. It follows that

$$ |~ x_n - x_m ~| \leq \epsilon,~~~ n,~ m \in~ ^*{\bf N} \setminus {\bf N} $$

\medskip
and since $\epsilon \in {\bf R},~ \epsilon > 0$ is arbitrary, the relation (3.9.2) is proved. \\ \\

{\bf Proposition 3.9.3} \\

The sequence $(~ x_n ~)_{n \in {\bf N}}$ is bounded, if and only if

\bigskip
(3.9.3) \quad $ x_n \in Gal ( 0 ),~~~ \mbox{for}~~ n \in~ ^*{\bf N} \setminus {\bf N} $ \\

{\bf Proof} \\

If the sequence $(~ x_n ~)_{n \in {\bf N}}$ is bounded, then for a certain $L \in {\bf N},~ L > 0$, the sentence

$$ ( \forall~ n )~ \Big [~ \underline{{\bf N}} < n > ~\longrightarrow~ |~ x_n ~| \leq L ~\Big ]  $$

\medskip
is true in ${\cal L}_{\cal R}$. Thus by transfer, the sentence

$$ ( \forall~ n )~ \Big [~ \underline{^*{\bf N}} < n > ~\longrightarrow~ |~ x_n ~| \leq L ~\Big ]  $$

\medskip
is true in ${\cal L}_{^*{\cal R}}$. It follows that

$$ x_n \in Gal ( 0 ),~~~ \mbox{for}~~ n \in~ ^*{\bf N} $$

\medskip
hence in particular, we have (3.9.3). \\

Conversely, if the sequence $(~ x_n ~)_{n \in {\bf N}}$ is not bounded, then there exists a Skolem function $\psi : {\bf N}
\longrightarrow {\bf N}$, such that $\psi ( n ) \geq n$, for $n \in {\bf N}$, and the sentence

$$ ( \forall~ n )~ \Big [~ \underline{{\bf N}} < n > ~\longrightarrow~ |~ x_{\underline{\psi(n)}} ~| \geq n ~\Big ]  $$

\medskip
is true in ${\cal L}_{\cal R}$. Then by transfer, the sentence

$$ ( \forall` n )~ \Big [~ \underline{^*{\bf N}} < n > ~\longrightarrow~ |~ x_{\underline{^*\psi(n)}} ~| \geq n ~\Big ]  $$

\medskip
is true in ${\cal L}_{^*{\cal R}}$. Thus in particular, we cannot have (3.9.3). \\ \\

{\bf Proposition 3.9.4} \\

A given real number $l \in {\bf R}$ is a usual limit point of the sequence $(~ x_n ~)_{n \in {\bf N}}$, if and only if

\bigskip
(3.9.4) \quad $ x_n ~\approx~ l,~~~ \mbox{for some}~~ n \in~ ^*{\bf N} \setminus {\bf N} $ \\

{\bf Proof} \\

Let us assume that $l \in {\bf R}$ is a usual limit point of the sequence $(~ x_n ~)_{n \in {\bf N}}$. Then

$$ \begin{array}{l}
                    \forall~ \epsilon \in {\bf R},~ \epsilon > 0 ~: \\
                    \forall~ m \in {\bf N} ~: \\
                    \exists~ n \in {\bf N},~ n \geq m ~: \\
                    ~~~ |~ l - x_n ~| \leq \epsilon
      \end{array} $$

\medskip
Hence there exists a Skolem function $\psi : {\bf R}_+ \times {\bf N} \longrightarrow {\bf N}$, where we denoted
${\bf R}_+ = ( 0, \infty )$, such that $\psi ( \epsilon, n ) \geq n$, and the sentence

$$ ( \forall~ \epsilon )~( \forall~ n )~ \Big [~ (~ \underline{{\bf R}_+} < \epsilon > \wedge~ \underline{{\bf N}} < n > ~)
                                  ~\longrightarrow~ |~ l - x_{\underline{\psi ( \epsilon, n )}} ~| \leq \epsilon ~\Big ] $$

\medskip
is true in ${\cal L}_{\cal R}$. Then by transfer, the sentence

$$ ( \forall~ \epsilon )~( \forall~ n )~ \Big [~ (~ \underline{^*(~{\bf R}_+ ~)} < \epsilon > \wedge~ \underline{^*{\bf N}} < n > ~)
                                  ~\longrightarrow~ |~ l - x_{\underline{^*\psi ( \epsilon, n )}} ~| \leq \epsilon ~\Big ] $$

\medskip
is true in ${\cal L}_{^*{\cal R}}$. Let us take now $\epsilon \in mon ( 0 ) \bigcap~ ^*(~ {\bf R}_+ ~)$ and $m \in~ ^*{\bf N}
\setminus {\bf N}$, and denote $n =~ ^*\psi ( \epsilon, m )$. Then $n \geq m$, hence $n \in~ ^*{\bf N} \setminus {\bf N}$,
and $|~ l - x_n ~| \leq \epsilon$, which means that (3.9.4) holds. \\

Conversely, if $l$ is not a limit point of the sequence $(~ x_n ~)_{n \in {\bf N}}$, then there exist $\epsilon \in {\bf R},~
\epsilon > 0$ and $m \in {\bf N}$, such that the sentence

$$ ( \forall~ n )~ \Big [~ (~ (~ \underline{{\bf N}} < n > ~) \wedge (~ n \geq m ~) ~)
                                                                \longrightarrow ~  |~ l - x_n ~| \geq \epsilon ~\Big ] $$

\medskip
is true in ${\cal L}_{\cal R}$. Then by transfer, the sentence

$$ ( \forall~ n )~ \Big [~ (~ (~ \underline{^*{\bf N}} < n > ~) \wedge (~ n \geq m ~) ~)
                                                                \longrightarrow ~  |~ l - x_n ~| \geq \epsilon ~\Big ] $$

\medskip
is true in ${\cal L}_{^*{\cal R}}$. And this clearly contradicts (3.9.4). \\ \\

{\bf Note 3.9.2} \\

As seen already in Proposition 3.9.3, the set $Gal ( 0 )$, which is {\it not} obtainable by transfer from any subset of
${\bf R}$, can be used to give a characterization of usual bounded sequences in ${\bf R}$. This is but one simple
example of the way in which nonstandard concepts which do not have a standard correspondent can nevertheless play
an important role in understanding and dealing with standard entities. \\ \\

{\bf 10. Proof of the Transfer Property in its Simple Version} \\ \\

\chapter{Short Review}

{\bf 1. Two Stages of Transfer} \\ \\

Let us review the essence of the constructions in chapters 2 and 3, in order to be better prepared to embark upon the
general approach based on Superstructures presented in chapter 5. As already mentioned, this essence has two
constituents, namely : \\

(TE)~:~~ the {\it transfer} of various {\it mathematical entities} from ${\bf R}$ to $^*{\bf R}$, \\
         \hspace*{1.45cm} the first and simplest case of it being the field homomorphism \\
         \hspace*{1.45cm} $^*(~)$ in (2.1.4) \\

(TS)~:~~ the {\it transfer} of {\it sentences} $\Phi$ from a given {\it language} ${\cal L}_{\bf R}$ about \\
         \hspace*{1.45cm} mathematical entities in ${\bf R}$ to corresponding sentences $^*\Phi$ in \\
         \hspace*{1.45cm} a language ${\cal L}_{^*{\bf R}}$ about mathematical entities in $^*{\bf R}$. \\

Furthermore, as already seen in chapter 3, the transfer of sentences is a {\it two way} process. Namely, when we want
to prove a standard property of ${\bf R}$, we formulate it as a sentence $\Phi$ in the language ${\cal L}_{\bf R}$,
then we transfer it to the sentence $^*\Phi$ in the language ${\cal L}_{^*{\bf R}}$, and prove this transferred
sentence $^*\Phi$ in $^*{\bf R}$, after which we deduce by reverse transfer that the initial sentence $\Phi$ is true
about ${\bf R}$. \\

Also, the stage (TS) of transfer is far more {it powerful} and {\it systematic} than (TE). The stage (TE), however, has the
advantage that it can be performed by employing usual mathematics. On the other hand, the stage (TS) needs a few
basic concepts and methods from Mathematical Logic.  Further, as seen in chapter 3, the most effective use of transfer
is by combining both (TE) and (TS), depending on the specific situation. \\

With respect to (TE), let us recall the following more important steps. First, as mentioned, the {\it field extension}
(2.1.4) is given by

\begin{math}
\setlength{\unitlength}{1cm}
\thicklines
\begin{picture}(15,2.5)
\put(2,1){${\bf R} \ni r$}
\put(3.2,1){\line(0,1){0.2}}
\put(3.2,1.1){\vector(1,0){2}}
\put(3.85,1.4){$^*(~)$}
\put(5.4,1){$^*r ~=~ [~ < r, r, r, ~.~.~.~ > ~] \in~ ^*{\bf R}$}
\put(0,1){$(4.1.1)$}
\end{picture}
\end{math}

Given now any n-ary {\it relation} $P$ on ${\bf R}$, its {\it transfer} $^*P$ is the n-ary relation on $^*{\bf R}$, defined as
follows. If $( s^1, ~.~.~.~ , s^n ) \in ( ^*{\bf R} )^n$ is any n-tuple, such that $s^i = [~ < s^i_1, s^i_2, s^i_3, ~.~.~.~ >~]$, with $1
\leq i \leq n$, then

\bigskip
(4.1.2) \quad $ \begin{array}{l}
                    ^*P < s^1, ~.~.~.~ , s^n > ~~~\Longleftrightarrow~~~ \\ \\
                         ~~~~~~~~~~\Longleftrightarrow~~~ \{~ j \in {\bf N} ~~|~~
                      P < s^1_j, ~.~.~.~ , s^n_j > ~\} \in {\cal U}
                 \end{array} $

\medskip
where ${\cal U}$ is the free ultrafilter on ${\bf N}$, used in the construction of $^*{\bf R}$, see (2.4.1). \\

As we have seen in chapters 2 and 3, once we know how to transfer arbitrary n-ary relations from ${\bf R}$ to $^*{\bf
R}$, we have in fact completed the task at (TE). \\

As for (TS), in chapter 3 we introduced the way to associate a {\it simple language} ${\cal L}_{\cal S}$ with any {\it
simple system} ${\cal S} = ( S, ( P_i ~|~ i \in I ), ( f_j ~|~ j \in J ) )$, see Definition 3.2.2. Then we applied this association
process to the simple systems of our interest, see (3.3.1) and (3.4.1)

\bigskip
(4.1.3) \quad $ {\cal R} ~=~ (~ {\bf R},~ {\cal P},~ {\cal F} ~),~~~ ^*{\cal R} ~=~ (~ ^*{\bf R},~ ^*{\cal P},~ ^*{\cal F} ~) $

\medskip
which coorespond to ${\bf R}$ and $^*{\bf R}$, respectively. In this way we obtained the corresponding simple
languages

\bigskip
(4.1.4) \quad $ {\cal L}_{\cal R},~~~ {\cal L}_{^*{\cal R}} $

\medskip
Now, the {\it transfer} on the level of {\it languages}, that is, the mapping which in the case of (4.1.4) corresponds to
(4.1.1), was given in 5 and 6 in Definition 3.4.1. And this transfer operates among {\it simple sentences}, be they atomic
or compund, according to

\begin{math}
\setlength{\unitlength}{1cm}
\thicklines
\begin{picture}(15,2.5)
\put(2,1){${\cal L}_{\cal R} \ni \Phi$}
\put(3.7,1){\line(0,1){0.2}}
\put(3.7,1.1){\vector(1,0){2}}
\put(4.35,1.4){$^*(~)$}
\put(5.9,1){$^*\Phi \in {\cal L}_{^*{\cal R}}$}
\put(0,1){$(4.1.5)$}
\end{picture}
\end{math}

Once we have {\it shifted} in (4.1.5) to the level of {\it languages} the transfer in (4.1.1), we obtain the corresponding
{\it simple} version of the Transfer Property, formulated in the Transfer Theorem, see section 4, in chapter 3, namely

\bigskip
(4.1.6) \quad $ \begin{array}{l}
                                        \forall~~ \Phi ~~\mbox{simple sentence in}~~ {\cal L}_{\cal R} ~: \\ \\
                                         ~~~ \Phi ~~\mbox{true in}~~ {\cal L}_{\cal R} ~~~\Longrightarrow~~~
                                                                                         ^*\Phi ~~\mbox{true in}~~ {\cal L}_{^*{\cal R}}
                         \end{array} $

\medskip
With this construction, the stage (TS) is completed, and as seen in chapter 3, this stage is far more powerful and
systematic than stage (TE), which in fact it contains as a particular case. Furthermore, as also seen, it is particularly
convenient to combine the use of both these two stages, according to the specifics of the situation. \\ \\

{\bf 2. Nonstandard Approach to General Mathematics} \\ \\

As mentioned, Nonstandard Analysis is not limited to the extension of the reals ${\bf R}$ into $^*{\bf R}$, and can be
successfully applied to a variety of mathematical theories, such as for instance, General Topology, Probablity Theory,
Functional Analysis, and so on. \\
In the next chapter, we shall present the basics of such a general apporach, by introducing Superstructures. \\

Let us recall that the main reason for using Superstructures is that we want to obtain the {\it full} version of the Transfer
Property, see (3.4.3). \\
The additional reason is that we want to create the possibility to apply Nonstandard Analysis to a variety of
mathematical theories. In chapter 6, we shall present such an application to Probability Theory, by introducing the Loeb
measures. \\

Concerning the versatility of Nonstandard Analysis, let us present here two more coments. \\

Suppose that we start the construction in chapter 2 with the rational numbers ${\bf Q}$, instead of the real numbers
${\bf R}$. Then similar with (2.1.4), we shall obtain

\begin{math}
\setlength{\unitlength}{1cm}
\thicklines
\begin{picture}(15,7)
\put(2,5){${\bf Q}$}
\put(2.7,5.15){\vector(1,0){6}}
\put(5.5,5.4){$^*(~)$}
\put(9,5){$^*{\bf Q} ~=~ {\bf Q}^{\bf N} ~/~ {\cal M}$}
\put(6,1.5){${\bf Q}^{\bf N}$}
\put(2.7,4.8){\vector(1,-1){3}}
\put(6.8,2){\vector(1,1){2.7}}
\put(0,3){$(4.2.1)$}
\end{picture}
\end{math}

where the {\it injective field homomorphism}, thus {\it field extension} $^*(~)$ is given by

\bigskip
(4.2.2) \quad $ \begin{array}{l}
                                    {\bf Q} \ni r ~~\longmapsto~~ ^* r ~=~ [~ < r, r, r, ~.~.~.~ > ~] ~=~ \\ \\
                                       ~~~~~~~~~~~~~ =~ < r, r, r, ~.~.~.~ > + {\cal M} \in ^*{\bf Q} ~=~ {\bf Q}^{\bf N} ~/~ {\cal M}
                     \end{array} $

\medskip
between the two fields ${\bf Q}$ and $~^*{\bf Q} ~=~ {\bf Q}^{\bf N} ~/~ {\cal M}$. \\

As in Definition 3.6.1, now in the case of $^*{\bf Q}$ as well, we can define the concepts of {\it infinitesimal, finite} and
{\it infinite} elements, and denote the first two classes with $mon_{\bf Q} ( 0 )$ and $Gal_{\bf Q} ( 0 )$, respectively. \\

However, since ${\bf Q}$ is {\it not} complete as a uniform topological space, we have to be careful when, as in 6 in
Proposition 3.6.1, we want to define the {\it standard part} mapping $st$, this time on $Gal_{\bf Q} ( 0 )$. Indeed, as it
turns out, we shall now have a {\it unique surjective algebra homomorphism}

\bigskip
(4.2.3) \quad $ st : Gal_{\bf Q} ( 0 ) \longrightarrow {\bf R} $

\medskip
such that $s - st ( s ) \in mon_{\bf Q} ( 0 )$, for $s \in Gal_{\bf Q} ( 0 )$. It follows therefore that we have the {\it field
isomorphism}

\bigskip
(4.2.4) \quad $ {\bf R} ~=~ Gal_{\bf Q} ( 0 ) / mon_{\bf Q} ( 0 ) $

\medskip
this being quite likely the most simple and direct way of constructing the reals ${\bf R}$ from the rationals ${\bf Q}$.\\

It is important to note, as a second comment, that the above (4.2.1) - (4.2.4) generalizes to spaces which have a
uniform topological structure. \\
Indeed, let $X$ be such a space. Then $X$ can be {\it extended} in two ways. Namely, it can be embedded into its {\it
uniform completion}, $X \subseteq X^{\#}$, or, it can be embedded into its {\it nonstandard} extension, $X \subseteq~
^*X$. \\

Now, under rather general conditions, see Fenstad \& Nyberg, one can also construct a {\it bounded} part $Gal_X$ of
$X$, which is similar to $Gal ( 0 )$ or $Gal_{\bf Q} ( 0 )$, in the case of ${\bf R}$ and ${\bf Q}$, respectively.
Furthermore, one can obtain a commutative diagram

\begin{math}
\setlength{\unitlength}{1cm}
\thicklines
\begin{picture}(15,6)
\put(2,2.5){$X$}
\put(2.7,2.8){\vector(1,1){2}}
\put(4.9,5){$Gal_X ~\subseteq~ ^*X$}
\put(2.7,2.4){\vector(1,-1){2}}
\put(5.2,0){$X^{\#}$}
\put(5.4,4.65){\vector(0,-1){4}}
\put(5.7,2.5){$\sigma$}
\put(0,2.5){$(4.2.5)$}
\end{picture}
\end{math} \\

\medskip
with the following property. Let $Y$ be any {\it complete} uniform space and let $f : X \longrightarrow Y$ a uniformly
continuous mapping which has a uniformly continuous extension $f^{\#} : X^{\#} \longrightarrow Y$. Then we have

\bigskip
(4.2.6) \quad $ st_Y ( ^*f ( x ) ) ~=~ f^{\#} ( \sigma ( x ) ),~~~ x \in Gal_X $

\medskip
where $st_Y : Gal_Y \longrightarrow Y$ is the standard part mapping in the nonstandard extension $Y \subseteq~ ^*Y$,
while $^*f :~ ^*X \longrightarrow~ ^*Y$ is the nonstandard extension of $f : X \longrightarrow Y$. \\
The relation (4.2.6) can be follwed on the diagram

\begin{math}
\setlength{\unitlength}{1cm}
\thicklines
\begin{picture}(15,12)
\put(2,8.5){$X$}
\put(2.7,8.8){\vector(1,1){2}}
\put(4.9,11){$Gal_X ~\subseteq~ ^*X$}
\put(2.7,8.4){\vector(1,-1){2}}
\put(5.2,6){$X^{\#}$}
\put(5.4,10.65){\vector(0,-1){4}}
\put(5.7,8.4){$\sigma$}
\put(2.1,3.5){$Y$}
\put(2.2,8.2){\vector(0,-1){4}}
\put(1.8,6){$f$}
\put(5,5.8){\vector(-1,-1){2}}
\put(4.2,4.6){$f^{\#}$}
\put(0,6){$(4.2.7)$}
\put(5.2,0.7){$Gal_Y ~\subseteq~ ^*Y$}
\put(7.1,10.55){\vector(0,-1){9.2}}
\put(7.4,6.1){$^*f$}
\put(2.9,3.2){\vector(1,-1){2}}
\put(4.5,0.9){\vector(-1,1){2}}
\put(2.8,1.6){$st_Y$}
\end{picture}
\end{math} \\

\medskip
The above can present a certain {\it explanation} for the frequent effectiveness of Nonstandard Analysis in various
branches of mathematics. Indeed, both extensions $X^{\#}$ and $^*X$ have advantages and disadvantages. \\

For instance, the topological extensions  $X^{\#}$ may perform poorly when $X$ have algebraic structures and we want
to extended them to $X^{\#}$. On the ohter hand, the nonstandard extensions $^*X$ will usually perform much better
when such algebraic extension are required. \\
Indeed, the mappings $\sigma : Gal_X \longrightarrow X^{\#}$ happen to collapse part of the algebraic structures on
$Gal_X$ or on $^*X$, the reason for that being that the topological completions $X^{\#}$ prove to be too small.
Certainly, it is precisely the {\it largeness} of the nonstandard extensions  $^*X$ which allows them the advantages
offered by Nonstandard Analysis. \\

As far as the nonstandard extensions $^*X$ are conscerned, they do not typically admit satisfactory topological
structures, see Zakon. The main reason for that is the presence of {\it infinitesimals}, as seen for instance in 8 in
Proposition 3.6.1. \\
Yet, topological type structures can be treated in a variety of other ways in Nonstandard Analysis, as seen for instance
with the few cases in section 9, in chapter 3. \\ \\

\chapter{Superstructures and their Languages}

With this chapter, we enter a {\it second} and, shall we say, full fledged stage in the study of Nonstandard Analysis.
Indeed, in chapters 2 and 3, with only the help of usual mathematics and of the rock bottom simple formal languages
employed there, we could not do but very little, when compared to the full power of Nonstandard Analysis. \\

Fortunately however, in order to be able to use that full power one need not venture too far in Mathematical Logic, since
one can remain within the framework of Superstructures and their formal languages, as shown by A Robinson and E
Zakon, back in 1969. \\

What we shall present in this chapter, based on Superstructures and their formal languages, contains as some rather
simple cases the nonstandard results of chapters 2 and 3. This fact can easily be seen upon a direct verification which
is useful to make step by step. We shall present a few of these steps at the end of section 4. \\ \\

{\bf 1. Superstructures} \\ \\

Without getting by that involved in complications, and in fact on the contrary and in a way which helps to clarify the
issues, we shall no longer restrict ourselves to studying only a nonstandard extension of the set ${\bf R}$ of real
numbers. Instead, we shall start with an arbitrary nonvoid set $X$ which contains the real numbers, that is

$$ {\bf R} ~\subseteq~ X $$

\medskip
and we shall present ways to construct nonstandard extensions of $X$. Such extensions are far from being unique, and
as seen related to the issue of {\it saturation}, for instance, the lack of uniqueness of nonstandard extensions can have
its own advantages. By the way, as we shall see, such a lack of uniqueness is already manifested when we construct
nonstandard extensions of the real numbers ${\bf R}$. \\

When doing usual mathematics related to an arbitrary set $X$, we may as well have to deal in addition to its elements,
also with its subsets, sets of such subsets, functions of one or several variables defined on $X$, sets of such
functions, and so on. \\

In this way, it is natural to think about a Superstructure $V ( X )$ associated with $X$ in such a way that it is the {\it
smallest} set which contains the above, and other similar, yet possibly more complex kind of objects related to $X$.
Therefore, we give the following  \\ \\

{\bf Definition 5.1.1} \\

Let $X$ be a nonvoid set, then we define inductively

\bigskip
(5.1.1) \quad $ V_0 ( X ) ~=~ X,~~~ V_{n + 1} ( X ) ~=~ V_n ( X ) \bigcup {\cal P} (~ V_n ( X ) ~),~~~ n \in {\bf N} $

\medskip
where for a set $Y$ we denote as usual by ${\cal P} ( Y )$ its power set, that is, the set of all the subsets of $Y$. \\

The {\it Superstructure} over $X$ is defined now by the {\it countable} union

\bigskip
(5.1.2) \quad $ V ( X ) ~=~ \bigcup_{n \in {\bf N}}~  V_n ( X ) $

\hfill $\Box$ \\ \\

Obviously, if $X \subseteq Y$ then $ V ( X ) \subseteq~V ( Y )$. Also we have

\bigskip
(5.1.3) \quad $ \begin{array}{l}
                                         X ~=~ V_0 ( X ) ~\subseteq~ V_1 ( X ) ~\subseteq~ ~.~.~.~ ~\subseteq~ V_n ( X )
                                                                                                       ~\subseteq~ ~.~.~.~  V ( X) \\ \\
                                         X ~=~ V_0 ( X ) ~\in~ V_1 ( X ) ~\in~ ~.~.~.~ ~\in V_n ( X )
                                                                                                       ~\in~ ~.~.~.~  V ( X) \\ \\
                                         X ~=~ V_0 ( X ),~ V_1 ( X ),~.~.~.~,~  V_n ( X ), ~.~.~.~ \in V ( X)
                        \end{array} $

\medskip
The elements of $X$ are called {\it individuals}, and it will be convenient to assume that they are {\it not} sets in
$V ( X )$, in other words

\bigskip
(5.1.4) \quad $ x \in X,~~ y \in V ( X ) ~~~\Longrightarrow ~~~ y \notin x $

\medskip
The elements of $V ( X )$ are called {\it entities}. This of course includes the individuals as a particular case, see
(5.1.3). \\
For $n \in {\bf N}$, the entities in $V_{n + 1} ( X ) \setminus V_n ( X )$ are called having {\it rank} n + 1, with the elements
of $X$, that is, the individuals, having by definition {\it rank} 0. \\ \\

{\bf Note 5.1.1} \\

1.~~ The assumption (5.1.4) is needed in order to avoid situations such as the following one. Suppose given $x \in X$
and $y \in V ( X )$, for which we have $y \in x$ in the specific form $x = \{~ y ~\}$, with $y \in X$. Then clearly $x \in
{\cal P} ( X )$. Thus $x \in X \bigcap {\cal P} ( X )$. \\
It follows that $x$ has rank 0, since $x \in X$, nevertheless at the same time $x \in {\cal P} ( X )$. \\

On the other hand, as we can note, with the sets usually encountered in mathematics, we have $X \bigcap {\cal P} ( X )
= \phi$. \\
For instance, we certainly prefer to consider that ${\bf R} \bigcap {\cal P} ( {\bf R} ) = \phi$, in other words,
subsets of ${\bf R}$ are {\it not} supposed to be real numbers. In this respect, for instance, we distinguish between the
real number $x \in {\bf R}$ and the one point set of real numbers $\{~ x ~\} \subseteq {\bf R}$, for which we have $\{~ x ~\}
\in {\cal P} ( {\bf R} )$. \\
Yet, here we have to exhibit a certain care. Indeed, according for instance to the Dedekind construction of ${\bf R}$,
every real number $x \in {\bf R}$ is a pair $( A, B )$, called a {\it cut}, of sets of rational numbers $A, B \subset
{\bf Q}$, such that $A \bigcup B = {\bf Q}$, and $a < b$, for every $a \in A$ and $b \in B$. Therefore, we can make the
identification $x = \{~ A, B ~\} \in {\cal P} ( {\cal P} ( {\bf R} ) )$, and thus obtain that ${\bf R} \bigcap {\cal P}
( {\cal P} ( {\bf R} ) ) \neq \phi$. \\
A similar situation occurs if we consider the Cauchy-Bolzano construction of the real numbers. \\
Such situations happen in various branches of mathematics. \\
Consequently, what is usually done is that, once ${\bf R}$, or other mathematical entities are constructed, they are
treated as sets which do satisfy (5.1.4). \\

2.~~ A third aspect of Nonstandard Analysis which may put off usual mathematicians is the need to work across an
{\it infinity}, even if only countable, of ranks, when dealing with Superstructures $V ( X )$. Indeed, in usual mathematics
one is not much concerned with ranks, let alone with infinitely many of them, since typically, one is limited to the first
three or at most four ranks only, as one deals with elements $x$, sets $X$, families of sets ${\cal X}$, and on occasion,
with families of such families of sets. \\
However, as seen in the sequel, dealing with the countable infinity of ranks involved in the Superstructures $V ( X )$,
will not lead to more complications than usual induction can hadle.

\hfill $\Box$ \\

Let us start by noting that, contrary to what may appear to be the case, Superstructures are {\it not} unduly large sets.
Indeed, although of course we have, as with any set, the inclusion $V ( X ) \subseteq V ( X )$, nevertheless, it follows
easily, see (5.1.7) below, that

\bigskip
(5.1.5) \quad $ V ( X ) \notin V ( X ) $

\medskip
and in general

\medskip
\begin{math}
\setlength{\unitlength}{1cm}
\thinlines
\begin{picture}(15,0.9)
\put(0,0.3){$(5.1.6)~~~~A ~\subseteq~ V ( X ) ~~~\Longrightarrow~~ A \in V ( X )$}
\put(0.1,0.3){$~~~~~~~~~~~~~~~~~~~~~~~~~~~~~~~~~/$}
\end{picture}
\end{math}

which means that $V ( X )$ is {\it not} closed under arbitrary subsets. Indeed, let $a \in X$, then obviously

$$ A ~=~ \{~ a,~ \{~ a ~\},~ \{~ \{~ a ~\} ~\},~ \{~ \{~ \{~ a ~\} ~\} ~\},~ .~.~.~ \} ~\subseteq~ V ( X ) $$

\medskip
since $~a \in X = V_0 ( X ),~ \{~ a ~\} \in V_1 ( X ),~ \{~ \{~ a ~\} ~\} \in V_2 ( X ),~  \{~ \{~ \{~ a ~\} ~\} ~\} \in V_3 ( X ),~.~.~ $, while in
view of (5.1.2), this infinite sequence of relationships implies that

$$ A \notin V ( X ) $$

\medskip
In this regard, it is easy to see that we have the following {\it characterization} of closure of $V ( X )$ under subsets,
given by a {\it bounded rank} condition, namely

\bigskip
(5.1.7) \quad $ \begin{array}{l}
                                \forall~ A \subseteq V ( X ) ~: \\ \\
                                ~~~ A \in V ( X ) ~~~\Longleftrightarrow~~~ A ~\subseteq~ V_n ( X ),~~~ \mbox{for some}~ n \in {\bf N}
                         \end{array} $

\medskip
Of course, we always have the converse implication in (5.1.6), in the form

\bigskip
(5.1.8) \quad $ A \in V ( X ) \setminus~ X ~~~\Longrightarrow~~~ A ~\subseteq~ V ( X ) $

\medskip
which follows directly from (5.1.7) and (5.1.3) \\

Also, in view of (5.1.7), we have the {\it transitivity} property

\bigskip
(5.1.9) \quad $ a \in A \in V ( X ) ~~~\Longrightarrow~~~ a \in V ( X ) $

\medskip
We respect to such transitivities, we can note that the following ones also hold

$$ \begin{array}{l}
                                a \in A \subseteq V ( X ) ~~~\Longrightarrow~~~ a \in V ( X ) \\ \\
                                B \subseteq A  \in V ( X ) ~~~\Longrightarrow~~~ B \in V ( X ) \\ \\
                                B \subseteq A  \subseteq V ( X ) ~~~\Longrightarrow~~~ B \subseteq V ( X )
      \end{array} $$

\medskip
And now the question arises : \\

What happens with {\it Cartesian products} related to $X$ ? \\

Indeed, such products do not appear explicitly in the definition (5.1.2) of Superstructures. Yet we obviously need them,
among others, in order to deal with relations, and in particular, functions connected with $X$. Fortunately, as seen next,
we do not have a problem in this regard. \\

Let us start with the case of a Cartesian product with 2 factors only. \\
Given $x,~ y \in V ( X )$, we can define the ordered pair

$$ < x, y > ~=~\{~ x, \{~ x, y ~\} ~\} $$

\medskip
In view of (5.1.2), we have $x,~ y \in V_n ( X )$, for a certain $n \in {\bf N}$. Therefore ~$< x, y >~ \in V_{n + 2} ( X )$, since
obviously $\{~ x ~\},~ \{~ x, y ~\} \in V_{n + 1} ( X )$. \\
The point in such a definition of a pair is that it is in terms of sets. Also, it is a correct definition, since clearly

$$ < x, y > ~=~ < u, v > ~~~\Longleftrightarrow~~~ x ~=~ u ~~~ \mbox{and}~~~ y ~=~ v $$

\medskip
It follows now that for $A,~ B \subseteq V_n ( X )$, we can define

\bigskip
(5.1.10) \quad $ A \times B ~=~ \{~ < x, y > ~~|~~ x \in A,~~ y \in B ~\} $

\medskip
and then we have

\bigskip
(5.1.11) \quad $ A \times B ~\subseteq~ V_{n + 2} ( X ),~~~ A \times B \in V_{n + 3} ( X ) $

\medskip
For m $\geq$ 2, we give the following definition for an m-tuple

$$ < x_1, ~.~.~.~ , x_m > ~=~ \{~ < 1, x_1 >, ~.~.~.~ , < m, x_m > ~\} $$

\medskip
where $x_1, ~.~.~.~ , x_m \in V ( X )$. Then again

$$ < x_1, ~.~.~.~ , x_m > ~=~ < y_1, ~.~.~.~ , y_m > ~~~\Longleftrightarrow~~~~ x_1 ~=~ y_1, ~.~.~.~ , x_m ~=~ y_m $$

\medskip
In this way, for $A_1, ~.~.~.~ , A_m \in V ( X )$, we can define the Cartesian product

$$ A_1 \times ~.~.~.~ \times A_m ~=~ \{~ < x_1, ~.~.~.~ , x_m > ~~|~~ x_1 \in A_1, ~.~.~.~ , x_m \in A_m ~\} \in V ( X ) $$

\medskip
Now, an m-ary {\it relation} $P$ on $A_1 \times ~.~.~.~ \times A_m$ is any subset $P$ of $A_1 \times ~.~.~.~ \times A_m$,
hence $P \in V ( X )$. \\

Given a binary relation $P \subseteq A_1 \times A_2$, then as usual, we define

$$ \begin{array}{l}
               domain~ P ~=~ \{~ a_1 \in A_1 ~~|~~ < a_1, a_2 >~ \in P~~~\mbox{for some}~ a_2 \in A_2 ~\} \\ \\
               range~ P ~=~ \{~ a_2 \in A_2 ~~|~~ < a_1, a_2 >~ \in P~~~\mbox{for some}~ a_1 \in A_1 ~\}
   \end{array} $$

\medskip
and for $B_1 \subseteq A_1,~ B_2 \subseteq A_2$, we define

$$ \begin{array}{l}
               P [ B_1 ] ~=~ \{~ b_2 \in A_2 ~~|~~ < b_1, b_2 > \in P~~~\mbox{for some}~ b_1 \in A_1 ~\} \\ \\
               P^{- 1} [ B_2 ] ~=~ \{~ b_1 \in A_1 ~~|~~ < b_1, b_2 > \in P~~~\mbox{for some}~ b_2 \in A_2 ~\}
   \end{array} $$

\medskip
Let us note that $domain~ P \subseteq A_1$, and $A_1 \in V_n ( X )$, for a certain \\ $n \in {\bf N}$, hence $domain~ P \in
V_{n + 1} ( X )$. Similarly, we have $range~ P,~ P [ B_1 ],~ \\ P^{- 1} [ B_2 ] \in V_{n + 1} ( X )$. \\

Functions being particular cases of binary relations, their treatment follows form the above. \\

For instance, a sequence $( x_n )_{n \in {\bf N}}$ of elements in $X$ can be seen as a function $x : {\bf N}
\longrightarrow X$, thus as a certain subset of ${\bf N} \times X$. In this way, from the above we obtain

$$ ( x_n )_{n \in {\bf N}} ~\subseteq~ V_2 ( X ),~~~ ( x_n )_{n \in {\bf N}} \in V_3 ( X ) $$

\medskip
hence $( x_n )_{n \in {\bf N}}$ has rank 3. \\

Cartesian products of an arbitrary number of factors can be defined as follows with the use of functions. Namely, we
start with the usual  general set theoretical formula

$$ \prod_{i \in I}~ A_i ~=~ \{~ f : I ~\longrightarrow~ (~ \bigcup_{i \in I}~A_i ~) ~~|~~ f (i ) \in A_i,~~ i\in I ~\}
                                                                                                ~~\subseteq~~ I ~\times~ (~ \bigcup_{i \in I}~A_i ~)$$

\medskip
where $I$ and $A_i$, with $i \in I$, are arbitrary sets. \\
Let us now assume that for a certain $n \in {\bf N}$, we have $A_i \subseteq V_n ( X )$, with $i \in I$, then (5.1.1) gives
that

$$ \bigcup_{i \in I}~A_i \subseteq V_n ( X ) $$

\medskip
which is valid even if $I \notin V ( X )$. \\

Now we further assume that we have as well $I \subseteq V_n ( X )$, for the same $n \in {\bf N}$ as above. Then in view
of (5.1.11), we have

\bigskip
(5.1.12) \quad $ \prod_{i \in I}~ A_i ~\subseteq~ V_{n + 2} ( X ),~~~ \prod_{i \in I}~ A_i \in V_{n + 3} ( X )$

\medskip
Here we should note that, according to (5.1.3), it is without loss of generality to assume that $I$, together with all the
$A_i$, with $i \in I$, are subsets of the same $V_n ( X )$, as long as we have already assumed $A_i \subseteq V_n
( X )$, with $i \in I$, and in addition, we also asked that $I \in V ( X )$. \\

The following easy to prove properties of $V ( X )$, for an arbitrary nonvoid $X$ and $n \in {\bf N}$, thus $n \geq 0$, will
be useful : \\

1.~~ $ V_{n + 1} ( X ) ~=~ X \bigcup {\cal P} (~ V_n ( X ) ~) $ \\

2.~~ $ a_1, ~.~.~.~ , a_m \in  V_n ( X ) ~~~\Longrightarrow~~~ \{~ a_1,~.~.~.~ , a_m ~\} \in V_{n + 1} ( X ) $ \\

3.~~ $ A_1, ~.~.~.~ , A_m \in V_n ( X ) \setminus X ~~~\Longrightarrow~~~ A_1 \bigcup  ~.~.~.~  \bigcup A_m\in V_n ( X ) $ \\

4.~~ $ A ~\subseteq~ B \in V_n ( X ) ~~~\Longrightarrow~~~ A \in V_{n + 1} ( X ) $ \\

5.~~ $ A \in V_n ( X ) \setminus X ~~~\Longrightarrow~~~ {\cal P} (~ A ~) \in V_{n + 2} ( X ) $ \\

6.~~ $ a \in A \in V_{n + 1} ( X ) ~~~\Longrightarrow~~~ a \in V_n ( X ) $ \\

7.~~ $ A \in V_{n + 1} ( X ),~ A \bigcap X ~=~ \phi  ~~~\Longrightarrow~~~ \bigcup_{B \in A}~ B \in V_n ( X ) $

\medskip
Finally, let $A_i \in V_n ( X ) \setminus X$, with $i \in I$, where $I$ is an arbitrary set, then \\

8.~~ $ \bigcup_{i \in I}~ A_i \in V_n ( X ) $ \\ \\

{\bf 2. Languages for Superstructures} \\ \\

Given a nonvoid set $X$ and the corresponding Superstructure $V ( X )$, we define a language ${\cal L}_X$ for $V ( X )$
as follows, see for comparison Definition 3.2.2. \\ \\

{\bf Definition 5.2.1} \\

The language ${\cal L}_X$ contains the following seven categories of {\it symbols} : \\

1. Five logical connectives : $\neg$, $\wedge$, $\vee$, $\rightarrow$ and $\leftrightarrow$ which are interpreted
as "not", "and", "or", "implies" and "if and only if". \\

2. Two quantifier symbols : $\forall$ and $\exists$ which are interpreted as "for all" and "there exists". \\

3. Six parentheses : [ , ] , ( , ) , $<$ , $>$ . \\

4. Variable symbols : a countable set x, y, ~.~.~.~ \\

5. Equality symbol : a binary relation ~=~ which is interpreted as "equals". \\

6. Predicate symbol : a binary relation $\in$ which is interpreted as "is an element of". \\

7. Constant symbols : at least one symbol $\underline{a}$ for every entity $a \in V ( X )$. \\

The {\it formulas} of ${\cal L}_X$ are defined inductively, as follows. \\

8. Atomic formulas :

$$ \begin{array}{l}
                               x ~\in~ y \\
                               x ~=~ y \\
                              < x_1, ~.~.~.~ , x_n >~ \in~ y \\
                              < x_1, ~.~.~.~ , x_n > ~= ~ y \\
                              < < x_1, ~.~.~.~ , x_n >, x >~ \in~ y \\
                              < < x_1, ~.~.~.~ , x_n >, x >~ =~ y
      \end{array} $$

\medskip
where $x_1, ~.~.~.~ , x_n, x,  y$ are variable or constant symbols. \\

9. If $\Phi,~ \Psi$ are formulas in ${\cal L}_X$, then

$$ \neg \Phi,~~~~ \Phi ~\wedge~ \Psi,~~~~ \Phi ~\vee~ \Psi,~~~~ \Phi ~\rightarrow~ \Psi,~~~~ \Phi ~\leftrightarrow~ \Psi $$

\medskip
are also formulas in ${\cal L}_X$. \\

10. If $\Phi$ is a formula in ${\cal L}_X$ and $x$ is a variable symbol, such that $\Phi$ does {\it not} contain
expressions of the form $(~ \forall~ x \in z ~)~ \Psi$ or $(~ \exists~ x \in z ~)~ \Psi$, for some formula $\Psi$ in
${\cal L}_X$, then

$$ (~ \forall~ x \in y ~)~ \Phi,~~~~ (~ \exists~ x \in y ~)~ \Phi $$

\medskip
are also formulas in ${\cal L}_X$, for every variable or constant symbol $y$. \\

Finally, the {\it sentences} of ${\cal L}_X$, which are particular cases of formulas in ${\cal L}_X$, are defined as
follows. \\

If $\Phi$ is a formula in ${\cal L}_X$ and $x$ is a variable symbol, we say that $x$ {\it occurs in} $\Phi$ {\it in the scope
of a quantifier}, if and only if $x$ appears in $\Phi$ in some expression of the form $(~ \forall~ x \in y ~)~ \Psi$, or $(~
\exists~ x \in y ~)~ \Psi$, where $\Psi$ is a formula in ${\cal L}_X$. \\
Obviously, a variable symbol $x$ can appear in a formula $\Phi$ both in the scope of a quantifier and outside of any
quantifier, such as for instance in the formula

$$ \Xi ~=~ [~ (~ \forall~ x \in y ~)~ (~ x > 1 ~) ~] \wedge (~ x \in z ~) $$

\medskip
where the first occurrence of $x$, that is, in $(~ \forall~ x \in y ~)~ (~ x > 1 ~)$ is in the scope of a quantifier, namely,
$\forall$, while the second one in $(~ x \in z ~)$ is not in the scope of any quantifier. \\

A variable symbol $x$ is called {\it bound} in a formula $\Phi$, if it occurs in it in the scope of a quantifier. Otherwise, a
variable symbol $x$ in a formula $\Phi$ is called {\it free}. \\
For instance, in the above formula $\Xi$, the first occurrence of $x$ is bound, while the second one is free. \\

And now : \\

11. A {\it sentence} of ${\cal L}_X$ is every formula $\Phi$ in ${\cal L}_X$ which does {\it not} contain free variable
symbols, that is, all variable symbols in it are bound.

\hfill $\Box$ \\ \\

The distinction between formulas in general, and their particular cases, the sentences, will be important - and will also
become obvious - in the sequel. Suffice it here to mention that a formula which is not a sentence, and thus contains at
least one free variable, makes a statement about the possible values - when interpreted - of its free variables. On the
other hand, a sentence makes a statement about the structure described by the language to which that statement may
happen to belong. \\

Let us note that within the {\it simple languages} given in Definition 3.2.2, the {\it atomic sentences} which contain
variable symbols are {\it not} sentences in the above sense, since they do not contain quantifiers, therefore, each
variable symbol which they contain is free. On the other hand, the {\it compound sentences} are always sentences in
the above sense, since all the variable symbols which they contain must by definition be bound. \\

Here we give a few simple examples which point out this difference between formulas in general, and on the other
hand, their particular cases, the sentences, and we do so from the point of view of the structure of the respective
formulas, as defined above in 1 - 11 in Definition 5.2.1. For instance, the expression

$$ (~ \forall~ x \in a ~)~(~ \exists~ y \in x ~)~(~ y \in b ~) $$

\medskip
where $a, b$ are constant symbols, is a formula in ${\cal L}_X$, which is also a sentence in ${\cal L}_X$, as both
variable symbols $x$ and $y$ which appear in it are {\it bound}. On the other hand, the expression

$$ (~ \exists~ y \in x ~)~(~ y \in b ~) $$

\medskip
is a formula in ${\cal L}_X$, which is {\it not} also a sentence in ${\cal L}_X$, since it contains the {\it free}
variable symbol $x$. \\
However, as noted before, a certain care is needed here. For instance, the expression

\bigskip
(5.2.1) \quad $ [~ (~ \exists~ x \in a ~)~(~ x = b ~) ~] ~\wedge~ (~ x = c ~) $

\medskip
where $a, b, c$ are constant symbols, is obviously a formula in ${\cal L}_X$. Yet this formula is {\it not} as well a
sentence in ${\cal L}_X$, since the second occurrence of the variable symbol $x$, namely, in $ (~ x = c ~)$ is free, and
thus it is not bound. Indeed, the scope of the quantifier $(~ \exists~ x \in a ~)$ does not extend over this second
occurrence of $x$ as well. \\

Let us note that a bound variable symbol in a formula can be replaced by any other variable symbol, and the result is a
formula which is identical. For instance, in (5.2.1) we can replace the bound occurrences of $x$ with $y$, and thus
obtain

$$ [~ (~ \exists~ y \in a ~)~(~ y = b ~) ~] ~\wedge~ (~ x = c ~) $$

\medskip
This is similar with replacing the index $i$ in the expression $\Sigma_{~1 \leq i \leq n}~ a_i$ with the index $j$, in which
case we obtain the identical expression $\Sigma_{~1 \leq j \leq n}~ a_j$. \\

In fact, the second occurrence of $x$ in (5.2.1), and which is free, namely in $(~ x = c ~)$, can also be replaced by any
other variable symbol, for instance, $y$, and we obtain the formula in ${\cal L}_X$, given by

$$  [~ (~ \exists~ x \in a ~)~(~ x = b ~) ~] ~\wedge~ (~ y = c ~) $$

\medskip
which, as we shall see in the next section, is {\it identical} with the one in (5.2.1), when subjected to interpretation. \\ \\

{\bf Note 5.2.1} \\

1~~~ We should note that according to 7 in Definition 5.2.1, the language ${\cal L}_X$ contains {\it uncountably} many
{\it constant symbols}, namely, at least one for every entity in $V ( X )$. \\

2~~~ The above languages ${\cal L}_X$ associated with Superstructures $V ( X )$ in Definition 5.2.1 are obviously more
complex than the simple languages in Definition 3.2.2, as this time more logical connectives are allowed, and also,
both quantifiers $\forall$ and $\exists$ are included. \\
There are however two exceptions. First, the simple languages contain {\it terms}, see 8 - 10 in the Definition 3.2.2,
while the languages ${\cal L}_X$ in Definition 5.2.1 do {\it not} contain terms. Second, the atomic formulas in 8 in
Definition 5.2.1 are few and rather simple, when compared with those in Definition 3.2.2. \\
This fact will have certain consequences. For instance, sometime formulas, and thus in particular, sentences in
${\cal L}_X$ will look more involved than their corresponding forms in simple languages, when such corresponding
forms may happen to exist. \\
Another difference with the languages ${\cal L}_X$ associated with Superstructures $V ( X )$ in Definition 5.2.1 is that,
they do {\it not} contain entities which are {\it not interpretable}, this being unlike with the case of simple languages in
Definition 3.2.2, where certain terms and atomic sentences may fail to be interpretable, see 2 and 3 in Definition
3.3.1. \\ \\

{\bf Convention 5.2.1} \\

As an extension of the spirit of the Convention Summary 3.4.1, we shall make the following further simplification. We
denote any given entity $a \in V ( X )$ by itself, that is, by $a$, as the corresponding {\it constant symbol} in the
language ${\cal L}_X$, see 7 in Definition 5.2.1. \\ \\

{\bf 3. Interpretations in Superstructures} \\ \\

As already mentioned in section 3, in chapter 3, the aim of a formal language is to obtain a rigorous and systematic
method for identifying and studying those of its {\it sentences} which are {\it true} in a certain given theory. And
the aim of the {\it interpretation} of a formal language is to give one of the possible ways to find out which
sentences are true, and which, on the contrary, are not true. \\

Here it should be recalled that the quality of being {\it true}, or on the contrary, of {\it not} being true, only
applies to the {\it sentences} of a formal language, and not also to other auxiliary entities of such a language. Such
auxiliary entities of a formal language and which are not sentences - like for instance, those under 1 - 10 in the
above definition of ${\cal L}_X$ - are only subjected to {\it interpretation}, without however acquiring by that the
quality of being true, or on the contrary, of not being true. \\

Let us now show how to {\it interpret} the language ${\cal L}_X$ in the Superstructure $V ( X )$, and let us recall that
in such languages ${\cal L}_X$ associated with Superstructures $V ( X )$ in Definition 5.2.1, {\it all} entities are
{\it interpretable}, see 2 in Note 5.2.1 above. \\
Here, we can take Definition 3.3.1 as a simple prototype for interpretation of languages. \\ \\

{\bf Definition 5.3.1} \\

Once again, {\it variable symbols} are {\it not} interpreted, more precisely, they remain the same with themselves
under interpretation, while any {\it constant symbol} $a$, see the above Convention 5.2.1, which corresponds to an
entity $a \in V ( X )$, is interpreted as this entity $a$. \\

Now, {\it formulas} in 8 - 10 in Definition 5.2.1 are interpreted in the obvious manner, namely, by interpreting their
constant symbols, and keeping the rest of their constituents unchanged. \\

At last, when we come to the {\it sentences} of the language ${\cal L}_X$ associated with the Superstructure $V ( X )$,
after their interpretation, we face the additional issue of whether they are {\it true}, or {\it false}.  In this
regard, we proceed as follows. \\

{\it Atomic formulas} which are sentences, see 8 and 11 in Definition 5.2.1, namely, $a \in b,~ < a_1, ~.~.~.~ , a_n >~
\in b,~ < < a_1, ~.~.~.~ , a_n >, c >~ \in b,~ a = b,~ < a_1, ~.~.~.~ , a_n >~ = b$ and $< < a_1, ~.~.~.~ , a_n >, c >~
= b$ are true in $V ( X )$, if and only if the respective relations hold in $V ( X )$, when the constant symbols
in these relations are replaced with their corresponding entities in $V ( X )$, see 7 in Definition 5.2.1. \\

Further, if $\Phi$ and $\Psi$ are sentences in ${\cal L}_X$, then \\

$ \neg \Phi ~~\mbox{is true in}~~ V ( X ) ~~~\Longleftrightarrow~~~ \Phi ~~\mbox{is not true in}~~ V ( X ) $ \\

$ \Phi \wedge \Psi ~~\mbox{is true in}~~ V ( X ) ~~~\Longleftrightarrow~~~ \Phi ~~\mbox{and}~~ \Psi ~~\mbox{are both
                 true in}~~ V ( X ) $ \\

$ \Phi \vee \Psi ~~\mbox{is true in}~~ V ( X ) ~~~\Longleftrightarrow~~~ \Phi ~~\mbox{or}~~ \Psi ~~\mbox{is true in}~~
                  V ( X ) $ \\

$ \Phi \rightarrow \Psi ~~\mbox{is true in}~~ V ( X ) ~~~\Longleftrightarrow~~~ \mbox{either}~~ \Psi ~~\mbox{is not
                   true in}~~ V ( X ), $ \\
$~~~~~~~~~~~~~~~~~~~~~~~~~~~~~~~~~~~~~~~~~~~~~~~\mbox{or}~~ \Psi ~~\mbox{is true in}~~ V ( X ) $ \\

$ \Phi \leftrightarrow \Psi ~~\mbox{is true in}~~ V ( X ) ~~~\Longleftrightarrow~~~ \mbox{either both}~~ \Phi
                   ~~\mbox{and}~~ \Psi ~~\mbox{are true in}~~ V ( X ), $ \\
$~~~~~~~~~~~~~~~~~~~~~~~~~~~~~~~~~~~~~~~~~~~~~~~~\mbox{or both} \Phi
                   ~~\mbox{and}~~ \Psi ~~\mbox{are false in}~~ V ( X ) $ \\

Let now $\Phi ( x )$ be a formula in ${\cal L}_X$ which contains $x$ as the only free variable, and given any two
constant symbols $a, b$. \\
Then the sentence $(~ \forall~ x \in b ~)~ \Phi ( x )$ is true in $V ( X )$, if and only if the sentence $\Phi ( a )$
is true in $V ( X )$, for every constant symbol $a \in b$. \\
On the other hand, the sentence $(~ \exists~ x \in b ~)~ \Phi ( x )$ is true in $V ( X )$, if and only if the sentence
$\Phi ( a )$ is true in $V ( X )$, for at least one constant symbol $a \in b$.

\hfill $\Box$ \\

For convenience, we shall often depart from the rigors in the Definition 5.2.1 concerning the notations in the languages
${\cal L}_X$ associated with Superstructures $V ( X )$, and use instead the more simple and familiar notations. However,
with minimal care, one can always reconstitute the respective rigorous notations. \\ \\

{\bf 4. Monomorphisms of Superstructures} \\ \\

Our aim, achieved constructively in the next section, is : \\

- to associate with each set $X$ its {\it nonstandard} extension $^*X$, \\

- to associate with the Superstructure $V ( X )$ corresponding to $X$, the \\
\hspace*{0.12cm}{ \it nonstandard} Superstructure $V ( ^*X )$ which corresponds to $^*X$, and \\
\hspace*{0.12cm} provide this association by a suitable mapping

\bigskip
(5.4.1) \quad $ ^*(~) : V ( X ) ~\longrightarrow V ( ^*X ) $

\medskip
called {\it monomorphism}, see Definition 5.4.2 below. \\

In other words, we want to construct a commutative diagram

\begin{math}
\setlength{\unitlength}{1cm}
\thicklines
\begin{picture}(15,6.5)
\put(0,3){$(5.4.2)$}
\put(3,5){$X$}
\put(6.2,5.4){$\subseteq$}
\put(3.8,5.1){\vector(1,0){5.2}}
\put(9.5,5){$^*X$}
\put(2.5,3){$\subseteq$}
\put(3.15,4.5){\vector(0,-1){3}}
\put(10.2,3){$\subseteq$}
\put(9.85,4.5){\vector(0,-1){3}}
\put(2.75,0.8){$V ( X )$}
\put(9.3,0.8){$V ( ^*X )$}
\put(4.2,0.9){\vector(1,0){4.6}}
\put(6.1,0.3){$^*(~)$}
\end{picture}
\end{math}

\medskip
with suitable properties to be specified in the sequel. \\

Here it is important to note two facts. First, the above inclusion $X \subseteq~ ^*X$ implies obviously the inclusion
$ V ( X ) ~\subseteq~ V ( ^*X )$, see (5.1.1), (5.1.2). However, this is inclusion, that is, the identity mapping
$V ( X ) \ni A ~\longmapsto~ A \in V ( ^*X )$, is {\it not} the relationship which is of main interest in Nonstandard
Analysis between $V ( X )$ and $V ( ^*X )$. Indeed, as we have seen in (2.7.8), the mapping

$$ {\bf R} ~\supseteq~ A ~~~\longmapsto~~~ ^*A ~\subseteq~ ^*{\bf R} $$

\medskip
is an identity, if and only if $A$ is a finite set. And we are interested in many more entities in $V ( {\bf R} )$, or for that
matter, in $V ( X )$, than finite subsets of ${\bf R}$. \\

Second, we are {\it not} interested in a diagram

\begin{math}
\setlength{\unitlength}{1cm}
\thicklines
\begin{picture}(15,6.5)
\put(0,3){$(5.4.2^*)$}
\put(3,5){$X$}
\put(3.8,5.1){\vector(1,0){5.2}}
\put(9.5,5){$^*X$}
\put(3.15,4.5){\vector(0,-1){3}}
\put(9.85,4.5){\vector(0,-1){3}}
\put(2.75,0.8){$V ( X )$}
\put(9.3,0.8){$^*V ( X )$}
\put(4.2,0.9){\vector(1,0){4.6}}
\end{picture}
\end{math}

Indeed, let us first of all note that, in terms of the mapping (5.4.1), $^*V ( X )$ {\it cannot} have any other meaning
but that of the {\it range} in $V ( ^*X )$ of the mapping (5.4.1), namely

\bigskip
(5.4.2$^{**}$) \quad $ ^*V ( X ) ~=~ \{~ ^*A ~~|~~ A \in V ( X ) ~\} ~=~ \bigcup_{n \in {\bf N}}~ ^*V_n ( X ) ~\subseteq~ V ( ^*X ) $

\medskip
since $V ( X ) \notin V ( X )$, see (5.1.5), therefore

\bigskip
(5.4.2$^{***}$) \quad $ ^*V ( X ) \notin V ( ^*X ) $

\medskip
In other words, $^*V ( X )$ is {\it not} an entity in, but {\it only} a subset of $V ( ^*X )$, see (5.1.6). \\

Furthermore, the nonstandard extension of $X$ is $^*X$. Therefore, just as the mathematics relevant to $X$ is supposed
to be described by the Superstructure $V ( X )$ associated with $X$, in the same way, the mathematics relevant to $^*X$
will be described by the Superstructure $V ( ^*X )$ associated with $^*X$, and {\it not} by $^*V ( X )$, which has and can
only have the meaning specified above. \\
Certainly, when we get involved in nonstandard extensions $^*X$, our primary aim is still the study of the mathematics
relevant to the initial sets $X$, that is, as described by $V ( X )$. In this way, the nonstandard approach is for us but a
{\it detour} into the richer and often more intuitive realms of $^*X$ and $V ( ^*X )$, realms from which, however, we
intend to come back to $X$ and $V ( X )$. This clearly points to the fact that we are interested in diagrams (5.4.2), and
{\it not} in diagrams (5.4.2$^*$). \\

In short, we are in fact interested in the mathematics described by

\begin{math}
\setlength{\unitlength}{1cm}
\thicklines
\begin{picture}(15,3)
\put(3,2){$X$}
\put(0,1.3){$( HOME )$}
\put(3.15,1.7){\vector(0,-1){0.8}}
\put(2.75,0.2){$V ( X )$}
\end{picture}
\end{math} \\

and we are only making a {\it detour} to the mathematics described by

\begin{math}
\setlength{\unitlength}{1cm}
\thicklines
\begin{picture}(15,3)
\put(4,2){$^*X$}
\put(0,1.3){$( DETOUR )$}
\put(4.3,1.7){\vector(0,-1){0.8}}
\put(3.75,0.2){$V ( ^*X )$}
\end{picture}
\end{math} \\

in order to obtain conjectures and/or proofs easier, after which we certainly intend to return $( HOME)$ with the results
we obtained in the detour. \\

As seen later in Theorem 5.7.1, $^*V ( X )$ is a {\it strict} subset of $V ( ^*X )$, and in fact, it is the set of all {\it internal}
sets in the vastly larger $V ( ^*X )$. \\

Now, with respect to the construction of the monomorphism (5.4.1), or equivalently, of the commutative diagram (5.4.2),
it is useful to approach the issue of monomorphisms in a more general setup, namely, between two arbitrary
Superstructures. \\

Let therefore $X$ and $Y$ be any two sets, such that

$$ {\bf R} ~\subseteq~ X ~\bigcap~ Y $$

\medskip
and let be given any {\it injective} mapping

\bigskip
(5.4.3) \quad $ ^*(~) : V ( X ) ~\longrightarrow~ V ( Y ) $

\medskip
In this general setting, we shall now define a corresponding mapping, called $^*$-{\it transform} and given by \\ \\

{\bf Definition 5.4.1} \\

The mapping

\bigskip
(5.4.4) \quad $ ^*(~) : {\cal L}_X ~\longrightarrow~ {\cal L}_Y $

\medskip
acts in the following very simple and direct manner. If $\Phi$ is a formula, or in particular, a sentence in ${\cal L}_X$,
then $^*\Phi$ is obtained from $\Phi$ by replacing every {\it constant symbol} $a$ in $\Phi$ with $^*a$.

\hfill $\Box$ \\

The deceptively simple way the $^*$-transform in (5.4.4) acts between the two respective languages corresponds in
fact to a powerful instrument in Nonstandard Analysis, as seen in the sequel. \\
Here we should point out that, with the above definition, the $^*$-{\it transform} in (5.4.4) does {\it not} make any of the
notational simplifications introduced in the Conventions in chapter 3. In other words, each and every constant symbol,
including those corresponding to the most familiar ones, be they certain numbers, relations, functions, operations, etc,
will now get the $^*$ in front of it, when we go from formulas in ${\cal L}_X$ to formulas in ${\cal L}_Y$. \\
The only exception, and simplification, will be that mentioned in Convention 5.4.1 below. \\

And now to the concept of {\it monomorphism} between two Superstructures. \\ \\

{\bf Definition 5.4.2} \\

The injective mapping (5.4.3) is a {\it monomorphism} between $V ( X )$ and $V ( Y )$, if and only if it satisfies the five
conditions : \\

1.~~~ $ ^*\phi ~=~ \phi $ \\

2.~~~ $ a \in X,~~ n \in {\bf N} ~~~\Longrightarrow~~ ^*a \in Y,~~ ^*n ~=~ n $ \\

3.~~~ for $n \in {\bf N}$, we have

$$ a \in V_{n + 1} ( X ) \setminus V_n ( X ) ~~~\Longrightarrow~~~ ^*a \in V_{n + 1} ( Y )  \setminus V_n ( Y ) $$

\medskip
4.~~~ for $n \in {\bf N}$, we have

$$ b \in a \in~ ^*V_{n + 1} ( X ) ~~~\Longrightarrow~~~ b \in~ ^*V_n ( X ) $$

\medskip
5.~~~ Transfer Property : for every sentence $\Phi$ in ${\cal L}_X$ we have

$$ \Phi ~~\mbox{true in}~~ {\cal L}_X ~~~\Longleftrightarrow~~~ ^*\Phi ~~\mbox{true in}~~ {\cal L}_Y $$ \\ \\

{\bf Convention 5.4.1} \\

1.~~ In view of 2 in Definition 5.4.2, and in the spirit of Convention 2.7.1, see (2.7.5), we shall assume that the
existence of a monomorphism (5.4.3) means that $X$ is a {\it subset} of $Y$, and furthermore, we shall make the
identification

\bigskip
(5.4.5) \quad $ X \ni a ~=~ ^*a \in Y $ \\

It follows that

\bigskip
(5.4.5$^*$) \quad $ X ~\subseteq~ Y,~~~~ V ( X ) ~\subseteq~ V ( Y ) $ \\

however, as noted earlier, the mapping $^*(~) : V ( X ) ~\longrightarrow~ V ( Y )$ in (5.4.3) is {\it not} supposed to be
merely the above inclusion, and instead, as seen in the sequel, it is a far less trivial mapping, encapsulating in fact
much of the power of Nonstandard Analysis. \\
Nevertheless, as follows from (5.4.5), (5.4.5$^*$), as well as 1 in Theorem 5.4.1 below, the mapping $^*(~) : V ( X )
~\longrightarrow~ V ( Y )$ in (5.4.3) reduces to an identity mapping for {\it finite} subsets $A \subset X$. \\

2.~~ Similar with Convention 3.4.2, and for the sake of simplicity, we shall say that a sentence is true in a language
${\cal L}_X$, if it is a sentence in that language, and in addition, it is true in the corresponding Superstructure
$V ( X )$. \\

Let us recall here that, throughout the treatment of Nonstandard Analysis in this book, we only deal with what in
Mathematical Logic is called {\it semantic} interpretation of truth.This means that the interpretation of the sentences
of a language as being true or not is done in terms of a given suitable mathematical structure which is outside of that
language. By contrast, a {\it syntactic} interpretation of the truth or falsity of sentences in a language is done
exclusively in terms of the structure of the respective sentences, thus it is done within the framework of the given
language, see Bell \& Slomson, or Marker. \\ \\

{\bf Note 5.4.1} \\

In the Transfer Property in 5 in Definition 5.4.2, the implication ~"$\Longleftarrow$"~ follows from the implication
~"$\Longrightarrow$"~, see Note 3.4.1. \\

Indeed, let $\Phi$ be a sentence in ${\cal L}_X$ which is {\it false} in $V ( X )$. Then in view of 9 in Definition 5.2.1,
$\neg \Phi$ is a sentence in ${\cal L}_X$, and according to Definition 5.3.1, $\neg \Phi$ is true in $V ( X )$. Now 5 in
Definition 5.4.2 implies that $^*( \neg \Phi )$ is true in ${\cal L}_Y$. However (5.4.4) obviously gives $^*( \neg \Phi ) =
\neg ^* \Phi$. \\
In this way it indeed follows that

$$ ^*\Phi ~~\mbox{true in}~~ {\cal L}_Y ~~~\Longrightarrow~~~ \Phi ~~\mbox{true in}~~ {\cal L}_X $$

\medskip
What prevents the simple languages ${\cal L}_{\cal S}$ in chapter 3 from having the full Transfer Property in 5 in
Definition 5.4.2 is the fact that they {\it lack} the above step

$$ \Phi ~~\mbox{is a sentence in}~~ {\cal L}_{\cal S} ~~~\Longrightarrow~~~  \neg \Phi ~~\mbox{is a sentence in}~~
                               {\cal L}_{\cal S} $$

\medskip
since those simple languages only allow one single quantifier, namely, $\forall$, see 12 in Definition 3.2.2

\hfill $\Box$ \\

By using the Transfer Property at 5 in Definition 5.4.2, we can easily obtain \\ \\

{\bf Theorem 5.4.1} \\

Given any monomorphism (5.4.3). Then for $a_1, ~.~.~.~ , a_n \in V ( X )$, we have \\

1.~~~ $^*\{~ a_1, ~.~.~.~ , a_n ~\} ~=~ \{~ ^*a_1, ~.~.~.~ ,~ ^*a_n ~\}$ \\

2.~~~ $^*<~ a_1, ~.~.~.~ , a_n ~> ~=~ <~ ^*a_1, ~.~.~.~ ,~ ^*a_n ~>$ \\

3.~~~ $^* ( \bigcup_{1 \leq i \leq n}~ a_i ) ~=~ \bigcup_{1 \leq i \leq n}~ ^*a_i,~~~~~~  ^* ( \bigcap_{1 \leq i \leq n}~ a_i ) ~=~
                        \bigcap_{1 \leq i \leq n}~ ^*a_i$ \\

4.~~~ $^* ( a_1 \times ~.~.~.~ \times a_n ) ~=~ ^*a_1 \times ~.~.~.~ \times ^*a_n$ \\

Further, for $a, b \in V ( X )$, we have \\

5.~~~ $a \in b ~~~\Longleftrightarrow~~~ ^*a \in~ ^*b$ \\

6.~~~ $a ~=~ b ~~~\Longleftrightarrow~~~ ^*a ~=~ ^*b$ \\

7.~~~$a ~\subseteq~ b ~~~\Longleftrightarrow~~~ ^*a ~\subseteq~ ^*b$ \\

If $P$ is a relation on $a_1 \times ~.~.~.~ \times a_n$, then \\

8.~~~ $^*P ~~\mbox{is a relation on}~~ ^*a_1 \times ~.~.~.~ \times~ ^*a_n $ \\

and if $n = 2$, then also \\

9.~~~ $^* dom P ~=~ dom~ ^*P,~~~ ^* range P ~=~ range~ ^*P $ \\

If we have a mapping $f : a \longrightarrow b$, then \\

10.~~ $^*f :~ ^*a \longrightarrow~ ^*b ~~~\mbox{and}~~ ^* ( f ( c ) ) ~=~ ^*f ( ^*c ), ~~\mbox{for}~~ c \in a $ \\

Also \\

11.~~ $ f ~~\mbox{is bijective} ~~~\Longleftrightarrow~~~ ^*f ~~\mbox{is bijective} $ \\

Obviously, in 2, 3 and 8 one has the restriction ~$a_1, ~.~.~.~ , a_n \in V ( X ) \setminus X$, while in 5, one has to assume
that $b \in V ( X ) \setminus X$, and at last, in 7 and 10, it assumed that $a, b \in V ( X ) \setminus X$. \\

Finally, if $a \in V ( X ) \setminus X$, then \\

12.~~ $ ^*{\cal P} ( a ) ~\subseteq~ {\cal P} ( ^*a ) $ \\

{\bf Proof} \\

1. Let us denote $b = \{ a_1, ~.~.~.~ , a_n \}$, then obviously $b \in V ( X )$. Now we can apply the $^*$-transform in
(5.4.4) to the $n + 1$ sentences in ${\cal L}_X$, which are clearly true in $V ( X )$, namely

$$ (~ \forall~ x \in b ~)~[~ x = a_1 \vee ~.~.~.~ \vee x = a_n ~],~~~ a_1 \in b, ~.~.~.~ , a_n \in b $$

\medskip
and obtain the $n + 1$ true sentences in ${\cal L}_Y$

$$ (~ \forall~ x \in~ ^*b ~)~[~ x =~ ^*a_1 \vee ~.~.~.~ \vee x =~ ^*a_n ~] $$
$$ ^*a_1 \in~ ^*b, ~.~.~.~ , ^*a_n \in~ ^*b $$

\medskip
The first sentence above obviously means the inclusion ~"$\subseteq$"~ in 1, while the other $n$ sentences mean the
converse inclusion. \\

2 - 7 follow by similar transfer arguments. \\

8 is a consequence of 7. \\

9. The sentence in ${\cal L}_X$

$$ (~ \forall~ x \in a ~)~[~ (~ x \in dom P \leftrightarrow (~ \exists~ y \in b ~)~[~ P < x, y > ~]~] $$

is true in $V ( X )$, since it is the definition of the $dom P$. By $^*$-transfer we obtain

$$ (~ \forall~ x \in~ ^*a ~)~[~ (~ x \in~ ^*dom P \leftrightarrow (~ \exists~ y \in~ ^*b ~)~[~ ^*f < x, y > ~]~] $$

which is a true sentence in ${\cal L}_Y$, and which by the definition of the $dom~ ^*P$ clearly gives $^* dom P ~=~
dom~ ^*P$. The other relation follows in a similar way. \\

10. The fact that $^*f :~ ^*a \longrightarrow~ ^*b$ results from 8. Let us now show that $^*f$ is indeed a function.
The sentence

$$ (~ \forall~ x \in a ~)~(~ \forall~ y \in b ~)~(~ \forall~ z \in b ~)~[~ (~ f < x, y > \wedge~ f < x, z > ~) \rightarrow (~ y = z ~) ~] $$

\medskip
is true in ${\cal L}_X$, since $f$ is a function. Thus by $^*$-transfer we obtain the true sentence in ${\cal L}_Y$

$$ (~ \forall~ x \in~ ^*a ~)~(~ \forall~ y \in~ ^*b ~)~(~ \forall~ z \in~ ^*b ~)[~ (~ ^*f < x, y > \wedge~ ^*f < x, z > ~)
                  \rightarrow (~ y = z ~) ~] $$

\medskip
which means that $^*f$ is also a function. \\
Let now be any $c \in a$, then $f ( c ) = d \in b$, hence we have the true sentence in ${\cal L}_X$

$$ f ( c ) = d $$

\medskip
thus by $^*$-transfer we get the true sentence in ${\cal L}_Y$

$$ ^*f ( ^*c ) =~ ^*d $$

\medskip
which means that indeed $^* ( f ( c ) ) ~=~ ^*d ~=~ ^*f ( ^*c )$. \\

11. We start with the implication ~"$\Longrightarrow$"~. The sentence

\bigskip
(5.4.6) \quad $ (~ \forall~ x \in a ~)~(~ \forall~ y \in a ~)~[~ (~ f ( x ) = f ( y ) ~) \rightarrow (~ x = y ~) ~] $

\medskip
is true in ${\cal L}_X$, since $f$ is injective. Thus by $^*$-transfer we obtain the true sentence in ${\cal L}_Y$

\bigskip
(5.4.7) \quad $ (~ \forall~ x \in~ ^*a ~)~(~ \forall~ y \in~ ^*a ~)~[~ (~ ^*f ( x ) =~ ^*f ( y ) ~) \rightarrow (~ x = y ~) ~] $

\medskip
which means that $^*f$ is also injective. The surjectivity of $f$ is given by the true sentence in ${\cal L}_X$

\bigskip
(5.4.8) \quad $ (~ \forall~ y \in b ~)~(~ \exists~ x \in a ~)~[~ f ( x ) = y ~] $

\medskip
thus by $^*$-transfer we obtain in ${\cal L}_Y$ the true sentence

\bigskip
(5.4.9) \quad $ (~ \forall~ y \in~ ^*b ~)~(~ \exists~ x \in~ ^*a ~)~[~ ^*f ( x ) = y ~] $

\medskip
therefore, the surjectivity of $^*f$. \\

Since the Transfer Property in 5 in Definition 5.4.2 is an equivalence, we also have in it the implication
~"$\Longleftarrow$"~. This means that (5.4.7) implies (5.4.6), while (5.4.9) implies (5.4.8). \\

12. The sentence

$$ (~\forall~ b \in {\cal P} ( a ) ~)~[~ b ~\subseteq~ a ~] $$

\medskip
is obviously true in ${\cal L}_X$, therefore by $^*$-transfer we obtain the sentence

$$ (~\forall~ b \in~ ^*{\cal P} ( a ) ~)~[~ b ~\subseteq~ ^*a ~] $$

\medskip
which will be true in ${\cal L}_{^*X}$. \\ \\

{\bf Note 5.4.2} \\

It is important to note that, in general

\bigskip
(5.4.10) \quad $ A \in V ( X ) \setminus X ~~~\not\Longrightarrow~~~ ^*A ~=~ \{~ ^*a ~~|~~ a \in A ~\} $

\medskip
and instead, we only have

\bigskip
(5.4.11) \quad $ A \in V ( X ) \setminus X ~~~\Longrightarrow~~~  \{~ ^*a ~~|~~ a \in A ~\} \subseteq~ ^*A $

\medskip
which follows from 5 in Theorem 5.4.1. Indeed, if we take for example $A = {\bf N}$, then in view of 2 in Definition 5.4.2,
we have $\{~ ^*n ~|~ n \in {\bf N} ~\} = {\bf N}$, thus we obtain an instance of (5.4.10), since in this particular case we have
${\bf N} \subseteq~ ^*{\bf N}$ and ${\bf N} \neq~ ^*{\bf N}$. \\

Let us also note that in view of (5.1.2), (5.4.1), we have

\bigskip
(5.4.12) \quad $ ^*V ( X ) ~=~ \bigcup_{n \in {\bf N}}~ ^*V_n ( X ) ~\subseteq~ V ( Y ) $ \\

\medskip
with the comment on (5.4.2$^{**}$) and (5.4.2$^{***}$), made about the meaning of $^*V ( X )$. \\

Related to the above, we introduce the following useful notation. For $A \in V ( X ) \setminus X$ we denote

\bigskip
(5.4.13) \quad $ ^*A_\infty ~=~ ^*A ~\setminus~ A $

\medskip
which we can see as the {\it nonstandard growth} of $A$. Let us note that, according to (5.4.11) and Convention 5.4.1,
we always have

\bigskip
(5.4.14) \quad $ A ~\subseteq~ ^*A $

\medskip
therefore (5.4.13) is meaningful. Now in view of 1 in Theorem 5.4.1, it follows that

\bigskip
(5.4.15) \quad $ A \in V ( X ) \setminus X,~~ A ~~\mbox{finite} ~~~\Longrightarrow~~~ ^*A ~=~ A $

\hfill $\Box$ \\

Let us indicate now the way the results in chapters 2 and 3 are included as corresponding to a particular case of
monomorphism of Superstructures. In this regard, let us start by returning to section 4, in chapter 3, and consider the
two simple systems, see (3.3.1), (3.3.2), (3.4.1), (3.4.2)

\bigskip
(5.4.16) \quad $ {\cal R} ~=~ (~ {\bf R},~ {\cal P},~ {\cal F} ~),~~~ ^*{\cal R} ~=~ (~ ^*{\bf R},~ ^*{\cal P},~ ^*{\cal F} ~) $

\medskip
In this case we have $X = {\bf R}$ and $^*X =~ ^*{\bf R}$. We also have the mapping

\bigskip
(5.4.17) \quad $ ^*(~) : X ~\longrightarrow~ ^*X $

\medskip
which can be seen as the restriction of (5.4.1) to $X$, and which is defined in (2.1.5). \\

Let us now indicate how to build up further the mapping (5.4.17), based on the constructions in chapters 2 and 3, so that
one may eventually get to the mapping in (5.4.1), and obtain it as a {\it monomorphism} between the two particular
Superstructures in (5.4.16). \\
Given on ${\bf R}$ any n-ary relation $P \in {\cal P}$, we obviously have $P \subseteq {\bf R}^n = X^n$. And then in view
of the way in section 1 Cartesian products were dealt with in Superstructures, we obviously have $P \subseteq V_4
( X )$, therefore $P \in V_5 ( X )$. \\
It follows that the transfers ${\cal P} \ni P ~\longmapsto~ ^*P \in~ ^*{\cal P}$ defined in section 7, in chapter 2, can
already indicate the way to extend the mapping (5.4.17) to a mapping

$$ ^*(~) : V_5 ( X ) ~\longrightarrow~ V_5 ( ^*X ) $$

\medskip
In order to obtain the full extension which gives the mapping in (5.4.1), one can make use of the simple version of the
Transfer Theorem in section 4, in chapter 3. At the same time while doing so, one can verify that, indeed, the mapping
(5.4.1) thus obtained satisfies the conditions of Definition 5.4.2, and thus it is a monomorphism between $V ( {\bf R} )$
and $V ( ^*{\bf R} )$. \\ \\

{\bf 5. Ultrapower Construction of Superstructures} \\ \\

Given any set $X$, such that

$$ {\bf R} ~\subseteq~ X $$

\medskip
we shall construct its {\it nonstandard} extension $^*X$ in a way which generalizes the earlier construction of
$^*{\bf R}$ from ${\bf R}$. Furthermore, we shall then construct the {\it monomorphism}

\bigskip
(5.5.1) \quad $ ^* : V ( X ) ~\longrightarrow V ( ^*X ) $

\medskip
For that purpose, we take any {\it infinite} index set $I$ and any {\it free ultrafilter} ${\cal U}$ on the set $I$. Here we
note that, earlier, in the construction of $^*{\bf R}$ from ${\bf R}$, we had $I = {\bf N}$. \\

We start with defining the mapping

\bigskip
(5.5.2) \quad $ V ( X ) \ni S ~~\longmapsto~~ \prod~ S ~=~ S^I ~=~ \{~ a ~|~ a : I ~\longrightarrow~ S ~\} $

\medskip
and note that in view of the assumption (5.1.4), if $S \in X$, then $\prod~ S = S^I = \{~ a ~|~ a : I ~\longrightarrow~ S ~\}$
has only one single element, namely, the {\it constant} function $I \ni i \longmapsto S$, which will be identified with
$S$. \\

We define now on $\prod~ S$ the equivalence relation $\approx_{\cal U}$ by

\bigskip
(5.5.3) \quad $ a ~\approx_{\cal U}~ b ~~~\Longleftrightarrow~~~ \{~ i \in I ~~|~~ a ( i ) ~=~ b ( i ) ~\} \in {\cal U} $

\medskip
for $a, b \in \prod~ S$. Further, we define quotient set of corresponding equivalence classes

\bigskip
(5.5.4) \quad $ \prod_{~{\cal U}}~ S ~=~ \prod~ S ~/~ \approx_{\cal U} $

\medskip
and call it the {\it ultrapower} of $S$. For every $a \in \prod~ S$, we shall denote by $[ a ]$ its equivalence class in
$\prod_{~{\cal U}}~ S = \prod~ S ~/~ \approx_{\cal U}$, thus we have the {\it surjective} mapping

\bigskip
$~~~~~~~~~~~~ \prod~ S \ni a ~\longmapsto~ [ a ] \in \prod_{~{\cal U}}~ S $

\medskip
Now we can define the {\it nonstandard} extension $^*X$ of $X$ as being given by

\bigskip
(5.5.5) \quad $ ^*X ~=~ \prod_{~{\cal U}}~ X $

\medskip
Clearly, this is an obvious generalization of the construction of $^*{\bf R}$ from ${\bf R}$ in (2.1.4). \\

Next we define the {\it bounded ultrapower} of $V ( X )$, namely

\bigskip
(5.5.6) \quad $ \prod^0_{~{\cal U}}~ V ( X ) ~=~  \bigcup_{n \in {\bf N}}~ \prod_{~{\cal U}}~ (~ V_n ( X )
                                                     \setminus V_{n - 1} ( X ) ~) $

\medskip
where we recall that $V_0 ( X ) = X$, while we denote from now on $V_{-1} ( X ) = \phi$. \\

The set $\prod^0_{~{\cal U}}~ V ( X )$ is called a bounded ultrapower, since for every element $[ a ]$ in it we obviously
have a certain $n \in {\bf N}$, such that $a ( i ) \in V_n ( X ) \setminus V_{n - 1} ( X )$, with $i \in I$. \\

Let us now define the mapping

\bigskip
(5.5.7) \quad $ e : V ( X ) ~~\longrightarrow~~ \prod^0_{~{\cal U}}~ V ( X ) $

\medskip
where for $a \in V ( X )$, we set $e ( a ) = [ \hat{a} ]$, with $\hat{a} : I \longrightarrow V ( X )$ given by $I \ni i
\longmapsto \hat{a} ( i ) = a \in V ( X )$. It follows that in case $a \in V_n ( X ) \setminus V_{n - 1} ( X )$, for a certain $n \in
{\bf N}$, then $e ( a ) \in \prod_{~{\cal U}}~ (~ V_n ( X ) \setminus V_{n - 1} (X )~)$. \\

Obviously, $e$ is a generalization of the mapping (2.1.5). \\

Further, we define the {\it Mostowski collapsing mapping}

\bigskip
(5.5.8) \quad $ M : \prod^0_{~{\cal U}}~ V ( X ) ~~\longrightarrow~~ V ( ^*X ) $

\medskip
in successive steps, as follows. First we note that $\prod_{~{\cal U}}~ (~ V_0 ( X ) \setminus V_{-1} ( X ) ~) =
\prod_{~{\cal U}}~ X =~ ^*X$. And we define $M$ on $^*X$ as the identity, that is

\bigskip
(5.5.9) \quad $ ^*X \ni [ a ] ~~\longmapsto~~ M ( [ a ] ) ~=~ [ a ] \in~ ^*X $

\medskip
Let us note here that $[ a ] \in~ ^*X$ means in view of (5.5.2) - (5.5.5) that $a : I \longrightarrow X$ and $[ a ] = \{~ b : I
\longrightarrow X  ~|~ b \approx_{\cal U} a ~\}$. \\

And now in the general situation, that is, for $[ a ] \in \prod_{~{\cal U}}~ (~ V_{n + 1} ( X ) \setminus V_n ( X ) ~)$, with $n \in
{\bf N}$, we define

\bigskip
(5.5.10) \quad $ M ( [ a ] ) ~=~ \left\{~ M ( [ b ] ) ~~
                        \begin{array}{|l}
                               ~1)~~ [ b ] \in_{\cal U} [ a ] \\ \\
                               ~2)~~ [ b ] \in  \prod_{~{\cal U}}~(~ V_{m} ( X ) \setminus V_{m - 1} ( X ) ~) \\ \\
                               ~~~~~~\mbox{for a certain}~~ 0 ~\leq~ m ~\leq~ n - 1
                        \end{array} ~\right\} $ \\

\medskip
where for $[ a ],~ [ b ] \in \prod^0_{~{\cal U}}~ V ( X )$ we denote $[ b ] \in_{\cal U} [ a ]$, if and only if
$\{ i \in I ~|~ b ( i ) \in a ( i ) ~\} \in {\cal U}$. \\

Before going further, let us illustrate in a few simpler cases the action of the mapping $M$. For that purpose, it
useful to return to the particular situation in chapters 2 and 3, where we had $X = {\bf R}$ and $I = {\bf N}$. \\
Now we shall follow, according to (5.5.6), the first two steps in the definition (5.5.10) of $M$. \\

In view of (5.5.9), the action of $M$ on $\prod_{~{\cal U}}~ (~ V_0 ( X ) \setminus V_{-1} ( X ) ~) =  \prod_{~{\cal U}}~ X =
\prod_{~{\cal U}}~ {\bf R} =~ ^*{\bf R}$ is obvious, being the identity mapping. \\

Let now be given any $[ a ] \in \prod_{~{\cal U}}~ (~ V_1 ( X ) \setminus V_{0} ( X ) ~) = \prod_{~{\cal U}}~ {\cal P} ( X )  =
\prod_{~{\cal U}}~ {\cal P} ( {\bf R} )$. Then (5.5.4) implies that $a : I \longrightarrow {\cal P} ( {\bf R} )$, thus $ I \ni i
\longmapsto a ( i ) = A_i \subseteq {\bf R}$. \\
Now (5.5.10) gives

\bigskip
$~~~~~~~~~~~ M ( [ a ] ) ~=~  \{~ M ( [ b ] ) ~~|~~ [ b ] \in_{\cal U} [ a ],~~~ [ b ] \in  \prod_{~{\cal U}}~ X ~=~ ^*{\bf R} ~\} ~=~ $ \\

$~~~~~~~~~~~~~~~~~~~~~=~ \{~ [ b ] \in~ ^*{\bf R} ~~|~~  \{~ i \in I ~|~ b ( i ) \in A_i ~\} \in {\cal U} ~\} $

\medskip
Let us take the particular case of $[ a ] \in \prod_{~{\cal U}}~ {\cal P} ( {\bf R} )$, when $a ( i ) = A$, with $i \in I$, for a
certain $A \subseteq {\bf R}$. This means according to (5.5.7) that $[ a ] = e ( A ) = [ \hat{A} ]$, thus the preceding
relations about $M$ give

\bigskip
$ M ( e ( A ) ) ~=~  \{~ [ b ] \in~ ^*{\bf R} ~~|~~  \{~ i \in I ~|~ b ( i ) \in A ~\} \in {\cal U} ~\} ~=~ ^*A $

\medskip
the last equality resulting from (2.7.6). \\

The above clearly illustrates the need for Mostowski's collapsing map, since obviously $e ( A ) = [ \hat{A} ]
\not\subseteq ~^*{\bf R}$. However, we have $M ( e ( A ) ) = ~^*A \subseteq~ ^*{\bf R}$. \\ \\

{\bf Proposition 5.5.1} \\

1.~~~ The mappings $e$ and $M$ are injective. \\

2.~~~ $e$ maps $X$ into $^*X$, while $M$ maps $^*X$ onto $^*X$. \\

3.~~~ $e$ maps $V_{n + 1} ( X ) \setminus V_n ( X )$ into $\prod_{~{\cal U}}~ (~ V_{n + 1} ( X ) \setminus V_n ( X ) ~)$, while
$M$ \\
\hspace*{0.8cm} maps $\prod_{~{\cal U}}~ (~ V_{n + 1} ( X ) \setminus V_n ( X ) ~)$ into $ V_{n + 1} ( ^*X ) \setminus
V_n ( ^*X )$. \\

4.~~~ For $a, b \in V ( X )$ we have $a \in b ~\Longleftrightarrow~ e ( a ) \in_{\cal U}~ e ( b )$, while for \\
\hspace*{0.8cm} $[ a ], [ b ] \in \prod^0_{~{\cal U}}~ V ( X )$, we have  $[ a ] \in_{\cal U}~ [ b ] ~\Longleftrightarrow~
M ( [ a ] ) \in M ( [ b ] )$. \\

5.~~~$  e ( X ) ~=~ [ \hat{X} ],~~~ M ( [ \hat{X} ] ) ~=~ M ( e ( X ) ) ~=~ ^*X  $ \\

6.~~~ Let $[ a ], [ b ] \in \prod^0_{~{\cal U}}~ V ( X )$, and define $c : I \longrightarrow V ( X )$ by \\
\hspace*{0.8cm} $I \ni i \longmapsto c_i = \{ a_i, b_i \}$. Then $[ c ] \in \prod^0_{~{\cal U}}~ V ( X )$ and $M ( [ c ] ) =
\hspace*{0.8cm} \{~ M ( [ a ] ),~ M ( [ b ] ) ~\}$. \\
\hspace*{0.8cm} Similar relations hold with ${~}$ replaced by $<~>$, and $=$ by $\in$. \\
\hspace*{0.8cm} Also, all such relations hold for three or more terms. \\

7.~~~ If $a \in V_n ( X ) \setminus V_{n - 1} ( X )$, for a certain $n \in {\bf N}$, and $[ b ] \in_{\cal U} e ( a )$, \\
\hspace*{0.8cm} then $[ b ] \in_{\cal U} e ( V_{n - 1} ( X ) )$. \\

{\bf Proof}?? \\

\hfill $\Box$ \\

Now we can return to defining the desired mapping $^*$ in (5.5.1) as the composition of $e$ and $M$, see (5.5.7),
(5.5.8), namely

\bigskip
(5.5.11) \quad $ ^* ~=~ e \circ M ~:~ V ( X ) \stackrel{e}{~\longrightarrow~~} \prod^0_{~{\cal U}}~ V ( X ) \stackrel{M}
                                                                           {~\longrightarrow~} V ( ^*X ) $

\medskip
and then prove that, indeed, it is a monomorphism in the sense of Definition 5.4.2. In this respect we need the following
basic result from Model Theory, see Bell \& Slomson, Marker, or Jech \\ \\

{\bf Theorem 5.5.1 ( {\L}\u{o}s, 1954 )} \\

Given a formula $\Psi ( x_1, ~.~.~.~ , x_n )$ in ${\cal L}_X$ in which $x_, ~.~.~.~ , x_n$ are the only free variables. Then for
any $[ a_1 ], ~.~.~.~ , [ a_n ] \in \prod^0_{~{\cal U}}~ V ( X )$, we have

\bigskip
(5.5.12) \quad $ \begin{array}{l}
                        ^*\Psi ( M ( [ a_1 ] ), ~.~.~.~ , M ( [ a_n ] ) ) ~~\mbox{true in}~~ {\cal L}_{^*X} ~~~\Longleftrightarrow~~ \\ \\

                          \Longleftrightarrow~~~  \{~ i \in I ~~|~~ \Psi ( a_1 ( i ), ~.~.~.~ , a_n ( i ) ) ~~\mbox{true in}~~ {\cal L}_{X} ~\} \in
                                           {\cal U}
                           \end{array} $ \\

{\bf Proof}?? \\

\hfill $\Box$ \\ \\

{\bf Theorem 5.5.2} \\

The mapping $^*$ in (5.5.11) is a monomorphism. \\

{\bf Proof}?? \\

\hfill $\Box$ \\ \\

{\bf 6. Genuine Extensions, Hyperfinite Sets, Standard and \\
        \hspace*{0.5cm} Nonstandard Entities, Concurrent Relations, and \\
        \hspace*{0.5cm} Exhausting Sets} \\ \\

Needless to say, our interest is to construct extensions $^*X$ of $X$, see (5.5.5), which are {\it genuine}, that is, for
which we have

\bigskip
(5.6.1) \quad $ X ~\subseteq~ ^*X,~~~ X ~\neq~ ^*X $

\medskip
Here we shall show that this is always possible when $X$ is {\it infinite}, if we choose the index set $I$ and the
ultrafilter ${\cal U}$ on $I$ in appropriate ways. And we recall that in our case $X$ is indeed infinite, since we have
assumed from the start that ${\bf R} \subseteq X$. \\

Let us now recall the relevant results related to filters and ultrafilters presented in sections 2 - 4, in chapter 2. First we
note that we are constrained to use {\it ultrafilters}, in view of the fact that we want to recover as a particular case the
construction of the nonstandard reals $^*{\bf R}$ done in chapter 2, a construction which was essentially based on
ultrafilters on ${\bf N}$. Therefore, the only issue is whether, in the general case studied here, we use {\it fixed} or {\it
free} ultrafilters on the sets $I$. \\

If ${\cal U}$ is a {\it fixed} ultrafilter on $I$, that is, ${\cal U} = \{~ J \subseteq I ~|~ i \in I ~\}$, for a certain given $i \in I$,
then it follows easily that, see (5.5.5)

\bigskip
$ ~~~~~~~~~~~~~~~^*X ~=~ \prod_{~{\cal U}}~ X ~=~ X $

\medskip
Therefore, in order to secure (5.6.1), we have to use {\it free} ultrafilters ${\cal U}$ on $I$. \\

Here however, in the general case, it is {\it not} sufficient to use any free ultrafilter ${\cal U}$ on any infinite set $I$, in
order to obtain (5.6.1). This is indeed {\it unlike} in the case of the construction of the nonstandard reals $^*{\bf R}$ in
chapter 2, where we had $I = {\bf N}$, and then we could use {\it any} free ultrafilter ${\cal U}$ on that particular $I$. \\

In order to clarify this issue, let us recall the way we showed in section 3, in chapter 2 that $^*{\bf R}$ is a genuine
extension of ${\bf R}$. \\
Namely, we took any free ultrafilter ${\cal U}$ on $I = {\bf N}$, and defined $\omega \in {\bf R}^{\bf N}$ by

\bigskip
$ ~~~~~~~~~~{\bf N} \ni n ~\longmapsto~ \omega ( n ) = n \in {\bf R} $

\medskip
Then it followed that $[ \omega ] \neq~ ^*r$, for any $r \in {\bf R}$. Indeed, if we assume that $[ \omega ] =~ ^*r$, for some
$r \in {\bf R}$, then we must have $\{~ n \in {\bf N} ~~|~~ n = r ~\} \in {\cal U}$. This however is obviously not possible for
any free ultrafilter ${\cal U}$ on $I = {\bf N}$, since the set $\{~ n \in {\bf N} ~~|~~ n = r ~\}$ is either void, or it contains one
single element. \\

Now, in the general case of (5.6.1), a similar proof would require to construct $\omega : I \longrightarrow X$, such that

\bigskip
$ ~~~~~~~~~~\{~ i \in I ~~|~~ \omega ( i ) ~=~ x ~\} \not\in {\cal U},~~~ \mbox{for}~~ x \in X $

\medskip
which is obviously not as trivial as in the case of $I = {\bf N}$, since for instance, the cardinal of $I$ could now happen
to be larger than that of $X$, and then every such function $\omega$ would fail to be injective, thus the sets $\{~ i \in I
~|~ \omega ( i ) = x ~\}$ could end up being rather large. \\

In order to obtain (5.6.1) for arbitrary sets $X$, we shall introduce the concept of {\it hyperfinite} sets. This concept will
also have an important interest on its own. \\

For any given set $Y$, let us denote by ${\cal P}_F ( Y )$ the set of all {\it finite} subsets of $Y$. Obviously, if $A \in
V_{n + 1} ( X ) \setminus X$, for some $n \in {\bf N}$, then ${\cal P}_F ( A ) \in V_{n + 2} ( X )$, thus in view of (5.5.1), we
can define the mapping

\bigskip
(5.6.2) \quad $ V ( X ) \setminus X \ni A ~\longmapsto~ ^*{\cal P}_F ( A ) \in V ( ^*X ) $ \\

If $A \in V ( X ) \setminus X$ is {\it finite}, then so is obviously ${\cal P}_F ( A ) = {\cal P} ( A )$, therefore, according to 1
in Theorem 5.4.1, the set

$$  ^*{\cal P}_F ( A ) ~=~ \{~ ^*B ~~|~~ B ~\subseteq~ A ~\} ~=~  ^*{\cal P} ( A ) $$

\medskip
is also finite. However, our interest will be mainly in the cases when $A$, and thus ${\cal P}_F ( A )$ and $^*{\cal P}_F
( A )$ are infinite, and then we have in general, see 7 and 12 in Theorem 5.4.1

$$  \{~ ^*B ~~|~~ B ~\subseteq~ A,~~~ B ~~\mbox{finite}~ ~\} ~\subseteq~  ^*{\cal P}_F ( A ) ~\subseteq~  ^*{\cal P} ( A )
                                       ~\subseteq~ {\cal P} ( ^*A ) $$

\medskip
and also

$$ B ~\subseteq~ A \in V ( X ) \setminus X ~~~\Longrightarrow~~~ ^*{\cal P}_F ( B ) ~\subseteq~ ^*{\cal P}_F ( A ) $$ \\ \\

{\bf Definition 5.6.1} \\

Given $A \in V ( X ) \setminus X$, then

\bigskip
(5.6.3) \quad $ ^*{\cal P}_F ( A ) $

\medskip
is called the set of {\it hyperfinite} subsets of $^*A$. In other words, by definition

\bigskip
(5.6.3$^*$) \quad $  ^*{\cal P}_F ( A ) ~=~ \{~ B ~\subseteq~ ^*A ~~|~~ B ~~\mbox{hyperfinite} ~\} $ \\ \\

{\bf Lemma 5.6.1} \\

The set of all hyperfinite subsets of $V ( ^*X )$ is given by

\bigskip
(5.6.4) \quad $ \begin{array}{l}
                              \bigcup_{A \in V ( X ) \setminus X}~ ^*{\cal P}_F ( A ) ~=~
                                                \bigcup_{n \in {\bf N}}~ ^*{\cal P}_F ( V_n ( X ) ) ~\subseteq~ \\ \\
                              ~\subseteq~ ^*V ( X ) ~=~ \bigcup_{n \in {\bf N}}~ ^*V_n ( X ) ~\subseteq~ V ( ^*X )
                        \end{array} $ \\

{\bf Proof} \\

If $A \in V ( X ) \setminus X$, then (5.1.2) gives $A \in V_{n + 1} ( X )$, for a certain $n \in {\bf N}$. If we take the smallest
such $n$, then (5.1.1) implies that $A \subseteq V_n ( X )$. Hence ${\cal P}_F ( A ) \subseteq {\cal P}_F ( V_n ( X ) )$.
Thus 7 in Theorem 5.4.1 results in $^*{\cal P}_F ( A ) \subseteq~ ^*{\cal P}_F ( V_n ( X ) )$. \\
Thus we have the inclusion $\subseteq$ in the first equality. The converse inclusion is obvious, since $V_n ( X ) \in
V ( X ) \setminus X$, for $n \in {\bf N}$, see (5.1.3). \\ \\

{\bf Note 5.6.1} \\

1.~~ It is obvious that every {\it finite}  $A \in V ( X ) \setminus X$ is {\it hyperfinite}, since in that case $A \in {\cal P} ( A )
= {\cal P}_F ( A )$, while in view of (5.4.15) we have $A =~ ^*A$, thus $A \in~ ^*{\cal P}_F ( A )$. \\

2.~~ Simple examples show that

\bigskip
(5.6.5) \quad $ A ~~\mbox{hyperfinite},~~ B ~\subseteq~ A ~~~\not\Longrightarrow~~~ B ~~\mbox{hyperfinite} $

\medskip
Indeed, as seen in Lemma 5.6.2 below, the sets $K = \{~ l \in~ ^*{\bf N} ~|~ l ~\leq~ k ~\}$, with $k \in~ ^*{\bf N}_\infty $, are
hyperfinite, and clearly, for each of them we have ${\bf N} \subseteq K$. \\
Yet ${\bf N}$ is {\it not} hyperfinite, see Proposition 5.6.1 below. \\

3.~~ In Lemma 5.6.2 next, it will be useful to note the following. Let $k \in~ ^*{\bf N}_\infty$. Then  the set $\{~ l \in~
^*{\bf N} ~|~ l ~\leq~ k ~\}$ can be written in the form $\{~ l \in~ ^*{\bf N} ~|~ l ~<~ k^\prime ~\}$, by taking $k^\prime = k + 1$.
Similarly, the set $\{~ l \in~ ^*{\bf N} ~|~ l ~<~ k ~\}$ can be written in the form $\{~ l \in~ ^*{\bf N} ~|~ l ~\leq~ k^\prime ~\}$, by
taking $k^\prime = k - 1$. \\
Sets of the from $\{~ l \in~ ^*{\bf N} ~|~ l ~\leq~ k ~\}$, with $k \in~ ^*{\bf N}_\infty $, are called {\it initial segments} in
$^*{\bf N}$.

\hfill $\Box$ \\

As we shall see, many properties of finite sets have, through $^*$-transfer, their correspondent in the case of
hyperfinite sets. Here we mention one such instance which, as seen in Proposition 5.6.1 below, leads to a useful {\it
representation}, and later also {\it characterisation}, of arbitrary hyperfinite sets. \\ \\

{\bf Lemma 5.6.2} \\

Let $k \in~ ^*{\bf N}_\infty$, then the sets $\{~ l \in~ ^*{\bf N} ~|~ l ~<~ k ~\}$ and $\{~ l \in~ ^*{\bf N} ~|~ l ~\leq~ k ~\}$ are both
hyperfinite and uncountable. \\

{\bf Proof} \\

Let ${\cal N}$ be the set of all subsets of ${\bf N}$ of the form $\{~ m \in {\bf N} ~|~ m \leq n ~\}$, where $n \in {\bf N}$, is
arbitrary. Then ${\cal N} \subseteq {\cal P}_F ( {\bf N} )$, hence $^*{\cal N} \subseteq~ ^*{\cal P}_F ( {\bf N} )$, see 7 in
Theorem 5.4.1. \\
Now we show that in view of 2 and 5 in Definition 5.4.2, we have $K = \{~ l \in~ ^*{\bf N} ~|~ l ~\leq~ k ~\} \in~ ^*{\cal N}$. \\
Indeed, the sentence

$$ (~ \forall~ A \in V ( X ) ~)~[~ A \in {\cal N} \leftrightarrow $$
$$~~~~~~(~ \exists~ n \in {\bf N} ~)~(~ \forall~ x \in X ~)~ [~ x \in A \leftrightarrow [~ (~ x \in {\bf N} ~) \wedge (~ x \leq n ~) ~]~]~] $$

\medskip
is true in ${\cal L}_X$, therefore by $^*$-transfer, we have in ${\cal L}_{^*X}$ the true sentence

$$ (~ \forall~ A \in~ ^*V ( X ) ~)~[~ A \in~ ^*{\cal N} \leftrightarrow $$
$$~~~~~~(~ \exists~ n \in~ ^*{\bf N} ~)~(~ \forall~ x \in~ ^*X ~)~ [~ x \in A \leftrightarrow
                                                              [~ (~ x \in~ ^*{\bf N} ~) \wedge (~ x \leq n ~) ~]~]~] $$

\medskip
But in view of (5.4.12), we have

$$ ^*V ( X ) ~=~ \bigcup_{n \in {\bf N}}~ ^*V_n ( X ) $$

\medskip
and clearly ${\cal N} \in V_2 ( X )$, hence $^*{\cal N} \in~  ^*V_2 ( X ) \subseteq~ ^*V ( X )$, see 5 in Theorem 5.4.1. \\
Thus from the second sentence above, a sentence in ${\cal L}_{^*X}$, we indeed have that $\{~ l \in~ ^*{\bf N} ~|~ l ~\leq~
k ~\} \in~ ^*{\cal N} \subseteq~ ^*{\cal P}_F ( {\bf N} )$. \\

The uncountability of $K$ follows from 4 in Corollary 3.6.1.

\hfill $\Box$ \\

It is useful to note that the hyperfinite {\it initial segments} in Lemma 5.6.2 above {\it characterize} all the infinite
hyperfinite sets. This characterization, however, involves the concept of {\it internal} sets, which will be introduced in
section 7. Here, as a part of that characterization, we present a {\it representation} of hyperfinite sets given by inital
segments. \\ \\

{\bf Proposition 5.6.1} \\

Let, according to (5.6.4), be given any hyperfinite set $B \in~ ^*{\cal P}_F ( V_n ( X ) )$, where~ $n \in {\bf N}$. Then there
exists a hyperfinite set  $K = \{~ l \in~ ^*{\bf N} ~|~ l ~\leq~ k ~\}$, for a suitable $k \in~ ^*{\bf N}_\infty$, and a {\it bijective}
mapping $\beta : K \longrightarrow B$, such that $\beta \in~ ^*V_{n + 3} ( X )$. It follows that we can write

\bigskip
(5.6.6) \quad $ B ~=~ \{~ b_0, b_1, b_2, ~.~.~.~ , b_k ~\} $

\medskip
where $b_l = \beta ( l )$, with $0 \leq l \leq k$. \\

{\bf Note 5.6.2} \\

Since the hyperfinite sets $B$ in (5.6.6) are {\it uncountable}, the respective dots $"~.~.~.~"$ in the right hand term of that
relation are indicating an {\it uncountable} amount of terms. This is unlike with the usual meaning of such dots in
standard mathematics where, typically, they are used for indicating a countable amount of terms only. However, since
the sets $K = \{~ l \in~ ^*{\bf N} ~|~ l ~\leq~ k ~\}$ are obviously totally ordered and the mappings $\beta : K \longrightarrow
B$ are bijective, the notation with dots in (5.6.6) is a natural extension of that in standard mathematics. \\

{\bf Proof} \\

Obviously,  we have in ${\cal L}_X$ the true sentence

$$ (~ \forall~ B \in {\cal P}_F ( V_n ( X ) ) ~)~(~ \exists~ k \in {\bf N} ~)~(~ \exists~ \beta \in V_{n + 3}
        ( X ) ~) $$
$$ [~ \beta : \{~ l \in {\bf N} ~~|~~ l ~\leq~ k ~\} \longrightarrow B ~~\mbox{is bijective} ~] $$

which for the sake of simplicity has its part within the brackets $[~~~]$ written informally. \\ \\
Now by $^*$-transfer, we obtain in ${\cal L}_{^*X }$ the true sentence

$$ (~ \forall~ B \in~ ^*{\cal P}_F ( V_n ( X ) ) ~)~(~ \exists~ k \in~ ^*{\bf N} ~)~(~ \exists~ \beta \in~ ^*V_{n + 3}
        ( X ) ~) $$
$$ [~ \beta : \{~ l \in~ ^*{\bf N} ~~|~~ l ~\leq~ k ~\} \longrightarrow B ~~\mbox{is bijective} ~] $$

\medskip
and thus (5.6.6). \\

Conversely, let $\beta : K \longrightarrow B$ be a bijective mapping, where $K = \{~ l \in~ ^*{\bf N} ~|~ l ~\leq~ k ~\}$, for a
suitable $k \in~ ^*{\bf N}_\infty$.

\hfill $\Box$ \\

The above Lemma 5.6.2 also gives an improvement of the result in Proposition 3.6.3, namely \\ \\

{\bf Proposition 5.6.2} \\

For every $s \in~ ^*{\bf R}_\infty =~ ^*{\bf R} \setminus {\bf R},~ s > 0$, there are uncountably many copies of $Gal ( 0 )$
in the disjoint union, see (3.6.9)

$$ ^*{\bf R} ~=~ {\bigcup}_{\lambda \in \Lambda}~ (~ s^\lambda + Gal ( 0 ) ~) $$

\medskip
which are between $Gal ( 0 )$ and $s + Gal ( 0 )$. \\

{\bf Proof} \\

Each copy of $Gal ( 0 )$ in the above disjoint union contains at most countably many elements from $^*{\bf N}$.

\hfill $\Box$ \\

Here we give two simple examples of the use of hyperfinite sets which, however, will be important in chapter 6,
connected with Loeb integration and measures. \\ \\

{\bf Example 5.6.1} \\

First we show how to sum up series with {\it uncountably} many terms. Let us define the function

$$ \Sigma : {\bf N} \times {\bf R}^{\bf N} ~~\longrightarrow~~ {\bf R} $$

\medskip
by

$$ \Sigma ( n , a ) ~=~ \sum_{~0 \leq i \leq n~} a ( i ),~~~~ n \in {\bf N},~~ a \in {\bf R}^{\bf N} $$

\medskip
in other words, for every given $n \in {\bf N}$, $\Sigma$ sums the first $n + 1$ terms in the sequence of real numbers
$a ( 0 ), a ( 1 ), a ( 2 ), ~.~.~.~ $. Then clearly, we have in ${\cal L}_X$ the true sentence

$$ (~\forall~ n \in {\bf N} ~)~(~ \forall~ a \in {\bf R}^{\bf N} ~)~(~ \exists~ s \in {\bf R} ~)~[~ \Sigma ( n, a ) ~=~ s ~] $$

\medskip
Therefore by $^*$-transfer, we have in ${\cal L}_{^*X}$ the true sentence

$$ (~\forall~ n \in~ ^*{\bf N} ~)~(~ \forall~ a \in~ ^*( {\bf R}^{\bf N} ) ~)~(~ \exists~ s \in~ ^*{\bf R} ~)~[~ ^*\Sigma ( n, a ) ~=~ s ~] $$

\medskip
Let us note that $^*\Sigma$ in the above sentence is again a function, according to 9 and 10 in Theorem 5.4.1, namely,
$^*\Sigma :~ ^*( {\bf N} \times {\bf R}^{\bf N} ) ~\longrightarrow~ ^*{\bf R} $. Further, in view of 4, and again 9 and 10 in the
same theorem, we have

$$ ^*( {\bf N} \times {\bf R}^{\bf N} ) ~=~ ^*{\bf N} \times~ ^*( {\bf R}^{\bf N} ) ~=~ ^*{\bf N} \times (^*{\bf R} )^{(^*{\bf N} )} $$

\medskip
hence in fact

$$ ^*\Sigma :~ ^*{\bf N} \times (^*{\bf R} )^{(^*{\bf N} )} ~~\longrightarrow~~ ^*{\bf R} $$

\medskip
Now the last sentence becomes

$$ (~\forall~ n \in~ ^*{\bf N} ~)~(~ \forall~ a \in~ (^*{\bf R})^{(^*{\bf N} )} ~)~(~ \exists~ s \in~ ^*{\bf R} ~)~
                           [~ ^*\Sigma ( n, a ) ~=~ s ~] $$

\medskip
which being true in ${\cal L}_{^*X}$, it means that, for {\it arbitrary} sequences $a \in {\bf R}^{\bf N}$, we can make {\it
hyperfinite} sums

$$ \sum_{~0 \leq i \leq n~} {^*a ( i )} \in~ ^*{\bf R} $$

\medskip
with {\it arbitrary hyperfinite} number of terms $n \in~ ^*{\bf N}_\infty =~ ^*{\bf N} \setminus {\bf N}$, where $^*a :~
^*{\bf N} ~\longrightarrow~ ^*{\bf R}$ is the nonstandard extension, or in other words, the $^*$-transfer of the
sequence $a : {\bf N} \longrightarrow {\bf R}$. \\

Let us recall here the nontrivial fact that, in view of Lemma 5.6.2, such {\it hyperfinite} sums always contain {\it
uncountably} many terms. \\ \\

{\bf Example 5.6.2} \\

Now we show how to make products with {\it uncountably} many factors. Let us define the function

$$ \prod : {\bf N} \times {\bf R}^{\bf N} ~~\longrightarrow~~ {\bf R} $$

\medskip
by

$$ \prod ( n , a ) ~=~ \prod_{~0 \leq i \leq n~} a ( i ),~~~~ n \in {\bf N},~~ a \in {\bf R}^{\bf N} $$

\medskip
in other words, for every given $n \in {\bf N}$, $\prod$ multiplies the first $n + 1$ terms in the sequence of real numbers
$a ( 0 ), a ( 1 ), a ( 2 ), ~.~.~.~ $. Then clearly, we have in ${\cal L}_X$ the true sentence

$$ (~\forall~ n \in {\bf N} ~)~(~ \forall~ a \in {\bf R}^{\bf N} ~)~(~ \exists~ p \in {\bf R} ~)~[~ \prod ( n, a ) ~=~ p ~] $$

\medskip
Therefore by $^*$-transfer, we have in ${\cal L}_{^*X}$ the true sentence

$$ (~\forall~ n \in~ ^*{\bf N} ~)~(~ \forall~ a \in~ ^*( {\bf R}^{\bf N} ) ~)~(~ \exists~ p \in~ ^*{\bf R} ~)~[~ ^*\prod ( n, a ) ~=~ p ~] $$

\medskip
Let us note that $^*\prod$ in the above sentence is again a function, according to 9 and 10 in Theorem 5.4.1, namely,
$^*\prod :~ ^*( {\bf N} \times {\bf R}^{\bf N} ) ~\longrightarrow~ ^*{\bf R} $. Further, in view of 4, and again 9 and 10 in the
same theorem, we have

$$ ^*( {\bf N} \times {\bf R}^{\bf N} ) ~=~ ^*{\bf N} \times~ ^*( {\bf R}^{\bf N} ) ~=~ ^*{\bf N} \times (^*{\bf R} )^{(^*{\bf N} )} $$

\medskip
hence in fact

$$ ^*\prod :~ ^*{\bf N} \times (^*{\bf R} )^{(^*{\bf N} )} ~~\longrightarrow~~ ^*{\bf R} $$

\medskip
Now the last sentence becomes

$$ (~\forall~ n \in~ ^*{\bf N} ~)~(~ \forall~ a \in~ (^*{\bf R})^{(^*{\bf N} )} ~)~(~ \exists~ p \in~ ^*{\bf R} ~)~
                           [~ ^*\prod ( n, a ) ~=~ p ~] $$

\medskip
which being true in ${\cal L}_{^*X}$, it means that, for {\it arbitrary} sequences $a \in {\bf R}^{\bf N}$, we can make {\it
hyperfinite} products

$$ \prod_{~0 \leq i \leq n~} {^*a ( i )} \in~ ^*{\bf R} $$

\medskip
with {\it arbitrary hyperfinite} number of terms $n \in~ ^*{\bf N}_\infty =~ ^*{\bf N} \setminus {\bf N}$, where $^*a :~
^*{\bf N} ~\longrightarrow~ ^*{\bf R}$ is the nonstandard extension, or in other words, the $^*$-transfer of the
sequence $a : {\bf N} \longrightarrow {\bf R}$. \\

In particular, for every {\it hyperfinite} number $n \in~ ^*{\bf N}_\infty =~ ^*{\bf N} \setminus {\bf N}$, we can compute

$$ n ! \in~ ^*{\bf N}_\infty =~ ^*{\bf N} \setminus {\bf N} $$

\medskip
Indeed, for that purpose we take $a : {\bf N} \longrightarrow {\bf N}$, given by $a ( i ) = i$, with $i \in {\bf N}$. In this case
we shall obtain $^*a :~ ^*{\bf N} ~\longrightarrow~ ^*{\bf N}$, given by $^*a ( i ) = i$, for $i \in~ ^*{\bf N}$, and thus the
above relation for $n !$, provided that in all the products we have $1 \leq i \leq n$, instead of $0 \leq i \leq n$. \\

Here again we have the nontrivial fact that, in view of Lemma 5.6.2, such {\it hyperfinite} products always contain {\it
uncountably} many factors. \\ \\

{\bf Definition 5.6.2} \\

Entities $a \in V ( X )$, or those in $V ( ^*X )$ which, see (5.4.1), are of the form $a =~ ^*b$, for some $b \in V
( X )$, are called {\it standard}. All other entities in $V ( ^*X )$ are called {\it nonstandard}. \\ \\

{\bf Note 5.6.3} \\

Let $A \subseteq X$ then

$$ \{~ a \in~ ^*A ~~|~~ a ~~\mbox{standard} ~\} ~\subseteq~ A $$

$$ a \in~ ^*A_\infty ~=~ ^*A \setminus A ~~~\Longrightarrow~~~ a ~~\mbox{nonstandard} $$

$$ ^*A_\infty ~=~ ^*A \setminus A ~\subseteq~ ^*X \setminus X $$

\medskip
Indeed, for the first relation we note that $a$ standard implies $a =~ ^*b$, for some $b \in V ( X )$. Hence $a \in~ ^*A$
gives $^*b \in~ ^*A$. And then by inverse $^*$-transfer we have $b \in A$, thus in view of (5.4.5), $^*b = b$. In this way
$a =~ ^*b = b \in A$. The converse follows directly from (5.4.5). \\
The second and thirs relations are obviously implied by the first one. \\ \\

{\bf Example 5.6.3} \\

Let ${\cal I}$ be the set of all closed bounded intervals $[ a, b ] \subset {\bf R}$. Then we have in ${\cal L}_{ X }$
the true sentences

$$ (~ \forall~ A \in {\cal I} ~)~(~ \exists~ a, b \in {\bf R} ~)~(~ \forall~ x \in {\bf R} ~)~[~ x \in A
\leftrightarrow a \leq x \leq b ~] $$
$$ (~ \forall~ a, b \in {\bf R} ~)~(~ \exists~ A \in {\cal I} ~)~(~ \forall~ x \in {\bf R} ~)~[~ ~ x \in A
\leftrightarrow a \leq x \leq b ~] $$

\medskip
By $^*$-transfer we obtain in ${\cal L}_{ ^*X }$ the true sentences

$$ (~ \forall~ A \in~ ^*{\cal I} ~)~(~ \exists~ a, b \in~ ^*{\bf R} ~)~(~ \forall~ x~ ^*\in {\bf R} ~)~[~ x \in A
\leftrightarrow a \leq x \leq b ~] $$
$$ (~ \forall~ a, b \in~ ^*{\bf R} ~)~(~ \exists~ A \in~ ^*{\cal I} ~)~(~ \forall~ x \in~ ^*{\bf R} ~)~[~ ~ x \in A
\leftrightarrow a \leq x \leq b ~] $$

\medskip
It follows by applying (5.4.1), that

$$ {\bf R} ~\supset~ [ a, b ] ~~~\longmapsto~~~ ^*[ a, b ] ~=~ [ a, b ] ~\subset~ ^*{\bf R} $$

\medskip
where the interval at the right is given by

$$ [ a, b ] ~=~ \{~ x \in~ ^*{\bf R} ~~|~~ a ~\leq~ x ~\leq~ b ~\} $$

\medskip
thus the nonstandard interval $[ a, b ] \subset~ ^*{\bf R}$ is a {\it standard} set, and when $a < b$, obviously it
contains {\it nonstandard} numbers, since for instance $mon ( ( a + b ) /2 ) \subset [ a, b ]$. \\ \\

{\bf Definition 5.6.3} \\

The Superstructure $V ( ^*X )$ is called an {\it extension}, if and only if for every set $A \in V ( X ) \setminus X$, there
exists a set $B \in~ ^*{\cal P}_F ( A )$, such that, see (5.4.10), (5.4.11)

\bigskip
(5.6.7) \quad $  a \in A ~~~\Longrightarrow~~~ ^*a \in B $

\hfill $\Box$ \\

Let us note that in view of Convention 5.4.1, the above condition (5.6.7) can equivalently be written as

\bigskip
(5.6.7$^*$) \quad $ A ~\subseteq~ B $

\medskip
The point in Definition 5.6.3 is as follows. In Proposition 5.6.1 we have seen that every hyperfinite set can be
represented by a hyperfinite set $\{~ l \in~ ^*{\bf N} ~|~ l \leq k ~\}$, with a suitable $k \in~ ^*{\bf N}_\infty$. And clearly

$$ {\bf N} ~\subset~ \{~ l \in~ ^*{\bf N} ~~|~~ l ~\leq~ k ~\} $$

\medskip
thus in particular $^*{\bf N}$ is a genuine extension of ${\bf N}$. \\
In this way, Definition 5.6.3 is a generalization of this genuine extension property to {\it every} set in $V ( X )$. \\

We shall now construct an infinite set $I$ and a {\it free} ultrafilter ${\cal U}$ on it, such that the corresponding
Superstructure $V ( ^*X )$ is an {\it extension}. \\
Let $X$ be an arbitrary set for which, as before, we have ${\bf R} \subseteq X$. Then we take

\bigskip
(5.6.8) \quad $ I ~=~ \{~ J ~\subseteq~ V ( X ) ~~|~~ J ~\neq~ \phi,~~\mbox{finite} ~\} $

\medskip
It follows that $I$ depends on $X$. Also we note that for $a \in V ( X )$, we have

\bigskip
(5.6.9) \quad $ a \in I ~~~\Longleftrightarrow~~ \exists~~ b \in V ( X ) \setminus X ~:~ a \in {\cal P}_F ( b ),~~ a \neq \phi $

\medskip
Further, for $a \in I$, we denote

\bigskip
(5.6.10) \quad $ I_a ~=~ \{~ b \in I ~~|~~ a ~\subseteq~ b ~\} $

\medskip
Then we have \\ \\

{\bf Lemma 5.6.3} \\

(5.6.11) \quad $ {\cal F} ~=~ \{~ J ~\subseteq~ I ~~|~~ \exists~~ a \in I ~:~ I_a ~\subseteq J ~\}$ \\

\medskip
is a filter on $I$, and it satisfies the condition \\

(5.6.12) \quad $ I ~\setminus~ \{~ a ~\} \in {\cal F},~~~ a \in I $ \\

{\bf Proof} \\

Let $A, B \in {\cal F}$, then there exist $a, b \in I$, such that $I_a \subseteq A, I_b \subseteq B$. But obviously

$$ A \bigcap B ~\supseteq~ I_a \bigcap I_b ~\subseteq~ I_{a ~\cup~ b} $$

\medskip
and $a \cup b \in I$. Thus ${\cal F}$ is indeed a filter on $I$. \\

Let us now show that (5.6.12) also holds. We take any $a \in I$. Then clearly there exists $b \in I$, with $a \cap b =
\phi$. Thus $a \notin I_b$, which means that $I_b \subseteq I \setminus \{~ a ~\}$, therefore $I \setminus \{~ a ~\}
\in {\cal F}$.

\hfill $\Box$ \\

Let ${\cal U}$ be any {\it ultrafilter} on $I$, such that ${\cal F} \subseteq{\cal U}$, see Appendix 1. Then ${\cal U}$ must
be {\it free}. Indeed, assume that for a certain $i \in I$, we have ${\cal U} = \{~ J \subseteq I ~|~ i \in J ~\}$. Let us take any
$F \in {\cal F}$, then $F \in {\cal U}$, hence $i \in F$, and (5.6.12) is contradicted. \\
Now, with $I$ given in (5.6.8) and with any free ultrafilter ${\cal U}$ on $I$ as chosen above, we have \\ \\

{\bf Theorem 5.6.1} \\

$V ( ^*X )$ is an extension of $V ( X )$. \\

{\bf Proof} \\

Let $A \in V ( X ) \setminus X$ and define a mapping $\gamma : I \longrightarrow {\cal P}_F ( A )$ by $\gamma ( a ) = a
\cap A$, for $a \in I$. \\
Now we choose $B = M ( [ \gamma ] )$. Then

??

\hfill $\Box$ \\

Let us note that once we have Theorem 5.6.1, we are no longer obliged to use the set $I$ in (5.6.8) and the particular
ultrafilters on it chosen above. \\
Indeed, as seen for instance when we consider the issue of {\it saturation}, it will be convenient to deal with arbitrary
extensions $V ( ^*X )$ of $V ( X )$, extensions which correspond to their respective sets $I$, and to the free ultrafilters
${\cal U}$ on them. \\

A Robinson himself first arrived to the existence of extensions in an alternative way which uses the concept of {\it
concurrent} relation. This concept has an interest of its own, in addition to having been basic to the original
development of Nonstandard Analysis. \\ \\

{\bf Definition 5.6.4} \\

Given a binary relation $P$ on a set $X$ and a subset $A \subseteq domain~ P$. Then $P$ is called {\it concurrent}, or
finitely satisfiable on $A$, if and only if for every finite set $\{ x_1, ~.~.~.~ , x_n \} \subseteq A$, there exists $y \in range~
P$, such that $< x_i, y > \in P$, for $1 \leq i \leq n$. \\
The binary relation $P$ is called {\it concurrent}, or finitely satisfiable, if and only if it is concurrent on $domain~ P$.

\hfill $\Box$ \\

Obviously, the binary relations $\leq$ and $<$ are concurrent on ${\bf N}$, as well as on ${\bf Z}$. The binary relation
$\subseteq$ is concurrent on ${\cal P} ( Y )$, for every set $Y$. Also, the binary relations $\subseteq$ and $\subset$
are concurrent on ${\cal P}_F ( {\bf N} )$ and ${\cal P}_F ( {\bf Z} )$. \\ \\

{\bf Theorem 5.6.2} \\

The following properties are equivalent \\

1.~~~ $V ( ^*X )$ is an extension of $V ( X )$. \\

2.~~~ For every concurrent relation $P \in V ( X )$, there exists $b \in range~ P$, such that $<~ ^*x, b >~ \in~ ^*P$, for $x \in
domain~ P$. \\

{\bf Proof} \\

??

\hfill $\Box$ \\ \\

{\bf Corollary 5.6.1} \\

Let $V ( ^*X )$ be an extension of $V ( X )$, and suppose given any infinite set $Y \in V ( X )$. Then $^*Y \in V ( ^*X)$
contains {\it nonstandard} entities. In particular, for every infinite subset $A \subseteq X$, we have $A \subseteq~ ^*A,~
A \neq~ ^*A$. \\

{\bf Proof} \\

??

\hfill $\Box$ \\

Let us now introduce another concept useful in a large range of applications. \\ \\

{\bf Definition 5.6.5} \\

A set ${\cal S}$ of subsets of a set $A \in V ( X )$ is called {\it exhausting} for $A$, if and only if

\bigskip
(5.6.13) \quad $ \begin{array}{l}
                            \forall~~ B \in {\cal P}_F ( A ) ~: \\ \\
                            \exists~~ S \in {\cal S} ~: \\ \\
                            ~~~ B ~\subseteq~ S
                 \end{array} $ \\ \\

{\bf Theorem 5.6.3} \\

If ${\cal S}$ is exhausting for $A \in V ( X )$, and we have an extension $V ( ^*X )$ of $V ( X )$, then

\bigskip
(5.6.12) \quad $ \begin{array}{l}
                           \exists~~ \Sigma \in~ ^*{\cal S} ~\subseteq~ V ( ^*X ) ~: \\ \\
                           ~~~ \{~ ^*a ~~|~~ a \in A ~\} ~\subseteq~ \Sigma
                 \end{array} $ \\

{\bf Proof} \\

??

\hfill $\Box$ \\

We shall indicate now a general pattern for the application of the result in Theorem 5.6.3. \\
Suppose we are given an {\it infinite} set $A$ that has a certain mathematical structure which can be {\it exhausted} by
a family ${\cal S}$ of its substructures. For instance, $A$ can be an infinite graph, as in Theorem 5.6.4 next, where
${\cal S}$ is the set of finite subgraphs of $A$. Or $A$ can be an infinite dimensional Hilbert space, with ${\cal S}$
being the set of all its finite dimensional Hilbert subspaces. \\
Further suppose that each $S \in {\cal S}$ has a certain property ${\cal P}$. \\
Then one may under certain conditions establish the same property ${\cal P}$ for $A$ itself. \\

?? \\ \\

{\bf Theorem 5.6.4 ( De Bruin - Erd\"{o}s )} \\

If every finite subgraph of an infinite graph is k-colourable, then the graph itself is k-colourable. \\ \\

{\bf Note 5.6.???} \\

By a graph we mean here a structure $( V, E )$, where $V$ is the set of {\it vertices} and $E \subseteq V \times V$ is the
set of {\it edges}. The graph $( V, E )$ is {\it infinite}, if and only if $V$ is infinite. A graph $( V^\prime, E^\prime )$ is a {\it
subgraph} of $( V, E )$, if and only if $V^\prime \subseteq V$ and $E^\prime = E \bigcap ( V^\prime \times V^\prime )$.
Finally, a graph $( V, E )$ is {\it k-colourable}, where $k \in {\bf N},~ k \geq 1$, if and only if there exists a mapping
$\kappa : V ~\longrightarrow~ \{~ 1, 2, 3, ~.~.~.~ , k ~\}$, such that $< a, b > \in E ~\Longrightarrow~ \kappa ( a ) \neq \kappa
( b )$. \\ \\

{\bf Proof} \\

??

{~} \\ \\

{\bf 7. Internal and External Entities, Comprehensiveness} \\

One of the most important uses of the Transfer Property at 5 in Definition 5.4.2, a use also seen in the respective
simpler form in chapter 3, is the following. We take a sentence $\Phi \in {\cal L}_X$, apply to it the transfer, and obtain
the sentence $^*\Phi \in {\cal L}_{^*X}$. Then we prove this transferred sentence $^*\Phi$ in ${\cal L}_{^*X}$, and finally,
we deduce from that, by a {\it reverse} use of the transfer, the truth of the initial sentence $\Phi$ in ${\cal L}_{^*X}$. \\

It follows that, in view of the mentioned reverse transfer stage, it is important to recognize when a given sentence
$\Psi \in {\cal L}_{^*X}$ is of the form

\bigskip
(5.7.1) \quad $ \Psi =~ ^*\Phi,~~~ \mbox{where}~~ \Phi \in {\cal L}_X $

\medskip
And needless to say, the majority of sentences $\Psi \in {\cal L}_{^*X}$ are {\it not} of that particular form, see Definition
5.4.1. \\

In this regard, from the start, we can easily note the following, based on (5.4.4). If  we are given any $\Psi \in
{\cal L}_{^*X}$ which is of the form (5.7.1), then $\Psi$ can only contain {\it constant symbols} $~^*a$ which correspond
through the transfer (5.4.4) to {\it constant symbols} $~a \in V ( X )$. Thus in view of Definition 5.6.2, $\Psi$ can only
contain constant symbols which are {\it standard}. In particular, $\Psi$ can only contain quantifiers of the form

$$ ( \forall~ x \in A ),~~~( \forall~ x \in~ ^*A ),~~~( \exists~ x \in A ),~~~( \exists~ x \in~ ^*A ) $$

\medskip
where $A \in V ( X )$. \\

In this way, in order to prove the truth in ${\cal L}_{^*X}$ of sentences of the form (5.7.1), we only have to consider
entities $b \in V ( ^*X )$ which satisfy $b \in~ ^*a$, for certain $a \in V ( X )$. Therefore, we are led to \\ \\

{\bf Definition 5.7.1} \\

Given an extension $V ( ^*X )$ of $V ( X )$, an entity $b \in V ( ^*X )$ is called {\it internal}, if and only if it satisfies $b
\in~ ^*a$, for a certain $a \in V ( X )$. Otherwise $b$ is called {\it external}. \\

A formula $\Psi \in {\cal L}_{^*X}$ is called {\it standard}, if and only if all the constant symbols which it contains are
standard. Otherwise $\Psi \in {\cal L}_{^*X}$ is called {\it nonstandard}. \\

Finally, a formula $\Psi \in {\cal L}_{^*X}$ is called {\it internal}, if and only if all the constant symbols which it
contains are internal. Otherwise $\Psi$ is called {\it external}. \\ \\

{\bf Note 5.7.1} \\

Let us recall for convenience the {\it difference} between the concepts of {\it standard} and {\it internal} entities in
$V ( ^*X )$, see Definitions 5.6.2 and 5.7.1, namely

$$ a \in V ( ^*X ) ~~\mbox{is standard} ~~~\Longleftrightarrow~~~ a ~=~ ^*b ~~\mbox{for some}~~ b \in V ( X ) $$
$$ a \in V ( ^*X ) ~~\mbox{is internal} ~~~\Longleftrightarrow~~~ a \in~ ^*b ~~\mbox{for some}~~ b \in V ( X ) $$ \\ \\

{\bf Proposition 5.7.1} \\

The concept of {\it internal} is more {\it general} than the concept of {\it standard}, as for every entity $a, b \in V ( ^*X )$,
we have the implication \\

\bigskip
(5.7.2) \quad $ a ~~\mbox{is standard} ~~~\Longrightarrow~~~ a ~~\mbox{is internal} $ \\

We also have the {\it transitivity} property

\bigskip
(5.7.3) \quad $ a ~~\mbox{is internal},~~ b \in a ~~~\Longrightarrow~~~ b ~~\mbox{is internal} $ \\

{\bf Proof} \\

If $a$ is standard then $a =~ ^*b$, for some $b \in V ( X )$. But obviously $b \in \{~ b ~\} \in V ( X )$, thus 5 in Theorem
5.4.1 gives $~^*b \in~ ^*\{~ b ~\}$, therefore $a \in~  ^*\{~ b ~\}$, and the proof of (5.7.2) is completed. \\

The assumption in (5.7.3) means that $b \in a \in~ ^*c$, for some $c \in V ( X )$, hence in view of  (5.4.3) and (5.1.2), we
have $b \in a \in~ ^*V_{n + 1} ( X )$, for a certain $n \in {\bf N}$. But then 4 in Definition 5.4.2 gives $b \in~ ^*V_n ( X )$,
and in view of (5.1.3), we have $V_n ( X ) \in V ( X )$, hence $b$ is internal. \\ \\

{\bf Important Note 5.7.2} \\

In view of (5.7.2), every standard entity is internal, and every standard formula $\Phi \in {\cal L}_{^*X}$ is internal. \\
Now in view of (5.7.3), it follows that in the case of any extension $V ( ^*X )$ of $V ( X )$, the Transfer Property at 5 in
Definition 5.4.2 refers {\it exclusively} to {\it standard} formulas $\Phi \in {\cal L}_{^*X}$ and {\it internal} entities $a \in~
^*V ( X )$ . \\
It follows that internal entities and standard formulas describe, correspondingly, the framework of the transfers in (TE)
and (TS), see for instance section 1, in chapter 4. \\

Such a limitation, however, does {\it not} at all mean that entities $a \in V ( ^*X )$, or formulas $\Phi \in {\cal L}_{^*X}$
which are not standard or internal are without interest in Nonstandard Analysis. Indeed, as indicated in the diagram in
chapter 0, such entities and formulas do play an important role, which is additional to the role played by the Transfer
Property. For instance, let us consider here just two simple examples among a considerable amount of similar ones. \\
First, the constant $mon ( 0 ) \in V ( ^*X )$, which is the set of {\it infinitesimals} in $^*{\bf R}$, and which ever since
Leibniz was an important, even if not rigorous, concept, is {\it not} internal, see Theorem 7.5.3 below. \\
It follows that $mon ( 0 )$ is {\it outside} of the ranges of the $^*$-transfer operator (5.4.3). \\

Second, let us consider the nonstandard characterization of bounded sequences of real numbers $(~ x_n ~)_{n \in
{\bf N}}$, given in Proposition 3.9.3, which according to (3.9.3), is expressed by the following sentence in
${\cal L}_{^*X}$

$$ (~ \forall~ n \in~ ^*{\bf N}_\infty ~)~[ x_n \in Gal ( 0 ) ~] $$

\medskip
Obviously, this sentence is {\it not} internal, since the constants $^*{\bf N}_\infty$ and $Gal ( 0 )$ which it contains are
not internal, see Theorem 7.5.3 below. Thus this sentence is {\it outside} of the ranges of the Transfer Property, even if
within its proof that property was used. \\

Nevertheless, many such entities and sentences can prove to be particularly useful from mathematical point of view,
since they can give a rigorous understanding which is not available in usual mathematics. \\ \\

{\bf Theorem 5.7.1} \\

The set of all {\it internal} entities in $V ( ^*X )$ is given by

\bigskip
(5.7.4) \quad $ ^*V ( X ) ~=~ \bigcup_{n \in {\bf N}}~ ^*V_n ( X ) $

\medskip
where we recall (5.4.2$^{**}$) and (5.4.2$^{***}$) for the meaning of $^*V ( X )$. \\

It follows that each {\it hyperfinite} subset of $^*X$ is {\it internal}. \\

On the other hand, for every $k \in~ ^*{\bf N}_\infty ~=~ ^*{\bf N} \setminus {\bf N}$, the set

$$ \{~ l \in~ ^*{\bf N} ~~|~~ l ~\leq~ k ~\} $$

\medskip
is {\it hyperfinite}, thus {\it internal}, but {\it not standard}. \\

{\bf Proof} \\

The relation (5.7.4) follows from (5.4.12). \\
Let now be given any internal entity $A \in V ( ^*X )$. Then by definition, we have $A \in~ ^*B$, for a certain $B \in V ( X )$.
It follows from (5.1.2) that $B \in V_{n + 1} ( X ) \setminus V_n ( X )$, for a certain $n \in {\bf N}$, which means
that $B \subseteq V_n ( X )$. Thus 7 in Theorem 5.4.1 gives $^*B \subseteq~ ^*V_n ( X )$, and then $A \in~ ^*V_n (
X )$. \\
Conversely, if $A \in~ ^*V ( X )$, then (5.7.4) gives $A \in~ ^*V_n ( X )$, for a certain $n \in {\bf N}$. But in view of (5.1.3),
we have $V_n ( X ) \in V ( X )$, thus indeed $A$ is internal. \\

In view of (5.6.4) in Lemma 5.6.1, hyperfinite sets are obviously internal. \\

Finally, the sets $\{~ l \in~ ^*{\bf N} ~~|~~ l ~\leq~ k ~\}$, with $k \in~ ^*{\bf N}_\infty ~=~ ^*{\bf N} \setminus {\bf N}$, are
hyperfinite, according to Lemma 5.6.2. But they are not standard. Indeed, assume $k \in~ ^*{\bf N}_\infty ~=~ ^*{\bf N}
\setminus {\bf N}$ and $A \in V ( X )$, such that

$$ \{~ l \in~ ^*{\bf N} ~~|~~ l ~\leq~ k ~\} ~=~ ^*A $$

\medskip
then (5.4.14) gives $A \subseteq~ ^*A = \{~ l \in~ ^*{\bf N} ~~|~~ l ~\leq~ k ~\} \subset~ ^*{\bf N} \subset~ ^*X $, hence $A \in
V_1 ( ^*X )$. And then 4 in Definition 5.4.2 together with (5.1.9) give

$$ a \in A ~~~\Longrightarrow~~~ a \in V_0 ( ^*X ) ~=~ ^*X \bigcap V ( X ) ~=~ X $$

\medskip
hence $A \subseteq X \cap~ ^*{\bf N} \subseteq {\bf N}$. Furthermore, we note that $A$ is infinite, see 1 in Theorem
5.4.1. Thus $A$ is an infinite subset of ${\bf N}$. \\

Since $A \subseteq X$, (5.5.9) - (5.5.11) give

$$ ^*A ~=~ M ( e ( A ) ) ~=~ M ( [ \hat{A} ] ) ~=~~~~~~~~~~~~~~~~~~~~~~~~~~~~~~~~~~~~~~~~~~~~~~~~~~~ $$
$$ ~=~\{~ M ( [ b ] ) ~~|~~ [ b ] \in \prod_{~{\cal U}}~ X ~=~ ^*X,~~~ [ b ] \in_{\cal U}~ [ \hat{A} ] ~\} ~=~
                                             \{~ [ b ] ~~|~~ [ b ] \in_{\cal U}~ [ \hat{A} ] ~\} $$

\medskip
which means that

$$ ^*A ~=~  \{~ [ b ] ~~|~~ b : I \longrightarrow V ( X ),~~~ \{~ i \in I ~~|~~ b ( i ) \in A ~\} \in {\cal U} ~\} $$

\medskip
Let us now recall that by assumption, $k \in~ ^*A$ is the largest element of $^*A$. Let us take any $b  : I \longrightarrow
V ( X )$, such that $k = [ b ]$. Then $ \{~ i \in I ~~|~~ b ( i ) \in A ~\} \in {\cal U}$. \\
Now we define $b : I \longrightarrow V ( X )$ as follows. Given $i \in I$, then

$$ c ( i ) ~=~ \begin{array}{|l} ~~ a ( i ) ~~~\mbox{if}~~ b ( i ) \in A \\ \\
                                              ~~ b ( i ) ~~~\mbox{if}~~ b ( i ) \notin A
                   \end{array} $$

\medskip
where $a ( i ) \in A,~ a ( i ) > b ( i )$, which is possible, since $A$ is infinite, thus unbounded in ${\bf N}$. But then,
according to the above relation giving $^*A$, it follows that $[ c ] \in~ ^*A$, while on the other hand, obviously, $[ c ] \in~
^*{\bf N}$ and $[ c ] > k$, thus contradicting that $k$ is the largest element of $^*A$.

\hfill $\Box$ \\

The above result identifies the set of all internal entities in $V ( ^*X )$. In addition to it, it is obviously particularly useful
to be able to identify as well individual internal entities in $V ( ^*X )$. This can be done as follows \\ \\

{\bf Theorem 5.7.2 ( Keisler's Internal Definition )} \\

Given an internal formula $\Phi ( x ) \in {\cal L}_{^*X}$ in which $x$ is the only free variable, and given an internal set
$A \in~ ^*V ( X )$, then

\bigskip
(5.7.5) \quad $ \{~ x \in A ~~|~~ \Phi ( x ) ~~\mbox{is true in}~~ {\cal L}_{^*X} ~\} ~~~\mbox{is internal} $ \\

{\bf Proof} \\

Let $c_1, ~.~.~.~ , c_n$ be all the constant symbols in $\Phi ( x )$, and then we shall write $\Phi ( x ) = \Phi ( c_1, ~.~.~.~ ,
c_n, x )$. According to Definition 5.7.1, all of $c_1, ~.~.~.~ , c_n$ are internal, hence we have $A, c_1, ~.~.~.~ , c_n \in~
^*V_m ( X )$, for a certain $m \in {\bf N}$. It follows that the sentence in ${\cal L}_X$

$$ (~ \forall~ x_1, ~.~.~.~ , x_n, y \in V_m ( X ) ~)~ (~ \exists~ z \in V_{m + 1} ( X ) ~)~ (~ \forall~ x \in V_m ( X ) ~) $$
$$ [~ x \in z \leftrightarrow [~ x \in y \wedge \Phi ( x_1, ~.~.~.~ , x_n, x ) ~] ~] $$

\medskip
is true in $V ( X )$. Thus the $^*$-transfer of it is true in $V ( ^*X )$. And in view of that sentence we have

$$  \{~ x \in A ~~|~~ \Phi ( x ) ~~\mbox{is true in}~~ {\cal L}_{^*X} ~\} \in~ ^*V_{m + 1} ( X ) $$.

\hfill $\Box$ \\

Several useful properties related to internal sets are collected in \\ \\

{\bf Theorem 5.7.3} \\

1.~~ If $A$ and $B$ are internal, then so are $A \bigcup B,~ A \bigcap B,~ A \setminus B$ and $A \times B$. \\

2.~~ If $A \in V ( X ) \setminus X$, then \\

(5.7.6) \quad $ \{~ a \in {\cal P} ( ^*A ) ~~|~~ a ~~\mbox{internal} ~\} ~=~ ^*{\cal P} ( A ) $ \\

3.~~ Let $K \subseteq~ ^*{\bf Z}$ be a nonvoid internal subset. If $K$ is bounded from below, then it has a least element,
and if it is bounded from above, then it has a largest element. \\

4.~~ The sets \\

(5.7.7) \quad $ \begin{array}{l}
                            {\bf N},~~~ {\bf Z},~~~ {\bf R},~~~ {\bf Q},~~ ^*{\bf N}_\infty,~~ ^*{\bf Z}_\infty,~~ ^*{\bf Q}_\infty,~~
                                          ^*{\bf R}_\infty, \\ \\
                             ^*{\bf R}^+ ~=~ \{~ s \in~ ^*{\bf R} ~~|~~ s > 0 ~\},~~ ^*{\bf R}^+_\infty,~~ mon ( 0 ),~~~ Gal ( 0 )
                         \end{array} $ \\

are external. \\

{\bf Proof} \\

?? \\ \\

{\bf Note 5.7.3} \\

Internal entities play an important role in Nonstandard Analysis, as pointed out for instance in Note 5.7.2. In view of
that the relation (5.7.4), namely

$$  ^*V ( X ) ~=~ \bigcup_{n \in {\bf N}}~ ^*V_n ( X ) $$

\medskip
which gives the set of all {\it internal} entities in $V ( ^*X )$, is particularly useful, see (5.4.2$^{**}$) and (5.4.2$^{***}$). \\

In the definition of Superstructures, see (5.1.1), (5.1.2), the set theoretic operation $Y ~\longmapsto~ {\cal P} ( Y )$,
which to each set $Y$ associates the set ${\cal P} ( Y )$ of all its subsets is fundamental. Therefore, it is important to
have as well an appropriate description of the $^*$-transfer operation

\bigskip
(5.7.8) \quad $ V ( X ) \setminus X \ni A ~~\longmapsto~~ ^*{\cal P} ( A ) \in V ( ^*X ) $

\medskip
This is done in (5.7.6) in the convenient terms of {\it internal} sets, namely

\bigskip
(5.7.9) \quad $ ^*{\cal P} ( A ) ~=~ \{~ B ~\subseteq~ ^*A ~~|~~ B ~~\mbox{internal} ~\},~~~ A \in V ( X ) \setminus X $

\medskip
and this relation further highlights the importance of {\it internal} sets in $V ( ^*X )$ as being precisely those which do
so simply and explicitly describe the above fundamental $^*$-transfer operation (5.7.8). \\

Here related to the above, we can recall that the {\it hyperfinite} sets are given by, see (5.6.3), (5.6.3$^*$), (5.6.4)

$$ \bigcup_{n \in {\bf N}}~ ^*{\cal P} ( V_n ( X ) ) $$

\medskip
while similar with the above expression in (5.7.9) for $^*{\cal P} ( A )$, we have

\bigskip
(5.7.10) \quad $ ^*{\cal P}_F ( A ) ~=~ \{~ B ~\subseteq~ ^*A ~~|~~ B ~~\mbox{hyperfinite} ~\},~~~A \in V ( X ) \setminus X  $

\medskip
a relation which further justifies the importance of {\it hyperfinite} sets in $V ( ^*X )$. \\

One can note the parallelism between (5.7.9) and (5.7.10), and respectively, the concept of {\it internal}, and on the other
hand, {\it finite} and {\it hyperfinite}.

\hfill $\Box$ \\

We return now to the characterization of infinite hyperfinite sets by initial segments, the first part of which was
presented in Proposition 5.6.1. We note that the bijection $\beta$ obtained there between the initial segment $K$ and
the hyperfinite set $B$ is {\it internal}, since we had $\beta \in~ ^*V_{n + 3} ( X )$, see (5.7.4). \\

The converse of that result is given now in \\ \\

{\bf Proposition 5.7.2} \\

Given any set $B \in V ( ^*X )$. If there exists an internal bijection $\beta : K \longrightarrow B$, where $K = \{~ 1, 2, 3,
~.~.~.~ , k ~\}$, for a certain $k \in~ ^*{\bf N}_\infty$, then $B$ is hyperfinite. \\

{\bf Proof} \\

?? \\ \\

{\bf Note 5.7.4} \\

The characterization of hyperfinite sets obtained in Propositions 5.6.1 and 5.7.2 is that, for any given $B \in V ( ^*X )$,
we have

$$ B ~~\mbox{hyperfinite} ~~~\Longleftrightarrow~~ \exists~~ \mbox{internal bijection}~~ \beta : K \longrightarrow B $$

\medskip
for a suitable initial segment $K$ in $^*{\bf N}_\infty$. \\
This however does {\it not} mean that being a hyperfinite set is merely a matter of cardinality. Indeed, the bijections in
the above characterization are {\it not} arbitrary, since they must be {\it internal} functions.

\hfill $\Box$ \\

Another useful concept is presented in \\ \\

{\bf Definition 5.7.2} \\

The monomorphism $^*(~) : V ( X ) \longrightarrow V ( ^*X )$ in (5.4.1) is called {\it comprehensive}, if and only if

\bigskip
(5.7.11) \quad $ \begin{array}{l}
                      \forall~ A, B \in V ( X ) \setminus X,~~ f : A \longrightarrow~ ^*B ~: \\ \\
                      \exists~ g :~ ^*A \longrightarrow~ ^*B,~~ g ~~\mbox{internal} ~: \\ \\
                      ~~~~ f ( a ) ~=~ g ( ^*a ), ~~\mbox{for}~~ a \in A
                \end{array} $

\medskip
The monomorphism $^*(~) : V ( X ) \longrightarrow V ( ^*X )$ in (5.4.1) is called {\it countably comprehensive}, if and
only if (5.7.11) only holds for countable $A$. \\ \\

{\bf Theorem 5.7.4} \\

The monomorphism $^*(~) : V ( X ) \longrightarrow V ( ^*X )$ constructed in (5.5.1) is comprehensive. \\

{\bf Proof} \\

In view of (5.5.8), we have

??

\hfill $\Box$ \\

We can illustrate the relationship between standard, hyperfinite and internal sets, as follows \\

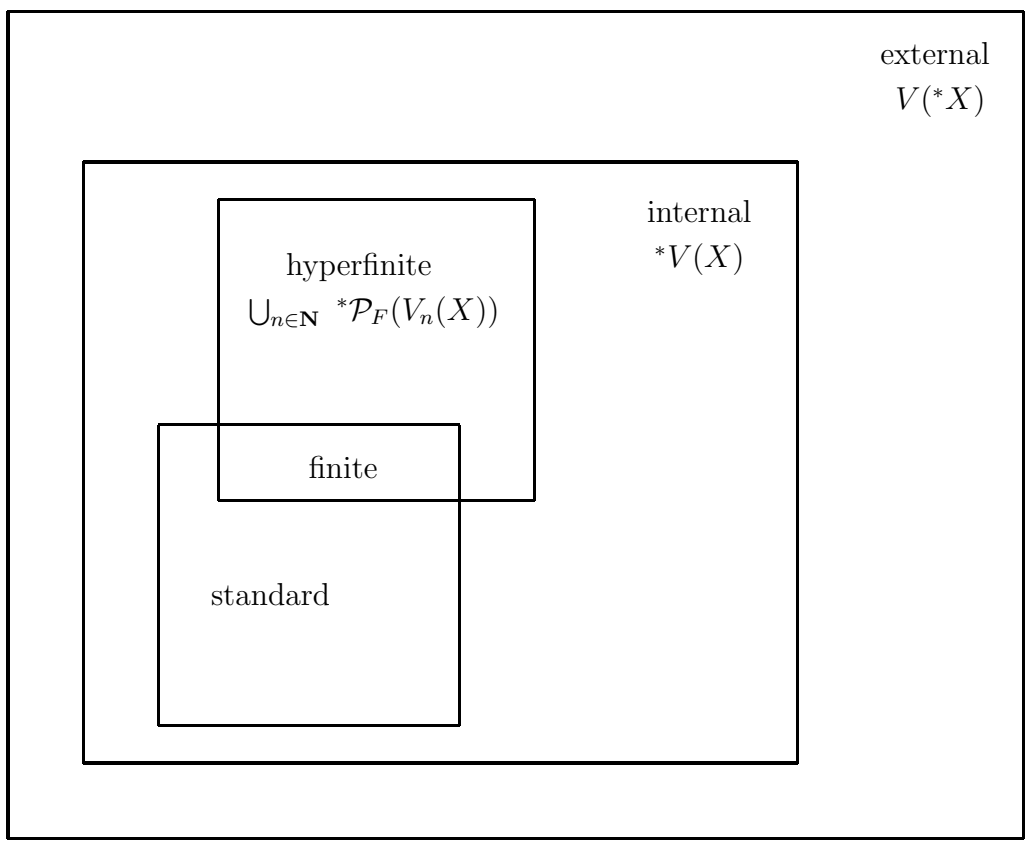
\begin{math}
\setlength{\unitlength}{1cm}
\thicklines
\begin{picture}(13,12)

\put(2.8,4.5){\line(0,1){4}}
\put(2.8,4.5){\line(1,0){4.2}}
\put(2.8,8.5){\line(1,0){4.2}}
\put(7,4.5){\line(0,1){4}}
\put(3.7,7.5){$\mbox{hyperfinite}$}
\put(3.2,6.9){$\bigcup_{n \in {\bf N}}~ ^*{\cal P}_F ( V_n ( X ) )$}

\put(2,1.5){\line(0,1){4}}
\put(2,1.5){\line(1,0){4}}
\put(2,5.5){\line(1,0){4}}
\put(6,1.5){\line(0,1){4}}
\put(2.7,3.1){$\mbox{standard}$}
\put(4,4.8){$\mbox{finite}$}

\put(1,1){\line(0,1){8}}
\put(1,1){\line(1,0){9.5}}
\put(1,9){\line(1,0){9.5}}
\put(10.5,1){\line(0,1){8}}
\put(8.5,8.2){$\mbox{internal}$}
\put(8.6,7.6){$^*V ( X )$}

\put(0,0){\line(0,1){11}}
\put(0,0){\line(1,0){13.5}}
\put(0,11){\line(1,0){13.5}}
\put(13.5,0){\line(0,1){11}}
\put(11.6,10.3){$\mbox{external}$}
\put(11.8,9.7){$V ( ^*X )$}
\end{picture}
\end{math} \\ \\

{\bf Note 5.7.5} \\

Let us show that, indeed, for every $A \in V ( ^*X )$, we have

$$ A ~~\mbox{standard and hyperfinite} ~~~\Longleftrightarrow~~~ A ~~\mbox{finite} $$

\hfill $\Box$ \\

Let us end this section by showing that the concept of {\it internal} subset $A$ of $^*{\bf R}$, given in (2.8.1), is indeed a
particular case of that in Definition 5.7.1. \\
We start by noting that the sequence of subsets $A_i \subseteq {\bf R}$, with $i \in {\bf N}$, in (2.8.1), can be identified
with the binary relation $P \subseteq {\bf N} \times {\bf R}$, defined by $P < i, x > ~\Longleftrightarrow~ x \in A_i$, where
$i \in {\bf N},~ x \in {\bf R}$. \\
Then according to Definition 2.7.1, for $ n =~ < n_0, n_1, n_2, ~.~.~.~ >~ \in {\bf N}^{\bf N},~ s =~ < s_0, s_1, s_2, ~.~.~.~ >~ \in
{\bf R}^{\bf N}$, we have

$$ ^*P < [ n ], [ s ] > ~~~\Longleftrightarrow~~~ \{~ i \in {\bf N} ~~|~~ s_i \in A_{n_i} ~\} \in {\cal U} $$

\medskip
Let us now take any $s =~ < s_0, s_1, s_2, ~.~.~.~ >~ \in {\bf R}^{\bf N}$, then (2.8.1) gives

$$ [ s ] \in A ~~~\Longleftrightarrow~~~  \{~ i \in {\bf N} ~~|~~ s_i \in A_i ~\} \in {\cal U} $$

\medskip
therefore

\bigskip
(5.7.12) \quad $ [ s ] \in A ~~~\Longleftrightarrow~~ ^*P < [ id ], [ s ] > ~~~\Longleftrightarrow~~~ < [ id ], [ s ] >~ \in~ ^*P $

\medskip
where we denoted $id =~ < 0, 1, 2, ~.~.~.~ >~ \in {\bf N}^{\bf N}$. It follows that, according to Definition 5.2.1, the formula
$\Phi ( x )$, in fact, an atomic one, given by

$$ < \eta, x >~ \in~ ^*P $$

\medskip
where $\eta = [ id ] \in~ ^*{\bf N}$ is a constant symbol, while $x$ is a variable symbol, belongs to the language ${\cal
L}_{^*X}$. But in view of Definition 5.7.1, this formula is internal, since $^*P$ and $^*{\bf N}$ are standard, and $\eta
\in~ ^*{\bf N}$. Therefore Theorem 5.7.2 implies that the set

$$ \{~ x \in~ ^*{\bf R} ~~|~~ \Phi ( x ) ~~\mbox{is true in}~~ {\cal L}_{^*X} ~\} $$

\medskip
is internal, since $^*{\bf R}$ being standard, it is also internal. However, in view of (5.7.12), the above set is precisely
$A$. \\ \\

{\bf 8. Permanence} \\ \\

{\bf Theorem 5.8.1} \\

Given any internal formula $\Phi ( x )$ in ${\cal L}_{^*X}$ in which $x$ is the only free variable. Then we have the {\it
overflow} properties : \\

1.~~ If $\Phi ( x )$ holds for every $x \in {\bf N}$, then there exists $k \in~ ^*{\bf N}_\infty$, such that $\Phi ( n )$ holds for
every $n \in~ ^*{\bf N},~~ n \leq k$. \\

2.~~ If $\Phi ( x )$ holds for every $x \in~ ^*{\bf R}^+ \setminus~ ^*{\bf R}_\infty$, then there exists $k \in~ ^*{\bf N}_\infty$,
such that $\Phi ( x )$ holds for every $x \in~ ^*{\bf R}^+,~~ x \leq k$. \\

Further, we have the {\it underflow} properties : \\

3.~~ If $\Phi ( x )$ holds for every $x \in~ ^*{\bf N}_\infty$, then there exists $k \in {\bf N}$, such that $\Phi ( n )$ holds for
every $x \in~ ^*{\bf N},~~ x \geq k$. \\

4.~~ If $\Phi ( x )$ holds for every $x \in~ ^*{\bf R}^+_\infty$, then there exists $k \in {\bf N}$, such that $\Phi ( x )$ holds
for every $x \in~ ^*{\bf R}^+,~~ x \geq k$. \\

Also, we have the {\it local overflow} property :  \\

5.~~ If $\Phi ( x )$ holds for every $x \in mon ( 0 )$, then there exists $r \in {\bf R},~~ r > 0$, such that $\Phi ( x )$ holds for
every $x \in~ ^*{\bf R},~~ | x | \leq r$. \\

{\bf Proof} \\

??

{\bf Corollary 5.8.1} \\

Let be given any internal subset $A \subseteq~ ^*{\bf R}$, then

\bigskip
(5.8.1) \quad $ \begin{array}{l}
                                 {\bf N} ~\subseteq~ A ~~~\Longrightarrow~~~ A ~\bigcap~ ^*{\bf N}_\infty ~\neq~ \phi \\ \\
                                 {\bf N}_\infty ~\subseteq~ A ~~~\Longrightarrow~~~ A ~\bigcap~ {\bf N} ~\neq~ \phi \\ \\
                                 mon ( 0 ) ~\bigcap~ {\bf R}^+ ~\subseteq~ A ~~~\Longrightarrow~~~ A ~\bigcap~ {\bf R}^+ ~\neq~ \phi
                        \end{array} $ \\

{\bf Proof} \\

{\bf Theorem 5.8.2 ( Robinson )} \\

Let $^*{\bf N} \ni n \longmapsto s_n \in~ ^*{\bf R}$ be an internal sequence such that

$$ s_n \in mon ( 0 ),~~~\mbox{for}~~ n \in {\bf N} $$

\medskip
Then there exists $\omega \in~ ^*{\bf N}_\infty$, such that

$$ s_n \in mon ( 0 ),~~~\mbox{for}~~ n \in~ ^*{\bf N},~~ n \leq \omega $$ \\

{\bf Proof} \\

??

\hfill $\Box$ \\

Here, according to A Robinson, we give an application of theorem 5.8.2 to the construction of {\it Banach limits of
bounded sequences} of real numbers. \\
Let us denote by $l_\infty$ the set of all the bounded sequences ${\bf N} \ni n \longmapsto s_n \in {\bf R}$, a set which
obviously contains the vector space $c$ of such convergent sequences. We shall consider a class of linear functionals
$L : l_\infty \longrightarrow {\bf R}$  which extend the usual limit $lim : c  \longrightarrow {\bf R}$, and do so in the
sense of \\ \\

{\bf Definition 5.8.1} \\

A linear functional $L : l_\infty \longrightarrow {\bf R}$ is called a {\it Banach limit}, if and only if for every sequence $s
\in l_\infty$ given by ${\bf N} \ni n \longmapsto s_n \in {\bf R}$ in $l_\infty$ we have

$$ \begin{array}{l}
                      \liminf_{n \in {\bf N}}~ s_n ~\leq~ L ( s ) ~\leq~ \limsup_{n \in {\bf N}} \\ \\
                      L ( s^+ )  ~=~ L ( s )
     \end{array} $$

\medskip
where $s^+ \in l_\infty$ is the shifted sequence with $s^+_n = s_{n + 1}$, for $n \in {\bf N}$.

\hfill $\Box$ \\

Now by using nonstandard methods, we shall {\it construct} such Banach limits

?? \\ \\

{\bf 9. Saturation} \\ \\

There are at least two aspects of the Superstructures $V ( ^*X )$ which lead to the concept of {\it saturation}. \\

First, as seen in section 5, for a given set $X$ one can construct more than one such  Superstructure $V ( ^*X )$,
depending on the choice of the sets $I$ and of the ultrafilters ${\cal U}$ on these sets $I$. And it is well known, Stroyan
\& Luxemburg, Keisler [2], Hurd \& Loeb, Goldblatt, that for a given set $X$ one can construct associated
Superstructures $V ( ^*X )$ of {\it arbitrary} large cardinals. In particular, one can have the set $^*{\bf R}$ of nonstandard
reals with arbitrary large cardinal. \\
This {\it lack} of uniqueness of Superstructures associated with a given set, or for that matter, of $^*{\bf R}$ itself may
in certain situation be undesirable. Saturation can avoid reestablish the uniqueness of Superstructures. \\

Second, internal sets prove to be very large, when they are not finite. This can prove to be convenient, and can be
further enhanced by saturation. \\

Let us start for illustration with a result about the latter phenomenon. \\ \\

{\bf Theorem 5.9.1} \\

??

\chapter{Loeb Integration and Measures}

Integration theory, and in particular, the study of measure and probability in a nonstandard setup did encounter early on
a crucial difficulty. Namely, the $^*$-transfer of a $\sigma$-additive measure is in general {\it not} again a $\sigma$-
additive measure. And in fact, the difficulties already arise with $\sigma$-algebras, when they are considered within a
nonstandard context, see Corollary 2.8.1. \\
These difficulties were overcome in the early 1970s by P Loeb. \\

For a better understanding of the issues involved, it is useful to return in short to the main reasons which led to the
development of modern integration theory in its standard form. One important reason why we could not remain with
the earlier Riemann integration, and instead, we had to develop the Lebesgue integration is related to the so called
{\it monotone limit} properties. \\

This issue of monotone limit properties can shortly be formulated as follows. Suppose on a certain interval $[ a, b ]
\subset {\bf R}$, we are given an {\it increasing} sequence of real valued functions $f_0 \leq f_1 \leq f_2 \leq ~.~.~.~$
which converges in some sense to a function $f$. In case all the functions in the sequence are integrable in a certain
sense, the following two {\it limit} questions arise : \\

- whether the limit function $f$ is also integrable, \\

and in case it is, \\

- whether in addition we have

$$ \lim_{n \to \infty}~ \int~ f_n ~=~ \int~ f $$

\medskip
As it happens, for the Riemann integral such monotone limit properties only hold under rather restrictive conditions.
The Lebesgue integral, on the other hand, has these two monotone limit properties in a significantly more general
setup. \\
Also the set of Lebesgue integrable functions is significantly larger than that of the Reimann integrable ones. \\

An alternative and more general approach than that of the Lebesgue integration was developed by P Daniel in the
1910s, based precisely on securing satisfactory monotone limit properties, and involving from the start partially order
structures on the real valued functions involved in the integration process. \\
According to that approach, one considers an arbitrary set $E$, a vector lattice ${\cal F}$ of real valued functions on
$E$, and a positive linear functional $I$, which is supposed to be the integral, defined on ${\cal F}$. \\
If we now assume that the respective {\it integration structure} $( {\cal F}, I )$ satisfies the {\it continuity} type condition

$$ \left (~~ f_n ~~\mbox{decreases to}~~ 0, ~~\mbox{when}~~ n \to \infty ~~\right ) ~~~\Longrightarrow~~~
                           \lim_{n \to \infty}~ I f_n  ~=~ 0 $$

\medskip
where $f_n \in {\cal F}$, then one can extend $( {\cal F}, I )$ to and integration structure $( {\cal F}^{\#}, I^{\#} )$
which will have the desired monotone limit properties. Furthermore, under suitable choices of $( {\cal F}, I )$, the
resulting extended integration structure $( {\cal F}^{\#}, I^{\#} )$ will yield the Lebesgue integral. \\

The Nonstandard Integration, introduced by P Loeb in the 1970s, see Hurd \& Loeb, follows the Daniel approach. An
advantage of such an approach is that it develops simultaneously a nonstandard treatment of both integration and
measure theory, and thus in particular, of probability as well. \\

One starts with an arbitrary set $X$ and an extension of it given by a monomorphism $^*(~) : V ( X ) \longrightarrow V
( ^*X )$. \\
Then one chooses any {\it internal} set $E \in~ ^*V ( X )$, and on $E$, one takes any {\it internal} integration structure
$( {\cal F}, I )$ on $E$, that is, with internal ${\cal F}, I \in~ ^*V ( X )$. \\
One of the remarkable aspects of this Nonstandard Integration procedure is that one can start with a structure $( {\cal
F}, I )$ on $E$ which need {\it not} satisfy any continuity type property, yet still be able to construct a standard
extension $( {\cal F}^{\#}, I^{\#} )$ on the same set $E$, an extension which will nevertheless have the monotone limit
 properties. \\
This process is called {\it standardization}, see section 6.2. \\
Needless to say, one can also recover as a particular case the Lebesgue integral. \\

One should note that, as long pointed out by the method of Daniel, integration is in essence about the {\it extension of
positive linear functionals}. Therefore, the integral is supposed to be defined upon an ordered vector space of functions,
which at the least should be a vector lattice. As far as the mentioned {\it monotonicity} properties of the integral are
concerned, they correspond to two things, namely, a sequential or $\sigma$ {\it completeness} of the vector lattice of
functions on which the integral is defined, and a corresponding {\it continuity} property of the integral. \\ \\

{\bf 1. Standard and Nonstandard Integration Structures} \\ \\

{\bf Definition 6.1.1} \\

Let be given any subset $E \subseteq X$ and any set ${\cal F} \subseteq {\bf R}^E$ of real valued functions on $E$.
Then ${\cal F}$ is called a {\it vector lattice} on $E$, if and only if \\

1.~~~ ${\cal F}$ is a vector space over ${\bf R}$, \\

2.~~~ $f \in {\cal F} ~~~\Longrightarrow~~~ | f | \in {\cal F}$ \\

We call ${\cal F}$ a {\it nonstandard vector lattice} on $E$, if an only if the functions in ${\cal F}$ can have nonstandard
values as well, that is, ${\cal F} \subseteq (^*{\bf R})^E$, and ${\cal F}$ is a vector space over $^*{\bf R}$.

\hfill $\Box$ \\

Given ${\cal F}$ a vector lattice, or a nonstandard vector lattice, we can define the operations min, max : ${\cal F} \times
{\cal F} \longrightarrow {\cal F}$ by

$$ \min~ ( f, g ) ~=~ ( f + g - | f - g |~ ) / 2,~~ \max~ ( f, g ) ~=~ ( f + g + | f - g |~ ) / 2,~~~ f, g \in
                        {\cal F} $$

\medskip
Often the alternative notation is also useful

$$ \min~ ( f, g ) ~=~ f \wedge g,~~~ \max~ ( f, g ) ~=~ f \vee  g $$

\medskip
Further, we can define the negative and positive part by

$$ f^{-} ~=~ ( - f ) \wedge 0,~~ f^{+} ~=~ f \wedge 0,~~~ f \in {\cal F} $$

\medskip
Finally, we can define a partial order on ${\cal F}$ by \\

$ f ~\leq~ g ~~~\Longleftrightarrow~~~ f \wedge g ~=~ f
                        ~~~\Longleftrightarrow~~~ f \vee g ~=~ g $ \\

It is easy to see that this partial order is compatible with the vector space structure of ${\cal F}$. \\ \\

{\bf Definition 6.1.2} \\

Let ${\cal F}$ be a vector lattice on $E$. Then a functional $I : {\cal F} \longrightarrow {\bf R}$ is called an {\it
integration} on ${\cal F}$, if and only if it is linear over ${\bf R}$, and it is positive, namely, $f \in {\cal F},~ f \geq 0
~~~\Longrightarrow~~~ I f \geq 0$ \\
In such a case $( {\cal F}, I )$ is called an {\it integration structure} on $E$. \\

A functional $I : {\cal F} \longrightarrow~ ^*{\bf R}$ is called a {\it nonstandard integration} on ${\cal F}$, if and
only if ${\cal F}$ is a nonstandard vector lattice, $I$ is linear over $^*{\bf R}$, and it is also positive. \\
Correspondingly, we call $( {\cal F}, I )$ a {\it nonstandard integration structure} on $E$. \\ \\

{\bf Example 6.1.1} \\

Let us recall in short Example 5.6.1, and relate it to the above definitions. We define $\Sigma : {\bf N} \times
{\bf R}^{\bf N} \longrightarrow {\bf R}$ by

$$ \Sigma ( n, f ) ~=~ \sum_{0 \leq i \leq n}~ f ( i ),~~~ n \in {\bf N},~ f \in {\bf R}^{\bf N} $$

\medskip
It follows that in ${\cal L}_X$ we have the true sentence

$$ (~ \forall~ n \in {\bf N} ~)~(~ \forall~ f \in {\bf R}^{\bf N} ~)~(~ \exists~ s \in {\bf R} ~)~[~ \Sigma ( n, f ) = s ~] $$

\medskip
Thus by $^*$-transfer, we have in ${\cal L}_{^*X}$ the true sentence

$$ (~ \forall~ n \in~ ^*{\bf N} ~)~(~ \forall~ f \in~ ( ^*{\bf R} )^{( ^*{\bf N} )} ~)~(~ \exists~ s \in~ ^*{\bf R} ~)~[~ ^*\Sigma ( n, f )
                              = s ~] $$

\medskip
which in particular means that for every $f :~ ^*{\bf N} \longrightarrow~ ^*{\bf R}$ and $k \in~ ^*{\bf N}_\infty =~ ^*{\bf N}
\setminus {\bf N}$, we can {\it always} effectuate the {\it uncountably} infinite sum

$$ \sum_{0 \leq l \leq k}~ f ( l ) \in~ ^*{\bf R}~~~~~~~~~~~~~~~~~~~~~~~~~~~~~~~~~~~~~~~~~~~~~~~~~~~~ $$

\medskip
and as a result, obtain a nonstandard real number. \\

Now, if we take $E = {\bf N},~ {\cal F} = {\bf R}^{\bf N}$, and for a given and fixed $n \in {\bf N}$, we take $I f = \Sigma
( n, f )$, with $f \in {\bf R}^{\bf N}$, then obviously $( {\cal F}, I )$ is an integration structure on $E$. \\

The nontrivial fact obtained above by $^*$-transfer is the following. Let us take this time $E =~ ^*{\bf N},~ {\cal F} =
( ^*{\bf R} )^{( ^*{\bf N} )}$, and fix any hyperfinite $k \in~ ^*{\bf N}_\infty =~ ^*{\bf N} \setminus {\bf N}$. Then as we have
seen, we can define $I : ( ^*{\bf R} )^{( ^*{\bf N} )} \longrightarrow~ ^*{\bf R}$ by

$$ I f ~=~ \sum_{0 \leq l \leq k}~ f ( l ),~~~ f \in ( ^*{\bf R} )^{( ^*{\bf N} )} $$

\medskip
and in this case $( {\cal F}, I )$ will be a nonstandard integration structure on $E$.

\hfill $\Box$ \\

The above example is in fact a particular case of \\ \\

{\bf Proposition 6.1.1} \\

If $( {\cal F}, I )$ is an integration structure on $E$, then $( ^*{\cal F},~ ^*I )$ is a nonstandard integration structure on
$^*E$. \\

{\bf Proof} \\

It follows easily by transfer.

\hfill $\Box$ \\

{\bf Definition 6.1.3} \\

Given two integration, or nonstandard integration structures $( {\cal F}, I )$ and $( {\cal G}, J )$ on the sets $E$
and $F$, respectively, the second one is an {\it extension} of the first one, if and only if $E \subseteq F$, ${\cal F}
\subseteq {\cal G}$ and $J|_{{\cal F}} = I$.

\hfill $\Box$ \\

In this way, in Proposition 6.1.1 it is clear that $( ^*{\cal F}, ^*I )$ is an extension of $( {\cal F}, I )$. \\

It is however important to note that {\it not} all nonstandard integration structures can be obtained by the procedure in
Proposition 6.1.1. In this regard we have \\ \\

{\bf Example 6.1.2} \\

Let $E = \{~ e_1, e_2, e_3, ~.~.~.~ , e_k ~\}$ be a hyperfinite set in $V ( ^*X )$, corresponding to a certain $k \in~
^*{\bf N}_\infty$. Let ${\cal F}$ be the set of all functions $f : E  ~\longrightarrow ~^*{\bf R}$. Finally, let $c_l \in~  ^*{\bf R},
~ c_l \geq 0$, with $1 \leq l \leq k$. Now as in Example 6.1.1, for every $f \in {\cal F}$ we can define the uncountable
sum

$$ I f ~=~ \sum_{1 \leq l \leq k}~ c_l f ( l ) ~\in~ ^*{\bf R} $$ \\

Then $( {\cal F}, I )$ is a nonstandard integration structure on $E$, and it {\it cannot} be obtained by the procedure in
Proposition 6.1.1. \\

Indeed, suppose given a set $F$ in $V ( X )$ and a standard integration structure $( {\cal G}, J )$ on $F$, such that $E =~
^*F,~ {\cal F} =~ ^*{\cal G}$ and $I =~ ^*J$. Now the relation $E =~ ^*F$ implies that $E$ is standard, see Definition 5.6.2,
which together with the hypothesis that $E$ is hyperfinite will result in $E$ being a finite set, see Note 5.7.4, while
clearly $E$ is by choice infinite.  \\ \\

{\bf 2. The Loeb Construction of Integration Structures \\
\hspace*{0.5cm} by Standardization} \\ \\

In this section, we shall show how to associate with every {\it internal nonstandard} integration structure $( {\cal F}, I )$
on an {\it internal} set $E$ in the Superstructure $V ( ^*X )$ a {\it standard} integration structure $( {\cal F}^{\#}, I^{\#} )$ on
the same set $E$. \\ \\

{\bf Definition 6.2.1} \\

Given a set $E$ in the Superstructure $V ( ^*X )$, we call $( {\cal F}, I )$ an {\it internal} nonstandard integration structure
on $E$, if and only if $( {\cal F}, I )$ is a nonstandard integration structure, and in addition $E,~ {\cal F}$ and $I$ are
internal. \\ \\

{\bf Note 6.2.1} \\

1.~~~If $( {\cal F}, I )$ is an internal nonstandard integration structure on an internal set $E$ in the Superstructure
$V ( ^*X )$ then every $f \in {\cal F}$ will also be internal, see (5.7.3). \\

2.~~ In Proposition 6.1.1, the resulting $( ^*{\cal F},~ ^*I )$ is obviously and internal nonstandard integration structure on
$^*E$, see Note 5.7.4. \\

3.~~ If in Example 6.1.2 we consider ${\cal F}$ constituted of all  {\it internal} functions $f : E \longrightarrow~ ^*{\bf R}$,
then ${\cal F}$ itself is internal, see ??. Therefore, $( {\cal F}, I )$ will be an internal nonstandard integration structure.

\hfill $\Box$ \\

In the sequel, we assume that the respective Superstructures $V ( ^*X )$ are {\it countably saturated}, see sections 2.8
and 5.9. \\

Also, we shall assume given an internal nonstandard integration structure $( {\cal F}, I )$ on an internal set $E$ in
the Superstructure $V ( ^*X )$. \\ \\

{\bf Definition 6.2.2} \\

We define the set ${\cal F}_0$ of {\it null} functions in the internal nonstandard integration
structure ${\cal F}$ as follows. If $g : E ~\longrightarrow~
^*{\bf R}$ is an internal or external function, then $g \in {\cal F}_0$, if and only if

\bigskip
(6.2.1) \quad $ \begin{array}{l}
                                     \forall~ \epsilon \in {\bf R},~ \epsilon > 0 ~: \\ \\
                                     \exists~ \phi \in {\cal F} ~: \\ \\
                                     ~~~ | g | ~\leq~ \phi,~~~ I \phi \in Gal ( 0 ),~~~
                                      st ( I \phi ) ~\leq~ \epsilon
                         \end{array} $

\medskip
We define now ${\cal F}^{\#}$ as the set of all functions $f : E ~\longrightarrow~ {\bf R}$ which admit a decomposition

\bigskip
(6.2.2) \quad $ f ~=~ \phi + g,~~~ \phi \in {\cal F},~~ I | \phi | \in Gal ( 0 ),~~ g \in {\cal F}_0 $ \\ \\

{\bf Note 6.2.2} \\

It is clear that in general we have

\bigskip
(6.2.3) \quad $ {\cal F}_0 \not \subseteq {\cal F},~~~ {\cal F}^{\#} \not \subseteq {\cal F} $ \\ \\

{\bf Proposition 6.2.1} \\

Given two decompositions

$$ f ~=~ \phi ~+~ g ~=~ \phi^\prime ~+~ g^\prime \in {\cal F}^{\#} $$

\medskip
with $\phi,~ \phi^\prime \in {\cal F},~ I | \phi | \in Gal ( 0 )$ and $g,~ g^\prime \in {\cal
F}_0$, then

$$ I | \phi^\prime | \in Gal ( 0 ),~~~ \phi^\prime ~-~ \phi \in {\cal F}_0,~~~ I \phi ~\approx~ I \phi^\prime $$

\medskip
Furthermore, given two decompositions

$$ f_i ~=~ \phi_i ~+~ g_i \in {\cal F}^{\#},~~~ 1 ~\leq~ i ~\leq~ 2 $$

\medskip
with $\phi_i \in {\cal F},~ I | \phi_i | \in Gal ( 0 ),~ g_i \in {\cal F}_0$, then

$$ ( f_1 \wedge f_2 ) ~-~ ( \phi_1 \wedge \phi_2 ),~~ ( f_1 \vee f_2 ) ~-~ ( \phi_1 \vee \phi_2 ) \in {\cal F}_0 $$ \\

{\bf Proof} \\

?? \\ \\

{\bf Theorem 6.2.1} \\

${\cal F}_0$ and ${\cal F}^{\#}$ are vector lattices on $E$. \\

{\bf Proof} \\

??

\hfill $\Box$ \\

In view of Proposition 6.2.1, we can give \\ \\

{\bf Definition 6.2.3} \\

The functional $I^{\#} : {\cal F}^{\#} ~\longrightarrow~ {\bf R}$ is defined by

$$ I^{\#} f ~=~ st ( I \phi ) $$

\medskip
where $f = \phi + g,~ \phi \in {\cal F},~ I | \phi | \in Gal ( 0 ),~ g \in {\cal F}_0$ \\ \\

{\bf Theorem 6.2.2} \\

$( {\cal F}^{\#}, I^{\#} )$ is an integration structure on $E$. \\

{\bf Proof} \\

?? \\ \\

{\bf Definition 6.2.4} \\

$( {\cal F}^{\#}, I^{\#} )$ is called the {\it standardization} of $( {\cal F}, I )$. \\ \\

{\bf Important Note 6.2.3} \\

1.~~~ The above construction of standardization

$$ ( {\cal F}, I ) ~~\longmapsto~~ ( {\cal F}^{\#}, I^{\#} ) $$

\medskip
which turns an internal nonstandard integration structure $( {\cal F}, I )$ on an internal set $E \in~ ^*V ( X )$ into an
integration structure $( {\cal F}^{\#}, I^{\#} )$ on the same set $E$, proves to be {\it natural}. Indeed, the crucial point of
that construction is the decomposition (6.2.2), namely

$$  f ~=~ \phi + g,~~~ \phi \in {\cal F},~~ I | \phi | \in Gal ( 0 ),~~ g \in {\cal F}_0 $$

\medskip
which by definition characterizes every function $f : E ~\longrightarrow~ {\bf R}$ that belongs to ${\cal F}^{\#}$. And in
this decomposition one should have $ \phi \in {\cal F}$, so that one may be able to compute $I \phi$, based on which
one intends to compute $I^{\#} f$. Further, one should also have $I | \phi | \in Gal ( 0 )$, in order to obtain usual real
values for $I^{\#}f$. Finally, one should also have $g \in {\cal F}_0$, so that in view of (6.2.1) one can disregard $g$ in
the integration process, and thus remain in $I^{\#}f$ only with the contribution of $I \phi$. \\
Needless to say, one must decompose $f$ in at least {\it two} parts as above, since $f$, and in fact, all the functions in
${\cal F}^{\#}$, are desired to have usual real values only, while the functions in ${\cal F}$, upon which the given
nonstandard integration $I$ is defined, may have nonstandard values as well. Therefore, in the above decomposition
$\phi$ and $g$ are to compensate one another in such a way that their sum is {\it always} a usual real number. \\

2.~~~ In view of (6.2.3), in general $( {\cal F}, I )$ is {\it not} an extension of $( {\cal F}^{\#}, I^{\#} )$. \\ \\

{\bf Theorem 6.2.3} \\

If $f : E ~\longrightarrow~ {\bf R}$, then $f \in {\cal F}^{\#}$, if and only if

\bigskip
(6.2.4) \quad $ \begin{array}{l}
                                  \forall~ \epsilon \in {\bf R},~ \epsilon > 0 ~: \\ \\
                                  \exists~ \alpha,~ \beta \in {\cal F} ~: \\ \\
                                  ~~~ \alpha ~\leq~ f ~\leq~ \beta,~~~ I | \alpha | \in Gal ( 0 ),~~~ I ( \beta ~-~ \alpha ) ~\leq~ \epsilon
                          \end{array} $

\medskip
in which case we also have

\bigskip
(6.2.5) \quad $ st ( I \alpha ) ~\leq~ I^{\#} f ~\leq~ st ( I \alpha ) ~+~ \epsilon $ \\

{\bf Proof} \\

?? \\ \\

{\bf Theorem 6.2.4 ( Monotone Convergence )} \\

Suppose $f_n \in {\cal F}^{\#}$, with $1 \leq n < \infty$ is monotone increasing, $f : E ~\longrightarrow~ {\bf R}$,
and the following two conditions are satisfied

\bigskip
(6.2.6) \quad $ \lim_{n \to \infty}~ f_n ( x ) ~=~ f ( x ) ~~\mbox{exists for}~~ x \in E $ \\

(6.2.7) \quad $ \lim_{n \to \infty}~ I^{\#} f_n ~<~ \infty $

\medskip
Then $f \in {\cal F}^{\#}$ and $I^{\#} f ~=~ \lim_{n \to \infty}~ I^{\#} f_n$. \\

{\bf Proof} \\

?? \\ \\

{\bf Theorem 6.2.5} \\

If $h \in {\cal F}$ only takes finite values, thus $h : E ~\longrightarrow~ Gal ( 0 )$, and there exists $\phi \in
{\cal F}$ with $\phi \geq 0,~ I \phi \in Gal ( 0 )$, such that

$$ \{~ x \in E ~~|~~ h ( x ) ~\neq~ 0 ~\} ~\subseteq~ \{~ x \in E ~~|~~ \phi ( x ) ~\geq~ 1 ~\} $$

\medskip
then

$$ h ~-~ st ( h ) \in {\cal F}_0,~~~ st ( h ) \in {\cal F}^{\#},~~~ I^{\#} st ( h ) ~=~ st ( I h ) $$ \\

{\bf Proof} \\

?? \\ \\

{\bf Theorem 6.2.6} \\

If $1 \in {\cal F}$ and $I 1 \in Gal ( 0 )$, then $1 \in {\cal F}^{\#}$. Furthermore, if $f : E ~\longrightarrow~
{\bf R}$ is of the form

$$ f ~=~ \phi ~+~ g,~~~ \phi \in {\cal F},~~~ g \in {\cal F}_0 $$

\medskip
then $I | \phi | \in Gal ( 0 )$, thus $f \in {\cal F}^{\#}$. \\

{\bf Proof} \\

?? \\ \\

{\bf 3. Loeb Measures} \\ \\

In view of Theorems 6.2.2 and 6.2.4, we can restrict ourselves to integration structures $( {\cal F}^{\#}, I^{\#} )$, given for
instance by standardization, see Definition 6.2.4, and which are described by the next definition. Here we recall that
such integration structures are defined on a set $E \in V ( ^*X )$, and ${\cal F}^{\#} \subseteq {\bf R}^E$, while $I^{\#} :
{\cal F}^{\#} \longrightarrow {\bf R}$, see Definition 6.1.2. \\ \\

{\bf Definition 6.3.1} \\

An integration structure $( {\cal F}^{\#}, I^{\#} )$ on a set $E$ in $V ( ^*X )$ is called {\it complete}, if and only if it has the
following property. Given any sequence of functions $f_n \in {\cal F}^{\#}$, with $n \in {\bf N}$, monotone increasing and
such that

\bigskip
(6.3.1) \quad $ \exists~ f : E \longrightarrow~ {\bf R} ~:~ \lim_{n \to \infty}~ f_n ( x ) ~=~ f ( x ),~~ x \in E $

\medskip
(6.3.2) \quad $ \lim_{n \to \infty}~ I^{\#} f_n ~<~ \infty $

\medskip
then

\bigskip
(6.3.3) \quad $ f \in {\cal F}^{\#},~~~  \lim_{n \to \infty}~ I^{\#} f_n ~=~ I^{\#} f $

\hfill $\Box$ \\

In this section the integration structures $( {\cal F}^{\#}, I^{\#} )$ will be considered complete. \\

Let us now associate with $( {\cal F}^{\#}, I^{\#} )$ a family of functions ${\cal M}^{\#} \subseteq (~ \overline{{\bf R}} ~)^E$
which contains ${\cal F}^{\#}$. Here as usual, we denoted $\overline{{\bf R}} = {\bf R} \cup \{~ - \infty,~ + \infty ~\}$.\\
These functions $f$ in ${\cal M}^{\#}$ play the role of {\it measurable} functions. Therefore, it may happen that the
integral $I^{\#} f$ is {\it not} defined for such functions, or may no longer be finite, taking the value $\pm \infty$. An
advantage of the set ${\cal M}^{\#}$ of measurable functions is that it is an {\it algebra}, unlike in general the set ${\cal
F}^{\#}$. \\

Further, we shall extend the integration $I^{\#}$ from ${\cal F}^{\#}$ to a subset ${\cal F}^{\#}_1 \subseteq {\cal M}^{\#}$,
obtaining in this way an extension $( {\cal F}^{\#}_1, I^{\#}_1 )$ of $( {\cal F}^{\#}, I^{\#} )$. \\
Finally, we shall call {\it measurable} those subsets $A \subseteq E$ whose characteristic function $\chi_A$ is in ${\cal
M}^{\#}$. \\

It should be noted that, as seen later, the measure and integration theory started from the above measurable sets is
equivalent with that started with integration. \\ \\

{\bf Definition 6.3.2} \\

We define

\bigskip
(6.3.4) \quad $ {\cal M}^{\#+} ~=~ \{~ f : E \longrightarrow [~ 0, + \infty ~] ~~|~~ \forall~ h \in {\cal F}^{\#} ~:~
                                            f \wedge h \in {\cal F}^{\#} ~\} $

\medskip
Then we define $J^{\#} : {\cal M}^{\#+} \longrightarrow [~ 0, + \infty ~]$ by

\bigskip
(6.3.5) \quad $ J^{\#} f ~=~ \sup~ \{~ I^{\#} ( f \wedge h ) ~~|~~ h \in {\cal F}^{\#} ~\},~~~ f \in {\cal M}^{\#+} $ \\

Further, we define

\bigskip
(6.3.6) \quad $ {\cal M}^{\#} ~=~ \{~ f : E \longrightarrow \overline{{\bf R}} ~~|~~ f^{-} = f \wedge 0,~~ f^{+} = f \vee 0 \in
                                        {\cal M}^{\#+} ~\} $

\medskip
Also, we define

\bigskip
(6.3.7) \quad $ {\cal M}^{\#}_1 ~=~ \{~ f : E \longrightarrow \overline{{\bf R}} ~~|~~ J^{\#} f^{+} ~~~\mbox{or}~~~ J^{\#} f^{-}
                                     ~~~\mbox{are finite} ~\} $

\medskip
Finally, we define $J^{\#} : {\cal M}^{\#}_1 \longrightarrow \overline{{\bf R}}$ by

\bigskip
(A6.3.8) \quad $ J^{\#} f ~=~ J^{\#} f^{+} ~-~ J^{\#} f^{-},~~~ f \in {\cal M}^{\#}_1 $

\hfill $\Box$ \\

With these preliminaries, we can now give \\ \\

{\bf Definition 6.3.3} \\

Let

\bigskip
(6.3.9) \quad $ {\cal F}^{\#}_1 ~=~ \{~ f \in {\cal M}^{\#}_1 ~~|~~ ~~~\mbox{both}~~~ J^{\#} f^{+} ~~~\mbox{and}~~~ J^{\#} f^{-}
                                     ~~~\mbox{are finite} ~\} $

\medskip
and let $I^{\#}_1 : {\cal F}^{\#}_1 \longrightarrow {\bf R}$ be given by

\bigskip
(6.3.10) \quad $ I^{\#}_1 f ~=~ J^{\#} f^{+} ~-~ J^{\#} f^{-},~~~ f \in {\cal F}^{\#}_1 $

\hfill $\Box$ \\

We note that

\bigskip
(6.3.11) \quad $ {\cal F}^{\#} ~\subseteq~ {\cal F}^{\#}_1 ~\subseteq~  {\cal M}^{\#}_1 ~\subseteq~ {\cal M}^{\#} $

\medskip
and $I^{\#}_1$ is the restriction of $J^{\#}$ to ${\cal F}^{\#}_1$.

\end{document}